\documentclass[10pt,twoside]{amsbook} 
\usepackage{amscd, euscript, amssymb}
\makeindex

\makeatletter
\def\chapter{\clearpage \thispagestyle{plain} \global\@topnum\z@
  \@afterindentfalse \secdef\@chapter\@schapter} 
 \@addtoreset{equation}{chapter}
\@addtoreset{figure}{chapter} \@addtoreset{footnote}{chapter}
\@addtoreset{section}{chapter}
 \makeatother

\newcommand{\itemm}[1]{\item[\textrm{#1)}]}
\newcommand{\bir}{\dasharrow}
\newcommand{\mt}[1]{{\operatorname{#1}}}
\newcommand{\compos}{\mathbin{\scriptstyle{\circ}}}
\newcommand{\red}{_{\mt{red}}}
\newcommand{\Supp}{\mt{Supp}}
\newcommand{\Const}{\mt{Const}}
\newcommand{\Pic}{\mt{Pic}}
\newcommand{\Sing}{\mt{Sing}}
\newcommand{\LCS}{\mt{LCS}}
\newcommand{\GL}{\mt{GL}}
\newcommand{\PGL}{\mt{PGL}}
\newcommand{\Bs}{\mt{Bs}}
\newcommand{\lcm}{\mt{lcm}}
\newcommand{\Weil}{\mt{Weil}}
\newcommand{\Weillin}{\Weil_{\mt{lin}}}
\newcommand{\Weilalg}{\Weil_{\mt{alg}}}
\newcommand{\rhonum}{\rho_{\mt{num}}}
\newcommand{\compl}{\mt{compl}}
\newcommand{\mult}{\mt{mult}}
\newcommand{\codim}{\mt{codim}}
\newcommand{\Diff}[1]{\mt{Diff}_{#1}}
\newcommand{\NE}{\overline{{NE}}}
\newcommand{\fr}[1]{\left\{ #1 \right\}}

\newcommand{\ov}[1]{\overline{#1}}
\newcommand{\down}[1]{\left\lfloor #1\right\rfloor}
\newcommand{\downn}[1]{\left[ #1\right]}

\newcommand{\wcir}[1]{\setlength{\unitlength}{2mm}
\hbox{\begin{picture}(2,5.5)(0,-5.0)
 \put(1,-4.3){\circle{2}}
 \put(1,-2.3){\makebox(0,0)[b]{$#1$}}
 \end{picture}}}

\newcommand{\wci}{\setlength{\unitlength}{2mm}
\hbox{\begin{picture}(2,2)(0,-5.0)
 \put(1,-4.3){\circle{2}}
 \end{picture}}}

\newcommand{\bci}{\setlength{\unitlength}{2mm}
\hbox{\begin{picture}(2,2)(0,-5.0)
 \put(1,-4.3){\circle*{2}}
 \end{picture}}}

\newcommand{\bcir}[1]{\setlength{\unitlength}{2mm}
\hbox{\begin{picture}(2,5.5)(0,-5.0)
 \put(1,-4.3){\circle*{2}}
 \put(1,-2.3){\makebox(0,0)[b]{$#1$}}
 \end{picture}}}

\newcommand{\cirmin}[1]{\setlength{\unitlength}{2mm}
\hbox{\begin{picture}(2,5.5)(0,-5.0)
 \put(1,-4.3){\circle{2}}
 \put(0.3,-4.29){\line(1,0){1.5}}
 \put(1,-2.3){\makebox(0,0)[b]{$#1$}}
 \end{picture}}}

\newcommand{\cirmi}{\setlength{\unitlength}{2mm}
\hbox{\begin{picture}(2,2)(0,-5.0)
 \put(1,-4.3){\circle{2}}
 \put(0.3,-4.29){\line(1,0){1.5}}
 \end{picture}}}

\newcommand{\vlin}{\setlength{\unitlength}{2mm}
\hbox{\begin{picture}(0.5,4.3)
 \put(0.25,0){\line(0,1){3.8}}
 \end{picture}}}

\newcommand{\lin}{\setlength{\unitlength}{2mm}
\hbox{\begin{picture}(4,2)(0,-5.0)
 \put(0,-4.3){\line(1,0){4}}
 \end{picture}}}

\newcommand{\up}[1]{\left\lceil #1\right\rceil}
\newcommand{\qq}{\mathbin{\sim_{\scriptscriptstyle{\mathbb Q}}}}
\newcommand{\ep}{\varepsilon}
\newcommand{\var}{\varphi}
\newcommand{\De}{\Delta}

\newcommand{\CC}{{\mathbb{C}}}
\newcommand{\EE}{{\mathbb{E}}}
\newcommand{\RR}{{\mathbb{R}}}
\newcommand{\ZZ}{{\mathbb{Z}}}
\newcommand{\DD}{{\mathbb{D}}}
\newcommand{\BA}{{\mathbb{A}}}
\newcommand{\QQ}{{\mathbb{Q}}}
\newcommand{\PP}{{\mathbb{P}}}
\newcommand{\NN}{{\mathbb{N}}}
\newcommand{\FF}{{\mathbb{F}}}
\newcommand{\OOO}{{\mathcal{O}}}

\newcommand{\EEE}{{\mathcal{E}}}
\newcommand{\LLL}{{\mathcal{L}}}
\newcommand{\FFF}{{\mathcal{F}}}
\newcommand{\KKK}{{\EuScript{K}}}

\newcommand{\NNN}{{\EuScript{N}}}
\newcommand{\MMM}{\Phi}
\newcommand{\Msm}{{\MMM}_{\mt{\mathbf{sm}}}}
\newcommand{\Mm}{{\MMM}_{\mt{\mathbf{m}}}}
\newcommand{\RRR}{{\EuScript{R}}}
\newcommand{\PPP}{{\EuScript{P}}}
\newcommand{\cyc}[1]{\mathbb{Z}_{#1}}
\renewcommand{\emptyset}{\varnothing}
\newcommand{\eref}[1]{\rm{(\ref{#1})}}

\newtheorem{theorem}[subsection]{Theorem}
\newtheorem{theorem1}[subsection]{Theorem}

\newtheorem{proposition}[subsection]{Proposition}
\newtheorem{proposition1}[subsection]{Proposition}
\newtheorem{conjecture}[subsection]{Conjecture}
\newtheorem{conjecture1}[subsection]{Conjecture}
\newtheorem{lemma}[subsection]{Lemma}
\newtheorem{lemma1}[subsection]{Lemma}
\newtheorem*{claim*}{Claim}
\newtheorem{corollary}[subsection]{Corollary}
\newtheorem{corollary1}[subsection]{Corollary}
\theoremstyle{definition}
\newtheorem*{definition*}{Definition}
\newtheorem{definition}[subsection]{Definition}
\newtheorem{definition1}[subsection]{Definition}
\newtheorem{example}[subsection]{Example}
\newtheorem{example1}[subsection]{Example}
\newtheorem{proposition-definition}[subsection]{Proposition-Definition}
\newtheorem{proposition-definition1}[subsection]{Proposition-Definition}
\newtheorem{corollary-definition}[subsection]{Corollary-Definition}
\newtheorem{corollary-definition1}[subsection]{Corollary-Definition}

\newtheorem{problem1}[subsection]{Problem}
\theoremstyle{remark} \theoremstyle{remark}
\newtheorem{remark}[subsection]{Remark}
\newtheorem{remark1}[subsection]{Remark}
\newtheorem*{remark*}{Remark}
\newtheorem{exercise1}[subsection]{Exercise}

\newenvironment{outline}{\begin{proof}[Sketch of proof]}{\end{proof}}
\newenvironment{hint}{\ \textit{Hint.}\ }{}

\begin{document}

\frontmatter

\vspace*{1in}
\thispagestyle{empty}
\begin{centering}
\rule{5in}{.04in}\vspace{.25in}
{\Huge Lectures on complements on log surfaces}\\
\vspace{.1in}
\rule{5in}{.04in}\vspace{1.5in}\\
\large by\\
\vspace{.2in}
\LARGE {\bf Yuri G. Prokhorov}\\ \vspace{0.9in}
\large
Department of Algebra, Faculty of Mathematics\\
Moscow State Lomonosov University\\ 
Moscow 117234, RUSSIA\\ \vspace{.1in}
\normalsize E-mail: prokhoro@mech.math.msu.su\\
\vspace{0.9in}
{\it AMS Mathematics subject classification}: Primary
14E30; Secondary 14E05, 14J17, 14J45, 14E15
\end{centering}

\tableofcontents
\setcounter{chapter}{-1}
\mainmatter

\chapter{Introduction}
These notes grew out of lectures I gave at Moscow University.
Basically we follow works of V.~V.~Shokurov \cite{Sh}, \cite{Sh1}.
These notes can help the reader to understand the main ideas of
the theory. In particular, they can be considered as an
introduction to \cite{Sh1}.
\par

One of the basic problems in the modern birational minimal model
program is the problem of constructing a divisor with rather
``good'' singularities in the anticanonical or multiple
anticanonical linear system. Probably this question arises for the
first time in the classification of Fano $3$-folds (see e.g.,
\cite{Isk}) and was solved by V.~V.~Shokurov \cite{Sh0}. Later
this technique was improved by Kawamata and others (see e.g.,
\cite{R}). However, this method only applies to linear systems of
Cartier divisors.

\par
Another approach to constructing ``good'' divisors in the
anticanonical linear system was proposed by Mori \cite{Mo} in the
proof of existence of three-dimensional flips. Unfortunately, now
this method was applied only in analytic situation and in
dimension $3$ (see \cite{KoM}, \cite{Pr} and also \cite{Ka}).

\par
The concept of complement was introduced by V.~V.~Shokurov in his
proof of the existence of good divisors near nontrivial fibers of
small contractions \cite{Sh}. Roughly speaking, an $n$-complement
of the canonical divisor $K_X$ is an element of the multiple
anticanonical linear system $D\in |{-}nK_X|$ such that the pair
$\left(X, \frac1{n}D\right)$ has only log canonical singularities
(for precise definitions we refer to \ref{ldef}). Thus the
question of the existence of complements can be posed for Fano or
Calabi-Yau varieties (i.e., with numerically effective
anticanonical divisor), and also for varieties with a fiber space
structure on varieties of such types. For example, if on a smooth
(or even with log canonical singularities) variety $X$ the
anticanonical divisor ${-}K_X$ is ample, then by Bertini's theorem
there is an $n$-complement for some $n\in\NN$. In case of
Calabi-Yau varieties of the property of the canonical divisor
$K_X$ to be $n$-complementary is equivalent to triviality of
$nK_X$. For example, it is known that for smooth surfaces such $n$
can be taken in $\{1,2,3,4,6\}$.
\par
These notes aim to do two things. The first is to give an
introduction to the theory of complements, with rigorous proofs
and motivating examples. For the first time the reader can be
confused by rather tricky definition of $n$-complements (see (iii)
of \ref{ldef}). However this definition is justified because of
their useful properties -- both birational and inductive (see
\ref{Birational-properties-of-complements} and
\ref{Inductive-properties-of-complements}). This shows that the
property of the existence of complements is very flexible: we can
change the birational model and then extend complements from the
reduced part of the boundary.
\par
The second is to prove several concrete results on the existence
of complements. Our results, presented mainly in dimension two,
illustrate main ideas and difficulties for higher-dimensional
generalizations (cf. \cite[Sect.7]{Sh1}, \cite{Pr1}, \cite{Pr2},
\cite{Pr3}).

The first type of results concerns with the inductive implication
(local) $\Longrightarrow$ (lower dimensional global). Although we
prove this result in dimension two only, it can be generalized for
the case of arbitrary dimension without difficulties (see
\cite{PSh}). The next type of results are inductive implications
(global) $\Longrightarrow$ (lower dimensional global). Roughly
speaking these results assert that if a log pair has singularities
worse than Kawamata log terminal, then complements can be induced
from lower dimensions. This result is not difficult for varieties
of log Calabi-Yau and log Fano type (see \cite{Ishii1},
\cite[Sect. 5]{PSh}) and becomes very hard for varieties of
intermediate numerical anticanonical dimension. We expect that
this result can be generalized in arbitrary dimension modulo log
MMP, groups of automorphisms of Calabi-Yau and Kodaira-type
formula. Our proof (which is a reworked version of \cite{Sh1}) is
based on detailed analysis of the reduced part of a (non-klt)
boundary.

The last group of discussed problems is the boundedness of
exceptional log del Pezzo surfaces (refer to \ref{s_except} for
the precise definition). The ideal result (which is expected to be
valid in any dimension) is that these surfaces have bounded
complements and they belong to a finite number of ``families''.
The proofs uses heavily Alexeev-Nikulin's theorems about
boundedness of log del Pezzo surfaces with $\epsilon$-log terminal
singularities and inductive results (see above). In recent work
\cite{KeM} similar techniques were used to study families of
rational curves on nonprojective surfaces. It is very important
problem here is to classify exceptional log del Pezzos. Some
results in this direction are obtained in \cite{KeM}, \cite{Sh1},
\cite{Abe} but in general the question is still open. We present,
following \cite{Sh1}, the classification of exceptional log Del
Pezzos with ``large'' boundary. This classification is interesting
for applications to three-dimensional birational geometry. For
example, there is a conjecture that any Mori contraction is
nonexceptional (see \ref{12346}, cf. \cite[6.4]{RYPG}). Adopting
this conjecture one can by using the inductive procedure construct
an $1$, $2$, $3$, $4$, or $6$-complements which leads to the
classification of all Mori contractions (see \cite[Sect. 7]{Sh1},
\cite{Pr2}, cf. \cite{KoM}).

As applications of the developed techniques, we reprove several
well known facts such as classifications of quotient singularities
and singular fibers of elliptic fibrations. This outlines ways for
higher-dimensional generalizations. Note, however, that these
expected generalizations cannot be straightforward. For example
the classification of all three-dimensional log-canonical (and
even canonical) singularities is not a reasonable problem (cf.
\cite{IP}).

The following diagram exhibit the logical structure of the results
(proved in dimension $2$ and conjectured in higher dimensions). We
omit some technical details and conditions (for example
restrictions on the coefficients of the boundaries).
\par\medskip

\begin{center}
\begin{tabular}{lr}
\multicolumn{2}{c}{\fbox{\parbox{12cm}{$(X/Z\ni o,D)$: a
$d$-dimensional log pair of local type (i.e. a contraction with
$\dim Z>0$) such that $-(K_X+D)$ is nef and big}}}
\\ & \\
\multicolumn{2}{c}{EXCEPTIONAL ?}\\ 
\multicolumn{1}{c}{YES}&\multicolumn{1}{c}{NO}\\ 
\\ 
\multicolumn{1}{c}{$\Downarrow$}&\multicolumn{1}{c}{$\Downarrow$}\\
\fbox{\parbox{5cm}{complements are induced from $d-1$-dimensional
log Fano pairs (possibly of global
type)}}&\fbox{\parbox{5cm}{complements are induced from
$d-2$-dimensional log Fano pairs \phantom{(possibly of global
type)}}}
\\ &\\&\\
\multicolumn{2}{c}{\fbox{\parbox{12cm}{$(X,D)$: a $d$-dimensional
log pair of global type (i.e. $X$ is compact) with nef and big
$-(K_X+D)$}}}
\\ & \\
\multicolumn{2}{c}{EXCEPTIONAL ?}\\
\multicolumn{1}{c}{YES}&\multicolumn{1}{c}{NO}\\ 
\\ 
\multicolumn{1}{c}{$\Downarrow$}&\multicolumn{1}{c}{$\Downarrow$}\\
\fbox{\parbox{5cm}{complements are bounded by a constant
$C(d)$}}&\fbox{\parbox{5cm}{complements are induced from
$d-1$-dimensional log Fano pairs}}\\
\end{tabular}
\end{center}

\par\bigskip
Our notes are organized as follows. We start with three
introductory chapters: Chapter~\ref{sect-1} gives precise
definitions of singularities of pairs and gives some known
constructions. Chapter~\ref{sect-2} contains some important tools:
Inversion of Adjunction and Connectedness Lemma. In
Chapter~\ref{sect-3} we collect facts about log terminal
modifications (these are consequences of the Minimal Model
Program, and therefore restricted to dimension two or three). In
Chapter~\ref{sect-4} we define complements and present some
properties. The following three Chapters deal with applications of
complements to surface geometry. The central part of these notes
is chapters~\ref{sect-8} and \ref{sect-9}. Roughly speaking, all
log surfaces can be divided into two classes: exceptional and
nonexceptional. If the log surface is nonexceptional, then under
some additional assumptions there is a regular (i.e., $1$, $2$,
$3$, $4$ or $6$-complement). This is done in Chapter~\ref{sect-8}.
The main result there is the Inductive
Theorem~\ref{Inductive-Theorem}. On the contrary, exceptional log
surfaces can have no regular complements. However, we show in
Theorems~\ref{main_Sh_A-1} and \ref{main_Sh_A-2} that they are
bounded (again under additional assumptions). In
Chapter~\ref{last-1} we discuss the problem of classification of
exceptional complements.
\par
Our presentation of the subject is rather elementary. The notes
contain many examples.
\par

\par\medskip\noindent
\textsc{Acknowledgements.} I am grateful to V.~V.~Shokurov for
numerous useful discussions. I am also grateful to the
participants of my seminar at Moscow State University, especially
R.~A.~Sarkisyan, for attention and useful remarks. In the spring
of 1999 I gave lectures on this subject at Waseda University. Deep
gratitude should be given to the participants of the Waseda
Tuesday seminar. I would like to express my thanks to Professors
S.~Ishii and V.~Ma\c{s}ek, who read preliminary versions of the
paper and gave me useful criticisms which served to improve the
text and to correct some mistakes. The final version of these
notes was prepared during my stay in Tokyo Institute of
Technology. I thank the staff of this institute for hospitality.
Finally, this revised version owes much to the referees' numerous
suggestions, both in mathematics and English, which are greatly
appreciated. This work also was partially supported by the grant
INTAS-OPEN-97-2072.

\chapter{Preliminary results}
\label{sect-1}
\section{Singularities of pairs}
\label{Singularities-of-pairs}
\subsection*{List of notations}
\quad
\newline
\begin{tabular}{ll}
 $\equiv$&numerical equivalence\\
 $\sim$&linear equivalence\\
 $\qq$&$\QQ$-linear equivalence\\
 $\KKK(X)$&function field of $X$\\
$D\approx D'$&$D$ and $D'$ gives the same valuation of $\KKK(X)$\\
$\rho(X)$\index{$\rho(X)$} &Picard number of $X$, rank of the
N\'eron-Severi group\\ $Z_1(X/Z)$&group of $1$-cycles on $X$ over
$Z$ (see \cite{KMM})\\
 $N_1(X/Z)$\index{$N_1(X/Z)$}&quotient of $Z_1(X/Z)$ modulo numerical
 equivalence (cf. )\\
$\NE(X/Z)$\index{$\NE(X/Z)$}&Mori cone (see \cite{KMM})\\
 $\Weil(X)$&group of Weil divisors,
i.e., the free abelian group \\& generated by prime divisors on
$X$\index{$\Weil(X)$}\\ \!\!\!$\left.\begin{array}{l}\Weillin(X)\\
\Weilalg(X)
\end{array}\right\}$ &\!\!\!\begin{tabular}{l}
quotients of $\Weil(X)$ modulo linear and algebraic\\ equivalence
respectively.\index{$\Weillin(X)$}\index{$\Weilalg(X)$}
\end{tabular}
\end{tabular}
\par\medskip\noindent
All varieties are assumed to be algebraic varieties defined over
the field $\CC$. By a \textit{contraction}\index{contraction} we
mean a projective morphism $f\colon X\to Z$ of normal varieties
such that $f_*\OOO_X=\OOO_Z$ (i.e., having connected fibers). We
call a birational contraction a \textit{blowdown}\index{blowdown}
or \textit{blowup}, \index{blowup} depending on our choice of
initial variety.
\par
A \textit{boundary}\index{boundary} on a variety $X$ is a
$\QQ$-Weil divisor $D=\sum d_iD_i$ with coefficients $0\le d_i\le
1$. If we have only $d_i\le 1$, we say that $D$ is a
\textit{subboundary}. \index{subboundary} All varieties are
usually considered supplied with boundary (or subboundary) as an
additional structure. If $D$ is a boundary, then we say that
$(X,D)$ is a \textit{log variety} or \textit{log pair}. \index{log
variety (log pair)} Moreover, if we have a contraction $f\colon
X\to Z$, then we say that $(X,D)$ is a \textit{log variety} over
$Z$ and denote it simply by $(X/Z,D)$. If $\dim Z>0$, we often
consider $Z$ as a germ near some point $o\in Z$. To specify this
we denote the corresponding log variety by $(X/Z\ni o,D)$.

Given a birational morphism $f\colon X\to Y$, the boundary $D_Y$
on $Y$ is usually considered as the image of the boundary $D_X$ on
$X$: $D_Y=f_*D_X$. The integral part of a $\QQ$-divisor $D=\sum
d_i D_i$ is defined in the usual way:\quad $\down{D}:=\sum
\down{d_i} D_i$,\index{$\down{\cdot}$} where $\down{d_i}$ is the
greatest integer such that $\down{d_i}\le d_i$. The (round up)
upper integral part $\up{D}$\index{$\up{\cdot}$} and the
fractional part $\fr{D}$\index{$\fr{\cdot}$} are similarly
defined.

\par
\textit{A log resolution}\index{log resolution} is a resolution
$f\colon \widetilde X\to X$ of singularities of $X$ such that the
union $\left(\bigcup \widetilde D_i\right)\cup\mt{Exc}(f)$ of
proper transforms\footnote{The proper transform is sometimes also
called the birational or strict transform.} of all the $D_i$ and
the exceptional locus the exceptional locus $\mt{Exc}(f)$ is a
divisor with simple normal crossings.
\par
Let $X$ be a normal variety, $D$ a $\QQ$-divisor on $X$, and
$f\colon\widetilde X\to X$ any projective birational morphism,
where $X$ and $\widetilde X$ are normal. Assume that $K_X+D$ is
$\QQ$-Cartier. Then we can write
\begin{equation}
\label{ddd}
K_{\widetilde X}+\widetilde D=f^*(K_X+D)+\sum_E a(E,D)E,
\end{equation}
where $\widetilde D$ is the proper transform $D$ and
$a(E,D)\in\QQ$. The numbers $a(E,D)$\index{$a(E,D)$} depends only
on $X$, $D$ and the discrete valuation of the field $\KKK(X)$
corresponding to $E$ (i.e., they do not depend on $f$). They are
called \textit{discrepancies} \index{discrepancy} or
\textit{discrepancy coefficients}. Define
\[
\mt{discrep}(X,D):=\inf\nolimits_E\{a(E,D) \mid E\ \text{is an
exceptional divisor over}\ X\}.
\]
\index{$\mt{discrep}(X,D)$} We also put for nonexceptional
divisors
\[
a(E,D):=\left\{
\begin{array}{ll}
 -d_i&\text{if}\ E=D_i;\\
 0&\text{otherwise}.\\
\end{array}
\right.
\]
Let us say that the pair $(X,D)$ has
\begin{itemize}
\item[]\textit{terminal} singularities,
\index{terminal singularities} if $\mt{discrep}(X,D)>0$;
\item[]\textit{canonical} singularities,
\index{canonical singularities} if $\mt{discrep}(X,D)\ge 0$;
\item[]\textit{Kawamata log terminal (klt)} singularities,
\index{Kawamata log terminal (klt)} if $\mt{discrep}(X,D)>-1$ and
$\down{D}\le 0$;
\item[]\textit{purely log terminal (plt)} singularities,
\index{purely log terminal (plt) singularities} if
$\mt{discrep}(X,D)>-1$;
\item[]\textit{log canonical (lc)} singularities,
\index{log canonical (lc) singularities} if $\mt{discrep}(X,D)\ge
-1$;
\item[]\textit{$\ep$-log terminal ($\ep$-lt)} singularities,
\index{$\ep$-log terminal ($\ep$-lt) singularities}\footnote{Note
that our definition of $\ep$-lt pairs is weaker than that given by
Alexeev \cite{A}: we do not claim that $-d_i>-1+\ep$.} if
$\mt{discrep}(X,D)>-1+\ep$;
\item[]\textit{$\ep$-log canonical ($\ep$-lc)} singularities,
\index{$\ep$-log canonical ($\ep$-lc) singularities} if
$\mt{discrep}(X,D)\ge-1+\ep$;
\item[]\textit{divisorial log terminal (dlt)} singularities,
\index{divisorial log terminal (dlt)} if $a(E_i,D)>-1$ for all
exceptional divisors $E_i$ of some log resolution $f$ whose
exceptional locus consists of divisors (or it is empty).
\end{itemize}
In these cases we also say simply that $K_X+D$ is lc (resp. klt,
etc.) We usually omit $D$ if it is trivial.

\subsection{}\label{remark-def-klt}
For the klt and lc properties the inequalities $a(E,D)>-1$
($\ge-1$) can be checked for exceptional divisors of some log
resolution (see \cite[0-2-12]{KMM}). The plt property of $(X,D)$
is equivalent to the existence of a log resolution $f\colon
\widetilde X\to X$ such that $a(E,D)>-1$ for all exceptional
divisors of $f$ and the proper transform $\widetilde{\down{D}}$ of
$\down{D}$ on $\widetilde X$ is smooth. It is easy to see that if
$K_X+D$ is lc, then $D$ is a subboundary. In the two-dimensional
case we can use Mumford's numerical pull back of any Weil divisor,
so all the above definitions can be given in this situation
numerically, without the $\QQ$-Cartier assumption (see e.g.,
\cite{Sakai} and \ref{MMP-in-dimension-2}).

\begin{example1}
Let $X$ be a smooth surface and $D=D_1+D_2$ a pair of smooth
curves intersecting transversally at one point. The identity map
is a log resolution, so $(X,D)$ is dlt. However, the blowup of the
point of intersection gives an exceptional divisor $E$ with
discrepancy $a(E,D)=-1$. Hence $(X, D)$ is not plt. If $D$ is an
irreducible curve with a node on a nonsingular surface $X$, then
the pair $(X,D)$ is not dlt. This shows that the dlt condition is
not local.
\end{example1}

\begin{example1}
Let $Q\subset\CC^4$ be a quadratic cone given by $xy=zt$ and $D$
its hyperplane section $\{ x=0\}$. Then $D=D_1+D_2$, where $D_1$
and $D_2$ are planes in $\CC^4$. There is a small resolution
$f\colon \widetilde{Q}\to Q$ with exceptional locus $\PP^1$. The
intersection of the proper transforms of the planes $D_1$ and
$D_2$ is a line on $\widetilde{Q}$. So $f$ is a log resolution.
However, $(Q, D)$ is not dlt because $f$ is small.
\end{example1}

\par
In \cite{Ut} the notion of weakly Kawamata log terminal
singularity was introduced. Later it was proved that this is
equivalent to the dlt property \cite{Sz}. The very close (but
wider) class of weakly log terminal pairs was considered in
\cite{KMM} and \cite{Sh}.
\par
Log varieties with dlt singularities form a convenient class of
varieties in which the log Minimal Model Program (log MMP) works
\cite{KMM}. In particular, these singularities are rational
\cite[1-3-1]{KMM}, \cite[5.22]{KM} and the Cone Theorem and
Contraction Theorem hold for these varieties \cite[4.2.1,
3-2-1]{KMM}. Log canonical singularities are not necessarily
rational. However, it was shown in \cite{Sh2} that reasonable log
MMP also works in this category.

More precisely, the property of a $\QQ$-factorial log variety to
have klt (resp. dlt) singularities is preserved under contractions
of extremal rays and flips, i.e., they form classes of log
varieties closed under the log MMP. We refer to \cite{KMM} for
technical details of this theory; see also \cite{Sh2} and (for two
dimensional case) Appendix~\ref{MMP-in-dimension-2}, \cite{A},
\cite{KoKov}. Note also that all distinctions between different
notions of log terminal singularities arise only if $D$ has
components with coefficient $1$.
\par
The following property can be obtained directly from the
definitions.

\begin{proposition}[{\cite{Sh}, \cite[2.17]{Ut}}]
\label{first_prop_sing}
Let $X$ be a normal variety and $D=\sum d_iD_i$ a subboundary on
$X$ such that $K_X+D$ is a $\QQ$-Cartier divisor.
\begin{enumerate}
\item
If $D'\le D$ and $K_X+D$ is lc (resp. $\ep$-lt, dlt, plt or klt)
and $K_X+D'$ is $\QQ$-Cartier, then $K_X+D'$ also is lc (resp.
$\ep$-lt, dlt, plt or klt);
\item
If $K_X+D$ is dlt, then there exists $\ep>0$ such that all
$\QQ$-Cartier divisors $K_X+D'$ also are dlt for all $D'=\sum
d_i'D_i$ with $d_i'\le\min \{d_i+\ep, 1 \}$;
\item
If $K_X+D$ is plt (resp. klt) and $K_X+D+D'$ is lc, then
$K_X+D+tD'$ is plt (resp. klt) for all $t<1$.
\end{enumerate}
\end{proposition}

\begin{remark}
\label{equiv-bir}
The formula \eref{ddd} can be written as
\begin{equation}
\label{eqv-bir-new}
K_{\widetilde X}+D'=f^*(K_X+D),
\end{equation}
where $D':=\widetilde D+\sum a_i E_i$. In particular, $D=f_*D'$.
Then the lc property of $K_X+D$ is equivalent to that $D'$ is a
subboundary. In this case, $K_{\widetilde X}+D'$ is called the
\textit{crepant pull back}\index{crepant pull back} of $K_X+D$.
\end{remark}
This trivial remark has the following useful generalization

\begin{proposition1}[\cite{Ko}]
\label{crepant-singul}
Let $f\colon Y\to X$ be a birational contraction and $D$ a
subboundary on $X$ such that $K_X+D$ is $\QQ$-Cartier. As in
\eref{eqv-bir-new} take the crepant pull back
\begin{equation}
\label{for-crepant-singul}
K_Y+D_Y=f^*(K_X+D),\qquad \text{with}\qquad D=f_*D_Y.
\end{equation}
Then
\begin{enumerate}
\item
$K_X+D$ is lc (resp. klt) $\Longleftrightarrow$ $K_Y+D_Y$ is lc
(resp. klt);
\item
$K_X+D$ is plt (resp. dlt) $\Longleftrightarrow$ $K_Y+D_Y$ is plt
(resp. dlt) and $f$ does not contract components of $D_Y$ with
coefficient $1$;
\item
$K_X+D$ is $\ep$-lt $\Longleftrightarrow$ $K_Y+D_Y$ is $\ep$-lt
and $f$ does not contract components of $D_Y$ with coefficient
$\ge 1-\ep$.
\end{enumerate}
\end{proposition1}

\begin{corollary1}[{\cite[2.18]{Ut}}]
\label{crepant-contraction}
Let $f\colon (Y,D_Y)\to (X,D)$ be a birational contraction, where
$D=f_*D_Y$. Assume that $K_X+D$ is $\QQ$-Cartier. If $K_Y+D_Y$ is
lc (resp. klt) and $f$-(numerically) nonpositive, then $K_X+D$ is
lc (resp. klt).
\end{corollary1}

\begin{example1}
Let $X$ be a normal toric variety and $D$ the reduced toric
boundary on $X$. Then $K_X+D$ is lc. This follows by
\ref{crepant-singul}, from the fact that $K_X+D\sim 0$ and from
the existence of toric resolutions.
\end{example1}

\section{Finite morphisms and singularities of pairs}
\label{finite}
Let $f\colon Y\to X$ be a finite surjective morphism of normal
varieties and $D=\sum d_iD_i$ a subboundary on $X$. We assume that
the ramification divisor is contained in $\Supp D$ (we allow $D$
to have coefficients $=0$). Define a $\QQ$-divisor $B$ on $Y$ by
the condition
\begin{equation}
\label{eqv-vetvl}
K_Y+B=f^*(K_X+D).
\end{equation}
Write $B$ as $B=\sum b_{i,j} B_{i, j}$, where $f(B_{i, j})=D_i$,
and $r_{i, j}$ the ramification index along $B_{i, j}$ (i.e., at
the general point of $B_{i, j}$). By the Hurwitz formula we have
\begin{equation}
\label{eqv-vetvl1_2}
b_{i, j}=1-r_{i, j}(1-d_i).
\end{equation}
Hence $B$ is also subboundary. Note however that $B$ may not be a
boundary even if $D$ is.

\begin{proposition1}[{\cite[\S 2]{Sh}, \cite[20.3]{Ut},
\cite{Ko}}]
\label{finite-plt}
Notation as above. Then $K_X+D$ is lc (resp. plt, klt) if and only
if $K_Y+B$ is lc (resp. is plt, klt).
\end{proposition1}

Propositions \ref{crepant-singul} and \ref{finite-plt} show that
the classes of klt, plt and lc singularities of pairs are very
natural and are closed under birational and finite morphisms. The
implication $\Longrightarrow$ also holds for dlt pairs if $f$ is
\'etale in codimension one \cite{Sz}.

\begin{outline}
Let $g\colon X'\to X$ be a birational morphism. Consider the base
change
\[
\begin{CD}
Y'@>f'>>X'\\ @VhVV @VgVV\\ Y@>f>>X\\
\end{CD}
\]
where $Y'$ is a dominant component of the normalization of
$Y\times_X X'$. As in \eref{for-crepant-singul}, write
\[
K_{X'}+D'= g^*(K_X+D)\qquad\text{and}\qquad K_{Y'}+B'=h^*(K_Y+B),
\]
and similar to \eref{eqv-vetvl} we can write
\[
K_{Y'}+B'={f'}^*(K_{X'}+D'),
\]
where by \eref{eqv-vetvl1_2} the coefficients of $B'$ are
\[
b'_{i, j}=1-r'_{i, j}(1-d_i').
\]
Let $E:=B_{i, j}'$ be an $h$-exceptional divisor and $F:=f'(E)$.
Then this formula can be rewritten as
\begin{equation}
\label{equiv-ramif}
a(E,B)+1=r'_{i,j}(a(F,D)+1), \qquad r'_{i,j}\le \deg f.
\end{equation}
This yields $a(E,B)\ge a(F,D)$ and all the implications
$\Longrightarrow$. The implications $\Longleftarrow$ follow by
\eref{equiv-ramif} and from the (nontrivial) fact that each
exceptional divisor $E$ over $Y$ can be obtained in the way
specified above (see \cite[3.17]{Ko}).
\end{outline}
Note that we have shown more:
\begin{equation}
\label{eqv-finite-new}
1+\mt{discrep}(X,D)\le 1+\mt{discrep}(Y,B) \le (\deg
f)(1+\mt{discrep}(X,D)).
\end{equation}

The following particular case of Proposition~\ref{finite-plt} is
very interesting for applications.

\begin{corollary1}
\label{corollary-finite}
If a morphism $f\colon Y\to X$ is \'etale in codimension one, then
$K_X+D$ is lc (resp. plt, klt) if and only if $K_Y+f^*D$ is lc
(resp. plt, klt).
\end{corollary1}

From Proposition~\ref{finite-plt} it is easy also to obtain the
following

\begin{corollary1}
Let $Y$ be a variety with at worst klt (resp. lc) singularities
and $Y\to X$ a finite surjective morphism. Then $X$ also has at
worst klt (resp. lc) singularities.
\end{corollary1}

In particular, all quotient singularities are klt. However, the
converse is true only in dimension two.

\begin{remark1}[{\cite[3.1]{Reid-canonical}}]
\label{Reid-canonical-ref}
Let $G\subset \GL_n(\CC)$ be a finite subgroup without
quasi-reflections. Then $\CC^n/G$ has only canonical singularities
if and only if for every element $g\in G$ of order $r$ and for any
primitive $r$th root of unity $\ep$, the diagonal form of the
action of $g$ is
\[
g\colon x_i\longrightarrow\ep^{a_i}x_i\quad\text{with}\quad 0\le
a_i<r,\quad \text{and}\quad \sum a_i\ge r.
\]
\end{remark1}

\begin{example1}[{\cite{Catanese}}]
Let $(Y\ni o)$ be a Du Val singularity and $f\colon(Y\ni
o)\to(X\ni P)$ a quotient by an involution. Write
$K_Y=f^*(K_X+\frac12\De)$, where $\De$ is the ramification
divisor. Then $(X,\frac12\De)$ is $(1/2)$-lt. There is an explicit
list of all such involutions and quotients \cite{Catanese}.
\end{example1}

\begin{example1}[{\cite{K}, \cite[7.2]{Mo}}]
Let $X\ni P $ be a germ of a three-dimensional terminal
singularity. By \cite{RYPG}, a general divisor $F\in |{-}K_X|$ has
only Du Val singularities. Then according to Inversion of
Adjunction \ref{inv} (see Example~\ref{YPG}) $K_X+F$ is plt. From
this, for general $S\in |{-}2K_X|$ the divisor $K_X+\frac12S$ is
also plt. Consider the double cover $f\colon Y\to X$ with
ramification divisor $S$. Then $Y$ has only klt singularities and
$K_Y\sim 0$. Hence the singularities of $Y$ are canonical of index
one.
\end{example1}

The existence of a good divisor $S\in |{-}2K_X|$ in the global
case or for extremal contractions $X\to Z$ is a much more
difficult problem. For example, it is sufficient for the existence
of three-dimensional flips \cite{K}.

\section{Log canonical covers}
\label{l-c-cover}
The following construction is well known (see e.g.,
\cite[2.4]{Sh}, \cite{Kawamata-models}, \cite[8.5]{K},
\cite{Ko1}). Let $X$ be a normal variety and $D=\sum d_iD_i$ a
boundary such that $m(K_X+D)\sim 0$. We take $m$ to be the least
positive integer satisfying this condition. Such $m$ is called the
\textit{index}\index{index} of $K_X+D$. Assume that $K_X+D$ is lc
and $d_i\in\{1-1/k\mid k\in\NN\cup\{\infty\}\}$ for all $i$ (i.e.
all the $d_i$ are standard, see \ref{def-SM}). Then the natural
map $\OOO_X(-m(K_X+D))\to\OOO_X$ defines an $\OOO_X$-algebra
structure on $\sum_{i=0}^{m-1}
\OOO_X\left(\down{-iK_X-iD}\right)$. Put
\[
Y:=\mt{Spec}\left(\sum_{i=0}^{m-1}
\OOO_X\left(\down{-iK_X-iD}\right)\right)
\]
and $\var\colon Y\to X$ the projection. Then $Y$ is irreducible,
$\var$ is a cyclic Galois $\ZZ_m$-cover. Put
$B:=\var^*(\down{D})$. Then the ramification divisor (i.e.
codimension one ramification locus) of $\var$ is
$\Supp(D-\down{D})$. Further, the ramification index along $D_i$
is $r_i$, where $d_i=1-1/r_i$. Therefore,
\[
\var^*(K_X+D)=K_Y+B\sim 0.
\]
By \ref{finite-plt}, $K_Y+B$ is lc. Moreover, $K_X+D$ is plt
(resp. klt) if and only if $K_Y+B$ is plt (resp. klt).

\begin{exercise1}
Let $(X\ni P)$ be a Du Val singularity of type $D_n$ (given by the
equation $x^2+y^2z+z^{n-1}=0$) and $H$ be a general hyperplane
section. Show that the double cover ramified along $H$ is a lc
singularity of index one. Write down the equation of this
singularity.
\end{exercise1}

Finally, we present some Bertini's type results.
\begin{proposition}
\label{Bertini}
Let $X$ be a normal variety and $D$ a $\QQ$-divisor on $X$. Let
$\mathcal{H}$ be a base point free linear system and $H\in
\mathcal{H}$ a general member.
\begin{enumerate}
\item
Then $K_X+D$ is lc if and only if $K_X+D+H$ is lc.
\item
Assume additionally that $\down D=0$. Then $K_X+D$ is klt if and
only if $K_X+D+H$ is plt.
\end{enumerate}
\end{proposition}
\begin{proof}
It is sufficient to show only the implications $\Longrightarrow$.
Let $f\colon Y\to X$ be a log resolution of $(X,D)$,
$E_1,\dots,E_r$ exceptional divisors and $D_Y$, $H_Y$ proper
transforms of $D$ and $H$, respectively. By Bertini's theorem,
$D_Y+H_Y$ is a simple normal crossing divisor, so $f$ is also a
log resolution of $(X,D+H)$. We can choose $H\in\mathcal{H}$ so
that $H$ does not contain $f(E_1),\dots,f(E_r)$. Thus we have
$a(E_i,D+H)= a(E_i,D)\ge -1$. This implies the first part of our
proposition. In the second part we can use remark in
\ref{remark-def-klt}.
\end{proof}

\begin{proposition}[{\cite[1.13]{Reid-canonical}}, {\cite[5.17]{KM}}]
\label{Bertini1}
Let $X$ be a normal variety and $D$ a $\QQ$-divisor on $X$. Let
$\mathcal{H}$ be a base point free linear system and $H\in
\mathcal{H}$ a general member. Assume that $K_X+D$ is klt (resp.\
plt, lc, canonical, terminal). Then $K_H+D|_H$ is klt (resp.\
plt, lc, canonical, terminal).
\end{proposition}

\chapter{Inversion of adjunction}
\label{sect-2}
\section{Two-dimensional toric singularities
and log canonical singularities with a reduced boundary}
\subsection{}
\label{frac}
If the cyclic group $\cyc{m}$ acts linearly on $\CC^n$ by
\[
x_1\to \ep^{a_1} x_1, \qquad x_2\to \ep^{a_2} x_2,\dots,x_n\to
\ep^{a_n} x_n,
\]
where $\ep$ is a chosen primitive root of degree $m$ of unity, we
call the integers $a_1,\dots,a_n$ the
\textit{weights}\index{weights} of the action. In this case, the
quotient is denoted by $\CC^n/\cyc{m}(a_1,\dots,a_n)$.
\index{$\CC^n/\cyc{m}(a_1,\dots,a_n)$} It is clear that the
weights are defined modulo $m$ and also depend on the choice of
the primitive root $\ep$.
\par
Let $(Z, Q)$ be a two-dimensional quotient singularity
$\CC^2/\cyc{m}(1, q)$, where $\gcd(q,m)=1$ (in particular, this
means that $\cyc{m}$ acts on $\CC^2$ freely in codimension one).
Then this singularity is toric, hence it is klt. The minimal
resolution is obtained as a sequence of \textit{weighted blowups}
(see \ref{w-blow-up}). The dual graph is a chain
\[
\begin{array}{ccccccccc}
\wcir{-a_1}&\lin&\wcir{-a_2}&\lin&\cdots&\lin&
\wcir{-a_{r-1}}&\lin& \wcir{-a_r},
\end{array}
\qquad\text{with}\quad a_i\ge 2,
\]
where the sequence $a_1, a_2,\dots,a_r$ is obtained from the
continued fraction decomposition of $m/q$ (see \cite{Hi} or
\cite{Brieskorn}):
\begin{equation}
\label{frac_chain}
\frac{m}{q}=a_1-\cfrac{1}{a_2-\cfrac{1}{\cdots\cfrac{1}{a_r}}}.
\end{equation}

Now we give the classification of two-dimensional log canonical
singularities with nonempty reduced boundary, following Kawamata
\cite{K}. Note that this is much easier than the classification of
all two-dimensional log canonical singularities.

\begin{theorem1}[{\cite[9.6]{K}, \cite[ch. 3]{Ut}}]
\label{kawamata_1}
Let $X\ni P$ be an analytic germ of a two-dimensional normal
singularity and $X\supset C$ a (possibly reducible) reduced curve.
Assume that $K_X+C$ is plt. Then
\[
(X, C)\simeq(\CC^2, \{x=0 \})/\cyc{m}(1, a),\quad\text{with}\quad
\gcd(a, m)=1.
\]
In particular, $C$ is irreducible and smooth. In this case,
$K_X+C$ has index $m$ and the graph of the minimal resolution of
$(X\supset C\ni P)$ is of type
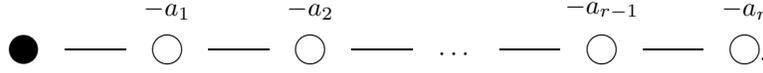
\[
\begin{array}{ccccccccccc}
\bcir{}&\lin& \wcir{-a_1}&\lin&\wcir{-a_2}&\lin&
\cdots&\lin&\wcir{-a_{r-1}}&\lin&\wcir{-a_r}.\\
\end{array}
\]
where the black vertex $\bci$ corresponds to the proper transform
of $C$, and the white ones $\wci$ correspond to exceptional
divisors. The numbers attached to white vertices are
self-intersection numbers.
\end{theorem1}

\begin{outline}
Let $m$ be the index of $C$ (i.e., $mC\sim 0$) and $\psi\colon
X'\to X$ the corresponding cyclic $m$-cover. Then
$C':=\psi^{-1}(C)_{\red}$ is a Cartier divisor and $K_{X'}+C'$ is
plt. This gives that $X'$ and $C'$ are smooth and $X'\simeq\CC^2$
up to analytic isomorphism.
\end{outline}

\begin{theorem1}[{\cite[9.6]{K}, \cite[ch. 3]{Ut}}]
\label{kawamata_2}
Let $X\ni P$ be an analytic germ of a two-dimensional normal
singularity and $X\supset C$ a (possibly reducible) reduced curve.
Assume that $K_X+C$ is lc but not plt. Then just one of the
following two possibilities holds:
\begin{enumerate}
\item
$C$ has two smooth components,
\[
(X, C)\simeq(\CC^2, \{xy=0 \})/\cyc{m}(1, a),
\quad\text{with}\quad\gcd(a, m)=1.
\]
The index of $K_X+C$ is
equal to $1$ (i.e., $K_X+C\sim 0$). The graph of the minimal
resolution of $(X\supset C\ni P)$ is of the form
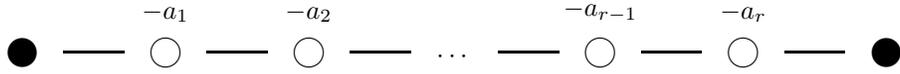
\[
\begin{array}{ccccccccccccc}
\bcir{}&\lin&\wcir{-a_1}&\lin&\wcir{-a_2}&\lin&\cdots
&\lin&\wcir{-a_{r-1}}&\lin&\wcir{-a_r}& \lin&\bcir{}.\\
\end{array}
\]
In this case, $K_X+C$ is not dlt for $m>1$ and dlt for $m=1$.
\item
The curve $C$ is smooth and irreducible,
\[
(X, C)\simeq(\CC^2, \{xy=0 \})/\DD_m,
\]
where $\DD_m\subset \GL_2(\CC)$ is a subgroup of dihedral type
without reflections (see \cite{Brieskorn} for precise description
of $\DD_m$). In this case, $K_X+C$ is not dlt and is of index two
(i.e., $2(K_X+C)\sim 0$), the log canonical cover is the
singularity from (i) and the graph of the minimal resolution of
$(X\supset C\ni P)$ is of type
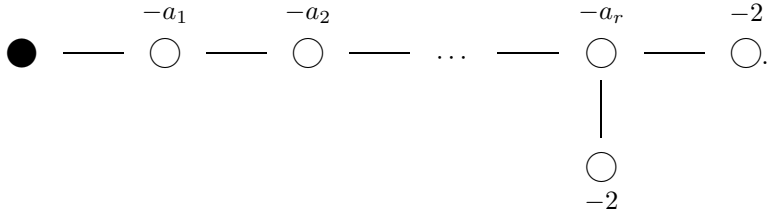
\[
\begin{array}{cccccccccccc}
\bcir{}&\lin&\wcir{-a_1}&\lin&\wcir{-a_2}&\lin&\cdots
&\lin&\wcir{-a_r}&\lin&\wcir{-2}.&\\&&&&&&&&\vlin&&&\\
&&&&&&&&\wci&&&\\&&&&&&&&{-2}&&&\\
\end{array}
\]
The degenerate case $r=1$ is included here (then $\DD_m$ is a
cyclic group).
\end{enumerate}
\end{theorem1}

\begin{corollary1}
Let $(X,D)$ be a log variety. Assume that $K_X+D$ is lc and
$W\subset X$ an irreducible subvariety of codimension two. Assume
that $W\subset\down{D}$. Then near a general point $w\in W$ there
is an analytic isomorphism between $(X,\down{D},W)$ and the
product of a surface singularity from \ref{kawamata_1} or (i)-(ii)
of \ref{kawamata_2} by $\CC^{\dim X-2}$.
\end{corollary1}

\begin{exercise1}[cf. {\ref{ex-Du Val}}]
\label{new}
Assume that in the conditions of the theorem above $C$ is a
Cartier divisor. Show without using the theorem that then
$(X,C)\simeq(\CC^2,\CC^1)$ or $X$ is a Du Val point of type $A_n$,
and $C$ is its general hyperplane section.
\end{exercise1}

\begin{exercise1}
Express in the form $\CC^2/\cyc{m}(1,q)$ the singularity with the
minimal resolution
\[
\begin{array}{ccccc}
\wcir{-2}&\lin&\wcir{-a}&\lin& \wcir{-2}\\&&\vlin&&\\&&\bci&&\\
\end{array}
\]
\end{exercise1}

\begin{example1}
\label{ex-Du Val}
Let $(Z,Q)$ be a Du Val singularity of type $A_n$ given by the
equation $x^2+y^2+z^{n+1}$ and $C$ the hyperplane section given by
$z=0$. Then $(Z,C)$ is a lc pair as in (i) of
Theorem~\ref{kawamata_2}. Similarly, for the case (ii) of
Theorem~\ref{kawamata_2} we can take $(Z,Q)$ of type $D_n$ given
by the equation $x^2+y^2z+z^{n-1}$, $n\ge 4$ and $C$ as
$\{z=0\}_{\red}$.
\end{example1}

\section{Adjunction}
\begin{example1}
\label{exam-difff}
Let $X=X_n\subset\PP^{n+1}$ be a two-dimensional projective cone
over a rational normal curve $C_n\subset\PP^n$ and $L\subset X$
its generator. The group of classes of Weil divisors modulo linear
equivalence is generated by the class of $L$:
$\Weillin(X)\simeq\ZZ\cdot L$ and $nL$ is the class of the
hyperplane section of $X$. Thus we have $L|_L=\frac1nP$, where $P$
is the class of a point on $L\simeq\PP^1$. It is also easy to
compute that $K_X\sim-(n+2)L$. This yields
\[
(K_X+L)|_L{-}K_L=-(n+1)L|_L+2P=(1-1/n)P.
\]
This is one instance where the \emph{standard coefficients} (see
\ref{def-SM}) arise naturally.
\end{example1}

This example shows that adjunction formula in its usual form fails
for the case of Weil divisors. This phenomenon was first observed
by M.~Reid and is called also \textit{subadjunction}. Shokurov
\cite[\S 3]{Sh} introduced the notion of different for the
difference $(K_X+L)|_L{-}K_L$ (see also \cite[5-1-9]{KMM},
\cite[ch. 16]{Ut}). The corresponding ideal sheaf sometimes is
called the \textit{conductor} ideal.
\par

The following construction is a codimension two construction, i.e.
the variety $X$ may always be replaced with any open subset
$X\setminus Z$, where $\codim_X Z\ge 3$.

\begin{proposition-definition1}
\label{def-deff}
Let $X$ be a normal variety and $S\subset X$ a reduced subscheme
of pure codimension one. For simplicity we assume that $K_X+S$ is
lc in codimension two. Then by Theorem~\ref{kawamata_1} and
Theorem~\ref{kawamata_2}, $S$ has only normal crossings in
codimension one. In particular, the scheme $S$ is Gorenstein in
codimension one. Then there exists naturally defined an effective
$\QQ$-Weil divisor $\Diff S(0)$, called the \emph{different},
\index{different} such that
\[
(K_X+S)|_S=K_S+\Diff S(0).
\]
Now let $B$ be a $\QQ$-divisor, which is $\QQ$-Cartier in
codimension two. Then the different for $K_X+S+B$ is defined by
the formula
\[
(K_X+S+B)|_S=K_S+\Diff S(B). \index{$\Diff S(B)$}
\]
In particular, if $B$ is a boundary and $K_X+S+B$ is lc in
codimension two, then by \ref{kawamata_1} and \ref{kawamata_2},
$B$ is $\QQ$-Cartier in codimension two. Moreover, none of the
components of $\Diff S(B)$ are contained in the singular locus of
$S$.
\end{proposition-definition1}

\begin{example1}
Let $Q\subset\PP^4$ be a quadratic cone over $xy=zt$ and $S\subset
Q$ a plane. Then $S|_S=0$ modulo codimension two subsets and
$(K_Q+S)|_S=K_S$. Therefore $\Diff S(0)=0$. This shows that
codimension three singularities are not essential for
\ref{def-deff}.
\end{example1}

The following theorem allows us to compute coefficients of the
different and shows that computations in Example~\ref{exam-difff}
are very general.

\begin{theorem1}[{\cite[3.9]{Sh}, \cite[16.6]{Ut}}]
\label{diff}
In the conditions of Theorem~\ref{kawamata_1} and
Theorem~\ref{kawamata_2} for the different $\Diff C(0)$ at $P$ we
have
\begin{enumerate}
\item
If $K_X+C$ is plt, then $\Diff C(0)=(1-1/m)P$, where $m$ is the
index of $K_X+C$ (see Theorem~\ref{kawamata_1}).
\item
If $(X\supset C\ni P)$ is as in (i) of Theorem~\ref{kawamata_2},
then $\Diff C(0)=0$.
\item
If $(X\supset C\ni P)$ is as in (ii) of Theorem~\ref{kawamata_2},
then $\Diff C(0)=P$.
\end{enumerate}
\end{theorem1}

\subsection{Notation}
\label{def-SM}
Put
\[
\Msm:=\{1-1/m\mid m\in\NN\cup\{\infty\}\}. \index{$\Msm$}
\]
We distinguish this set because it naturally appears in the
Adjunction Formula. Latter we will see that the class of
boundaries with coefficients $\in\Msm$ is closed under finite
Galois morphisms (see~\ref{finite}) and Adjunction Formula
(Corollary~\ref{coeff-Diff-1}, cf. \cite{Ko1}). We say that the
boundary $D=\sum d_iD_i$ has \textit{standard
coefficients},\index{standard coefficient} if $d_i\in\Msm$ for all
$i$. Unfortunately the property $\in\Msm$ is not closed under
crepant birational transformations (see \ref{crepant-singul}) to
avoid this difficulty Shokurov considered the class of boundaries
with coefficients from the set
\[
\Mm:=\Msm\cup [6/7,1]. \index{$\Mm$}
\]

The following result is very important for applications and is
called \textit{Inversion of Adjunction}. \index{Inversion of
Adjunction} \index{Adjunction}

\begin{theorem}[{\cite[3.3, 3.12, 5.13]{Sh}},
{\cite[17.6]{Ut}}]
\label{inv}
Notation as in \ref{def-deff}. Let $B$ be an effective
$\QQ$-divisor on $X$ such that $S$ and $B$ have no common
components and $\down{B}=0$. Assume that $K_X+S+B$ is
$\QQ$-Cartier. Then $K_X+S+B$ is plt near $S$ if and only if $S$
is normal and $K_S+\Diff S(B)$ is klt.
\end{theorem}

\begin{corollary1}[{\cite[17.7]{Ut}}]
\label{inv1}
Let $X$ be a normal variety, $S$ an irreducible divisor and let
$B$, $B'$ be effective $\QQ$-divisors such that $S$ and $B+B'$
have no common components and $\down{B}=0$. Assume that $K_X+S+B$
and $B'$ are $\QQ$-Cartier and $K_X+S+B$ is plt. Then $K_X+S+B+B'$
is lc near $S$ if and only if so is $K_S+\Diff S(B+B')$.
\end{corollary1}

\begin{corollary1}[{\cite[3.10]{Sh}}]
\label{coeff-Diff}
Let $X$ be a normal variety, $S$ a reduced Weil divisor on $X$ and
$B=\sum b_iB_i$ a boundary on $X$ such that $S$ and $B$ have no
common components. Assume that $K_X+S+B$ is plt. Then $\Diff S(B)$
has the form
\[
\Diff S(B)=\sum_{P_i}\left(\frac{m_i-1}{m_i}+\sum_j
\frac{b_jn_{i,j}}{m_i}\right)P_i,
\]
where each $P_i$ is a prime divisor on $S$, $m_i$ is the index of
$S$ at a general point of $P_i$, and $n_{i,j}\in\NN$. Moreover,
assume that $B$ has only standard coefficients. Then so has $\Diff
S(B)$. More precisely, if $B=\sum (1-1/r_i)B_i$ and $P$ is a prime
divisor on $S$, then $P$ is contained in at most one component,
say $B_i$, of $B$ and the coefficient of $\Diff S(B)$ along $P$ is
equal to $1-\frac1{m_ir_i}$.
\end{corollary1}
\begin{corollary1}
\label{coeff-Diff-1}
Notation as in \ref{coeff-Diff}. Then
\[
\begin{array}{lll}
 B\in\Msm\quad&\Longrightarrow\quad&
 \Diff S(B)\in\Msm,\\
 B\in\Mm\quad&\Longrightarrow\quad&
 \Diff S(B)\in\Mm.
\end{array}
\]
\end{corollary1}

The example below shows that Inversion of Adjunction fails in the
case of noneffective divisors.
\begin{example1}
Consider the following smooth curves on $\CC^2$: $C:=\{x=0\}$,
$B_1:=\{y=x^2\}$, $B_2:=\{y=2x^2\}$, $B_3:=\{y=x\}$ and consider
the subboundary $B:=bB_1+bB_2+(\frac32-3b)B_3$, where
$\frac12<b<1$. Then $\Diff C(B)=(\frac32-b)(\mt{pt})$ because $C$
intersects transversally $B_1$, $B_2$, $B_3$. Hence $K_C+\Diff
C(B)$ is klt. On the other hand $K_{\CC^2}+C+B$ is not lc. Indeed,
a log resolution of $(\CC^2,C+B)$ can be obtained by two blowing
ups:
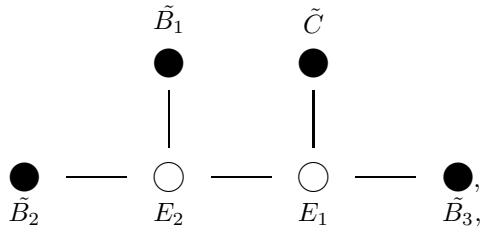
\[
\begin{array}{cccccccc}
&&\bcir{\tilde{B_1}} &&\bcir{\tilde{C}}&&\\ &&\vlin&&\vlin&&\\
\bci&\lin& \wci&\lin&\wci& \lin&\bci,\\ \tilde {B_2}&&
E_2&&E_1&&\tilde{B_3},\\
\end{array}
\]
where $E_1$ is a $-2$-curve and $E_2$ is a $-1$-curve. It is easy
to compute
\[
a(E_1,C+B)=-\frac32+b>-1,\qquad a(E_2,C+B)=-\frac12-b<-1.
\]
Therefore $K_{\CC^2}+C+B$ is not lc at the origin.
\end{example1}

For dlt singularities there are only weaker results, which use
generalizations of the definition of dlt singularities to the case
of nonnormal varieties (cf. {\cite[3.2.3, 3.6, 3.8]{Sh}},
{\cite[17.5, 16.9]{Ut}}):

\begin{proposition1}[\cite{Sz}]
\label{s2}
Let $(X,S+B)$ be a log variety, where $S$ is reduced, $\down{B}=0$
and $S$, $B$ have no common components. Assume that $K_X+S+B$ is
dlt. Then $K_S+\Diff S(B)$ is generalized divisorial log terminal.
\end{proposition1}

\begin{example1}
\label{YPG}
Let $(X\ni P)$ be a germ of three-dimensional terminal
singularity. Then by \cite{RYPG} a general divisor $F\in |{-}K_X|$
has only Du Val singularity at $P$. Hence by Theorem~\ref{inv},
$K_X+F$ is plt (and even canonical, because $K_X+F$ is Cartier).
\end{example1}

\begin{exercise1}
Let $H$ be a general hyperplane section of the canonical quotient
singularity $X:=\CC^3/\cyc{3}(1,1,1)$. Show that $K_X+H$ is not
plt.
\end{exercise1}

\begin{example1}[cf. {\cite[Sect. 1]{Pagoda}}]
\label{YPG_G}
Let $(X\ni P)$ be a normal three-dimensional $\QQ$-Gorenstein
singularity and $H\ni P$ a hyperplane section. Assume that the
singularity $(H\ni P)$ is Du Val. Then by Inversion of Adjunction,
$(X\ni P)$ is canonical. Moreover, if $(X\ni P)$ is isolated, then
it is terminal.
\end{example1}

\begin{exercise1}
Let $(X\ni o, D)$ be a germ of a normal singularity. Assume that
$K_X+D$ is lc and each component of $\down{D}$ is $\QQ$-Cartier.
Prove that $\down{D}$ has at most $\dim X$ components.
\end{exercise1}

We also have a more general results:

\begin{example1}
Let $(X\ni o, D=\sum d_iD_i)$ be a germ of a normal singularity of
dimension $\le 3$. Assume that $K_X+D$ is lc and each component of
$D$ is $\QQ$-Cartier. Then $\sum d_i\le\dim X$. Indeed, by taking
cyclic covers \'etale in codimension one we obtain $(X'\ni o',
D'=\sum d_i'D_i')$ such that $K_{X'}+D'$ is lc and each component
of $D'$ is Cartier. Obviously, $\sum d_i'\ge\sum d_i$. It is known
that in dimension $\le 3$ there exists a divisor $E$ over $X'\ni
o'$ such that $a(E,0)\le\dim X-1$ (see Kawamata's appendix to
\cite{Sh} and \cite{M}). Then $-1\le a(E,D')\le a(E,0)-\sum d_i'$.
This yields $\sum d_i\le\dim X$. Moreover, if the equality holds,
then $X'$ is smooth. In this case, $X'\to X$ gives the universal
cover of the smooth locus of $X$. Therefore $X'\to X$ is a
quotient by a finite group, say $G$, which acts freely in
codimension one. Then we have also $\sum d_i'=\sum d_i$. Hence $G$
does not permute components of $D'$. So $G$ must be abelian. This
shows that the equality $\sum d_i=\dim X$ implies that $X$ is
analytically isomorphic to a toric singularity and $\down{D}$ is
contained in the toric boundary.
\end{example1}
Actually, the above result can be proved in any dimension without
using \cite{M}:
\begin{theorem}[{\cite[18.22]{Ut}}, {\cite{A1}}]
\label{teor-toric}
Let $(X\ni o,D=\sum d_iD_i)$ be a germ of a log variety such that
$K_X+D$ is lc. Assume that all the $D_i$ are $\QQ$-Cartier at $o$.
Then $\sum d_i\le \dim X$. Moreover, the equality holds only if
$X\ni o$ is a cyclic quotient singularity.
\end{theorem}

Let $(X/Z\ni o,D)$ be a log variety. Then $(X/Z\ni o,\down{D})$ is
said to be a \textit{toric pair}\index{toric pair} if there are
analytic isomorphisms $\pi\colon X\to X^T$, $Z\to Z^T$ and the
commutative diagram
\[
\begin{CD}
X@>\pi>>X^T\\ @VVV @VVV\\ Z@>>>Z^T\\
\end{CD}
\]
where $X^T\to Z^T$ is an algebraic toric contraction and
$\down{\pi(D)}$ is the toric boundary (i.e., $\down{\pi(D)}$ is
contained in the set $X^T\setminus \{\text{open orbit}\}$).

Shokurov proposed the following generalization of
\ref{teor-toric}.
\begin{conjecture1}[{\cite{Sh1}}]
\label{conjecture-toric}
Let $(X/Z\ni o,D=\sum d_iD_i)$ be a log variety such that $K_X+D$
is lc and $-(K_X+D)$ is nef over $Z$. Then\footnote{Shokurov
pointed out that the stronger version of inequality
\eqref{eq-conj-toric} should be $\Weil(X)/ \Weil_0 \ge\sum
d_i-\dim X$, where $\Weil_0\subset \Weil(X)$ is the subgroup of
all numerically trivial over $Z$ (and $\QQ$-Cartier) divisors.}
\begin{equation}
\label{eq-conj-toric}
\mt{rk}\Weilalg(X)\ge\sum d_i-\dim X.
\end{equation}
If $X$ is $\QQ$-factorial, then
\begin{equation}
\label{eq-conj-toric-1}
\rho(X/Z)\ge\sum d_i-\dim X.
\end{equation}
Moreover, equalities hold only if $(X/Z\ni o,\down{D})$ is a toric
pair.
\end{conjecture1}

In the case when $Z$ is a point and $\rho(X)=1$ the inequality
\eref{eq-conj-toric-1} was proved in \cite[18.24]{Ut}, see also
\cite{A1}. Shokurov \cite{Sh1} proved this conjecture in dimension
two; see theorems~\ref{th-toric} and \ref{th-toric-1}.

Note that inequality \eref{eq-conj-toric-1} is stronger than
\eref{eq-conj-toric}:
\begin{proposition1}
Notation as in \ref{conjecture-toric}. Assume that the pair
$(X,D)$ has at least one minimal log terminal ($\QQ$-factorial)
modification $f\colon (\tilde X,\tilde D=\sum \tilde d_i\tilde
D_i)\to (X,D)$ (see \ref{minimal_lt}). Then
\begin{multline*}
\mt{rk}\Weilalg(X)-\sum d_i+\dim X\ge\\ \mt{rk}\Weilalg(\tilde
X)-\sum\tilde d_i+\dim \tilde X\ge \rho(\tilde X/Z)-\sum \tilde
d_i+\dim \tilde X.
\end{multline*}
\end{proposition1}
\begin{proof}
Let us prove, for example, the first inequality. Write $\tilde
D=\sum d_i B_i+\sum_{j=1}^r E_j$, where each $B_i$ is the proper
transform of $D_i$ and $\sum_{j=1}^r E_j$ is the (reduced)
exceptional divisor. Then $\sum \tilde d_i=r+\sum d_i$. From the
exact sequence
\[
\bigoplus_{i=1}^r \ZZ\cdot E_i\longrightarrow \Weilalg(\tilde
X)\longrightarrow \Weilalg\left(\tilde X\setminus \sum E_i\right)
\longrightarrow 0
\]
(cf. \cite[Ch. II, 6.5]{Ha}) we have
\[
\mt{rk}\Weilalg(\tilde X)\le \Weilalg\left(\tilde X\setminus \sum
E_i\right)+r= \mt{rk}\Weilalg(X)+r.
\]
\end{proof}

\begin{example1}
Let $X=Z$ be the hypersurface singularity given in $\CC^4$ by the
equation $xy=zt$. Consider four planes $D_1:=\{x=z=0\}$,
$D_2:=\{x=z=0\}$, $D_3:=\{y=z=0\}$, $D_4:=\{y=t=0\}$. Let $D:=\sum
D_i$. It is easy to check that $K_X+D$ is lc. The group
$\Weilalg(X)$ is generated by $D_1$ and $D_2$. However,
$D_1+D_2\sim 0$ (because it is Cartier). Thus
$\mt{rk}\Weilalg(X)=1$ and $\sum d_i=4$. We have equality in
\eqref{eq-conj-toric} and the pair $(X,D)$ is toric.
\end{example1}

\section{Connectedness Lemma}
The most essential part of the proof of Theorem~\ref{inv} is the
following result which was proved firstly by Shokurov
\cite[5.7]{Sh} in dimension two and latter by Koll\'ar
\cite[17.4]{Ut}, \cite[7.4]{Ko} in arbitrary dimension.

\begin{theorem1}[Connectedness Lemma]
\label{connect} \index{Connectedness Lemma}
Let $f\colon X\to Z$ be a contraction and $D=\sum d_iD_i$ an
effective $\QQ$-divisor on $X$ such that $K_X+D$ is $\QQ$-Cartier.
Assume that $-(K_X+D)$ is $f$-nef and $f$-big. Let
\[
h\colon Y\stackrel{g}{\longrightarrow} X
\stackrel{f}{\longrightarrow} Z
\]
be a log resolution. Write
\[
K_Y=g^*(K_X+D)+E^{(+)}-E^{(-)},
\]
where $E^{(-)}\ge 0$ and the coefficients of $E^{(+)}$ are $>-1$,
and the coefficients of $E^{(-)}$ are $\ge 1$. Then $\Supp
E^{(-)}$ is connected in a neighborhood of any fiber of $h$.
\end{theorem1}
Note that in the case when $f$ is birational, the big condition
holds automatically.

\begin{proof}
We have
\[
\up{E^{(+)}}-\down{E^{(-)}}=
K_X-g^*(K_X+D)+\fr{-E^{(+)}}+\fr{E^{(-)}}.
\]
From this by Kawamata-Viehweg Vanishing Theorem \cite[1-2-3]{KMM},
\[
R^1f_*\OOO_Y\left(\up{E^{(+)}}-\down{E^{(-)}}\right)=0.
\]
Applying $f_*$ to an exact sequence
\begin{multline*}
0\longrightarrow\OOO_Y\left(\up{E^{(+)}}-\down{E^{(-)}}\right)
\longrightarrow\OOO_Y\left(\up{E^{(+)}}\right)\\
\longrightarrow\OOO_{\down{E^{(-)}}}\left(\up{E^{(+)}}\right)
\longrightarrow 0
\end{multline*}
we get the surjectivity of the map
\[
h_*\OOO_Y\left(\up{E^{(+)}}\right) \longrightarrow
h_*\OOO_{\down{E^{(-)}}}\left(\up{E^{(+)}}\right).
\]
Let $E_i$ be a component $\up{E^{(+)}}$. Then either $E_i$ is
$g$-exceptional or $E_i$ is the proper transform of some $D_i$
whose coefficient $d_i<1$. Thus $\up{E^{(+)}}$ is $g$-exceptional
and
\[
h_*\OOO_Y\left(\up{E^{(+)}}\right)= f_*
\left(\OOO_X\left(g_*\left(\up{E^{(+)}}\right)\right)\right)=
\OOO_Z.
\]
Assume that in a neighborhood of some fiber $ h^{-1}(z)$, $z\in Z$
the set $\down{E^{(-)}}$ has two connected components $F_1$ and
$F_2$. Then
\[
h_*\OOO_{\down{E^{(-)}}}\left(\up{E^{(+)}}\right)_{(z)}=
h_*\OOO_{F_1}\left(\up{E^{(+)}}\right)_{(z)}+
h_*\OOO_{F_2}\left(\up{E^{(+)}}\right)_{(z)},
\]
and both terms do not vanish. Hence
$h_*\OOO_{\down{E^{(-)}}}\left(\up{E^{(+)}}\right)_{(z)}$ cannot
be a quotient of the cyclic module $\OOO_{z, Z}\simeq
h_*\OOO_Y\left(\up{E^{(+)}}\right)_{(z)}$.
\end{proof}

\begin{definition1}[{\cite[3.14]{Sh}}]
Let $X$ be a normal variety and $D=\sum d_iD_i$ a $\QQ$-divisor on
$X$ such that $K_X+D$ is $\QQ$-Cartier. We say that a subvariety
$W\subset X$ is a \textit{center of log canonical
singularities},\index{center of log canonical singularities} if
there exists a birational contraction $f\colon Y\to X$ and a
divisor $E$ (not necessarily $f$-exceptional) with discrepancy
$a(E, D, X)\le-1$ such that $f(E)=W$. The union of all centers of
lc singularities is called \textit{the locus of log canonical
singularities}\index{locus of log canonical singularities} of
$(X,D)$ and is denoted by $\LCS(X,D)$.\index{$\LCS(X,D)$}
\end{definition1}

\begin{corollary1}
\label{connect1}
Notation as in Theorem~\ref{connect}. Then the set $\LCS(X,D)$ is
connected in a neighborhood of any fiber of $f$.
\end{corollary1}

\chapter{Log terminal modifications}
\label{sect-3}
\section{Log terminal modifications}
Many results of this chapter holds in arbitrary dimension modulo
log MMP.
\begin{proposition-definition}[cf. {\cite[9.1]{Sh}};
see also {\cite[6.16]{Ut}}]
\label{blow-up}
Let $(X,D)$ be a log variety of dimension $\le 3$. Assume that
$K_X+D$ is lc. Then there exists a \emph{log terminal
modification} \index{log terminal modification} of $(X,D)$; that
is, a birational contraction $g\colon X'\to X$ and a boundary $D'$
on $X'$ such that
\begin{enumerate}
\item
$K_{X'}+D'\equiv g^*(K_X+D)$;
\item
$K_{X'}+D'$ is dlt;
\item
$X'$ is $\QQ$-factorial.
\end{enumerate}
Moreover, if $\dim X=2$, it is possible to choose $X'$ smooth.
\end{proposition-definition}

\begin{proof}
Consider a log resolution $h\colon Y\to X$. We have
\begin{equation}
\label{ooo}
K_Y+D_Y=h^*(K_X+D)+E^{(+)}-E^{(-)},
\end{equation}
where $D_Y$ is the proper transform $D$ on $Y$ and $E^{(+)}$,
$E^{(-)}$ are effective exceptional $\QQ$-divisors without common
components. Then $D_Y+E^{(-)}$ is a boundary and $K_Y+D_Y+E^{(-)}$
is dlt. Apply log MMP to $(Y,D_Y+E^{(-)})$ over $X$. We get a
birational contraction $g\colon X'\to X$ from a normal
$\QQ$-factorial variety $X'$. Denote by $D'$ the proper transform
of $D_Y+E^{(-)}$ on $X'$. Then $K_{X'}+D'$ is dlt and $g$-nef. It
is also obvious that $g_*D'=D$. We prove (i). Since the inverse to
the birational map $h\colon Y\bir X'$ does not contract divisors,
\[\begin{array}{l}
K_{X'}+D'=h_*\left(K_Y+D_Y+E^{(-)}\right)=
h_*\left(f^*(K_X+D)+E^{(+)}\right)=\\ \qquad
g^*(K_X+D)+h_*E^{(+)}.
\end{array}
\]
On the other hand, by numerical properties of contractions (see
e.g., \cite[1.1]{Sh}) in the last formula all the coefficients of
$h_*E^{(+)}$ should be nonpositive, i.e., all of them are equal to
zero.
\par
Finally, we consider the case $\dim X=2$. If $E^{(+)}\ne 0$, then
$(E^{(+)})^2<0$. From this $E^{(+)}\cdot E<0$ for some $E$. Then
$E^2<0$ and by \eref{ooo}, $K_Y\cdot E<0$ and
$(K_Y+D_Y+E^{(-)})\cdot E<0$. Hence $E$ is a $-1$-curve and steps
of log MMP over $X$ are contractions of such curves. Continuing
the process, we obtain a smooth surface $X'$. This proves the
statement.
\end{proof}

\begin{proposition}[{\cite[3.1]{Sh2}}, {\cite[21.6.1]{Ut}}]
\label{per-gen}
Notation as in Theorem~\ref{blow-up}. Let $ h\colon Y\to X$ be any
log resolution. Consider a set $\EEE=\{E_i\}$ of exceptional
divisors on $Y$ such that
\begin{enumerate}
\itemm{a}
if $a(E_i, D)=-1$, then $E_i\in\EEE$;
\itemm{b}
if $E_i\in\EEE$, then $a(E_i, D)\le 0$.
\end{enumerate}
Then there exists a blowup $g\colon X'\to X$ and a boundary $D'$
on $X'$ such that (i), (ii) and (iii) of \ref{blow-up} holds and
moreover,
\begin{enumerate}
\item[\textrm{(iv)}]
the exceptional divisors of $g$ are exactly the elements of $\EEE$
(i.e., they give the same discrete valuations of the field
$\KKK(X)$).
\end{enumerate}
\end{proposition}
\begin{proof}
Take a sufficiently small $\ep>0$ and put
\[
d_i=\left\{
\begin{array}{ll}
-a(E_i,D)\quad&\text{if}\quad E_i\in\EEE,\\ \max \{-a(E_i,D)+\ep,
0\}\quad&\text{otherwise.}\\
\end{array}
\right.
\]
Then
\begin{equation}
\label{eq1}
K_Y+D_Y+\sum d_iE_i\equiv
h^*(K_X+D)+\sum_{E_j\notin\EEE}(d_j+a(E_j, D))E_j.
\end{equation}
Next, run $(K_Y+D_Y+\sum d_iE_i)$-MMP over $X$. By \eref{eq1},
each extremal ray is negative with respect to the proper transform
of $\sum_{E_j\notin\EEE}(d_j+a(E_j, D))E_j$, an effective divisor.
Such a divisor can be nef only if it is trivial. Hence the process
terminates when the proper transform of
$\sum_{E_j\notin\EEE}(d_j+a(E_j, D))E_j$ becomes zero.
\end{proof}

\begin{corollary-definition}[{\cite[3.1]{Sh2}}, {\cite[21.6.1]{Ut}}]
\label{log-term}
\label{minimal_lt}
Notation as in Theorem~\ref{blow-up}. Then there exists a blowup
$g\colon X'\to X$ and a boundary $D'$ on $X'$ such that (i), (ii)
and (iii) of \ref{blow-up} holds and moreover,
\begin{enumerate}
\item[\textrm{(iv)}${}'$]
if $K_X+D$ is dlt, then $f$ can be taken small; if $K_X+D$ is not
dlt, then for all exceptional divisors $E_i$ of $g$ we have
$a(E_i,D)=-1$.
\end{enumerate}
We call $g\colon X'\to X$ a \emph{minimal log terminal
modification}.\index{log terminal modification!minimal}
\end{corollary-definition}

We generalize slightly the last result:
\begin{proposition}
\label{plt}
Let $X$ be a normal $\QQ$-factorial variety of dimension $\le 3$
and $D$ a boundary on $X$ such that $K_X+D$ is lc, but is not plt.
Assume also that $X$ has only klt singularities. Then there exists
a blowup $f\colon Y\to X$ such that
\begin{enumerate}
\item
$Y$ is $\QQ$-factorial, $\rho(Y/X)=1$ and the exceptional locus of
$f$ is an irreducible divisor, say $E$;
\item
$K_Y+E+D_Y=f^*(K_X+D)$ is lc, where $D_Y$ is the proper transform
of $D$;
\item
$K_Y+E+(1-\ep)D_Y$ is plt and is negative over $X$ for any
$\ep>0$.
\end{enumerate}
\end{proposition}

\begin{proof}
As a first approximation to $Y$ we take a minimal log terminal
modification $g\colon X'\to X$ as in \ref{minimal_lt}. Then
$g^*(K_X+D)=K_{X'}+E'+D'$, where $E'\ne 0$ is an integral reduced
divisor and $D'$ is the proper transform of $D$. In particular,
$\rho(X'/X)$ is the number of components of $E'$. Since $X$ has
only klt singularities, $K_{X'}=g^*K_X+\sum a_iE'_i$, where the
$E'_i$ are components of $E'$ and $a_i>-1$ for all $i$. Therefore
$K_{X'}+E'=g^*K_X+\sum(a_i+1)E'_i$ cannot be $g$-nef by numerical
properties of contractions \cite[1.1]{Sh}. Run $(K_{X'}+E')$-MMP
over $X$. At each step, as above, $K+E$ cannot be nef over $X$.
Hence at the last step we get a divisorial extremal contraction
$f\colon Y\to X$, negative with respect to $K_{Y}+E$ and such that
$\rho(Y/X)=1$. Since at each step $K+E+D$ is numerically trivial
over $X$, the log divisor $f^*(K_X+D)=K_{Y}+E+D_{Y}$ is lc (see
Corollary~\ref{crepant-contraction}). Obviously, $E$ is
irreducible and $\rho(Y/X)=1$. Then $K_Y+E$ is plt and $K_Y+E+D_Y$
is lc. By Proposition~\ref{first_prop_sing}, $K_Y+E+(1-\ep)D$ is
plt.
\end{proof}

\begin{definition}
\label{def_plt_blowup}
Let $X$ be a normal variety and $f\colon Y\to X$ a blowup such
that the exceptional locus of $f$ contains only one irreducible
divisor, say $S$. Assume that $K_Y+S$ is plt and $-(K_Y+S)$ is
$f$-ample. Then $f\colon (Y\supset S)\to X$ is called a
\textit{purely log terminal (plt) blowup} of $X$. \index{purely
log terminal (plt) blowup} The blowup $f\colon(Y, E+D_Y)\to(X,D)$
constructed in \ref{plt} is called an \textit{inductive
blowup}\index{inductive blowup} of $(X,D)$. Note that it not
necessarily unique (cf. \ref{unique}).
\end{definition}

\begin{example1}
\label{ex_finite}
Let $f\colon (Y,S)\to X$ be a plt blowup and $\var\colon X'\to X$
a finite \'etale in codimension one cover. Consider the following
commutative diagram
\begin{equation}
\label{exampl_plt_fin}
\begin{CD}
Y'@>\psi>>Y\\ @Vf'VV@VfVV\\ X'@>\var>>X,\\
\end{CD}
\end{equation}
where $Y'$ is the normalization of $X'\times_XY$. Then $f'\colon
Y'\to X'$ is a plt blowup. Indeed, $f'$ is a contraction and
$\psi$ is a finite morphism, \'etale in codimension one outside of
$S$. Set $S':=\psi^{-1}(S)$. The ramification formula
\eref{eqv-vetvl1_2}, gives $K_{Y'}+S'=\psi^*(K_Y+S)$. Therefore
$K_{Y'}+S'$ is plt (see Proposition~\ref{finite-plt}). By
Connectedness Lemma and Theorem~\ref{inv}, $S'$ is irreducible.
\end{example1}

\begin{example1}
Let $f'\colon Y'\to X'$ be a plt blowup and $G\times Y'\to Y'$ an
equivariant action of a finite group such that the induced action
$G\times X'\to X'$ is free in codimension one. Put $Y:=Y'/G$,
$X:=X'/G$ and consider the commutative diagram
\eref{exampl_plt_fin}. As in \ref{ex_finite} we obtain that
$f\colon Y\to X$ is a plt blowup.
\end{example1}

\begin{proposition-definition}[{\cite{Kawamata-flip},
\cite{Ut}}]
\label{terminal-mod}
Let $(X,D)$ be a klt log variety of dimension $\le 3$. Then there
exists a blowup $f\colon X^t\to X$ and a boundary $D^t$ on $Y$
such that
\begin{enumerate}
\item
$X^t$ is $\QQ$-factorial;
\item
$K_{X^t}+D^t=f^*(K_X+D)$;
\item
$K_{X^t}+D^t$ is terminal.
\end{enumerate}
This blowup is called a \textit{terminal blowup}\index{terminal
blowup} of $(X,D)$.
\end{proposition-definition}
The proof uses Proposition~\ref{per-gen} and the following simple
lemma.

\begin{lemma1}[{\cite{Sh-nonvan}}, {\cite[2.12.2]{Ut}}]
\label{finite-term}
Let $(X,D)$ be a klt log variety. Then the number of divisors $E$
of the function field $\KKK(X)$ with $a(E,D)\le 0$ is finite.
\end{lemma1}
\begin{outline}
Let $f\colon Y\to X$ be a log resolution and $B$ a crepant pull
back of $D$:
\[
f^*(K_X+D)=K_Y+B,\quad\text{with}\quad f_*B=D.
\]
Write $B=B^{(+)}-B^{(-)}$, where $B^{(+)}, B^{(-)}$ are effective
and have no common components. We assume that $\Supp B^{(+)}$
contains also all $f$-exceptional divisors $E$ with discrepancy
$a(E,D)=0$. By Hironaka it is sufficient to construct $f$ so that
all components of $B^{(+)}$ are disjoint. Let $B_i$, $B_j$ be two
components of $B$ such that $B_i\cap B_j\ne \emptyset$. Put
$b_{i,j}:=a(B_i,D)+a(B_j,D)$. We want to change $f\colon Y\to X$
so that $b_{i,j}>0$ whenever $B_i\cap B_j\ne \emptyset$. Let
$a:=1+\inf_E\{a(E,D)\}$. Since $(X,D)$ is klt, $a>0$. By
blowing-up $B_i\cap B_j$ we obtain a new log resolution such that
the proper transforms of $B_i$ and $B_j$ are disjoint and new
exceptional divisor $B_k$ has the discrepancy
\[
a(B_k,D)=1+a(B_i,D)+a(B_j,D).
\]
Then we have
\[
b_{i,k}=a(B_i,D)+1+a(B_i,D)+a(B_j,D)\ge b_{i,j}+a.
\]
Similarly, $b_{j,k}\ge b_{i,j}+a$. Thus after a finite number of
such blowing ups we get the situation when $b_{i,k}>0$ whenever
$B_i\cap B_k\ne \emptyset$. In particular, all components of
$B^{(+)}$ are disjoint.
\end{outline}

One can see that if $(X^t_1,D^t_1)$ is another terminal blowup,
then the induced map $X^t\bir X^t_1$ is an isomorphism in
codimension one. In particular, the terminal blowup is unique in
the surface case.

\section{Weighted blowups}
\label{w-blow-up}
Consider a cyclic quotient singularity $X:=\CC^n/\cyc{m}(a_1,
\dots,a_n)$, where $ a_i\in\NN$ and $\gcd(a_1,\dots,a_n)=1$ (the
case $m=1$, i.e., $X\simeq\CC^n$, is also possible). Let $x_1,
\dots, x_n$ be eigen-coordinates in $\CC^n$, for $\cyc{m}$.
\textit{The weighted blowup}\index{weighted blowup} of $X$ with
weights $a_1,\dots,a_n$ is a projective birational morphism
$f\colon Y\to X$ such that $Y$ is covered by affine charts
$U_1,\dots,U_n$, where
\[
\begin{array}{ccc}
U_i=\CC^n_{y_1,\dots,y_n}/\cyc{a_i}(-a_1,\dots,&m,&\dots,-a_n).\\
&\uparrow&\\&i&\\
\end{array}
\]
The coordinates in $X$ and in $U_i$ are related by
\[
x_i=y_i^{a_i/m}, \qquad x_j=y_jy_i^{a_j/m},\quad j\ne i.
\]
The exceptional locus $E$ of $f$ is an irreducible divisor and
$E\cap U_i=\{y_i=0 \}/\cyc{a_i}$. The morphism $f\colon Y\to X$ is
toric, i.e., there is an equivariant natural action of
$(\CC^*)^n$. It is easy to show that $E$ is the \textit{weighted
projective space}\index{weighted
projective space} $\PP(a_1,\dots, a_n)$\index{$\PP(a_1,\dots, a_n)$} and
$\OOO_E(bE)=\OOO_{\PP}(-mb)$, if $b$ is divisible by
$\lcm(a_1,\cdots,a_n)$ (and then $bE$ is a Cartier divisor).

Note that the blowup constructed above depends on a choice of
numbers $a_1,\dots,a_n$, and not just on their values $\mod m$.

\begin{lemma1}
\label{discr-tor}
In the above conditions we have
\begin{enumerate}
\item
$K_Y=f^*K_X+(-1+\sum a_i/m)E$;
\item
if $D=\{x_i=0 \}/\cyc{m}$ and $D_Y$ is the proper transform of
$D$, then $D_Y=f^*D-\frac{a_i}{m}E$.
\end{enumerate}
\end{lemma1}
\begin{proof}
The relation in (i) follows from the equality
\[
d x_1\wedge\dots\wedge dx_n= y_1^{(\sum a_i/m-1)}
dy_1\wedge\dots\wedge dy_n.
\]
The assertion (ii) can be proved similarly.
\end{proof}

Any weighted blowup $f\colon Y\to \CC^n/\cyc{m}$ of a cyclic
quotient singularity is a plt blowup.

\begin{example1}
Let $X\subset\CC^3$ be a Du Val singularity given by one of the
equations
\[
\begin{array}{ll}
 D_n(n\ge 4):&x^2+y^2z+z^{n-1}=0,\\
 E_6:&x^2+y^3+z^4=0,\\
 E_7:&x^2+y^3+yz^3=0,\\
 E_8:&x^2+y^3+z^5=0.\\
\end{array}
\]
Let $f\colon Y\to X$ be the weighted blowup with weights
$(n-1,n-2,2)$, $(6,4,3)$, $(9,6,4)$, $(15,10,6)$ in cases $D_n$,
$E_6$, $E_7$, $E_8$, respectively. Then $f$ is a plt blowup. We
will see below that it is unique. In the case $D_n$ any weighted
blowup with weights $(w+1,w,2)$, where $w=1,\dots,n-2$ gives a
blowup with irreducible exceptional divisor $C$ such that $K_Y+C$
is lc. It is proved in \cite{IP} that for any hypersurface
canonical singularity given in $\CC^n$ by a nondegenerate
function, there exists a weighted blowup which gives a plt blowup.
\end{example1}

\section{Generalizations of Connectedness Lemma}
Now we generalize Connectedness Lemma to the case nef
anticanonical divisor. We prove them in dimension two. However
there are similar results in arbitrary dimension (modulo log MMP)
\cite{F}.
\begin{proposition1}[{\cite[6.9]{Sh}}]
\label{dipol1}
Let $(X/Z\ni o,D)$ be a log surface, where $Z$ is a curve. Assume
that $K_X+D$ is lc and $-(K_X+D)$ is nef over $Z$. Then in a
neighborhood of the fiber over $o$ the locus of lc singularities
of $(X,D)$ has at most two connected components. Moreover, if
$\LCS(X,D)$ has exactly two connected components, then $(X,D)$ is
plt and $\LCS(X,D)=\down{D}$ is a disjoint union of two sections
(and a general fiber of $X\to Z$ is $\PP^1$).
\end{proposition1}
\begin{proof}
Let $(Y,D_Y)\to(X,D)$ be a minimal log terminal modification (see
\ref{log-term}). Since the fibers of $g$ are connected and
$\LCS(X,D)=g(\LCS(Y,D_Y))$, it is sufficient to prove the
assertion for $(Y,D_Y)$. Let $h\colon Y\to Z$ be the composition
map. Set $C_Y:=\down{D_Y}$ and $B_Y:=\fr{D_Y}$. Then
$\LCS(Y,D_Y)=C_Y$. Assume that $C_Y$ is disconnected. Run
$(K_Y+B_Y)$-MMP over $Z$. If the fiber $h^{-1}(o)$ is reducible,
then there is its component $F\not\subset C_Y$ meeting $C_Y$. Then
$F$ is an extremal curve. Let $Y\to Y_1$ be its contraction. Since
$F\cdot(K_Y+D_Y)=0$ and $F\not\subset C_Y$, the dlt property of
$K_Y+D_Y$ is preserved (see \ref{crepant-contraction}). On the
other hand, by Connectedness Lemma \ref{connect}, $F$ meets only
one connected component of $C_Y$. Hence the number of connected
components of $C_Y$ remains the same. Continuing the process, we
obtain a contraction $\ov h\colon\ov Y\to Z$ with irreducible
fiber $\ov h^{-1}(o)$. Since $K_{\ov Y}+D_{\ov Y}$ is nef, for a
general fiber $L$ of $\ov h$ we have $L\cdot C_{\ov Y}\le L\cdot
D_{\ov Y}={-}K_{\ov Y}\cdot L\le 2$. By our assumption, the fiber
$\ov h^{-1}(o)$ does not contain $C_{\ov Y}$. Hence $C_{\ov Y}$
has exactly two connected components, which are sections of $\ov
Y\to Z$. It is also clear that $(K_{\ov Y}, C_{\ov Y})$ is plt.
The components of $C_Y$ cannot be contractible over $Z$. Hence
$Y\to X$ is the identity map. This proves the statement.
\end{proof}

Similarly we have

\begin{proposition1}[{\cite[6.9]{Sh}}]
\label{dipol2}
Let $(X,D)$ be a projective log surface such that $K_X+D$ is lc
and $-(K_X+D)$ is nef. Then the locus of log canonical
singularities of $(X,D)$ has at most two connected components.
Moreover, if $\LCS(X,D)$ has exactly two connected components and
$K_X+D$ is dlt, then $(X,D)$ is plt and there exists a contraction
$f\colon X\to Z$ with a general fiber $\PP^1$ onto a curve $Z$ of
genus $0$ or $1$ such that $\LCS(X,D)=\down{D}$ is a disjoint
union of two sections.
\end{proposition1}

\begin{proof}
As in the proof of \ref{dipol1}, $g\colon(Y,D_Y)\to(X,D)$ a
minimal log terminal modification. Again set $C_Y:=\down{D_Y}$ and
$B_Y:=\fr{D_Y}$. Assume that $C_Y$ is disconnected. Run
$(K_Y+B_Y)$-MMP. All intermediate contractions is
$(K_Y+D_Y)$-nonpositive. Therefore the log canonical property of
$K_Y+D_Y$ is preserved (see \ref{crepant-singul}). Since at each
step $K_Y+B_Y$ is klt, $K_Y+C_Y+B_Y$ is klt outside of $C_Y$ and
$\LCS(Y,D_Y)=C_Y$. By Connectedness Lemma \ref{connect}, each
contractible curve meets only one connected component of $C_Y$.
Therefore the number of connected components of $\LCS(Y,D_Y)$ is
preserved. At the last step there are two possibilities:
\begin{enumerate}
\itemm{1}
$\rho(\ov Y)=1$, then irreducible components of $\LCS(\ov Y,
D_{\ov Y})$ are intersect each other and gives only one connected
component of $\LCS(X,D)$;
\itemm{2}
$\rho(\ov Y)=2$ and there is a nonbirational contraction $\ov
h\colon\ov Y\to Z$ onto a curve. Here we can apply
Proposition~\ref{dipol1}.
\end{enumerate}
\end{proof}

Log surfaces $(X/Z,D)$, such that $K_X+D$ is lc and numerically
trivial are called \textit{monopoles} if $\LCS(X,D)$ is connected
and \textit{dipoles} if $\LCS(X,D)$ has two connected components.
From \ref{dipol2} we can see that dipoles have a simpler
structure.

\begin{example1}
Let $Z$ be a rational or elliptic curve and
$X:=\PP(\OOO_Z\oplus\FFF)$, where $\FFF$ is an invertible sheaf of
degree $d\ge 0$. There are two nonintersecting sections $C_1$,
$C_2$. If $g(Z)=1$, then $K_X+C_1+C_2=0$ (see \cite[Ch.~5, \S
2]{Ha}) and $\LCS(X,C_1+C_2)$ has two connected components, i.e.
$(X, C_1+C_2)$ is a dipole. Similarly, in the case of a rational
curve $Z$, we can take the log divisor $K_X+C_1+C_2+\sum b_i F_i$,
where $F_i$ are different fibers of $X\to Z$, $\sum b_i=2$,
$b_i<1$, $\forall i$. Then $(X, C_1+C_2+\sum b_i F_i)$ is also a
dipole. We may construct many examples of dipoles by blowing up
points on $C_i$ or blowing down the negative section of $X\to Z$.
For example, we can take a cone over a projectively normal
elliptic curve $Z_d\subset\PP^{d-1}$, and its general hyperplane
section as boundary.
\end{example1}
\begin{example1}[see {\cite{Bl}}, cf. {\cite{Um}}]
Let $X$ be a log Enriques surface (i.e., $K_X$ is lc and
numerically trivial, see \ref{def-log-del-Pezzo}). Then $X$ has at
most two nonklt points. Moreover, if $X$ has exactly two nonklt
points, then they are simple elliptic singularities (see
Theorem~\ref{class_lc}).
\end{example1}

\chapter{Definition of complements and elementary properties}
\label{sect-4}
\label{s3}
\section{Introduction}
The following conjecture is called \textit{Reid's general elephant
conjecture}
\begin{conjecture}
\label{elephant}
Let $f\colon X\to Z\ni o$ be a $K_X$-negative contraction from a
threefold with only terminal singularities. Then near the fiber
over $o$ the linear system ${-}K_X$ contains a divisor having only
Du Val singularities.
\end{conjecture}
At the moment it is known that this conjecture is true (only in
analytic situation) in the following cases:
\begin{itemize}
\item
$X=Z\ni o$ is an isolated singularity \cite{RYPG}, moreover, this
is equivalent to the classification of three-dimensional terminal
singularities;
\item
$f\colon X\to Z$ is an extremal flipping or divisorial small
contraction \cite{Mo}, \cite{KoM}, this is a sufficient condition
for the existence of flips \cite{K}.
\end{itemize}
Some particular results are known in the case when $f\colon X\to
Z$ is an extremal contraction to a surface \cite{Pr}. This case is
interesting for applications to rationality problem of conic
bundles.
\par
However, at the moment it is not so clear how one can prove Reid's
conjecture in the algebraic situation. Moreover, it fails for the
case $Z=\mt{pt}$ (there are examples of $\QQ$-Fano threefolds with
empty $|{-}K_X|$). Shokurov proposed the notion of complements,
which is weaker then ``general elephant'' but much more easier to
work with.

\begin{definition}
Let $(X,D)$ be a log pair, where $D$ is a subboundary. Then a
\textit{$\QQ$-complement}\index{$\QQ$-complement} of $K_X+D$ is a
log divisor $K_X+D'$ such that $D'\ge D$, $K_X+D'$ is lc and
$n(K_X+D')\sim 0$ for some $n\in\NN$.
\end{definition}

\begin{definition}[\cite{Sh}]
\label{ldef}
Let $X$ be a normal variety and $D=S+B$ a subboundary on $X$, such
that $B$ and $S$ have no common components, $S$ is an effective
integral divisor and $\down{B}\le 0$. Then we say that $K_X+D$ is
\textit{$n$-complementary}, if there is a $\QQ$-divisor $D^+$ such
that
\begin{enumerate}
\item
$n(K_X+D^+)\sim 0$ (in particular, $nD^+$ is integral divisor);
\item
$K_X+D^+$ is lc;
\item
$nD^+\ge nS+\down{(n+1)B}$.
\end{enumerate}
In this situation the
\textit{$n$-complement}\index{$n$-complement}\index{complement} of
$K_X+D$ is $K_X+D^+$. If moreover $K_X+D^+$ is plt, then we say
that $K_X+D$ is \textit{strongly $n$-complementary}.\index{strong
$n$-complement}
\end{definition}
Note that an $n$-complement is not necessarily a $\QQ$-complement
because of condition (iii). We need this condition for technical
reasons (see \ref{prodolj}). If $B=0$, then (iii) holds
automatically. In applications this is the most interesting case.
We give also a generalization of this definition for the case of
nodal curves.

\begin{definition1}
Let $X$ be a reduced (not necessarily irreducible) curve. Then $X$
is said to be \textit{nodal}\index{nodal curve} if all its
singularities are normal crossing points. A subboundary $D=\sum
d_iD_i$ on a nodal curve is said to be \textit{semilog
canonical}\index{semilog canonical (slc) singularities}
(\textit{slc}) if $\Supp D\cap\Sing X=\emptyset$ and $d_i\le 1$
for all $i$.
\par
Let $X$ be a nodal curve and $D=S+B$ a subboundary on $X$, such
that $B$ and $S$ have no common components, $S$ is an effective
integral divisor and $\down{B}\le 0$. Assume that $\Supp
D\cap\Sing X=\emptyset$. Then an
\textit{$n$-semicomplement}\index{$n$-semicomplement} of $K_X+D$
is a log divisor $K_X+D^+$ such that conditions (i), (iii) of
\ref{ldef} and the following (ii${}'$) below holds.
\begin{enumerate}
\item[(ii${}'$)]
$K_X+D^+$ is slc.
\end{enumerate}
\end{definition1}
The last definition can be generalized to the higher dimensional
case (see~\cite{Ut}).

\begin{remark1}
Assume that on a variety $X$ the canonical divisor $K_X$ is
strongly $1$-complementary. Let $K_X+B$ be this complement. Then
$B$ is an integral divisor, $B\in |{-}K_X|$ and $K_X+B$ is plt
(and even canonical because $K_X+B\sim 0$). By \ref{diff}, $\Diff
B(0)=0$ and by Inversion of Adjunction, $K_B$ is klt. Since
$K_B\sim 0$, $B$ has only canonical Gorenstein singularities. This
shows that $K_X$ is strongly $1$-complementary if and only if
Reid's general elephant conjecture holds for $X/Z$.
\end{remark1}

The following conjecture seems to be more realistic than
Conjecture~\ref{elephant}:

\begin{conjecture1}
\label{12346}
Let $f\colon X\to Z\ni o$ be a contraction from a threefold with
only terminal singularities such that ${-}K_X$ is $f$-nef and
$f$-big. Then near the fiber over $o$ the canonical divisor $K_X$
is $1$, $2$, $3$, $4$, or $6$-complementary.
\end{conjecture1}

Note that the condition that $K_X+D$ is $n$-complementary implies
the existence of an integral effective divisor
\begin{equation}
\label{compl-prop}
\ov{D}\in\left|{-}nK_X-nS-\down{(n+1)B}\right|
\end{equation}
related to $D^+$ by the equality
\begin{equation}
\label{compl-prop1}
D^+:=S+\frac1{n}\left(\down{(n+1)B}+\ov{D}\right).
\end{equation}
It is also easy to see that if $D$ is a boundary, then so is
$D^+$. As an immediate consequence of the definition we have

\begin{proposition}
\label{simple-prp}
Let $X$ be a normal variety.
\begin{enumerate}
\item
Fix $n\in\NN$. Let $D=\sum d_iD_i$ and $D'=\sum d_i'D_i$ be
subboundaries on $X$ such that the following conditions hold
\begin{enumerate}
\itemm{a}
$d_i'\ge d_i-\ep$ for $0<\ep\ll 1$;
\itemm{b}
$d_i'\ge d_i$ whenever $(n+1)d_i$ is an integer $\le n$.
\end{enumerate}
Assume that $K_X+D'$ is $n$-complementary. Then $K_X+D$ is
$n$-complementary.
\item
Fix $n\in\NN$. Let $D=\sum d_iD_i$ be a subboundary. Assume that
$K_X+D$ is $n$-complementary. Then $K_X+D'$ is $n$-complementary
for any subboundary $D'=\sum d_i'D_i$ such that
$|d_i'-d_i|<\frac1{(n+1)q_i}$, $\forall i$, where $q_i\ge 1$ is
the denominator of $d_i$.
\end{enumerate}
\end{proposition}
\begin{outline}
We show, for example, (i). Let ${D'}^+=\sum {d'_i}^+D_i$ be an
$n$-complement of $K_X+D'$. Put $D^+:={D'}^+$. It is sufficient to
verify (iii) of \ref{ldef}, i.e.
\begin{equation}
\label{eq-prop-d}
{d'_i}^+\ge \left\{
\begin{array}{ll}
\frac1n\down{(n+1)d_i}\quad& \text{if}\ d_i<1\\
 1& \text{otherwise}.
\end{array}
\right.
\end{equation}
On the other hand, we have
\[
{d'_i}^+\ge \left\{
\begin{array}{ll}
\frac1n\down{(n+1)d_i'}\quad& \text{if}\ d_i'<1\\
 1&
\text{otherwise}.
\end{array}
\right.
\]
If $n+1$ is a denominator of $d_i$, then $d_i\le d_i'$ and
\eref{eq-prop-d} is obvious. If $n+1$ is not a denominator of
$d_i$, then $d_i'\ge d_i-\ep$. Hence
$\down{(n+1)d_i}=\down{(n+1)d_i'}$ for small positive $\ep$. Again
we obtain \eref{eq-prop-d}. Finally, if $d_i=1$, then ${d'_i}^+\ge
\frac1n\down{(n+1)(1-\ep)}=1$ for $\ep<1/(n+1)$.
\end{outline}

\begin{corollary1}
\label{simple-prp-1}
Let $X$ be a normal variety and $D=\sum d_iD_i$ a subboundary on
$X$. Fix $n\in\NN$. Let $D'=\sum d_i'D_i$, where $d_i'\ge
\min\{d_i,\frac{n}{n+1}\}$. Assume that $K_X+D'$ is
$n$-complementary. Then so is $K_X+D$.
\end{corollary1}

\begin{example}
\label{toric}
\begin{enumerate}
\item
Let $X$ be a toric variety and $S=\sum S_i$ be the toric boundary.
Then $K_X+S\sim 0$ and $K_X+S$ is lc. Hence $K_X+S$ is
$1$-complementary.
\item
Let $(X\ni P)$ be an analytic germ of a three-dimensional terminal
singularity. Then $K_X$ is strongly $1$-complementary (see
\ref{YPG}). Conversely, if there is a strong nontrivial
$1$-complement near an isolated three-dimensional $\QQ$-Gorenstein
singularity $(X\ni P)$, then $(X\ni P)$ is terminal. A
three-dimensional Gorenstein canonical singularity is
(nontrivially) $1$-complementary if and only if it is cDV (see
\ref{YPG_G}).
\item
Consider the cyclic quotient singularity $X:=\CC^3/\ZZ_9(1,4,7)$.
By \ref{Reid-canonical-ref} it is canonical. Since $X$ is not
terminal, $K_X$ is not strongly $1$-complementary. However, $K_X$
is strongly $2$-complementary. Indeed, $(x^2y+y^2z+z^2x)(dx\wedge
dy\wedge dz)^{-2}$ is an invariant form. Hence this gives us a
member of $|-2K_X|$. It is easy to check that $x^2y+y^2z+z^2x=0$
defines a log canonical singularity $F\subset \CC^3$. Let
$F':=F/\ZZ_9\subset X$. By Corollary \ref{corollary-finite}
$K_X+\frac12F'$ is klt. Note that $K_X$ is $1$-complementary in
this case (see (i)).
\item
According to \cite{MMM} there are four-dimensional terminal cyclic
quotient singularities which have no strong $1$ or
$2$-complements. However, it is expected that there are only a
finite number of such exceptions. For example, the singularity
$\CC^4/\ZZ_{83}(3,14,23,44)$ has no strong $1$ or $2$-complements.
As above, the invariant $(dx\wedge dy\wedge dz)^{-3}$ gives us a
strong $3$-complement.
\item
Let $f\colon X\to Z\ni o$ be an analytic germ of a
three-dimensional flipping extremal contraction. Then $K_X$ is
strongly $1$-complementary \cite{Mo}, \cite{KoM}.
\item
Let $X$ be a Fano threefold with Gorenstein canonical
singularities. Then $K_X$ is strongly $1$-complementary
\cite{Sh0}, \cite{R}.
\item
Let $X$ be a variety with log canonical singularities and
numerically trivial canonical divisor $K_X$. Then $K_X$ is
$n$-complementary if and only if there exists $n$ such that
$nK_X\sim 0$. For example, in the case of a smooth surface of
Kodaira dimension $0$ the canonical divisor is either $1$, $2$,
$3$, $4$ or $6$-complementary (see e.g. \cite{BPV}).
\item
\label{relatively-minimal-elliptic}
Let $g\colon X\to \PP^1$ be a relatively minimal elliptic
fibration, where $X$ is a smooth surface of Kodaira dimension
$\kappa(X)\le 0$. Then $X$ is $n$-complementary for some $n\in\{1,
2, 3, 4, 6\}$. Indeed, we have the canonical bundle formula (see
e.g. \cite[Ch. V, \S 12]{BPV})
\[
K_X\sim \bigl(\chi(\OOO_X)-2\bigr)L+\sum_{i=1}^s (r_i-1)E_i,
\]
where $E_i$ are multiple fibers of multiplicities $r_i$ and $L$ is
a general fiber. Consider, for example, case $\kappa(X)=-\infty$.
Then $X$ is a ruled surface over an elliptic curve. Let $F$ be a
general fiber of the rulling and denote $\delta:=L\cdot F$.
Clearly, $K_X\cdot F=-2$ and $E_i\cdot F=\delta/r_i$. This gives
us
\[
-2=-2\delta+\sum_{i=1}^s (r_i-1)\delta/r_i, \qquad
2-2/\delta=\sum_{i=1}^s(1-1/r_i).
\]
It is easy to see that $2\le s\le 3$ and $K_X\sim (s-1)L+\sum
E_i$. There are only the following possibilities:
\par
\medskip
\begin{center}
\begin{tabular}{|c||c|c|c|c|c|}
\hline $s$&$2$&$3$&$3$&$3$&$3$\\ \hline
$(r_1,\dots,r_s)$&$(r_1,r_2)$&$(2,2,r)$&$(2,3,3)$&
$(2,3,4)$&$(2,3,5)$\\ \hline
$\delta$&$2\,\left(r_1^{-1}+r_2^{-1}\right)^{-1}
$&$2r$&$12$&$24$&$60$\\\hline $n$&$1$&$2$&$3$&$4$&$6$\\ \hline
$n$-complement&$E_1+E_2$&$E_3$&$\frac13E_1$&$\frac14E_2$&$\frac16E_3$\\
\hline
\end{tabular}
\end{center}
\par\medskip\noindent
This shows that $n\in\{1,2,3,4,6\}$.
\end{enumerate}
\end{example}

Further, in the two-dimensional case $1$, $2$, $3$, $4$ and
$6$-complement we call \textit{regular}\index{regular complement}
and define
\[
\RRR_2:=\{1,2,3,4,6\}.\index{$\RRR_2$}
\]
In the higher-dimensional case we should replace the set $\RRR_2$
with bigger one $\RRR_n$ (see \cite{PSh}).

A very important question is:
\begin{quote}
when does some $n$-complement of $K_X+D$ exist?
\end{quote}
Obviously, these exist for some $n\gg 0$ when $-(K_X+D)$ is ample
(or even semiample) \cite[5.5]{Sh}. By Base Point Free Theorem
(see \cite[3-12]{KMM}), $n$-complements exist for some $n\gg 0$ if
$K_X+D$ is klt and $-(K_X+D)$ is nef and big. It is expected also
that we can remove klt condition on lc and $D\in\Msm$ (see
Proposition~\ref{E-C}).

In general, only the nef condition is not sufficient for the
existence of complements (see Example~\ref{Atiyah}).

\begin{theorem}[see {\cite[5.2]{Sh}}, {\cite[19.4]{Ut}}]
\label{1}
Let $X$ be a nodal connected (but not necessarily compact) curve.
Let $D$ be a boundary on $X$ contained in the smooth and compact
part of $X$. Assume that the degree of $-(K_X+D)$ is nonnegative
(on each compact component of $X$). Then
\begin{enumerate}
\item
$K_X+D$ is $n$-semicomplementary for $n\in\RRR_2$;
\item
if $K_X+D$ is not $1$ or $2$-semicomplementary, then
$X\simeq\PP^1$ and $\down{D}=\down{D^+}=0$;
\item
if $X$ contains a noncomplete component and $K_X+D$ is not
$1$-semicomplementary, then the compact components of $X$ form a
chain $\sum_{i=1}^r X_i$, a (unique) noncomplete component $X'$
intersect an end $X_1$ of $\sum X_i$, $\Supp D$ is contained in
another one $X_r$ and $D=1/2P_1+1/2P_2$ (the case $r=1$, $X_1=X_r$
is also possible).
\end{enumerate}
\end{theorem}

For each log pair $(X/Z\ni o,D)$ we define the minimal
complementary number by
\begin{equation}
\label{eq-def-compl}
\compl(X,D):=\min\{m\mid K_X+D\ \text{is $m$-complementary}\}.
\index{$\compl(X,D)$}
\end{equation}
This is an invariant to ``measure'' how singular a log pair is. We
also define
\begin{equation}
\label{eq-def-compl-1}
\compl'(X,D):=\min\{m\mid\exists\ \text{$m$-complement of}\ K_X+D\
\text{which is not klt}\}. \index{$\compl'(X,D)$}
\end{equation}

By definition, $\compl(X,D), \compl'(X,D)\in\NN\cup\{\infty\}$.
Take a subset $\MMM\subset [0,1]$. For example, consider cases
$\MMM=\Msm$ (see \ref{def-SM}), $\MMM=[0,1]$, or $\MMM=\Mm$. For a
$\QQ$-divisor $D$ we write simply $D\in\MMM$ if all the
coefficients of $D$ belong to $\MMM$.

Define the set of natural numbers $\NNN_n(\MMM)$ by
\begin{multline*}
\NNN_n(\MMM):=\{m\in\NN\mid\exists\ \text{a log Fano variety} \
(X,D)\ \text{of dimension $n$}\\ \text{such that} \ D\in\MMM\
\text{and}\ \compl(X,D)=m\}.
\end{multline*}
\index{$\NNN_n(\MMM)$} Thus Theorem~\ref{1} and
Corollary~\ref{m_i} below give us
$\NNN_1([0,1])=\RRR_2:=\{1,2,3,4,6\}$. Obviously,
$\NNN_n(\MMM')\subset\NNN_n(\MMM'')$ if $\MMM'\subset\MMM''$.
Theorem~\ref{1} and Corollary~\ref{m_i} show that
$\NNN_1(\Msm)=\NNN_1([0,1])=\RRR_2$. We will see below that
$\NNN_2(\Msm)$ is bounded (Theorem~\ref{main_Sh_A}).

\begin{corollary1}
\label{m_i}
Notation as in Theorem~\ref{1}. Assume that $X\simeq\PP^1$,
$\down{D}=0$, $-(K_X+D)$ is ample and $D\in\Msm$. Then
$D=\sum_{i=1}^r (1-1/m_i)D_i$, where for $(m_1,\dots,m_r)$ there
is only one of the following possibilities (up to permutations):
\begin{description}
\item[$A_n$]
$(m)$, $(m_1,m_2)$, $K_X+D$ is $1$-complementary;
\item[$D_n$]
$(2,2,m)$, $K_X+D$ is $2$-complementary;
\item[$E_6$]
$(2,3,3)$, $K_X+D$ is $3$-complementary;
\item[$E_7$]
$(2,3,4)$, $K_X+D$ is $4$-complementary;
\item[$E_8$]
$(2,3,5)$, $K_X+D$ is $6$-complementary.
\end{description}
\end{corollary1}

Relations between our notation $A_n$, $D_n$, $E_6$, $E_7$, $E_8$
and two-dimensional Du Val singularities will be explained in
Ch.~\ref{s4}.

\begin{exercise1}
\label{=}
Let $X\simeq\PP^1$ and $D\in\Msm\cap \left[0,1\right)$. Assume
also that $\deg D=2$. Show that $D=\sum_{i=1}^r (1-1/m_i)D_i$,
where for $(m_1,\dots,m_r)$ there is only one of the following
possibilities:
\begin{description}
\item[$\widetilde D_4$]
$(2,2,2,2)$, $K_X+D$ is $2$-complementary;
\item[$\widetilde E_6$]
$(3,3,3)$, $K_X+D$ is $3$-complementary;
\item[$\widetilde E_7$]
$(2,4,4)$, $K_X+D$ is $4$-complementary;
\item[$\widetilde E_8$]
$(2,3,6)$, $K_X+D$ is $6$-complementary.
\end{description}
\end{exercise1}

\section{Monotonicity}
We noticed that the inequality $D^+\ge D$ does not hold in
general. However, under some additional restrictions on
coefficients we can expect $D^+\ge D$ to be true.

\subsection{}\label{def_PPP_n_sect}
Fix $n\in\NN$ and define the set $\PPP_n$ by
\[
\label{def_PPP_n}
\alpha\in\PPP_n\quad\Longleftrightarrow\quad 0\le\alpha\le
1\quad\text{and}\quad \down{(n+1)\alpha}\ge n\alpha.
\index{$\PPP_n$}
\]
\begin{corollary1}
\label{monot-1}
Let $(X,D)$ be a log pair such that $D\in\PPP_n$ and $K_X+D^+$ any
$n$-complement. Then $D^+\ge D$.
\end{corollary1}

\begin{lemma1}[Monotonicity of the integral part]
\label{monot}
\begin{enumerate}
\item
Let $r\in\QQ$ such that $r<1$ and $nr\in\ZZ$. Then
\[
\down{(n+1)r}\le nr.
\]
\item
Let $r=1-1/m$, $n\in\NN$. Then for any $n\in\NN$
\[
\down{(n+1)r}\ge nr.
\]
\end{enumerate}
\end{lemma1}
\begin{proof}
Let us proof, for example, (ii). Write $nr=q+k/m$, where
$q=\down{nr}$ and $k/m=\fr{nr}$, $k\in\ZZ$, $0\le k\le m-1$. Then
\[
\down{(n+1)r}=\down{q+k/m+1-1/m}=\left\{
\begin{array}{ll}
 q\qquad&\text{if}\ k=0,\\
 q+1\qquad&\text{otherwise.}
\end{array}
\right.
\]
In both cases $\down{(n+1)r}\ge q+k/m=nr$.
\end{proof}

By Monotonicity Lemma, $\PPP_n\supset\Msm$ for any $n\in\NN$.
Moreover, we have
\begin{corollary1}
\[
\Msm=\bigcap_{n\in\NN}\PPP_n.
\]
\end{corollary1}
\begin{proof}
Let $\alpha\notin\Msm$. Then $1-1/m<\alpha<1-1/(m+1)$ for some
$m\in\NN$. This yields $\down{(m+1)\alpha}\le m-1$ and
$m\alpha>m-1$. Hence $\alpha\notin\PPP_m$.
\end{proof}

\begin{example1}
Let $(X,D)$ be a log variety with a standard boundary (i.e.
$D\in\Msm$). Assume that $K_X+D$ is numerically trivial. Let
$K_X+D^+$ be some $n$-complement. Then $D^+\ge D$ and $D^+\equiv
D$. If $X$ is projective, then this yields $D^+=D$. In this case,
$n$ is any natural such that $nD$ is an integral divisor. In
general case, we say that a complement $K_X+D^+$ of $K_X+D$ is
\textit{trivial} if $D=D^+$.\index{trivial complement}
\end{example1}

It is easy to check that
\begin{multline*}
\PPP_n=\left\{0\right\}\bigcup \downn{\frac1{n+1},\frac1{n}}
\bigcup \downn{\frac2{n+1},\frac2{n}}\bigcup\dots\\
 \bigcup\downn{\frac k{n+1},\frac k{n}}\bigcup\dots\bigcup
\downn{\frac{n}{n+1},1}.
\end{multline*}
This gives
\begin{lemma1}
\begin{enumerate}
\item
If $\alpha_1,\alpha_2\in\PPP_n$, then either
$\alpha_1+\alpha_2\in\PPP_n$ or $\alpha_1+\alpha_2>1$.
\item
Let $m\in\NN$, $k_j\in\NN\cup\{0\}$ and $b_j\in\PPP_n$,
$j=1,\dots,r$. Assume that
\[
\alpha:=\frac{m-1}{m}+\frac1{m}\sum_{j=1}^rk_jb_j\le 1.
\]
Then $\alpha\in\PPP_n$.
\end{enumerate}
\end{lemma1}
\begin{proof}
(i) is trivial. As for (ii) we notice that $\sum k_jb_j\le 1$.
Hence by (i), $\sum k_jb_j\in\PPP_n$ and we may assume that $r=1$
and $k_1=1$. Put $b:=b_1$. It is sufficient to show that there
exists $q\in\NN$ such that
\[
\frac{q}{n+1}\le \frac{m-1+b}{m}\le\frac{q}{n}.
\]
This is equivalent to
\[
n(m-1+b)\le mq\le (n+1)(m-1+b).
\]
Taking into account that $b\in\PPP_n$, we have $l/(n+1)\le b\le
l/n$ for some $l\in\NN$. So there exists $q\in\NN$ such that
\begin{multline*}
n(m-1+b)\le n(m-1)+l\le mq\le\\ (n+1)(m-1)+l\le (n+1)(m-1+b).
\end{multline*}
This proves the lemma.
\end{proof}
From Corollary~\ref{coeff-Diff} we have
\begin{corollary1}
\label{coeff-Diff-nonst}
Let $(X,S+B)$ be a log variety, where $S$ is reduced and $B$ is
effective. Assume that $K_X+S+B$ is plt and $B\in\PPP_n$. Then
$\Diff S(B)\in\PPP_n$.
\end{corollary1}

\begin{remark1}
\label{important_compl}
It is easy to see that $\Mm\subset \PPP_1\cup\PPP_2\cup
\PPP_3\cup\PPP_4 \cup\PPP_6$. Therefore if $D\in\Mm$ and $K_X+D^+$
is an $1$, $2$, $3$, $4$ or $6$-complement of $K_X+D$, then
$D^+\ge D$.
\end{remark1}

\section{Birational properties of complements}
\label{Birational-properties-of-complements}
Now we will see that complements have good birational properties.
\begin{proposition1}[{\cite{Sh}}]
\label{bir-prop-a}
Let $f\colon X\to Y$ be a birational contraction and $D$ a
subboundary on $X$. Assume that $K_X+D$ is $n$-complementary for
some $n\in\NN$. Then $K_Y+f(D)$ is also $n$-complementary.
\end{proposition1}
\begin{proof}
Take $f(D)^+:=f_*(D^+)$ and apply \ref{crepant-singul}.
\end{proof}

Under additional assumptions we have the inverse implication:

\begin{proposition1}[{\cite[2.13]{Sh1}}]
\label{bir-prop}
Fix $n\in\NN$. Let $f\colon Y\to X$ be a birational contraction
and $D$ a subboundary on $Y$ such that
\begin{enumerate}
\item
$K_Y+D$ is nef over $X$;
\item
$f(D)\in\PPP_n$ (in particular, $f(D)$ is a boundary).
\end{enumerate}
Assume that $K_X+f(D)$ is $n$-complementary. Then $K_Y+D$ is also
$n$-complementary.
\end{proposition1}

\begin{proof}
Consider the crepant pull back
\[
K_Y+D'=f^*(K_X+f(D)^+),\quad \text{with}\quad f_*D'=f(D)^+.
\]
Write $D'=S'+B'$, where $S'$ is reduced, $S'$, $B'$ have no common
components, and $\down{B'}\le 0$. We claim that $K_Y+D'$ is an
$n$-complement of $K_Y+D$. The only thing we need to check is that
$nB'\ge \down{(n+1)\fr{D}}$. From (ii) we have $f(D)^+\ge f(D)$.
This gives that $D'\ge D$ (because $D-D'$ is $f$-nef; see
\cite[1.1]{Sh} or \cite[3.39]{KM}). Finally, by Monotonicity
Lemma~\ref{monot} and because $nD'$ is an integral divisor, we
have
\[
nD'\ge nS'+\down{(n+1)B'}\ge n\down{D}+\down{(n+1)\fr{D}}.
\]
\end{proof}

\begin{remark1}
\label{bir-prop-comment}
\begin{enumerate}
\item
By Monotonicity Lemma~\ref{monot}, the condition (ii) holds if all
the coefficients of $f(D)$ are standard, i.e., $f(D)\in\Msm$. By
\ref{important_compl} (ii) also holds if $n\in\RRR_2$ and
$f(D)\in\Mm$.
\item
The above proof shows that the proposition holds under the
following weaker assumption instead of (ii):
\end{enumerate}
\begin{enumerate}
\item[(ii)${}'$]
for each nonexceptional component $D_i$ of $D=\sum d_iD_i$ meeting
the exceptional divisor of $f$ we have $d_i\in\PPP_n$.
\end{enumerate}
\end{remark1}

\section{Inductive properties of complements}
\label{Inductive-properties-of-complements}
\begin{proposition1}[{cf. \cite[Proof of 5.6]{Sh}},
{\cite[19.6]{Ut}}]
\label{prodolj}
Let $(X/Z\ni o,D=S+B)$ be a log variety. Set $S:=\down{D}$ and
$B:=\fr{D}$. Assume that
\begin{enumerate}
\item
$K_X+D$ is plt;
\item
$-(K_X+D)$ is nef and big over $Z$;
\item
$S\ne 0$ near $f^{-1}(o)$;
\item
$D\in\PPP_n$ for some $n\in\NN$.
\end{enumerate}
Further, assume that near $f^{-1}(o)\cap S$ there exists an
$n$-complement $K_S+\Diff S(B)^+$ of $K_S+\Diff S(B)$. Then near
$f^{-1}(o)$ there exists an $n$-complement $K_X+S+B^+$ of
$K_X+S+B$ such that $\Diff S(B)^+=\Diff S(B^+)$.
\end{proposition1}

This proposition should be true in the case when $K_X+D$ is dlt.
We need only good definitions of complements on nonnormal
varieties (see \cite{Ut}).

\begin{proof}
Let $g\colon Y\to X$ be a log resolution. Write
$K_Y+S_Y+A=g^*(K_X+S+B)$, where $S_Y$ is the proper transform of
$S$ on $Y$ and $\down{A}\le 0$. By Inversion of Adjunction, $S$ is
normal and $K_S+\Diff S(B)$ is plt. In particular, $g_S\colon
S_Y\to S$ is a birational contraction. Therefore we have
\[
K_{S_Y}+\Diff{S_Y}(A)=g_S^*(K_S+\Diff S(B)).
\]
Note that $\Diff{S_Y}(A)=A|_{S_Y}$, because $Y$ is smooth. By
Corollary~\ref{coeff-Diff-nonst} we see that $\Diff
S(B)\in\PPP_n$. So we can apply Proposition~\ref{bir-prop} to
$g_S$. We get an $n$-complement $K_{S_Y}+\Diff{S_Y}(A)^+$ of
$K_{S_Y}+\Diff{S_Y}(A)$. In particular, by \eref{compl-prop},
there exists
\[
\Theta\in \left|{-}nK_{S_Y}-\down{(n+1) \Diff{S_Y}(A)}\right|
\]
such that
\[
n\Diff{S_Y}(A)^+= \down{(n+1)\Diff{S_Y}(A)}+\Theta.
\]
By Kawamata-Viehweg Vanishing,
\begin{multline*}
R^1h_*\left(\OOO_{Y}({-}nK_Y-(n+1)S_Y-\down{(n+1)A})\right)=\\
R^1h_*\left(\OOO_{Y}(K_Y+\up{-(n+1)(K_Y+S_Y+A)})\right)=0.
\end{multline*}
From the exact sequence
\begin{multline*}
0\longrightarrow\OOO_{Y}({-}nK_Y-(n+1)S_Y-\down{(n+1)A})\\
\longrightarrow\OOO_{Y}({-}nK_Y-nS_Y-\down{(n+1)A})\\
\longrightarrow\OOO_{S_Y}({-}nK_{S_Y}-\down{(n+1)A}|_{S_Y})
\longrightarrow 0
\end{multline*}
we get surjectivity of the restriction map
\begin{multline*}
H^0(Y,\OOO_{Y}({-}nK_Y-nS_Y-\down{(n+1)A})) \longrightarrow\\
H^0(S_Y,\OOO_{S_Y}({-}nK_{S_Y}-\down{(n+1)A}|_{S_Y})).
\end{multline*}
Therefore there exists a divisor
\[
\Xi\in\left|{-}nK_Y-nS_Y-\down{(n+1)A}\right|
\]
such that $\Xi|_{S_Y}=\Theta$. Set
\[
A^+:=\frac1{n}(\down{(n+1)A}+\Xi).
\]
Then $n(K_Y+S_Y+A^+)\sim 0$ and $(K_Y+S_Y+A^+)|_{S_Y}=
K_{S_Y}+\Diff{S_Y}(A)^+$. Note that we cannot apply Inversion of
Adjunction on $Y$ because $A^+$ can have negative coefficients. So
we put $B^+:=g_*A^+$. Again we have $n(K_X+S+B^+)\sim 0$ and
$(K_X+S+B^+)|_S=K_S+\Diff S(B)^+$. We have to show only that
$K_X+S+B^+$ is lc. Assume that $K_X+S+B^+$ is not lc. Then
$K_X+S+B+\alpha(B^+-B)$ is also not lc for some $\alpha<1$. It is
clear that $-(K_X+S+B+\alpha(B^+-B))$ is nef and big over $Z$. By
Inversion of Adjunction, $K_X+S+B+\alpha(B^+-B)$ is plt near
$S\cap f^{-1}(o)$. Hence $\LCS(X,B+\alpha(B^+-B))=S$ near $S\cap
f^{-1}(o)$. On the other hand, by Connectedness Lemma,
$\LCS(X,B+\alpha(B^+-B))$ is connected near $f^{-1}(o)$. Thus
$K_X+S+B+\alpha(B^+-B)$ is plt. This contradiction proves the
proposition.
\end{proof}

\begin{remark1}
\label{prodolj-1}
It follows from the proof that we can replace (iv) in
Proposition~\ref{prodolj} with
\begin{enumerate}
\item[(iv)${}'$]
$\Diff S(B)\in\PPP_n$ for some $n$.
\end{enumerate}
\end{remark1}

In the two-dimensional case we have a stronger result.

\begin{proposition1}[{cf. \cite[Proof of 5.6]{Sh}},
{\cite[19.6]{Ut}}]
\label{prodoljenie}
Let $(X/Z\ni o,D=S+B)$ be a log surface such that
\begin{enumerate}
\item
$K_X+D$ is dlt;
\item
$-(K_X+D)$ is nef and big over $Z$;
\item
$S:=\down{D}\ne 0$ near $f^{-1}(o)$.
\end{enumerate}
Assume that near $f^{-1}(o)\cap S$ there exists an
$n$-semicomplement $K_S+\Diff S(B)^+$ of $K_S+\Diff S(B)$. Then
near $f^{-1}(o)$ there exists an $n$-complement $K_X+S+B^+$ of
$K_X+S+B$ such that $\Diff S(B)^+=\Diff S(B^+)$.
\end{proposition1}
\begin{proof}
Similar to the proof of \ref{prodolj}. By Propositions
\ref{kawamata_1} and \ref{kawamata_2}, the curve $S$ is nodal.
Further, we can take a log resolution $g\colon Y\to X$ so that
$S_Y\simeq S$.
\end{proof}

\begin{exercise1}[\cite{Sh1}]
\label{exer-odna-comp}
Let $(X\ni P)$ be a germ of a two-dimensional normal singularity,
let $C\ne 0$ be a reduced divisor on $X$, and $B=\sum b_iB_i\ne 0$
a boundary on $X$ such that $K_X+C+B$ is plt. Assume that $b_i\ge
1/2$ for all $i$. Show that $K_X+C+\up{B}$ is lc and $\Supp B$ is
irreducible. Moreover,
\[
\Diff C(B)=\left(1-\frac1n+\frac{b_1}n\right)P,\quad\text{where}
\quad (X\ni P)\simeq \CC^2/\ZZ_n(1,q).
\]
If $B\in\Msm$ (i.e., $B=(1-1/m_1)B_1$, $m_1\in \NN$), then
\[
(K_X+C+B)|_C=(1-1/m)P, \qquad m=m_1n.
\]
\begin{hint}
Show that $K_X+C+B$ is $1$-complementary (using \ref{1} and
\ref{prodoljenie}).
\end{hint}
\end{exercise1}

\section{Exceptionality}
\label{s_except}
\begin{definition1}
\label{def_except}
Let $(X/Z\ni o,D)$ be a log variety such that there is at least
one $\QQ$-complement of $K_X+D$ near the fiber over $o$.
\begin{itemize}
\item
Assume that $Z$ is not a point (local case). Then $(X/Z\ni o,D)$
is said to be \textit{exceptional}\index{exceptional log variety}
over $o$ if for any $\QQ$-complement of $K_X+D$ near the fiber
over $o$ there exists at most one (not necessarily exceptional)
divisor $S$ such that $a(S,D)=-1$.
\item
Assume that $Z$ is a point (global case). Then $(X,D)$ is said to
be \textit{exceptional} if every $\QQ$-complement of $K_X+D$ is
klt.
\end{itemize}
\end{definition1}

The main advantage of this definition is Shokurov's conjecture
that exceptional log varieties are bounded in some sense (see
\ref{ex-1-ex}, \ref{except-bound}, \ref{rem-lc-exc} (ii),
\ref{except-bound-conic}, \ref{elliptic-model-classification},
\ref{main_Sh_A-1}, \ref{main_Sh_A-2}, \cite[\S 7]{Sh1},
\cite{KeM}, \cite{Pr2}). On the contrary, nonexceptional ones has
``regular'' complements (i.e., $n$-complements with small $n$).
This phenomena was discovered by Shokurov \cite{Sh1}. In
\cite{KeM} exceptional log del Pezzo surfaces are called
\textit{del Pezzo surfaces without tiger}. Studying of such
surfaces is closely related to the uniruledness of affine surfaces
\cite[6.1]{KeM}.

\begin{example1}
\label{ex-1-ex}
Let $D$ be a boundary on a curve $X$. If $(X,D)$ is
nonexceptional, then by Theorem~\ref{1}, $K_X+D$ is $1$ or
$2$-complementary. Assume additionally that $D\in\Msm$ and
$X=\PP^1$. By \ref{m_i} and \ref{=}, $(X,D)$ is exceptional only
in the following cases: $E_6$, $E_7$, $E_8$, $\widetilde D_4$,
$\widetilde E_6$, $\widetilde E_7$, $\widetilde E_8$.
\end{example1}

We discuss two-dimensional generalizations of this fact in
Ch.~\ref{sect-8} and \ref{sect-9}.

\begin{example1}
\label{ex-2-ex}
A log canonical singularity $(X,o)$ is exceptional if and only if
for every boundary $B$ such that $K_X+B$ is lc there exists at
most one divisor $S$ (not necessarily exceptional) such that
$a(S,B)=-1$. We see in Ch.~\ref{s4} that a two-dimensional klt
singularity is exceptional if and only if it is of type $\EE_6$,
$\EE_7$ or $\EE_8$. Note that they are bounded. In contrary,
nonexceptional klt singularities belong to two infinite series
$\BA_n$ and $\DD_n$. Refer to \cite{Ishii}, \cite{MP}, \cite{IP}
for generalizations of this observation.

An isolated log canonical nonklt singularity $(X,o)$ is
exceptional if and only if there is exactly one divisor with
discrepancy $a(\cdot,0)=-1$. Under the assumption that $X$ is
Gorenstein such singularities are called \textit{simple
elliptic}\index{simple elliptic singularity} in dimension two and
\textit{simple $K3$}\index{simple $K3$ singularity} in higher
dimensions \cite{IW}.
\end{example1}

\begin{example1}
Let $C\subset\CC^2$ be a curve given by $x^2=y^3$. Then
$K_{\CC^2}+\frac{5}{6}C$ is lc and not klt. Simple computations
show that there exists only one divisor with discrepancy $-1$.
Therefore $\left(\CC^2,\frac{5}{6}C\right)$ is exceptional.
\end{example1}

The following proposition gives a nice relationship between local
and global exceptional objects.

\begin{proposition1}[{\cite[Theorem~5]{Pr1}}]
Let $(X\ni P)$ be a klt singularity and $f\colon (Y,S)\to X$ a plt
blowup of $P$. Then the following are equivalent:
\begin{enumerate}
\item
$(X\ni P)$ is exceptional;
\item
$f(S)=P$ and $(S,\Diff S(0))$ is exceptional;
\item
$(S,\Diff S(0))$ is exceptional.
\end{enumerate}
\end{proposition1}

\begin{proposition1}[\cite{PSh}, see also {\cite{MP}, \cite{IP}, \cite{Pr2}}]
Let $(X/Z\ni o,D)$ be an exceptional log variety of local type.
Then there exists a divisor $S$ of $\KKK(X)$ such that
$a(S,D^+)=-1$ for any nonklt $\QQ$-complement of $K_X+D$ (i.e. $S$
does not depend on the choice of $D^+$).
\end{proposition1}

\begin{corollary1}[{\cite{Pr1}}]
Let $(X\ni P)$ be a $\QQ$-factorial exceptional lc singularity.
Then a plt blowup is unique up to isomorphisms.
\end{corollary1}

\begin{example1}[\cite{MP}, cf. {\cite[Theorem 5]{Pr1}}]
Let $G$ be a finite group acting on $\CC^n$ freely in codimension
one. Then the quotient singularity $\CC^n/G$ is exceptional if and
only if so is the log Fano $(\PP^{n-1}/G,D)$, where $D$ is given
by the formula \eref{eqv-vetvl}. In dimension two there are
exactly three types of exceptional groups: tetrahedral, octahedral
and icosahedral (up to conjugation and scalar multiplication). In
dimension three there are four types of them: $F$, $G$, $I$, $J$,
in the classical notation.
\end{example1}

\chapter{Log del Pezzo surfaces}
\label{sect-5}
In the present chapter we discuss some properties of log del Pezzo
surfaces.
\section{Definitions and examples}
\begin{definition1}
\label{def-log-del-Pezzo}
A projective log surface $(X,D)$ is called
\begin{itemize}
\item
a \textit{log del Pezzo surface}\index{log del Pezzo surface} if
$K_X+D$ is lc and $-(K_X+D)$ is nef and big;
\item
a \textit{log Enriques surface}\index{log Enriques surface} if
$K_X+D$ is lc and $K_X+D\equiv 0$.
\end{itemize}
Higher dimensional analogs of these are called \textit{log
Fano}\index{log Fano variety} and \textit{log
Calabi-Yau}\index{log Calabi-Yau variety} varieties, respectively.
Usually we omit $D$ if $D=0$.
\end{definition1}

If $(X,D)$ is a log del Pezzo, then by Proposition~\ref{E-C} there
exists some $\QQ$-complement $K_X+D^+$ of $K_X+D$. The pair
$(X,D^+)$ is a log Enriques surface.

Examples of log del Pezzo surfaces are the classical ones,
weighted projective planes $\PP(a_1,a_2,a_3)$ with boundary
$D=\sum d_iD_i$, where $D_i:=\{x_i=0\}$ and $\sum d_i<3$,
Hirzebruch surfaces $\FF_n$\index{$\FF_n$} with boundary
$\alpha\Sigma_0$, where $\Sigma_0$\index{$\Sigma_0$} is the
negative section and $(n-2)/n\le\alpha\le 1$.
\par
Let $f\colon(X',D')\to(X,D)$ be a birational log crepant morphism;
that is,
\[
K_{X'}+D'=f^*(K_X+D),\quad\text{with}\quad f_*D'=D.
\]
Then
$(X,D)$ is a log del Pezzo if and only if so is $(X',D')$ (see
\ref{equiv-bir}). Conversely, if $f\colon X'\to X$ is a birational
morphism and $(X',D')$ is a log del Pezzo then so is $(X,f_*D')$.
Many examples can also be obtained by taking finite quotients; see
\ref{finite}.

\begin{example1}
Let $G\subset \PGL_2(\CC)$ be a finite subgroup, $X:=\PP^2/G$ and
$f\colon\PP^2\to X$ the natural projection. As in \ref{finite}, we
define a boundary $D$ on $X$ by the condition
$K_{\PP^2}=f^*(K_X+D)$, where $D=\sum(1-1/r_i)D_i$, all the $D_i$
are images of lines on $\PP^2$, and $r_i$ is the ramification
index over $D_i$. For example, if $G$ is the symmetric group
$\mathfrak{S}_3$, acting on $\PP^2$ by permutations of
coordinates, $X$ is the weighted projective plane
$\PP(1,2,3)=\mt{Proj}\CC[\sigma_1,\sigma_2,\sigma_3]$, where the
$\sigma_i$ are the symmetric functions on coordinates on $\PP^2$.
The divisor $D$ has exactly one component $D_1$ with coefficient
$1/2$, where $D_1$ is determined by the equation
\[
\sigma_1^2\sigma_2^2-4\sigma_2^3-
4\sigma_1^3\sigma_3-27\sigma_3^2+18\sigma_1\sigma_2\sigma_3 =0
\]
(the equation of the discriminant). The surface $X$ has exactly
two singular points which are Du Val of types $A_1$ and $A_2$.
Therefore $X$ is a Gorenstein del Pezzo surface of degree
$K_X^2=6$. The curve $D_1$ is contained in the smooth locus and
has a unique singularity at $(1,1/3,1/27)$, which is a cusp.
\end{example1}

\begin{lemma}
\label{log-del-Pezzo-Pic}
Let $(X,D)$ be a log del Pezzo surface. Assume additionally either
$K_X+D$ is klt or $K_X+D$ is dlt and $-(K_X+D)$ is ample. Then
\begin{enumerate}
\item
$\Pic(X)$ is a finitely generated free abelian group;
\item
the numerical equivalence in $\Pic(X)$ coincides with linear one;
\item
group of classes of Weil divisors $\Weillin(X)$ is finitely
generated.
\end{enumerate}
\end{lemma}
Recall that two-dimensional log terminal singularities are
automatically $\QQ$-factorial.
\begin{outline}
From the exponential sequence and Kawamata-Viehweg vanishing we
have $\Pic(X)\simeq H^2(X,\ZZ)$. Assume that $D\in\Pic(X)$ is
$n$-torsion. Then again Kawamata-Viehweg vanishing and by
Riemann-Roch, $|D|\ne \emptyset$. Therefore $D\sim 0$. (iii)
follows by \cite[Lemma~1.1]{K}.
\end{outline}

\section{Boundedness of log del Pezzos}
It is well known that the degree $K^2_X$ of classical del Pezzo
surfaces is bounded by $9$. This however is not true for log del
Pezzo surfaces. Indeed, we can take $(\FF_n,(1-2/n)\Sigma_0)$,
where $\FF_n$ is the Hirzebruch surface and $\Sigma_0$ is the
negative section. Then $-(K+(1-2/n)\Sigma_0)$ is nef and big. It
is easy to compute that
\[
(K+(1-2/n)\Sigma_0)^2=n+4+4/n,
\]
so it is unbounded. However, results of Alexeev and Nikulin (see
\cite{A}) show that the degree $(K_X+D)^2$ of log del Pezzo
surfaces is bounded by a constant $\Const(\ep)$ if $K_X+D$ is
$\ep$-lt. More precisely we have

\begin{theorem1}[\cite{A}, see also {\cite[Sect.~9]{KeM}}]
\label{Al}
Fix $\ep>0$. Let $(X,D=\sum d_iD_i)$ be a projective log surface
such that $-(K_X+D)$ is nef, $K_X+D$ is $\ep$-lt and $d_i<1-\ep$.
Then the class $\{X\}$ is bounded in the algebraic moduli sense
except for the case when $D=0$ and $K_X\sim 0$.
\end{theorem1}
In the case when $D=0$ and ${-}K_X$ is ample we have more
effective Nikulin's results.
\begin{theorem1}[\cite{N3}]
\label{Nikulin}
Let $X$ be a projective surface with only klt singularities such
that ${-}K_X$ is ample and $X^{\min}\to X$ a minimal resolution.
Then
\[
\rho(X^{\min})<\left\{
\begin{array}{ll}
 3141&\text{if}\ e=2,\\
 5317&\text{if}\ e=3,\\
 17735&\text{if}\ e=4,\\
 \frac{192e(e-3)(6e-27)}{\ep(e)}+1536e^2(e-3)+1820e+69
 &\text{if}\ e\ge 5,\\
\end{array}
\right.
\]
where $e$ is the maximal multiplicity of the singularities of $X$
and $\ep(e)$ is some function of $e$.
\end{theorem1}
Nikulin \cite{N1} obtained also a better bound (linear in $e$)
\[
\rho(X^{\min})\le 352e+1284
\]
in the case when all the discrepancies of $X$ satisfy the
inequalities
\[
a(E,0)\le -1/2\qquad \text{or}\qquad a(E,0)\ge 0.
\]
Proofs of theorems~\ref{Al}, \ref{Nikulin} use weighted graph
technique and Nikulin's diagram method. In the case $\rho(X)=1$
Theorem~\ref{Nikulin} was proved also in \cite[Sect. 9]{KeM} by
using nonnegativity of $\hat c_2(\hat\Omega^1_X(\log D))$ and
Bogomolov type inequality $\hat c_1^2(\hat\Omega^1_X(\log D))\le
3\hat c_2(\hat\Omega^1_X(\log D))$ (see also \cite[Ch. 10]{Ut} and
\cite{Kawamata-bound}). As an easy consequence of this theory we
have the following

\begin{theorem1}[{\cite[9.2]{KeM}}]
\label{Keel-McCern-}
Let $(X,C)$ be a log surface with $\rho(X)=1$ such that $K_X+C$ is
lc, $K_X$ is klt and $C=\sum C_i$ is reduced. Then
\begin{equation*}
\sum_{P\in (X\setminus C)}\frac{m_P-1}{m_P}\le
\chi_{top}(X)-\chi_{top}(C)
\end{equation*}
where $m_P$ is the order of the local fundamental group
$\pi_1(U_P\setminus\{P\})$ ($U_P$ is a sufficiently small
neighborhood of $P$). If $X$ is rational and $p_a(C)=0$, then
\begin{equation}
\label{eq-KeMcC}
\sum_{P\in (X\setminus C)}\frac{m_P-1}{m_P}\left\{
\begin{array}{cc}
\le 3 & \text{if}\ C=0 \\ \le 1 & \text{if}\ \#\{C_i\}=1
\\ =0 & \text{if}\ \#\{C_i\}=2.
\end{array}
\right.
\end{equation}
\end{theorem1}
Using this fact one can easily show the following:

\begin{corollary1}
Let $X$ be a log del Pezzo surface with $\rho(X)=1$ such that
$K_X$ is klt. Then the number of singular points of $X$ is at most
$5$.
\end{corollary1}
\begin{proof}
By Theorem~\ref{Keel-McCern-} the number of singular points is
$\le 6$. Assume that $X$ has exactly six singular points
$P_1,\dots,P_6$. Then by inequality \eqref{eq-KeMcC} we have
$m_{P_1}=\cdots=m_{P_6}=2$. This means that $P_1,\dots,P_6$ are
ordinary double points. In particular, $K_X$ is Cartier. Applying
Noether's formula to the minimal resolution $\tilde X$ of $X$, we
obtain $K_X^2=K_{\tilde X}^2=10-\rho(\tilde X)=10-1-6=3$. Let
$\tilde L\subset \tilde X$ be a $-1$-curve and $L\subset X$ its
image. Then ${-}K_X\cdot L={-}K_{\tilde X}\cdot \tilde L=1$. Since
$\rho(X)=1$, we have $L\equiv - \frac13K_X$, so $L^2=\frac13$. On
the other hand, $2L$ is Cartier, a contradiction.
\end{proof}

\section{On the existence of regular complements}

\begin{proposition1}[Inductive Theorem, Weak Form \cite{Sh1}]
\label{Inductive-Theorem-Weak-Form}
Let $(X,D)$ be a log del Pezzo surface. If $K_X+D$ is not klt,
then there exists a regular complement of $K_X+D$ (i.e.
$n$-complement with $n\in\RRR_2$). Moreover, if $K_X+D$ is not $1$
or $2$-complementary, then there is at most one divisor of
$\KKK(X)$ with discrepancy $a(\cdot,D)=-1$.
\end{proposition1}

In \cite{KeM} such a log divisor $K_X+D$ was called a
\textit{tiger}. This is a sort of antithesis to Reid's general
elephant (see \ref{elephant}).
\begin{proof}
Replacing $(X,D)$ with a log terminal modification, we may assume
that $K_X+D$ is dlt. Then $\down{D}\ne 0$. In this situation we
can apply Proposition~\ref{prodoljenie} and \ref{1}. The last
statement follows by Connectedness Lemma.
\end{proof}

\begin{corollary1}
\label{down-compl}
Let $(X,D)$ be a log del Pezzo surface with $D\in\Mm$. Write, as
usual, $D=C+B$, where $C:=\down{D}$, $B:=\fr{D}$. Assume that
$|{-}nK_X-nC-\down{(n+1)B}|\ne \emptyset$ for some $n\in\RRR_2$.
Then $K_X+D$ has a regular complement.
\end{corollary1}
Note that the inverse implication follows by \eref{compl-prop}.
\begin{definition1}[{\cite[18.2]{Ut}}]
\label{def-maximally-lc}
A log divisor $K_X+D+\sum b_iB_i$ is said to be \textit{maximally
log canonical}\index{maximally log canonical} if $K_X+D+\sum
b_i'B_i$ is not lc, where $b_i'\ge b_i$ with inequality holding
for at least one index $i$. Note that this definition depends on
the decomposition $D+\sum b_iB_i$, not only on the sum $D+\sum
b_iB_i$.
\end{definition1}
\begin{proof}
If $K_X+D$ is not klt, the assertion follows by
Proposition~\ref{Inductive-Theorem-Weak-Form}. Thus we may assume
that $K_X+D$ is klt (in particular, $C=0$). Let
\[
\ov{D}\in\left|{-}nK_X-\down{(n+1)B}\right|
\]
As in \eref{compl-prop} and \eref{compl-prop1} put
\[
D':=\frac1{n}\left(\down{(n+1)B}+\ov{D}\right).
\]
Note that $D'\ge D$ (because $D\in\Mm$, see
\ref{important_compl}). If $K_X+D'$ is lc, then this is a regular
complement. Assume that $K_X+D'$ is not lc. Take $\alpha$ so that
$K_X+D+\alpha(D'-D)$ is maximally lc. It is clear that
$0<\alpha<1$ and $-(K_X+D+\alpha(D'-D))$ is nef and big. Now we
can apply Proposition~\ref{Inductive-Theorem-Weak-Form} to
$K_X+D+\alpha(D'-D)$ to get the desired regular complement.
\end{proof}

\begin{corollary1}[{\cite{Sh1}}, cf. {Corollary~\ref{cor-=4}}]
\label{>4}
Let $(X,D)$ be a log del Pezzo surface. If $(K_X+D)^2>4$, then it
is nonexceptional. In this case, there exists a regular complement
of $K_X+D$. Moreover, there exists such a complement which is not
klt.
\end{corollary1}
\begin{proof}
Riemann-Roch gives that $\dim|-n(K_X+D)|$ is sufficiently large,
where $n$ is divisible enough and $n\gg 0$. Then standard
arguments show that $K_X+D+\alpha H$ is not klt for some $H\in
|-n(K_X+D)|$ and $\alpha<1/n$ (see e.g. \cite[Lemma 6.1]{Ko} or
the proof of Corollary~\ref{cor-=4}).
\end{proof}

\begin{corollary1}
Let $X$ be a log del Pezzo surface. Assume that $K_X^2>4$. Then
$|{-}nK_X|\ne\emptyset$ for some $n\in\RRR_2$.
\end{corollary1}

\begin{exercise1}
Let $G\subset \PGL_2(\CC)$ be a finite subgroup. Then $G$ acts
naturally on $\PP^1\times\PP^1$ so that the action is free in
codimension one. Prove that the quotient $X:=(\PP^1\times\PP^1)/G$
is a log del Pezzo and $K_X$ is $1$, $2$, $3$, $4$, or
$6$-complementary.
\begin{hint}
Apply Proposition~\ref{Inductive-Theorem-Weak-Form} to $K_X+\De$,
where $\De$ is the image of the diagonal.
\end{hint}
\end{exercise1}

\subsection{}
\label{types}
Log del Pezzo surfaces can be classified in terms of complements.
A log del Pezzo (or log Enriques) surface $(X,D)$ is said to be
\textit{regular}\index{regular log surface} if $K_X+D$ is
$r$-complementary for some $r\in\RRR_2$. Shokurov \cite{Sh1}
proposed the following rough classification of them.
\par
Let $(X,D)$ be a log surface having a regular $r$-complement
$K_X+G$ (i.e., with $r\in\RRR_2$). Let $\EuScript{A}(X,G)$ be the
set of divisors with discrepancy $a(\cdot,G)=-1$.

We say that $K_X+G$ is of type
\begin{description}
\item[$\BA^n_m$]
if $r=1$ and $\EuScript{A}(X,G)$ is infinite;
\item[$\EE 1^n_m$]
if $r=1$ and $\EuScript{A}(X,G)$ is finite;
\item[$\DD^n_m$]
if $r=2$ and $\EuScript{A}(X,G)$ is infinite;
\item[$\EE 2^n_m$]
if $r=2$ and $\EuScript{A}(X,G)$ is finite;
\item[$\EE 3^n_m$, $\EE 4^n_m$, $\EE 6^n_m$] if $r=3,4,6$, respectively,
\end{description}
where $n$ is the number of components of $\down G$ and $m$ is the
number of exceptional divisors with $a(\cdot,G)=-1$ on a minimal
log terminal modification. Thus $n+m=0$ if and only if $K_X+G$ is
klt. In cases $\EE 1$-$6^n_m$ we always have $n+m\le 2$. Moreover,
$n+m=2$ only in the dipole case. For example, a weighted
projective plane $\PP(a,b,c)$ has a natural structure of toric
$1$-complement of type $\BA^3_m$. More general, any toric surface
has a complement of type $\BA^n_m$. Note that this division into
cases gives is very rough classification, more delicate invariant
of a nonexceptional log variety is the simplicial topological
space introduced in \cite[Sect. 7]{Sh1}, see also \cite{Ishii}.

\begin{exercise1}[\cite{Sh1}, cf. \ref{conjecture-toric}]
Let $(X,D)$ be a log del Pezzo surface such that $K_X+D$ is dlt
and $-(K_X+D)$ is ample. Prove that $\down{D}$ has at most two
components. Moreover, if $\down{D}$ has exactly two components,
then $K_X+D$ is $1$ or $2$-complementary. If $K_X+D$ is
$1$-complementary, then $(X,\down{D^+})$ is a toric pair (see
\ref{conjecture-toric}). \begin{hint} Use Adjunction and
\ref{prodoljenie}.\end{hint}
\end{exercise1}

\section{Nonrational log del Pezzo surfaces}
Of course we cannot expect to get a reasonable classification of
all log del Pezzos. Below we describe nonrational ones. Results of
Ch.~\ref{last} shows that exceptional log del Pezzos (see
\ref{def_except}) at least in principle can be classified. By
Kawamata-Viehweg vanishing we have
\begin{lemma1}
\label{lemma_non_rat}
Let $(X,D)$ be a log del Pezzo surface. Assume that $K_X+D$ is
dlt. Then $X$ is rational if and only if $H^1(X,\OOO_X)=0$.
Moreover, $X$ is rational if either $K_X+D$ is klt or $-(K_X+D)$
is ample.
\end{lemma1}

\begin{proposition1}
\label{exer2}
Let $(X,D)$ be a log del Pezzo surface such that $K_X+D$ is dlt.
Put $C:=\down{D}$ and $B:=\fr{D}$. Assume that $X$ is nonrational.
Then
\begin{enumerate}
\item\label{exer2_1}
$\rho(X)\ge 2$;

\item\label{exer2_2}
if $\rho(X)=2$, then $X$ is smooth, $X\simeq\PP(\EEE)$, where
$\EEE$ is a rank two vector bundle on an elliptic curve, all
components of $D$ are horizontal and $C=\Sigma_0$ is the negative
section, which does not intersect other components of $D$;

\item\label{exer2_3}
if $\rho(X)\ge 3$, then there exists a contraction with rational
fibers $g\colon X\to X'$ onto a log del Pezzo $(X',D'=C'+B')$ with
$\rho(X')=2$. Moreover, $C'$ is a smooth elliptic curve contained
in the smooth locus of $X$ (and then it is a section of the
composition map $f\colon X\to X'\to Z$ and it does not intersects
other components of $D$);

\item\label{exer2_4}
$H^1(X,\OOO_X)\simeq\CC$.
\end{enumerate}
\end{proposition1}

\begin{proof}
First note that $C\ne 0$, because $K_X+D$ is not klt and $C$ is
connected by \ref{connect}.
\par

The assertion (i) follows by Lemma~\ref{lemma_non_rat}.
\par

To prove (ii) we note that there exists an extremal
$(K_X+D)$-negative contraction $f\colon X\to Z$. By (i), $f$ is
not birational. Hence $Z$ is a curve of genus $g(Z)\ge 1$ and
fibers of $f$ are irreducible. If $C$ is contained in fibers of
$f$, then $-(K_X+D-\ep C)$ is nef and big. By
Lemma~\ref{lemma_non_rat}, $X$ is rational in this case. So we
assume that there is a component $C'\subset C$ such that
$f(C')=Z$. Thus $p_a(C')\ge 1$. Let $F$ be a general fiber. Then
$C'\cdot F\le D\cdot F<2$. It follows that $C'\cdot F=1$. Further,
\begin{multline*}
0\le 2p_a(C')-2+\deg\Diff{C'}(0)=(K_X+C')\cdot C'\\
\le(K_X+D)\cdot C'\le 0.
\end{multline*}
From this we have that $C'$ is smooth elliptic curve, $g(Z)=1$,
$C'$ is contained in the smooth part of $X$ and does not intersect
other components of $D$. In particular, $C'=C$. Since $C$ is the
section, $f\colon X\to Z$ has no multiple fibers. Therefore $X$ is
smooth and $X\simeq\PP(\EEE)$, where $\EEE$ is a rank two vector
bundle on $Z$. From $K_X^2=0$ and $(K_X+C)^2>0$ we have $C^2<0$.

As for (iii), run $(K_X+D)$-MMP $g\colon X\to X'$. At the end we
obtain $(X',D')$ with an extremal contraction $X'\to Z$ as in
(ii). By Lemma~\ref{lemma_non_rat} $C':=g(C)\ne 0$ cannot be
contracted to a point on $X'$. As in (ii), by Adjunction we have
\[
0\ge (K_X+D)\cdot C\ge (K_X+C)\cdot C=2p_a(C)-2+\deg\Diff C(0)\ge
0.
\]
This yields $p_a(C)=p_a(Z)=1$ and $\Diff C(0)=0$. Hence $C$ is a
smooth elliptic curve, $X$ is smooth along $C$ and $C$ does not
intersect other components of $D$.

Finally, (iv) follows by (iii) because $R^1f_*\OOO_X=0$.
\end{proof}

\begin{corollary1}
\label{cor-log-del-Pezzo}
Let $(X,D)$ be a log del Pezzo surface. Then $X$ is rational or is
birationally isomorphic to a ruled surface over an elliptic curve.
\end{corollary1}

\begin{corollary1}
\label{Del-231} Let $(X,D=C+B)$ be a log del Pezzo surface.
Assume that $X$ is nonrational. If $\rho(X)=1$, then $X$ is a
generalized cone over an elliptic curve (the contraction of the
negative section on $\PP(\EEE)$; see Proposition~\ref{exer2}).
\end{corollary1}

\begin{exercise1}
\label{cone} \label{Del-cone} Let $(X,D)$ be a log del Pezzo
surface. Assume that $X$ is nonrational. Prove that $K_X+D$ is
$1$-complementary.
\begin{hint}
Apply Proposition~\ref{prodoljenie}.
\end{hint}
\end{exercise1}

\chapter[Birational contractions]{Birational contractions and
two-dimensional log canonical singularities}
\label{sect-6} \label{s4}
\begin{theorem}
\label{local}
Let $(X/Z\ni o,D)$ be a log surface of local type, where $f\colon
X\to Z\ni o$ is a contraction. Assume that $K_X+D$ is lc and
$-(K_X+D)$ is $f$-nef and $f$-big. Then there exists an $1$, $2$,
$3$, $4$, or $6$-complement of $K_X+D$ which is not klt near
$f^{-1}(o)$. Moreover, if there are no nonklt $1$ or
$2$-complements, then $(X/Z\ni o,D)$ is exceptional. These
complements $K_X+D^+$ can be taken so that $a(E,D)=-1$ implies
$a(E,D^+)=-1$ for any divisor $E$ of $\KKK(X)$.
\end{theorem}
\begin{proof}
Let $H$ be an effective Cartier divisor on $Z$ containing $o$ and
let $F:=f^*H$. First we take the $c\in\QQ$ such that $K_X+D+cF$ is
maximally lc (see \ref{def-maximally-lc}) and replace $D$ with
$D+cF$. This gives that $\LCS(X,D)\ne\emptyset$. Next we replace
$(X,D)$ with a log terminal modification. So we may assume that
$K_X+D$ is dlt and $\down{D}\ne 0$. Then
Proposition~\ref{prodoljenie} and Theorem~\ref{1} give us that
there exists a regular complement $K_X+D^+$ of $K_X+D$. By
construction, $\down{D^+}\ge\down{D}$. If $K_X+D$ is not
exceptional, then there exists a $\QQ$-complement $K_X+D'$ of
$K_X+D$ and at least two divisors with discrepancy
$a(\cdot,D')=-1$. Then we can replace $D$ with $D'$. Taking a log
terminal blowup, we obtain that $\down{D}$ is reducible. The rest
follows by Theorem~\ref{1}.
\end{proof}

\begin{corollary1}
\label{index-lc-2}
Let $(Z, Q)$ be a lc, but not klt two-dimensional singularity.
Then the index of $(Z,Q)$ is $1$, $2$, $3$, $4$, or $6$.
\end{corollary1}
This fact has three-dimensional generalizations \cite{Ishii}.
\begin{proof}
Apply Theorem~\ref{local} to $f=\mt{id}$ and $K_Z$. We get an
$n$-complement $K_Z+D$ with $n\in\{1,2,3,4,6 \}$. Then $K_Z+D$ is
lc and $n(K_Z+D)\sim 0$. But if $D\ne 0$, $K_Z$ is klt (because
$Q\in\Supp D$).
\end{proof}

\begin{corollary1}
\label{lc-threshold}
Let $(X\ni P)$ be a normal surface germ. Let $D$ be a boundary
such that $D\in\Mm$ and $C$ a reduced Weil divisor on $X$. Assume
that $D$ and $C$ have no common components. Then one of the
following holds:
\begin{enumerate}
\item
$K_X+D+C$ is lc; or
\item
$K_X+D+\alpha C$ is not lc for any $\alpha\ge 6/7$.
\end{enumerate}
\end{corollary1}
Actually, we have more precise result \ref{lc-threshold-1}. See
\cite{Ko1} for three-dimensional generalizations.
\begin{proof}
Assume that $K_X+D+\alpha C$ is lc for some $\alpha\ge 6/7$. By
Theorem~\ref{local} there is a regular complement
$K_X+D^++\alpha^+ C$ near $P$. Since $D\in\Mm$, $D^+\ge D$. By the
definition of complements, $\alpha^+=1$. Hence $K_X+D+C$ is lc.
\end{proof}

Let $(X\ni o)$ be a klt singularity and $D$ an effective Weil
divisor on $X$. Assume that $D$ is $\QQ$-Cartier. The \textit{log
canonical threshold}\index{log canonical threshold} is defined as
follows
\[
c_o(X,D):=\sup\{ c \mid K_X+cD\ \text{is lc} \}.
\]
\index{$c_o(X,D)$}

\begin{corollary1}[cf. {\ref{lc-thresholds-3}}]
\label{lc-threshold-1}
Let $(X\ni P)$ be a normal klt surface germ. Let $D$ be a reduced
Weil divisor on $X$. Assume that $c_P(X,D)\ge 2/3$, then
\[
c_P(X,D)\in\mathcal{S}:=\left\{\frac23,\frac7{10},
\frac34,\frac56,1\right\}.
\]
\end{corollary1}

\begin{proof}
Put $c:= c_P(X,D)$. Assume that $2/3<c<1$. Clearly, $D$ is
reduced. Let $f\colon (Y,C)\to X$ be an inductive blow up of
$(X,cD)$ (see \ref{def_plt_blowup}). Then we can write
$f^*(K_X+cD)=K_Y+C+cD_Y$, where $D_Y$ is the proper transform of
$D$. If $K_Y+C+cD_Y$ is not plt, then by Theorem~\ref{local},
$K_X+cD$ is $1$ or $2$-complementary. Since $c\ge 2/3$, this gives
us that $(cD)^+\ge D$, $K_X+D$ is lc and $c=1$. Hence, we may
assume that $K_Y+C+cD_Y$ is plt. By \ref{exer-odna-comp}, $C$
intersects $\Supp D_Y$ transversally and
\begin{multline}
\label{diff-plt-threshold}
\Diff C(cD)=\sum_{i=1}^s\frac{n_i-1+c}{n_i}P_i+\sum_{j=1}^q
\frac{r_j-1}{r_j}Q_j,\\ \text{where}\quad
\{P_1,\dots,P_s\}=C\cap\Supp D,\quad n_i,r_j\in\NN,\\
\{Q_1,\dots,Q_q\}=\Sing Y\setminus \Supp D.
\end{multline}
Since $C\simeq\PP^1$, $\deg\Diff C(cD)=2$. If $s\ge 3$, then $s=3$
and in \eref{diff-plt-threshold}, $\frac{n_i-1+c}{n_i}=\frac23$
for $i=1,2,3$, a contradiction with $c>2/3$. Assume that $s=2$.
Then in \eref{diff-plt-threshold} we have
$2>\sum\frac{n_i-1+c}{n_i}>\frac43$ and $0<\sum
\frac{r_j-1}{r_j}<2/3$. Hence $\sum \frac{r_j-1}{r_j}=\frac12$ and
$\sum_{i=1}^2\frac{n_i-1+c}{n_i}=\frac32$. This yields
\[
\frac12=\frac{1-c}{n_1}+\frac{1-c}{n_2}<
\frac1{3n_1}+\frac1{3n_2}.
\]
Therefore $n_1=n_2=1$ and $c=3/4$. Finally, assume that $s=1$.
Similarly, in \eref{diff-plt-threshold} we have $1<\sum_{j=1}^q
\frac{r_j-1}{r_j}<4/3$. From this $q=2$ and up to permutations
$(r_1,r_2)$ is one of the following: $(2,3)$, $(2,4)$, $(2,5)$.
Thus $\frac{n_1-1+c}{n_1}=\frac56,\frac34$, or $\frac7{10}$. In
all cases $n_1=1$, so $c\in\mathcal{S}$.
\end{proof}

\section[Two-dimensional log canonical
singularities]
{Classification of two-dimensional log canonical
singularities}
\label{cll}
Two-dimensional log terminal singularities (=quotient
singularities) were classified for the first time by Brieskorn
\cite{Brieskorn} (see also \cite{Iliev}, \cite[ch. 3]{Ut}). We
reprove this classification in terms of plt blowups. It is
expected that this method has higher-dimensional generalizations,
cf. \cite{Sh1}, \cite{Pr1}. Recall that two-dimensional log
terminal singularities are exactly quotient singularities (see
\cite{K}).
\par
Let $(Z, Q)$ be a two-dimensional klt singularity. If $K_Z$ is
$1$-complementary, then by \ref{kawamata_2} $(Z, Q)$ is
analytically isomorphic to a cyclic quotient singularity. Assume
further that $K_Z$ is not $1$-complementary. By Lemma~\ref{plt}
there exists a plt blowup $f\colon(X,C)\to Z$ (where $C$ is the
exceptional divisor of $f$). Further, we classify these blowups
and propose the method to construct the minimal resolution. This
method also allows us to describe klt singularities as quotients
(see Proposition~\ref{m4}).

\begin{lemma1}
\label{mMmM}
Let $f\colon X\to Z\ni o$ be a plt blowup of a surface klt
singularity and $C$ the (irreducible) exceptional divisor. Then
\begin{enumerate}
\item
$X$ has at most three singular points on $C$;
\item
near each singular point the pair $C\subset X$ is analytically
isomorphic to $(\{x=0\}\subset\CC^2)/\cyc{m_i}(1, a_i)$, where
$\gcd(a_i, m_i)=1$;
\item
if $X$ has one or two singular points on $C$, then $K_X+C$ is
$1$-complementary;
\item
if $X$ has three singular points on $C$, then $(m_1, m_2,
m_3)=(2,2,m)$, $(2,3,3)$, $(2,3,4)$ or $(2,3,5)$ and $K_X+C$ is
respectively $2$, $3$, $4$, or $6$-complementary in these cases.
\end{enumerate}
\end{lemma1}
\begin{proof}
By Proposition~\ref{kawamata_1} we get that all singular points
$P_1,\dots,P_r\in X$ are cyclic quotients:
\[
(X\supset C\ni P_i) \simeq(\CC^2\supset \{x=0 \}\ni 0)/
\cyc{m_i}(1, a_i),\quad \gcd(a_i, m_i)=1,
\]
where the action of $\cyc{m_i}$ on $\CC^2$ is free outside of $
0$. Therefore $C\simeq\PP^1$, $\Diff C(0)=\sum(1-1/m_i)P_i$, where
$K_C+\Diff C(0)=(K_X+C)|_C$ is negative on $C$. From this it is
easy to see that for $(m_1, \dots,m_r)$ there are only the
possibilities $(m)$, $(m_1, m_2)$, $(2,2, m)$, $(2,3,3)$,
$(2,3,4)$ and $(2,3,5)$. Since $-(K_X+C)$ is $f$-ample,
$n$-complements for $K_C+\Diff C(0)$ can be extended to
$n$-complements of $K_X+C$. By \ref{1} we have the desired
$n$-complements. This proves (iii) and (iv).
\end{proof}

By our assumptions, $K_Z$ is not $1$-complementary and by
Lemma~\ref{mMmM} we have exactly three singular points on $C$.
Consider now the minimal resolution $g\colon Y\to X$ and put
$h:=f\compos g\colon Y\to Z$. Then on this resolution $C$
corresponds to some curve, say $C'$, and the singular points
$P_i$, $i=1,2,3$ correspond to ``tails'' meeting $C'$ and
consisting of smooth rational curves (see Fig.~\ref{figure1})
\begin{figure}
\vspace{3cm}
\begin{picture}(120,80)(80,0)
\unitlength=0.8mm \put(56,70){$-p$} \put(60,60){\circle{6}}
\put(63,60){\line(1,0){32}} \put(57,60){\line(-1,0){32}}
\put(60,57){\line(0,-1){22}} \put(110,60){\oval(30,10)}
\put(105,58){$\CC^2/\cyc{m_3}$} \put(10,60){\oval(30,10)}
\put(5,58){$\CC^2/\cyc{m_2}$} \put(60,20){\oval(10,30)}
\put(69,20){$\CC^2/\cyc{m_1}$}
\end{picture}
\caption{}
\label{figure1}
\end{figure}
\par
Here each oval is a chain for the minimal resolution of
$\CC^2/\cyc{m}$ (see \ref{frac}). Thus for $m\le 5$ it is one of
the following:
\begin{equation}
\label{uchit}
\begin{array}{lcl}
m=2&\qquad&\wcir{-2}\\&&\\ m=3&&\left\{\begin{array}{c}
\wcir{-3}\\ \wcir{-2}\lin\wcir{-2}\\
\end{array}\right.\\
&&\\ m=4&&\left\{\begin{array}{c} \wcir{-4}\\
\wcir{-2}\lin\wcir{-2}\lin\wcir{-2}\\
\end{array}\right.\\
&&\\ m=5&&\left\{\begin{array}{c} \wcir{-5}\\
\wcir{-2}\lin\wcir{-3}\\ \wcir{-2}\lin\wcir{-2}\lin\wcir{-2}
\lin\wcir{-2}\\
\end{array}\right.\\
\end{array}
\end{equation}

\begin{figure}[!ht]
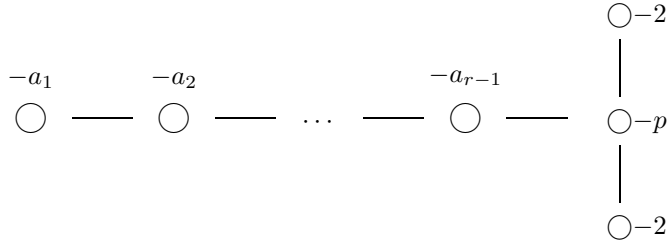

\vspace{2cm} $$
\begin{array}{ccccccccc}
 \wcir{-a_1}&\lin&\wcir{-a_2}&\lin&
\cdots&\lin&\wcir{-a_{r-1}}&\lin&
\begin{array}{c}
 \bigcirc\lefteqn{{-2}}\\
 \vlin\\
 \bigcirc\lefteqn{{-p}}\\
 \vlin\\
 \bigcirc\lefteqn{{-2}}\\
\end{array}
\end{array}
$$ \caption{Case $(2,2,m)$}
\label{figure2}
\end{figure}

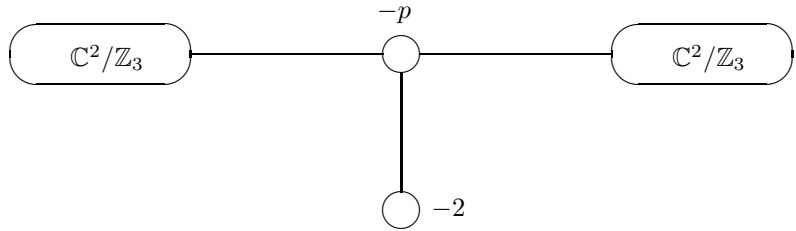
\begin{figure}[!ht]
\vspace{2.5cm} \qquad\qquad
\begin{picture}(120,50)(70,0)
\unitlength=0.8mm \put(56,36){${-p}$} \put(60,30){\circle{6}}
\put(63,30){\line(1,0){32}} \put(57,30){\line(-1,0){32}}
\put(60,27){\line(0,-1){20}} \put(110,30){\oval(30,10)}
\put(105,28){$\CC^2/\cyc{3}$} \put(10,30){\oval(30,10)}
\put(5,28){$\CC^2/\cyc{3}$} \put(60,4){\circle{6}}
\put(65,3){${-2}$}
\end{picture}
\caption{Case $(2,3,3)$}
\label{figure3}
\end{figure}

\begin{figure}[!ht]
\vspace{2.5cm} \qquad\qquad
\begin{picture}(120,50)(70,0)
\unitlength=0.8mm \put(56,36){${-p}$} \put(60,30){\circle{6}}
\put(63,30){\line(1,0){32}} \put(57,30){\line(-1,0){32}}
\put(60,27){\line(0,-1){20}} \put(110,30){\oval(30,10)}
\put(105,28){$\CC^2/\cyc{4}$} \put(10,30){\oval(30,10)}
\put(5,28){$\CC^2/\cyc{3}$} \put(60,4){\circle{6}}
\put(65,3){${-2}$}
\end{picture}
\caption{Case $(2,3,4)$}
\label{figure4}
\end{figure}
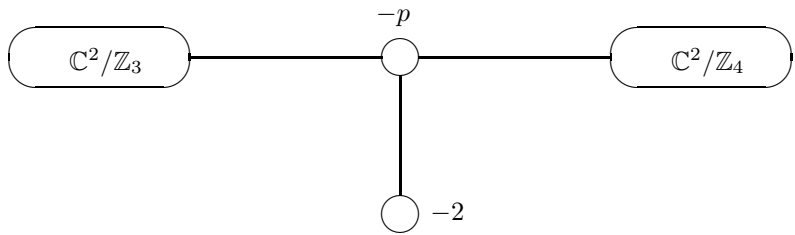

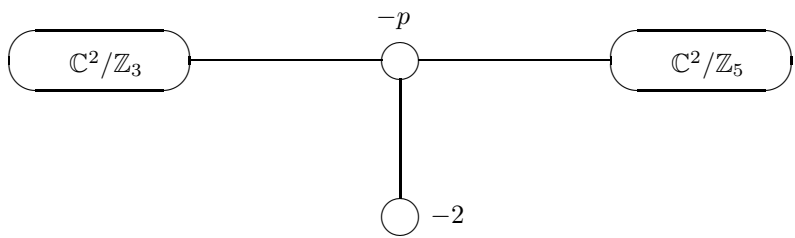
\begin{figure}[!ht]
\vspace{2.5cm} \qquad\qquad
\begin{picture}(120,50)(70,0)
\unitlength=0.8mm \put(56,36){${-p}$} \put(60,30){\circle{6}}
\put(63,30){\line(1,0){32}} \put(57,30){\line(-1,0){32}}
\put(60,27){\line(0,-1){20}} \put(110,30){\oval(30,10)}
\put(105,28){$\CC^2/\cyc{5}$} \put(10,30){\oval(30,10)}
\put(5,28){$\CC^2/\cyc{3}$} \put(60,4){\circle{6}}
\put(65,3){${-2}$}
\end{picture}
\caption{Case $(2,3,5)$}
\label{figure5}
\end{figure}

From Corollary~\ref{m_i} we obtain cases for $(m_1, m_2, m_3)$ in
figures \ref{figure2}-\ref{figure5}. Since the intersection matrix
of exceptional divisors is negative definite, $p\ge 2$. Then
taking into account \eref{uchit} it is easy to get the complete
list of klt singularities (see \cite{Brieskorn}, \cite{Iliev},
\cite[ch. 3]{Ut}).

\par

\begin{theorem}
\label{class-lt}
Let $(Z,Q)$ be a two-dimensional log terminal singularity. Then
one of the following holds:
\begin{enumerate}
\item
$(Z,Q)$ is nonexceptional and then it is either cyclic quotient
(case $\BA_n$ see \ref{frac}) or the dual graph of its minimal
resolution is as in Fig.~\ref{figure2} (case $\DD_n$);
\item
$(Z,Q)$ is exceptional and then the dual graph of its minimal
resolution is as in Fig.~\ref{figure3}-\ref{figure5} (cases
$\EE_6$, $\EE_7$, $\EE_8$).
\end{enumerate}
\end{theorem}

\begin{remark1}
\label{unique}
\begin{enumerate}
\item
Note that our classification uses only the \emph{numerical}
definition of log terminal singularities (by using numerical pull
backs \cite{Sakai}, see \cite{K}).
\item
In all cases of noncyclic quotient singularities (i.e.
Fig.~\ref{figure2}-\ref{figure5}) the plt blowup is unique.
\end{enumerate}
\end{remark1}

\begin{corollary1}
\label{except-bound}
Fix $\ep>0$. There is only a finite number of two-dimensional
exceptional $\ep$-lt singularities (up to analytic isomorphisms).
\end{corollary1}
\begin{proof}
Let $E_0$ be the ``central'' exceptional divisor of the minimal
resolution and $E_1$, $E_2$, $E_3$ exceptional divisors adjacent
to $E_0$. Write $K_Y=h^*K_Z+\sum a_iE_i$. Intersecting both sides
with $E_0$, we obtain
\[
p-2=-pa_0+a_1+a_2+a_3.
\]
This yields
\[
p\ep<p(1+a_0)\le 2,\qquad p<2/\ep.
\]
\end{proof}

\begin{exercise1}
Classify two-dimensional singularities $(Z,Q)$ with klt $K_Z+D$,
where
\[
D=(1-1/m_1)D_1+(1-1/m_2)D_2+(1-1/m_3)D_3,\quad Q\in D_i.
\]
\end{exercise1}

\begin{theorem}[{\cite{Sakai1}, \cite{K}, \cite[Ch.~3]{Ut}}]
\label{class_lc}
Let $(Z,Q)$ be a two-dimensional lc, but not klt singularity, and
let $f\colon X\to Z$ be the minimal resolution. Write
$K_X+D=f^*K_Z$ and put $C:=\down{D}$, $B:=\fr{D}$. Then one of the
following holds:
\begin{description}
\item[$Ell-\widetilde A_n$]
$B=0$, $p_a(C)=1$ and $C$ is either
\begin{description}
\item[$Ell$]
a smooth elliptic curve (\emph{simple elliptic} singularity),
\item[$\widetilde A_n$]
a rational curve with a node, or a wheel of smooth rational curves
(\emph{cusp} singularity);
\end{description}
\item[$\widetilde D_n$]
the dual graph of $f^{-1}(Q)$ is given by
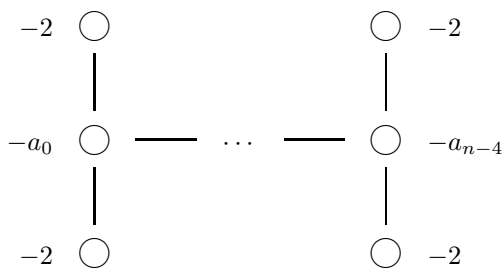
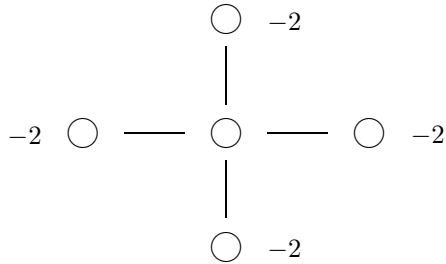
\[
\begin{array}{rcccccl}
{-2}&\wci&&&&\wci&{-2}\\ &\vlin&&&&\vlin&\\
{-a_0}&\wci&\lin&\cdots&\lin&\wci&{-a_{n-4}}\\ &\vlin&&&&\vlin&\\
{-2}&\wci&&&&\wci&{-2}
\end{array}
\]
or (when $n=4$)
\[
\begin{array}{cccclcc}
&&&\wci&-2&&\\&&&\vlin&&&\\ {-2}&\wci&\lin&\wci&\lin&\wci&{-2}\\
&&&\vlin&&\\&&&\wci&-2&&
\end{array}
\]
\item[$Exc$]
the dual graph of $f^{-1}(Q)$ is as in Fig.~\ref{figure1}, where
$(m_1,m_2,m_3)$ is one of the following: $(3,3,3)$ (case
$\widetilde E_6$), $(2,4,4)$ (case $\widetilde E_7$), or $(2,3,6)$
(case $\widetilde E_8$), cf. \ref{=}.
\end{description}
\end{theorem}

\begin{outline}
Similar to the proof of Theorem~\ref{class-lt}. Instead of plt
blowup we can use a minimal log terminal modification $f\colon
X\to Z$. Let $K_X+C=f^*K_Z$ be the crepant pull back. Then $C$ is
a reduced divisor and $K_X+C$ is dlt. If $C$ is reducible, we can
use Lemmas~\ref{wheel} and \ref{chain} below. If $C$ is
irreducible, then either $C$ is a smooth elliptic curve (and $X$
is also smooth) or $C\simeq\PP^1$. In the second case as in the
proof of we Theorem~\ref{class-lt} have cases according to
\ref{=}. We need to check only that $p\ge 2$ in
Fig.~\ref{figure1}. This follows by the fact that the intersection
matrix is negative definite.
\end{outline}

\begin{lemma1}\label{wheel}
Let $(X/Z,D)$ be a log surface such that $K_X+D$ is dlt and
$-(K_X+D)$ is nef over $Z$, $C:=\down{D}\ne\emptyset$ and
$B:=\fr{D}$. Assume that $C$ is compact, connected and not a tree
of smooth rational curves. Then $X$ is smooth along $C$,
$C\cap\Supp B=\emptyset$, $p_a(C)=1$ and $C$ is either a smooth
elliptic curve or a wheel of smooth rational curves.
\end{lemma1}
\begin{proof}
One has
\[
0\le 2p_a(C)-2+\deg\Diff C(B)= (K_X+C+B)\cdot C\le 0.
\]
This yields $p_a(C)=0$ and $\Diff S(B)=0$. In particular, $C$ is
contained in the smooth locus of $X$ and $C\cap\Supp B=\emptyset$.
Moreover, it follows from $(K_X+C+B)\cdot C=0$ that
$(K_X+C+B)\cdot C_i=0$ for any component $C_i\subset C$. Similarly
we can write
\[
0\le 2p_a(C_i)-2+\deg\Diff{C_i}(C-C_i)= (K_X+C)\cdot C_i=0.
\]
If $C=C_i$ is irreducible, then $p_a(C)=1$ and $C$ is a smooth
elliptic curve (because $K_X+C$ is dlt). If $C_i\subsetneq C$,
then $p_a(C_i)=0$, $C_i\simeq\PP^1$ and $\deg\Diff{C_i}(C-C_i)=2$.
Since $K_X+C$ is dlt, $C_i$ intersects $C-C_i$ transversally at
two points. The only possibility is when $C$ is a wheel of smooth
rational curves.
\end{proof}

\begin{remark1}
\label{wheel-rem}
Assuming that $K_X+D$ is only analytically dlt, we have
additionally the case when $C$ is a rational curve with a node.
\end{remark1}

Similar to \ref{wheel} one can prove the following

\begin{lemma1}
\label{chain}
Let $(X/Z,D)$ be a log surface such that $K_X+D$ is dlt and
numerically trivial over $Z$, $C:=\down{D}\ne\emptyset$ and
$B:=\fr{D}$. Assume that $B\in\Mm$, $C$ is compact and it is a
(connected and reducible) tree of smooth rational curves. Then $C$
is a chain. Further, write $C=\sum_{i=1}^r C_i$ where $C_1$, $C_r$
are ends. Then
\begin{enumerate}
\item
$\Diff C(B)=\frac12P^1_1+\frac12P^1_2+ \frac12P^r_1+\frac12P^r_2$,
where $P^1_1,P^1_2\in C_1$, $P^r_1,P^r_2\in C_r$ are smooth points
of $C$;
\item
$C\cap(\Sing X\cup \Supp B)\subset \{P^1_1, P^1_2, P^r_1,
P^r_2\}$;
\item
for each $P^i_j$, $(i,j)\in\{(1,1),(1,2), (r,1), (r,2)\}$ we have
one of the following:
\begin{enumerate}
\item[\textrm{(a)}]
$X$ is smooth at $P^i_j$ and there is exactly one component $B_k$
of $B$ passing through $P^i_j$, in this case $C_i$ intersects
$B_k$ transversally and the coefficient of $B_k$ is equal to
$1/2$;
\item[\textrm{(b)}]
$X$ has at $P^i_j$ Du Val point of type $A_1$ and no components of
$B$ pass through $P^i_j$.
\end{enumerate}
\end{enumerate}
In particular, if $B=0$, then $(X,C)$ looks like that on
Fig.~\ref{figure_log_can}, $X$ is singular only at $P^1_1,
P^1_2\in C_1$, $P^r_1, P^r_2\in C_r$ and these singularities are
Du Val of type $A_1$.
\par
\begin{figure}[htb]
\vspace{5cm} \qquad\qquad
\begin{picture}(160,50)(100,50)
\unitlength 0.7mm \linethickness{0.4pt}
\multiput(9.00,60.67)(0.08,0.74){3}{\line(0,1){0.74}}
\multiput(9.24,62.87)(0.11,0.73){3}{\line(0,1){0.73}}
\multiput(9.57,65.08)(0.10,0.55){4}{\line(0,1){0.55}}
\multiput(9.97,67.27)(0.10,0.44){5}{\line(0,1){0.44}}
\multiput(10.46,69.46)(0.11,0.44){5}{\line(0,1){0.44}}
\multiput(11.04,71.64)(0.11,0.36){6}{\line(0,1){0.36}}
\multiput(11.69,73.81)(0.11,0.31){7}{\line(0,1){0.31}}
\multiput(12.43,75.98)(0.12,0.31){7}{\line(0,1){0.31}}
\multiput(13.25,78.14)(0.11,0.27){8}{\line(0,1){0.27}}
\multiput(14.15,80.29)(0.11,0.24){9}{\line(0,1){0.24}}
\multiput(15.13,82.44)(0.12,0.24){9}{\line(0,1){0.24}}
\multiput(16.20,84.58)(0.11,0.21){10}{\line(0,1){0.21}}
\multiput(17.35,86.71)(0.11,0.19){11}{\line(0,1){0.19}}
\multiput(18.58,88.84)(0.12,0.19){11}{\line(0,1){0.19}}
\multiput(19.90,90.96)(0.12,0.18){12}{\line(0,1){0.18}}
\multiput(21.29,93.07)(0.11,0.16){13}{\line(0,1){0.16}}
\multiput(22.77,95.18)(0.11,0.15){14}{\line(0,1){0.15}}
\multiput(24.34,97.28)(0.12,0.15){14}{\line(0,1){0.15}}
\multiput(25.98,99.37)(0.12,0.14){15}{\line(0,1){0.14}}
\multiput(27.71,101.46)(0.12,0.13){14}{\line(0,1){0.13}}
\multiput(20.33,102.00)(0.20,-0.11){9}{\line(1,0){0.20}}
\multiput(22.15,101.03)(0.20,-0.12){9}{\line(1,0){0.20}}
\multiput(23.93,99.96)(0.17,-0.12){10}{\line(1,0){0.17}}
\multiput(25.68,98.78)(0.16,-0.12){11}{\line(1,0){0.16}}
\multiput(27.40,97.50)(0.14,-0.12){12}{\line(1,0){0.14}}
\multiput(29.08,96.12)(0.13,-0.11){13}{\line(1,0){0.13}}
\multiput(30.72,94.63)(0.12,-0.11){14}{\line(1,0){0.12}}
\multiput(32.33,93.04)(0.11,-0.12){14}{\line(0,-1){0.12}}
\multiput(33.91,91.35)(0.12,-0.14){13}{\line(0,-1){0.14}}
\multiput(35.45,89.55)(0.12,-0.15){13}{\line(0,-1){0.15}}
\multiput(36.96,87.65)(0.11,-0.15){13}{\line(0,-1){0.15}}
\multiput(38.43,85.64)(0.12,-0.18){12}{\line(0,-1){0.18}}
\multiput(39.87,83.53)(0.12,-0.18){12}{\line(0,-1){0.18}}
\multiput(41.27,81.32)(0.11,-0.19){12}{\line(0,-1){0.19}}
\multiput(42.64,79.00)(0.11,-0.20){12}{\line(0,-1){0.20}}
\multiput(43.97,76.58)(0.12,-0.23){11}{\line(0,-1){0.23}}
\multiput(45.27,74.06)(0.12,-0.24){11}{\line(0,-1){0.24}}
\multiput(46.54,71.43)(0.11,-0.25){11}{\line(0,-1){0.25}}
\multiput(47.77,68.70)(0.12,-0.28){10}{\line(0,-1){0.28}}
\multiput(48.96,65.87)(0.12,-0.29){10}{\line(0,-1){0.29}}
\multiput(50.13,62.93)(0.11,-0.30){11}{\line(0,-1){0.30}}
\multiput(139.00,59.67)(-0.12,0.47){7}{\line(0,1){0.47}}
\multiput(138.19,62.99)(-0.11,0.40){8}{\line(0,1){0.40}}
\multiput(137.34,66.18)(-0.11,0.38){8}{\line(0,1){0.38}}
\multiput(136.44,69.24)(-0.12,0.37){8}{\line(0,1){0.37}}
\multiput(135.50,72.17)(-0.11,0.31){9}{\line(0,1){0.31}}
\multiput(134.51,74.98)(-0.12,0.30){9}{\line(0,1){0.30}}
\multiput(133.47,77.65)(-0.11,0.25){10}{\line(0,1){0.25}}
\multiput(132.39,80.19)(-0.11,0.24){10}{\line(0,1){0.24}}
\multiput(131.27,82.61)(-0.12,0.23){10}{\line(0,1){0.23}}
\multiput(130.09,84.89)(-0.11,0.20){11}{\line(0,1){0.20}}
\multiput(128.88,87.04)(-0.11,0.18){11}{\line(0,1){0.18}}
\multiput(127.61,89.07)(-0.12,0.17){11}{\line(0,1){0.17}}
\multiput(126.30,90.96)(-0.11,0.15){12}{\line(0,1){0.15}}
\multiput(124.95,92.72)(-0.12,0.14){12}{\line(0,1){0.14}}
\multiput(123.54,94.36)(-0.11,0.12){13}{\line(0,1){0.12}}
\multiput(122.10,95.86)(-0.12,0.11){12}{\line(-1,0){0.12}}
\multiput(120.61,97.24)(-0.14,0.11){11}{\line(-1,0){0.14}}
\multiput(119.07,98.48)(-0.16,0.11){10}{\line(-1,0){0.16}}
\multiput(117.48,99.60)(-0.18,0.11){9}{\line(-1,0){0.18}}
\multiput(115.85,100.58)(-0.21,0.11){8}{\line(-1,0){0.21}}
\multiput(114.18,101.44)(-0.27,0.11){8}{\line(-1,0){0.27}}
\multiput(128.00,102.00)(-0.19,-0.11){8}{\line(-1,0){0.19}}
\multiput(126.47,101.13)(-0.17,-0.11){9}{\line(-1,0){0.17}}
\multiput(124.96,100.15)(-0.15,-0.11){10}{\line(-1,0){0.15}}
\multiput(123.48,99.04)(-0.13,-0.11){11}{\line(-1,0){0.13}}
\multiput(122.02,97.82)(-0.12,-0.11){12}{\line(-1,0){0.12}}
\multiput(120.60,96.47)(-0.12,-0.12){12}{\line(0,-1){0.12}}
\multiput(119.20,95.01)(-0.11,-0.13){12}{\line(0,-1){0.13}}
\multiput(117.82,93.42)(-0.11,-0.14){12}{\line(0,-1){0.14}}
\multiput(116.47,91.72)(-0.11,-0.15){12}{\line(0,-1){0.15}}
\multiput(115.15,89.89)(-0.12,-0.18){11}{\line(0,-1){0.18}}
\multiput(113.86,87.95)(-0.12,-0.19){11}{\line(0,-1){0.19}}
\multiput(112.59,85.88)(-0.11,-0.20){11}{\line(0,-1){0.20}}
\multiput(111.35,83.70)(-0.11,-0.21){11}{\line(0,-1){0.21}}
\multiput(110.13,81.39)(-0.12,-0.24){10}{\line(0,-1){0.24}}
\multiput(108.94,78.97)(-0.12,-0.25){10}{\line(0,-1){0.25}}
\multiput(107.78,76.42)(-0.11,-0.27){10}{\line(0,-1){0.27}}
\multiput(106.64,73.76)(-0.11,-0.28){10}{\line(0,-1){0.28}}
\multiput(105.53,70.98)(-0.11,-0.29){10}{\line(0,-1){0.29}}
\multiput(104.45,68.07)(-0.12,-0.34){9}{\line(0,-1){0.34}}
\multiput(103.39,65.05)(-0.11,-0.36){15}{\line(0,-1){0.36}}
\multiput(110.33,59.00)(-0.12,0.14){16}{\line(0,1){0.14}}
\multiput(108.45,61.31)(-0.12,0.15){15}{\line(0,1){0.15}}
\multiput(106.66,63.61)(-0.11,0.15){15}{\line(0,1){0.15}}
\multiput(104.96,65.90)(-0.12,0.16){14}{\line(0,1){0.16}}
\multiput(103.34,68.17)(-0.12,0.17){13}{\line(0,1){0.17}}
\multiput(101.81,70.42)(-0.11,0.17){13}{\line(0,1){0.17}}
\multiput(100.37,72.66)(-0.11,0.19){12}{\line(0,1){0.19}}
\multiput(99.01,74.89)(-0.12,0.20){11}{\line(0,1){0.20}}
\multiput(97.74,77.10)(-0.12,0.22){10}{\line(0,1){0.22}}
\multiput(96.56,79.29)(-0.11,0.22){10}{\line(0,1){0.22}}
\multiput(95.47,81.47)(-0.11,0.24){9}{\line(0,1){0.24}}
\multiput(94.46,83.64)(-0.11,0.27){8}{\line(0,1){0.27}}
\multiput(93.55,85.79)(-0.12,0.31){7}{\line(0,1){0.31}}
\multiput(92.72,87.93)(-0.11,0.30){7}{\line(0,1){0.30}}
\multiput(91.97,90.05)(-0.11,0.35){6}{\line(0,1){0.35}}
\multiput(91.32,92.15)(-0.11,0.42){5}{\line(0,1){0.42}}
\multiput(90.75,94.25)(-0.10,0.42){5}{\line(0,1){0.42}}
\multiput(90.27,96.32)(-0.10,0.52){4}{\line(0,1){0.52}}
\multiput(89.87,98.38)(-0.10,0.68){3}{\line(0,1){0.68}}
\multiput(89.57,100.43)(-0.12,1.12){2}{\line(0,1){1.12}}
\multiput(43.67,58.67)(0.11,0.15){16}{\line(0,1){0.15}}
\multiput(45.48,61.14)(0.11,0.16){15}{\line(0,1){0.16}}
\multiput(47.20,63.59)(0.12,0.17){14}{\line(0,1){0.17}}
\multiput(48.83,66.00)(0.12,0.18){13}{\line(0,1){0.18}}
\multiput(50.37,68.39)(0.11,0.18){13}{\line(0,1){0.18}}
\multiput(51.82,70.75)(0.11,0.19){12}{\line(0,1){0.19}}
\multiput(53.18,73.09)(0.12,0.21){11}{\line(0,1){0.21}}
\multiput(54.44,75.39)(0.12,0.23){10}{\line(0,1){0.23}}
\multiput(55.62,77.67)(0.11,0.22){10}{\line(0,1){0.22}}
\multiput(56.71,79.92)(0.11,0.25){9}{\line(0,1){0.25}}
\multiput(57.70,82.14)(0.11,0.27){8}{\line(0,1){0.27}}
\multiput(58.61,84.33)(0.12,0.31){7}{\line(0,1){0.31}}
\multiput(59.42,86.50)(0.10,0.31){7}{\line(0,1){0.31}}
\multiput(60.14,88.63)(0.11,0.35){6}{\line(0,1){0.35}}
\multiput(60.78,90.74)(0.11,0.42){5}{\line(0,1){0.42}}
\multiput(61.32,92.82)(0.11,0.51){4}{\line(0,1){0.51}}
\multiput(61.77,94.88)(0.12,0.67){3}{\line(0,1){0.67}}
\multiput(62.13,96.90)(0.09,0.67){3}{\line(0,1){0.67}}
\multiput(62.40,98.90)(0.09,0.98){2}{\line(0,1){0.98}}
\put(62.58,100.87){\line(0,1){2.13}} \put(69.67,101.67){. . . .}
\put(130.67,83.05){\circle*{1.89}}
\put(137.33,68.38){\circle*{1.89}}
\put(10.67,69.05){\circle*{1.33}}
\put(17.33,86.38){\circle*{1.33}}
\put(134,81.05){\makebox(0,0)[lb]{$Q_1$}}
\put(143,67){\makebox(0,0)[lb]{$Q_2$}}
\put(12,86){\makebox(0,0)[rb]{$P_1$}}
\put(14,66.05){\makebox(0,0)[lb]{$P_2$}}
\put(10,48){\makebox(0,0)[rb]{$C_1$}}
\put(138,48){\makebox(0,0)[lb]{$C_r$}}
\end{picture}
\caption{}
\label{figure_log_can}
\end{figure}
\end{lemma1}

\begin{remark1}
\label{rem-lc-exc}
\begin{enumerate}
\item
We have $nK_Z\sim 0$, where $n=1,2,3,4,6$ in cases $Ell-\widetilde
A_n$, $\widetilde D_n$, $\widetilde E_6$, $\widetilde E_7$, and
$\widetilde E_8$, respectively (see Corollary~\ref{index-lc-2}).
This gives that any two-dimensional lc but not klt singularity is
a quotient of a singularity of type $Ell-\widetilde A_n$ by a
cyclic group of order $1$, $2$, $3$, $4$, or $6$.
\item
The singularity $(Z,Q)$ is exceptional exactly in cases $Ell$,
$\widetilde D_4$, $\widetilde E_6$, $\widetilde E_7$, and
$\widetilde E_8$.
\end{enumerate}
\end{remark1}

Recall that a normal surface singularity $Z\ni Q$ is said to be
\textit{rational}\index{rational singularity} \cite{Ar} (resp.
\textit{elliptic}) \index{elliptic singularity} if
$R^1f_*\OOO_X=0$ (resp. $R^1f_*\OOO_X$ is one-dimensional) for any
resolution $f\colon X\to Z$.

\begin{corollary1}[\cite{K}]
\label{cor-rat-nonrat}
Let $(Z\ni Q)$ be a two-dimensional lc singularity and $f\colon
X\to Z$ its minimal resolution. Write $K_X+D=f^*K_Z$. Then one of
the following holds:
\begin{enumerate}
\item
$\fr{D}\ne 0$ and $Z\ni Q$ is a rational singularity;
\item
$\fr{D}=0$ and $D$ is either a smooth elliptic curve (type $Ell$),
a rational curve with a node or a wheel of smooth rational curves
(type $\widetilde A_n$). In this case, $(Z\ni Q)$ is a Gorenstein
elliptic singularity.
\end{enumerate}
\end{corollary1}

Note that exceptional log canonical singularities are rational
except for the case $Ell$.

\begin{exercise1}
Prove that the following hypersurface singularities are lc but not
klt:
\begin{enumerate}
\item[]
$x^3+y^3+z^3+axyz=0, \quad a^3+27\ne 0$;
\item[]
$x^2+y^4+z^4+ay^2z^2=0, \quad a^2\ne 4$;
\item[]
$x^2+y^3+z^6+ay^2z^2=0, \quad 4a^3+27\ne 0$.
\end{enumerate}
\end{exercise1}

\section[Log terminal singularities as quotients]
{Two-dimensional log terminal singularities as quotients}
Now we discuss the relation between two-dimensional log terminal
singularities and quotient singularities. We use the following
standard notation:
\par
\begin{tabular}{ll}
 $\mathfrak{S}_n$&symmetric group;\\
 $\mathfrak{A}_n$&alternating group;\\
$\mathfrak{D}_n=\left\langle \alpha, \beta\mid \alpha^n=\beta^2=1,
\beta\alpha\beta=\alpha^{-1} \right\rangle$&dihedral group of
order $2n$.\\
\end{tabular}
\index{$\mathfrak{S}_n$}\index{$\mathfrak{A}_n$}\index{$\mathfrak{D}_n$}

\begin{proposition1}
\label{m}
Notation as in Lemma~\ref{mMmM}. Then we have
\begin{enumerate}
\item
if $X$ has three singular points on $C$, $f\colon X\to Z$ is the
quotient of the minimal resolution of the cyclic quotient
singularity $\CC^2/\cyc{r}(1,1)$ by the group $\mathfrak{D}_m$,
$\mathfrak{A}_4$, $\mathfrak{S}_4$ and $\mathfrak{A}_5$, in cases
$(m_1, m_2, m_3)=(2,2, m)$, $(2,3,3)$, $(2,3,4)$ or $(2,3,5)$,
respectively;
\item
if moreover ${-}K_X$ is $f$-ample, then $X$ has at most two
singular points on $C$ (and $K_X+C$ is $1$-complementary).
\end{enumerate}
\end{proposition1}
\begin{proof}
Consider $X$ as a small analytic neighborhood of $C\simeq\PP^1$.
We calculate the fundamental group of $X\setminus \{P_1, P_2, P_3
\}$. Denote by $\Gamma(m_1, m_2, m_3)$ the group generated by
$\alpha_1, \alpha_2, \alpha_3$ with relations
\[
\alpha_1^{m_1}=\alpha_2^{m_2}=\alpha_3^{m_3}=
\alpha_1\alpha_2\alpha_3=1.
\]

\begin{lemma1}[cf. {\cite[0.4.13.3]{Mo}}]
\[
\pi_1(X\setminus \{P_1, P_2, P_3 \})\simeq \Gamma(m_1, m_2, m_3).
\]
\end{lemma1}

\begin{proof}
Let $U_i\subset X$ be a small neighborhood of $P_i$ and $
U_i^o:=U_i\setminus \{P_i \}$. From Theorem~\ref{kawamata_1} we
have $\pi_1(U_i^o)\simeq \cyc{m_i}$. Denote by $\alpha_i$ the
generators of these groups. The set $X\setminus \{P_1, P_2, P_3
\}$ is homotopically equivalent to $\PP^1\setminus \{P_1, P_2, P_3
\}$ glued along $\alpha_1, \alpha_2, \alpha_3$ with sets $ U_1^0$,
$ U_2^0$, $ U_3^0$. Denote loops around $P_i$ (with the
appropriate orientation) also by $\alpha_i$. Then
$\pi_1(\PP^1\setminus \{P_1, P_2, P_3 \})\simeq\langle\alpha_1,
\alpha_2, \alpha_3\mid \alpha_1\alpha_2\alpha_3=1\rangle$. From
the description of points \ref{kawamata_1} it follows also that
the map
\[
\pi_1(C\cap U_i^o)\simeq\ZZ\to\pi_1(U_i^o)\simeq\cyc{m_i}
\]
is surjective. Now the lemma follows by Van Kampen's theorem.
\end{proof}

Now for (i) we notice that the groups $\Gamma(2,2,m)$,
$\Gamma(2,3,3)$, $\Gamma(2,3,4)$ and $\Gamma(2,3,5)$ have finite
quotient groups isomorphic to $\mathfrak{D}_m$, $\mathfrak{A}_4$,
$\mathfrak{S}_4 $ and $\mathfrak{A}_5$, respectively, such that
the images of the elements $\alpha_i$ have orders $m_i$. This
follows from the fact that there exist actions of
$\mathfrak{D}_m$, $\mathfrak{A}_4$, $\mathfrak{S}_4$ and
$\mathfrak{A}_5$ on $\PP^1$ with ramification points of orders
$(m_1, m_2, m_3)$. Then this finite group determines a finite
cover $\widehat X\to X$ unramified outside of $P_1, P_2, P_3$,
where $\widehat X$ is smooth. The Stein factorization gives a
contraction $\widehat X\to\widehat Z$ of an irreducible curve
$\PP^1\simeq\widehat{C}\subset\widehat X$. If $\widehat{C}^2=-r$,
then this contraction is the minimal resolution of the singularity
$\CC^2/\cyc{r}(1,1)$. Finally, if ${-}K_X$ is ample, then so is
${-}K_{\widehat X}$. Thus $r=1$, i.e., $Z\ni o $ is a smooth
point. But the groups $\mathfrak{D}_m$, $\mathfrak{A}_4$,
$\mathfrak{S}_4$ and $\mathfrak{A}_5$ cannot act on $(Z\ni
o)\simeq(\CC^2,0)$ freely in codimension one. This proves (ii).
\end{proof}

\begin{corollary1}[{\cite[0-2-17]{KMM}}]
Any two-dimensional klt singularity is a quotient singularity.
\end{corollary1}

\begin{corollary1}[{\cite{Brieskorn}}]
Let $(Z,Q)$ be a two-dimensional klt singularity. Then
$\pi_1(Z\setminus\{Q\})$ is finite.
\end{corollary1}

\begin{example1}
Let $a, b, m\in\NN$, $\gcd(a, b)=1$. Consider a cyclic quotient
singularity $ o\ni Z=\CC^2/\cyc{m}(a, b)$ (the case $ m=1$ is not
excluded). Any weighted blowup $f\colon X\to Z$ with weights $(a,
b)$ is an extremal contraction with exceptional divisor
$C\simeq\PP^1$. By Lemma~\ref{discr-tor},
$K_X=f^*K_Z+((a+b)/m-1)C$. Hence for $a+b>m$ the divisor ${-}K_X$
is $f$-ample.
\end{example1}

\begin{proposition1}[cf. {Conjecture~\ref{conjecture-toric}}]
\label{m4}
Let $f\colon X\to Z$ be a birational contraction of normal
surfaces. Assume that $f$ contracts an irreducible curve $C$ (i.e.
$\rho(X/Z)=1$) and $K_X+C$ is a plt and $f$-antiample (i.e., $f$
is a plt blowup; see \ref{plt}). Assume also that $X$ has at most
two singular points on $C$. Then $f$ is a weighted blowup.
\end{proposition1}

\begin{remark1}
The condition of the antiampleness of $K_X+C$ is equivalent to the
klt property of $f(C)\in Z$. The condition that $X$ has $\le 2$
singular points is equivalent to that $Z\ni f(C)$ is a cyclic
quotient singularity (or smooth).
\end{remark1}

\begin{proof}
By Proposition~\ref{m}, $K_X+C$ and $K_Z$ are $1$-complementary.
Therefore there are two curves $C_1$, $C_2$ such that
$K_X+C+C_1+C_2$ is lc and linearly trivial over $Z$. Moreover, by
Theorem~\ref{kawamata_2} up to analytic isomorphisms we may assume
that $(Z, f(C_1)+f(C_2))$ is a toric pair. For example, assume
that $X$ has exactly two singular points. Consider the minimal
resolution $\mu\colon X'\to X$ and $f'\colon X'\to Z$ the
composition. It is sufficient to show that the morphism $f'$ is
toric. By \ref{kawamata_2}, in a fiber over $o\in Z$ we have the
following configuration of curves:
\[
\begin{array}{ccccccccccccccc}
\cirmi&\lin&\wci&
&\cdots&\wci&\lin&\bci& \lin&\wci&
&\cdots&\wci&\lin&\cirmi,
\end{array}
\]
where the black vertex corresponds to a fiber (and has
self-intersection number $a\le-1$), white vertices correspond to
exceptional divisors and have self-intersection numbers
$b_i\le-2$, and the vertices $\cirmi$ correspond to the curves
$C_1$, $C_2$. If $a<-1$, $f$ is the minimal resolution of a cyclic
quotient singularity $o\in Z$ and in this case the morphism $f'$
is toric. If $a=-1$, then $f'\colon X'\to Z$ factors through the
minimal resolution $g\colon Y\to Z$ of the singularity $o\in Z$
(which is a toric morphism) and $X'\to Y $ is a composition of
blowups with centers at points of intersections of curves. Such
blowups preserve the action of the two-dimensional torus, hence
$f'$ is a toric morphism.
\end{proof}

\begin{example1}[\cite{Morrison}]
Let $f\colon X\to Z\ni o$ be a $K_X$-negative extremal birational
contraction of surfaces. Assume that $X$ has only Du Val
singularities. Then
\begin{enumerate}
\item
$Z$ is smooth;
\item
$f$ is a weighted blowup (see \ref{w-blow-up}) with weights $(1,
q)$ (and then $X$ contains only one singular point, which is of
type $A_{q-1}$).
\end{enumerate}
\end{example1}

\begin{exercise1}[cf. {\ref{conjecture-toric}}]
\label{exe-tor-1}
Let $f\colon X\to Z\ni o$ be a birational two-dimensional
contraction and $D$ a boundary on $X$ such that $K_X+D$ is lc and
$-(K_X+D)$ is nef over $Z$. Prove that
\begin{equation*}
\rhonum(X/Z)+2\ge\sum d_i,
\end{equation*}
where $\rhonum(X/Z)$ is the rank of the quotient of $\Weil(X)$
modulo numerical equivalence.\index{$\rhonum(X/Z)$} Moreover, the
equality holds only if $(X/Z\ni o,\down{D})$ is a toric pair.
\end{exercise1}

\chapter{Contractions onto curves}
\label{sect-7}
In this chapter we discuss complements on log surfaces over
curves. The main result is Theorem~\ref{bir-nef}. From
Theorem~\ref{local} we have
\begin{corollary}
\label{two}
Let $f\colon X\to Z\ni o$ be a contraction from a normal surface
$X$ onto a smooth curve $Z$. Let $D$ be a boundary on $X$. Assume
that $K_X+D$ is lc and $-(K_X+D)$ is $f$-nef and $f$-big. Then
there exists a nonklt $1$, $2$, $3$, $4$, or $6$-complement of
$K_X+D$ near $f^{-1}(o)$. Moreover, if there are no nonklt $1$ or
$2$-complements of $K_X+D$, then $f\colon X\to Z\ni o$ is
exceptional.
\end{corollary}

Below we give generalization of this result for the case when
$K_X+D\equiv 0$ and classify two-dimensional log conic bundles.

\section{Log conic bundles}
\subsection{Assumptions}
\label{notations}
Let $(X\supset C)$ be a germ of normal surface $X$ with only klt
singularities along a reduced curve $C$, and $(Z\ni o)$ a smooth
curve germ. Let $f\colon(X,C)\to(Z,o)$ be a $K_X$-negative
contraction such that $f^{-1}(o)_{\mt{red}}=C$. Then it is easy to
prove that $p_a(C)=0$ and each irreducible component of $C$ is
isomorphic to $\PP^1$. Everywhere in this paragraph if we do not
specify the opposite, we assume that $C$ is irreducible (or,
equivalently, $\rho(X/Z)=1$, i.e., $f$ is extremal). Let
$X_{\min}\to X$ be the minimal resolution. Since the composition
map $f_{\min}\colon X_{\min}\to Z$ is flat, the fiber of
$f_{\min}^{-1}(o)$ is a tree of rational curves. Therefore it is
possible to define the dual graph of $f_{\min}^{-1}(o)$. We draw
it in the following way: $\bci$ denotes the proper transform of
$C$, while $\wci$ denotes the exceptional curve. We attach the
selfintersection number to the corresponding vertex. By
construction, the proper transform of $C$ is the only $-1$-curve
in $f_0^{-1}(o)$, so we usually omit $-1$ over $\bci$.

\begin{example1}
\label{ex2}\label{ex2_1}
Let $\PP^1\times\CC^1\to\CC^1$ be the natural projection. Consider
the following action of $\cyc{m}$ on $\PP^1_{x,y}\times\CC^1_{u}$:
\[
(x,y;u)\longrightarrow (x,\ep^q y;\ep u),\qquad \ep=\exp 2\pi i/m,
\qquad \gcd(m,q)=1.
\]
Then the morphism $f\colon
X=(\PP^1\times\CC^1)/\cyc{m}\to\CC^1/\cyc{m}$ satisfies the
conditions above. The surface $X$ has exactly two singular points
which are of types $\frac1{m}(1,q)$ and $\frac1{m}(1,-q)$. The
morphism $f$ is toric, so $K_X$ is $1$-complementary. One can
check that the minimal resolution of $X$ has the dual graph
\[
\begin{array}{ccccccccccccc}
\wcir{-b_1}&\lin&\cdots&\lin&\wcir{-b_s}&\lin&
\bci&\lin&\wcir{-a_r}&\lin&\cdots&\lin& \wcir{-a_1},
\end{array}
\]
where $(b_1,\dots,b_s)$ and $(a_r,\dots,a_1)$ are defined by
\eref{frac_chain}.
\end{example1}

\begin{proposition1}[see also {\cite[(11.5.12)]{KeM}}]
\label{new1}
Let $f\colon(X,C)\to(Z,o)$ be a contraction as in \ref{notations},
but not necessarily extremal (i.e., $C$ may be reducible). Assume
that $X$ singular and has only Du Val singularities. Then $X$ is
analytically isomorphic to a surface in
$\PP^2_{x,y,z}\times\CC^1_t$ which is given by one of the
following equations:
\begin{enumerate}
\item
$x^2+y^2+t^nz^2=0$, then the central fiber is a reducible conic
and $X$ has only one singular point, which is of type $A_{n-1}$;
\item
$x^2+ty^2+tz^2=0$, then the central fiber is a nonreduced conic
and $X$ has exactly two singular points, which are of type $A_1$;
\item
$x^2+ty^2+t^2z^2=0$, then the central fiber is a nonreduced conic
and $X$ has only one singular point, which is of type $A_3$;
\item
$x^2+ty^2+t^nz^2=0$, $t\ge 3$ then the central fiber is a
nonreduced conic and $X$ has only one singular point, which is of
type $D_{n+1}$.
\end{enumerate}
\end{proposition1}
\begin{outline}
One can show that the linear system $|{-}K_X|$ is very ample and
determines an embedding $X\subset\PP^2\times Z$. Then $X$ must be
given by the equation $x^2+t^ky^2+t^nz^2=0$.
\end{outline}

\subsection{Construction}
\label{construction}
Notation and assumptions as in \ref{notations}. Let $d$ be the
index of $C$ on $X$, i. e. the smallest positive integer such that
$dC\sim 0$. If $d=1$, then $C$ is a Cartier divisor and $X$ must
be smooth along $C$, because so is $C$. If $d>1$, then there
exists the following commutative diagram:
\[
\begin{CD}
\widehat X@>g>>X\\ @V\widehat{f}VV@VfVV\\ \widehat Z@>h>>Z,\\
\end{CD}
\]
where $\widehat X\to X$ is a cyclic \'etale outside $\Sing X$
cover of degree $d$ defined by $C$ and $\widehat X\to \widehat
Z\to Z$ is the Stein factorization. Then $\widehat{f}\colon
\widehat X\to \widehat Z$ is also a $K_{\widehat X}$-negative
contraction but not necessarily extremal. By construction, the
central fiber $\widehat{C}:=\widehat{f}^{-1}(\widehat{o})$ is a
reducible Cartier divisor. Note that $p_a(\widehat{C})=0$.
Therefore $\widehat X$ is smooth outside $\Sing \widehat{C}$. We
distinguish two cases.
\subsection{Case: $\widehat{C}$ is irreducible}
\label{(A)}
Then $\widehat X$ is smooth and $\widehat X\simeq
\PP^1\times\widehat Z$. Thus $f\colon X\to Z$ is analytically
isomorphic to the contraction from Example~\ref{ex2_1}.
\subsection{Case: $\widehat{C}$ is reducible}
\label{(B)}
Then the group $\ZZ_d$ permutes components of $\widehat{C}$
transitively. Since $p_a(\widehat{C})=0$, this gives that all the
components of $\widehat{C}$ passes through one point, say
$\widehat{P}$, and they do not intersect each other elsewhere. The
surface $\widehat X$ is smooth outside $\widehat{P}$. Note that in
this case $K_X+C$ is not plt, because neither is $K_{\widehat
X}+\widehat{C}$.
\par

\begin{corollary1}
\label{l1}
Notation as in \ref{construction}. Then $X$ has at most two
singular points on $C$.
\end{corollary1}

\begin{proof}
In Case~\ref{(B)} any nontrivial element $a\in\ZZ_d$ have
$\widehat{P}$ as a fixed point. It can have at most one more fixed
point $\widehat{P}_i$ on each component
$\widehat{C}_i\subset\widehat{C}$. Moreover, $\ZZ_d$ permutes
points $\widehat{P}_1,\dots$. Then $X$ can be singular only at
images of $\widehat{P}$ and $\widehat{P}_1,\dots$.
\end{proof}

\subsection{Additional notation}
\label{addd}
In Case~\ref{(B)} we denote $P:=g(\widehat{P})$. If $X$ has two
singular points, let $Q$ be another singular point. To distinguish
exceptional divisors over $P$ and $Q$ in the corresponding Dynkin
graph we reserve the notation $\wci$ for exceptional divisors over
$P$ and $\cirmi$ for exceptional divisors over $Q$.

\begin{corollary1}
\label{l2}
In the above conditions, $K_X+C$ is plt outside of $P$.
\end{corollary1}

\begin{lemma1}
\label{max}
Notation as in \ref{notations}, \ref{construction} and \ref{addd}.
Let $X'\to X$ be a finite \'etale in codimension one cover. Then
there exists the decomposition $\widehat X\to X'\to X$. In
particular, $X'\to X$ is cyclic and the preimage of $P$ on $X'$
consists of one point.
\end{lemma1}
\begin{proof}
Let $X''$ be the normalization of $X'\times_X\widehat X$. Consider
the Stein factorization $X''\to Z'' \to Z$. Then $X''\to Z''$ is
flat and a generically $\PP^1$-bundle. Therefore for the central
fiber $C''$ one has $({-}K_{X''}\cdot C'')=2$, where $C''$ is
reduced and it is the preimage of $\widehat{C}$. On the other
hand,
\[
({-}K_{X''}\cdot C'')=n({-}K_{\widehat X}\cdot \widehat{C})=2n,
\]
where $n$ is the degree of $X''\to \widehat X$. Whence $n=1$,
$X''\simeq\widehat X$. This proves the assertion.
\end{proof}

\begin{lemma1}
\label{m1}
Let $f\colon X\to(Z\ni o)$ be an extremal contraction as in
\ref{notations} (with irreducible $C$). Assume that $K_X+C$ is
plt. Then
\begin{enumerate}
\item
$f\colon X\to(Z\ni o)$ is analytically isomorphic to the
contraction from Example \ref{ex2} (so it is toroidal). In
particular, $X$ has exactly two singular points on $C$ which are
of types $\frac1{m}(1,q)$ and $\frac1{m}(1,-q)$;
\item
$K_X+C$ is $1$-complementary.
\end{enumerate}
\end{lemma1}

\begin{proof}
In the construction \ref{construction} we have Case~\ref{(A)}.
Then
\[
\Diff C(0)=(1-1/d)P_1+(1-1/d)P_2,
\]
where $P_1$, $P_2$ are singular points of $X$ and $d$ is the index
of $C$. By Corollary~\ref{m_i} and by
Proposition~\ref{prodoljenie}, $K_X+C$ is $1$-complementary.
\end{proof}

The following result gives the classification of surface log
terminal contractions of relative dimension one. For applications
to three-dimensional case and generalizations we refer to
\cite{Pr2}, \cite{Pr3}.

\begin{theorem1}[{\cite{Pr3}}]
\label{Log-conic-bundles}
Let $f\colon(X\supset C)\to(Z\ni o)$ be an extremal contraction as
in \ref{notations} (with irreducible $C$). Then $K_X$ is $1$, $2$
or $3$-complementary. Moreover, there are the following cases:
\begin{description}
\item[{\sl Case} $A^*$]
$K_X+C$ is plt, then $K_X+C$ is $1$-complementary and $f$ is
toroidal (see Example~\ref{ex2}, cf.
Conjecture~\ref{conjecture-toric});

\item[{\sl Case} $D^*$]
$K_X+C$ is lc, but not plt, then $K_X+C$ is $2$-complementary and
$f$ is a quotient of a conic bundle of type \textrm{(i)} of
Proposition \ref{new1} by a cyclic group $\cyc{2m}$ which permutes
components of the central fiber and acts on $X$ freely in
codimension one. The minimal resolution of $X$ is
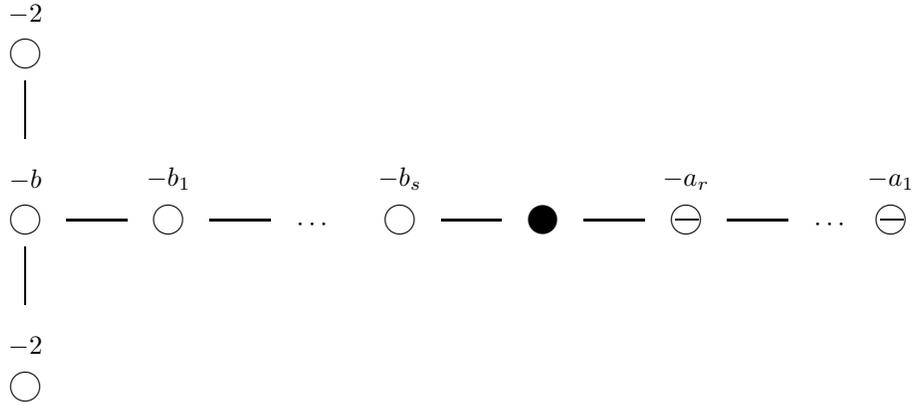
\[
\begin{array}{cccccccccccccc}
\wcir{-2}&&&&&&&&&&&&&\\ \vlin&&&&&&&&&&&&&\\
\wcir{-b}&\lin&\wcir{-b_1}&\lin&\cdots&&\wcir{-b_s}&
\lin&\bci&\lin&\cirmin{-a_r}&\lin&\cdots& \cirmin{-a_1}\\
\vlin&&&&&&&&&&&&&\\ \wcir{-2}&&&&&&&&&&&&&\\
\end{array}
\]
where $s, r\ge 0$ (recall that $X$ can be smooth outside $P$, so
$r=0$ is also possible).
\item[{\sl Case} $A^{**}$]
$K_X$ is $1$-complementary, but $K_X+C$ is not lc. The minimal
resolution of $X$ is
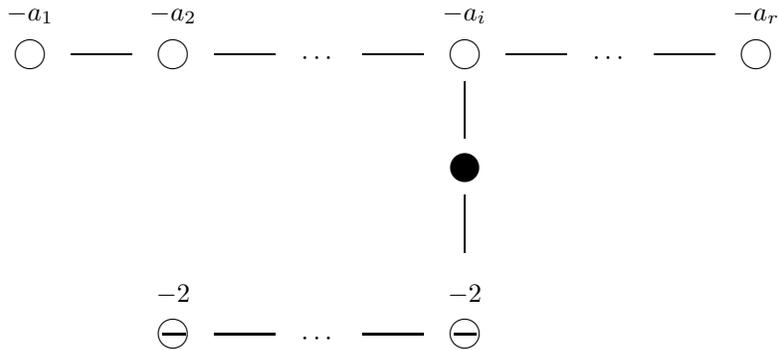
\[
\begin{array}{ccccccccccc}
\wcir{-a_1}&\lin&\wcir{-a_2}&\lin&\cdots&\lin&
\wcir{-a_i}&\lin&\cdots&\lin&\wcir{-a_r}\\
&&&&&&\vlin&&&&\\&&&&&&\bci&&&&\\&&&&&&\vlin&&&&\\
&&\cirmin{-2}&\lin&\cdots&\lin& \cirmin{-2}&&&&\\
\end{array}
\]
where $r\ge 4$, $i\ne 1, r$.
\item[{\sl Case} $D^{**}$]
$K_X$ is $2$-complementary, but not $1$-complementary and $K_X+C$
is not lc. The minimal resolution of $X$ is
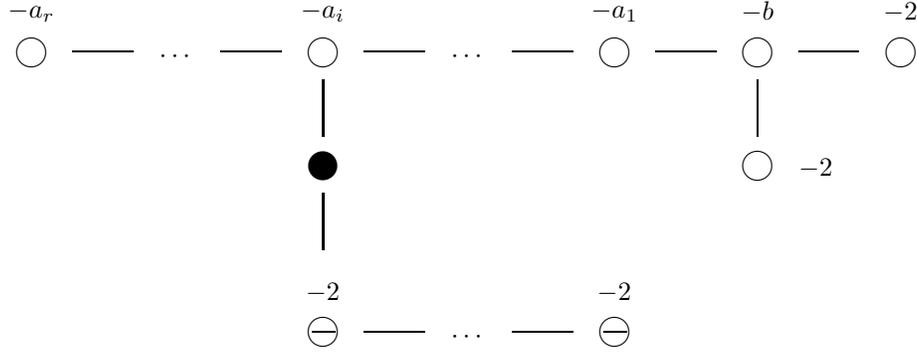
\[
\begin{array}{ccccccccccclc}
\wcir{-a_r}&\lin&\cdots&\lin& \wcir{-a_i}&\lin&\cdots&\lin
&\wcir{-a_1}&\lin&
\wcir{-b}&\lin&\wcir{-2}\\&&&&\vlin&&&&&&\vlin&&\\
&&&&\bci&&&&&&\wci&-2&\\&&&&\vlin&&&&&&&&\\ &&&&\cirmin{-2}&\lin&
\cdots&\lin&\cirmin{-2}&&&&\\
\end{array}
\]
where $r\ge 2$, $i\ne r$.
\item[{\sl Case} $E^*_6$ \textrm{(exceptional case)}]
$K_X$ is $3$-complementary, but not $1$- or $2$-complementary. The
minimal resolution of $X$ is
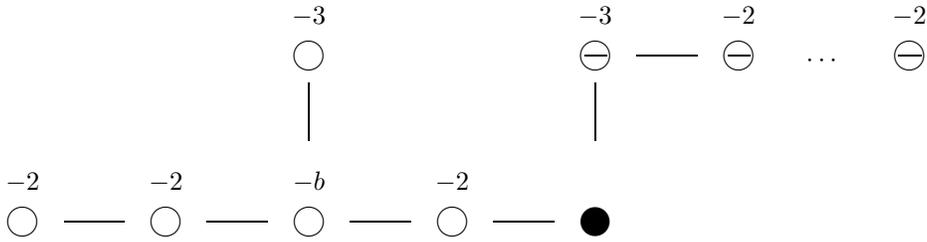
\[
\begin{array}{ccccccccccccccc}
&&&&\wcir{-3}&&&&\cirmin{-3}&\lin&\cirmin{-2}&
&\cdots&&\cirmin{-2}\\ &&&&\vlin&&&&\vlin&&&&&&\\
\wcir{-2}&\lin&\wcir{-2}&\lin&\wcir{-b}
&\lin&\wcir{-2}&\lin&\bcir{}&&&&&&\\
\end{array}
\]
Here the number of $\cirmi$-vertices is $b-2$ (it is possible that
$b=2$ and $Q\in X$ is smooth).
\end{description}
\end{theorem1}

\begin{remark1}
\begin{enumerate}
\item
In the case $D^*$ the canonical divisor $K_X$ can be
$1$-complementary:
\begin{enumerate}
\itemm{a}
if $P\in X$ is Du Val (see \ref{new1} \textrm{(ii)}), or
\itemm{b}
if $s=0$, $a_1=\cdots=a_r=2$, $b=r+2$.
\end{enumerate}
\item
In cases $D^*$, $A^{**}$ and $D^{**}$ there are additional
restrictions on the graph of the minimal resolution. For example,
in the case $A^{**}$ one easily can check that
\[
\left(\sum_{j=1}^{i-1}a_j\right)-(i-1)=\left(\sum_{j=i+1}^{r}a_j\right)-(r-i)
\]
and
\[
a_i=(\text{number}\ \circleddash\text{-vertices})+2.
\]
\end{enumerate}
\end{remark1}

\begin{proof}
If $K_X+C$ is plt, then by Lemma~\ref{m1} we have Case $A^*$. Thus
we may assume that $K_X+C$ is not plt.
\par
We claim that $K_X$ is $1$, $2$ or $3$-complementary. Assume that
$K_X$ is not $1$-complementary. For some $\alpha\le 1$ the log
divisor $K_X+\alpha C$ is lc, but not plt (so, $K_X+\alpha C$ is
maximally lc). Consider a minimal log terminal modification
$\var\colon(\check X,\sum E_i+\alpha\check{C})\to(X,\alpha C)$,
where $\sum E_i$ is the reduced exceptional divisor, $\check{C}$
is the proper transform of $C$ and $K_{\check X}+\sum
E_i+\alpha\check{C}=\var^*(K_X+\alpha C)$ is dlt. As in \ref{plt},
applying the $(K_{\check X}+\sum E_i)$-MMP to $\check X$ at the
last step we obtain the blowup $\sigma\colon\widetilde X\to X$
with irreducible exceptional divisor $E$. Moreover,
$\sigma^*(K_X+\alpha C)=K_{\widetilde X}+E+\alpha\widetilde C$ is
lc, where $\widetilde C$ is the proper transform of $C$ and
$K_{\widetilde X}+E$ is plt and negative over $X$. Since
$K_{\widetilde X}+E+(\alpha-\ep)\widetilde C$ is antiample for
$0<\ep\ll 1$, the curve $\widetilde C$ can be contracted in the
appropriate log MMP over $Z$ and this gives a contraction $(\ov
X,\ov E)\to Z$ with purely log terminal $K_{\ov X}+\ov E$. By
Lemma~\ref{m1} $(\ov X,\ov E)\to Z$ is as in Example \ref{ex2}. If
$K_{\widetilde X}+E$ in nonnegative on $\widetilde C$, then by
Proposition~\ref{bir-prop} we can pull back $1$-complements from
$\ov X$ on $\widetilde X$ and then push-down them on $X$ (see
\ref{bir-prop-a}). Thus we obtain $1$-complement of $K_X$, a
contradiction. From now on we assume that $-(K_{\widetilde X}+E)$
is ample over $Z$. Then by Proposition~\ref{prodoljenie}
complements for $K_E+\Diff E(0)$ can be extended on $\widetilde
X$. According to \ref{m_i}, $\Diff E(0)=\sum_{i=1}^3(1-1/m_i)P_i$,
where for $(m_1,m_2,m_3)$ there are the following possibilities:
\[
(2,2,m),\ (2,3,3),\ (2,3,4),\ (2,3,5).
\]
Further, $\ov X$ has exactly two singular points and these are of
type $\frac1{m}(1,q)$ and $\frac1{m}(1,-q)$, respectively (see
Lemma~\ref{m1}). Since $\widetilde C$ intersects $E$ at only one
point, this point must be singular and there are two more points
with $m_i=m_j$. We get two cases:
\subsection{}
\label{skobka-1}\quad
$(2,2,m)$, $\widetilde C\cap E=\{P_3\}$, there is a
$2$-complement;
\subsection{}
\label{skobka-2}\quad
$(2,3,3)$, $\widetilde C\cap E=\{P_1\}$, there is a
$3$-complement.
\par
This proves the claim.
\par
If $K_X+C$ is lc (but not plt), then in Construction
\ref{construction} $K_{\widehat X}+\widehat{C}$ is also lc but not
plt (see Proposition~\ref{finite-plt}). Since $\widehat{C}$ is a
Cartier divisor, $K_{\widehat X}$ is canonical. Hence
$\widehat{f}$ is as in Proposition \ref{new1}, (i). We get the
case $D^*$.

To prove that note that $\alpha=1$ and $K_{\widetilde
X}+E+\widetilde C$ is lc. Hence $f\colon X\to Z$ is not
exceptional and $K_X$ is $1$- or $2$-complementary by
Corollary~\ref{two}.
\par
Assume that $K_X$ is $1$-complementary, but $K_X+C$ is not lc.
Then there exists a reduced divisor $D$ such that $K_X+D$ is lc
and linearly trivial. By our assumption and by
Propositions~\ref{kawamata_1} and \ref{kawamata_2}, $C\not\subset
D$. Let $P\in X$ be a point of index $>1$. Then $P\in C\cap D$ and
again by Propositions~\ref{kawamata_1} and \ref{kawamata_2} there
are two components $D_1, D_2\subset D$ passing through $P$. But
since $D\cdot L=2$, where $L$ is a generic fiber of $f$,
$D=D_1+D_2$, $P\in D_1\cap D_2$ and $P$ is the only point of index
$>1$ on $X$.
\par
Now assume that $K_X$ is $2$-complementary, but not
$1$-complementary and $K_X+C$ is not lc. Then we are in the case
\ref{skobka-1}. Therefore
\begin{gather*}
(\widetilde X\ni P_1)\simeq(\widetilde X\ni P_2) \simeq(\CC^2,0)/
\cyc{2}(1,1),\\
 (\widetilde X\ni P_3)\simeq(\CC^2,0)/\cyc{m}(1,q),\quad \gcd(m,q)=1.
\end{gather*}
Take the minimal resolution $X_{\min}\to \widetilde X$ of $P_1,
P_2, P_3\in\widetilde X$. Over $P_1$ and $P_2$ we have only single
$-2$-curves and over $P_3$ we have a chain which must intersect
the proper transform of $\widetilde C$, because $\widetilde C$
passes through $P_3$. Since the fiber of $\widetilde X_{\min}\to
Z$ over $o$ is a tree of rational curves, there are no three of
them passing through one point. Whence proper transforms of $E$
and $\widetilde C$ on $\widetilde X_{\min}$ are disjoint.
Moreover, the proper transform of $E$ cannot be a $-1$-curve.
Indeed, otherwise contracting it we get three components of the
fiber over $o\in Z$ passing through one point. It gives that
$\widetilde X_{\min}$ coincides with the minimal resolution
$X_{\min}$ of $X$. Therefore configuration of curves on $X_{\min}$
looks like that in Case $D^{**}$. We have to show only that all
the curves in the down part have selfintersections $-2$. Indeed,
contracting $-1$-curves over $Z$ we obtain a $\PP^1$-bundle. Each
time, we contract a $-1$-curve, we have the configuration of the
same type. If there is a vertex with selfintersection $<-2$, then
at some step we get the configuration
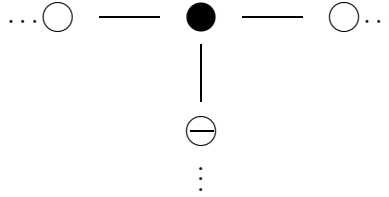
\[
\begin{array}{ccccccc}
\cdots\wci&\lin&\bci&\lin&\wci\cdots\\&&\vlin&&\\&&\cirmi&&\\
&&\vdots&&\\
\end{array}
\]
It is easy to see that this configuration cannot be contracted to
a smooth point over $o\in Z$, because contraction of the central
$-1$-curve gives configuration curves which is not a tree. This
completes Case $D^{**}$.
\par
Case $E^*_6$ is very similar to $D^{**}$. We omit it.
\end{proof}
From Corollary~\ref{except-bound} we have
\begin{corollary1}[cf. {\cite{Pr2}}]
\label{except-bound-conic}
Fix $\ep>0$. There is only a finite number of exceptional (i.e.,
of type $E^*_6$) log conic bundles $f\colon X\to Z$ as in
Theorem~\ref{Log-conic-bundles} with $\ep$-lt $X$.
\end{corollary1}

\begin{exercise1}[cf. {\ref{conjecture-toric}}, {\ref{exe-tor-1}}]
Let $f\colon X\to Z\ni o$ be a contraction from a surface onto a
curve and $D=\sum d_iD_i$ a boundary on $X$ such that $K_X+D$ is
lc and $-(K_X+D)$ is nef over $Z$. Prove that
\begin{equation*}
\rhonum(X/Z)+2\ge\sum d_i.
\end{equation*}
Moreover, the equality holds only if $(X/Z\ni o,\down{D})$ is a
toric pair.
\end{exercise1}

\section{Elliptic fibrations}
As an application of complements we obtain Kodaira's
classification of degenerate of elliptic fibers (see \cite{Sh1}).

\begin{definition1}
\textit{An elliptic fibration}\index{elliptic fibration} is a
contraction from a surface to a curve such that its general fiber
is a smooth elliptic curve. An elliptic fibration $f\colon X\to Z$
is said to be \textit{minimal} \index{elliptic fibration!minimal}
if $X$ is smooth and $K_X\equiv 0$ over $Z$.
\end{definition1}

\begin{remark1}
\label{Du Val}
\textrm{(i)} Note that any elliptic fibration obtained from
minimal one by contracting curves in fibers has only Du Val
singularities.
\par
\textrm{(ii)} Let $K_X+B$ be a $\QQ$-complement on an elliptic
fibration $f\colon X\to Z\ni o$ with $K_X\equiv 0$. Then $B\equiv
0$. By Zariski's lemma, $pB\sim qf^*o$ for some $p,q\in\NN$. In
particular, there exists exactly one complement $K_X+B$ which is
not klt.
\end{remark1}
Recall also that a minimal model is unique up to isomorphisms.

\begin{proposition-definition1}
\label{elliptic-model}
Let $f\colon X\to Z\ni o$ be a minimal elliptic fibration. Then
there exists a birational model $\ov f\colon\ov X\to Z$ such that
$K_{\ov X}+\ov F$ is dlt and numerically trivial near ${\ov
f}^{-1}(o)$, where $\ov F:={\ov f}^{-1}(o)_{\red}$. Such a model
is called a \emph{dlt model}\index{dlt model of elliptic
fibration} of $f$. Moreover, if $K_{\ov X}+\ov F$ is
$n$-complementary, then $K_X$ is $n$-complementary. More
precisely, $\compl'(X)\le \compl(\ov X,\ov F)$. If $(X/Z\ni o)$ is
exceptional, then $\ov F$ is irreducible, $K_{\ov X}+\ov F$ is plt
and a dlt model is unique.
\end{proposition-definition1}

\begin{proof}
First take the maximal $c\in\QQ$ such that $K_X+cf^*o$ is lc. Put
$B:=cf^*o$. Next we consider a minimal log terminal modification
$g\colon Y\to X$ of $(X,B)$ (if $K_X+B$ is dlt, we put $Y=X$).
Thus we can write $g^*(K_X+B)=K_Y+C+B_Y\equiv 0$, where $C$ is
reduced and nonempty, $\down{B_Y}=0$ and $\Supp(C+B_Y)$ is
contained in the fiber over $o$. Run $(K_Y+C+(1+\ep)B_Y)$-MMP over
$Z$:
\begin{equation}
\label{eq-dlt-model}
\begin{array}{c}
 \mbox{
 \begin{picture}(215,80)
 \put(65,70){\vector(-2,-1) {40}}
 \put(25,30){\vector(2,-1){40}}
 \put(110,70){\vector(2,-1){40}}
 \put(150,30){\vector(-2,-1) {40}}
 \put(90,75){\makebox(0,0){$Y$}}
 \put(10,40){\makebox(0,0){$X$}}
 \put(170,40){\makebox(0,0){$\ov X$}}
 \put(90,3){\makebox(0,0){$Z$}}
 \put(40,66){\makebox(0,0)[c]{\scriptsize $g$}}
 \put(40,13){\makebox(0,0)[c]{\scriptsize $f$}}
 \put(135,14){\makebox(0,0)[c]{\scriptsize $\ov f$}}
 \put(135,66){\makebox(0,0)[c]{\scriptsize $\ov g$}}
 \end{picture}}
 \end{array}
\end{equation}
If $B_Y\ne 0$, then $B_Y^2<0$ and we can contract a component of
$B_Y$. At the end we get the situation when $B_Y=0$. Taking $\ov
F:=\ov g(C)$ we see the first part of the proposition. The second
part follows by \ref{bir-prop} and the fact that all contractions
$Y\to \ov X$ are positive with respect to $K+C$.

Finally, assume that $(X/Z\ni o)$ is exceptional. Then by
Remark~\ref{Du Val}, there is exactly one nonklt complement
$K_X+B$ (where $B=cf^*o$). Clearly, $C$ is irreducible in this
case. Contractions $g$ and $\ov g$ are crepant with respect to
$K_Y+C+B_Y$. By Proposition~\ref{crepant-singul} $K_{\ov X}+\ov F$
is plt. Assume that there are two dlt models $(\ov X/Z\ni o,\ov
F)$ and $(\ov X'/Z\ni o,\ov F')$. Consider the diagram
\[
 \mbox{
 \begin{picture}(215,80)
 \put(65,70){\vector(-2,-1) {40}}
 \put(25,30){\vector(2,-1){40}}
 \put(110,70){\vector(2,-1){40}}
 \put(150,30){\vector(-2,-1) {40}}
 \put(90,75){\makebox(0,0){$X_{\min}$}}
 \put(10,40){\makebox(0,0){$X$}}
 \put(170,40){\makebox(0,0){$\ov X$}}
 \put(90,3){\makebox(0,0){$Z$}}
 \put(40,66){\makebox(0,0)[c]{\scriptsize $h$}}
 \put(40,13){\makebox(0,0)[c]{\scriptsize $f$}}
 \put(135,14){\makebox(0,0)[c]{\scriptsize $\ov f$}}
 \put(135,66){\makebox(0,0)[c]{\scriptsize $\ov h$}}
 \end{picture}
 }
\]
where $\ov h\colon X_{\min}\to\ov X$ is the minimal resolution and
$h\colon X_{\min}\to X$ is a composition of contractions of
$-1$-curves. Let $K_{\ov X}+\ov F+\ov D$ be a $\QQ$-complement and
\[
K_{X_{\min}}+F_{\min}+D_{\min}=\ov h^*(K_{\ov X}+\ov F+\ov D)
\]
the crepant pull back, where $F_{\min}$ is the proper transform of
$\ov F$ and $D_{\min}$ is a boundary. Clearly,
\[
-1=a(F_{\min}^{i},F_{\min}+D_{\min})=a(F_{\min}^{i},h_*(F_{\min}+D_{\min}))
\]
for any irreducible component $F_{\min}^{i}$ of $F_{\min}$. Hence
$K_X+h_*(F_{\min}+D_{\min})$ is a nonklt $\QQ$-complement, so
$h_*(F_{\min}+D_{\min})=B$ and $a(\ov F^{i},B)=-1$. Similarly, we
get $a(\ov {F'}^j,B)=-1$. By exceptionality, $\ov F$ and $\ov F'$
are irreducible and $\ov F\approx\ov F'$ (as discrete valuations
of $\KKK(X)$). Then $\ov X\bir\ov X'$ is an isomorphism in
codimension one, hence it is an isomorphism.
\end{proof}

\begin{remark1}
\label{dlt-model-complement}
Let $\ov f\colon (\ov X,\ov F)\to Z\ni o$ be a dlt model of an
elliptic fibration and $K_{\ov X}+\ov F+\ov B$ a $\QQ$-complement.
As in Remark~\ref{Du Val} we have $\Supp \ov B\subset\ov F$, hence
$\ov B=0$.
\end{remark1}

\begin{corollary1}
Under notation of \ref{elliptic-model} the following are
equivalent:
\begin{enumerate}
\item
$(X/Z\ni o)$ is exceptional;
\item
$(\ov X/Z\ni o,\ov F)$ is exceptional;
\item
$K_{\ov X}+\ov F$ is plt.
\end{enumerate}
\end{corollary1}
\begin{proof}
The implication (ii)$\Longrightarrow$(iii) is obvious (because
$\ov F$ is reduced, see \ref{inv}). If $K_{\ov X}+\ov F$ is plt,
then by Remark~\ref{dlt-model-complement} $K_{\ov X}+\ov F$ is the
only nonklt complement and $\ov F$ is the only divisor with $a(\ov
F,\ov F)=-1$. This shows (iii)$\Longrightarrow$(ii). (i)
$\Longrightarrow$ (ii) follows by \ref{elliptic-model}.

Let us prove the implication (ii) $\Longrightarrow$ (i). Assume
that $(X/Z\ni o)$ is nonexceptional. By Remark~\ref{Du Val} there
are two different divisors $E_1$, $E_2$ such that
$a(E_1,B)=a(E_2,B)=-1$. Then in \eref{eq-dlt-model} we have
$a(E_1,C+B_Y)=a(E_2,C+B_Y)=-1$. Since $K_Y+C+B_Y\equiv 0$,
$a(E_1,\ov F)=a(E_2,\ov F)=-1$, i.e., $(\ov X/Z\ni o,\ov F)$ is
nonexceptional.
\end{proof}

Similar to Theorem~\ref{class_lc} we have the following

\begin{proposition1}
\label{elliptic-model-classification}
Let $\ov f\colon \ov X\to Z\ni o$ be dlt model of an elliptic
fibration and $\ov{F}:={\ov f}^{-1}(o)_{\red}$. Then one of the
following holds:
\begin{description}
\item[$Ell$-$\widetilde A_n$]
$p_a(\ov{F})=1$, $\ov X$ is smooth and $\ov{F}$ is either
\begin{description}
\item[$Ell$]
a smooth elliptic curve, or
\item[$\widetilde A_n$]
a wheel of smooth rational curves;
\end{description}
\item[$\widetilde D_n$, $n\ge 5$]
$\ov{F}$ is a chain of smooth rational curves, and it is as in
Lemma~\ref{chain} and Fig.~\ref{figure_log_can} (here $n-3$ is the
number of components of $\ov F$);
\item[$Exc$]
$K_{\ov X}+\ov{F}$ is plt (therefore it is exceptional), then
$\Diff{\ov{F}}(0)=\sum_{i=1}^r (1-1/m_i)$ where for
$(m_1,\dots,m_r)$ there are possibilities as in \ref{=}:
\begin{description}
\item[$\widetilde D_4$]
$(2,2,2,2)$;
\item[$\widetilde E_6$]
$(3,3,3)$;
\item[$\widetilde E_7$]
$(2,4,4)$;
\item[$\widetilde E_8$]
$(2,3,6)$.
\end{description}
\end{description}
\end{proposition1}
\begin{proof}
Follows by \ref{wheel} and \ref{chain}.
\end{proof}

\begin{corollary1}
\label{elliptic-model-classification-compl}
Notation as in Proposition~\ref{elliptic-model-classification}.
Then the index of $K_{\ov X}+\ov{F}$ is equal to $1$, $2$, $3$,
$4$, or $6$, in cases $\widetilde A_n$ (and $Ell$), $\widetilde
D_n$ ($n\ge 4$), $\widetilde E_6$, $\widetilde E_7$ and
$\widetilde E_8$, respectively.
\end{corollary1}
\begin{outline}
Applying Zariski's lemma on the minimal resolution we get $\ov
F\qq 0$. Let $r$ be the index of $\ov F$, i.e., the smallest
positive integer such that $r\ov F\sim 0$. By taking the
corresponding cyclic cover (cf. \ref{l-c-cover})
\[
X':=\operatorname{Spec} \left(\bigoplus_{i=0}^{r-1} \OOO_{\ov
X}(-i\ov F)\right)\to \ov X
\]
we obtain an elliptic fibration $f'\colon X'\to Z'\ni o'$ such
that $F'$ is linearly trivial and log canonical. Since $\ov X$ is
smooth at singular points of $\ov F$, we have that $K_{X'}+F'$ is
dlt (see Theorem~\ref{kawamata_2} or \cite{Sz}). Again by
Theorem~\ref{kawamata_2} $X'$ is smooth along $F'$ (because $F'$
is Cartier). Hence the elliptic fibration $f'\colon X'\to Z'\ni
o'$ must be of type $Ell$ or $\widetilde A_k$. By the canonical
bundle formula, $K_{X'}+F'\sim 0$ (see e.g. \cite[Ch. V, \S
12]{BPV}). Therefore, $m(K_{\ov X}+\ov F)\sim 0$ for some $m$.
Again let $m$ be the index of $K_{\ov X}+\ov F$. Now we consider
the log canonical cover (see \ref{l-c-cover})
\[
X'':=\operatorname{Spec} \left(\bigoplus_{i=0}^{m-1} \OOO_{\ov
X}(-iK_{\ov X}-i\ov F)\right)\to \ov X
\]
As above, $K_{X''}+F''$ is dlt and the elliptic fibration
$f''\colon X''\to Z''\ni o''$ is of type $Ell$ or $\widetilde
A_k$.

If $f''$ is of type $\widetilde A_k$, then the group
$\mt{Gal}(X''/\ov X)$ acts on $F''$ so that the stabilazer of
every singular point is trivial. If $m>1$, then the only
possibility is $m=2$ and $f$ is of type $\widetilde D_{n}$, $n\ge
5$.

Assume that $f''$ is of type $Ell$. Note that $\mt{Gal}(X''/\ov
X)$ contains no subgroups $G$ acting freely on $F''$ (otherwise
the quotient $X''/G\to Z''/G$ is again of type $Ell$). In
particular, $\mt{Gal}(X''/\ov X)\subset \mt{Aut}(F'')$ and this
group contains no translations of the elliptic curve $F''$. It is
well known (see e.g., \cite{Ha}) that, in this situation, the
order of $\mt{Gal}(X''/\ov X)$ can be $2$, $3$, $4$ or $6$.
Moreover, it is easy to see that the ramification indices are such
as in $\widetilde D_4$, $\widetilde E_6$, $\widetilde E_7$, or
$\widetilde E_8$ of \ref{=}.
\end{outline}

\begin{corollary1}
Notation as in Proposition~\ref{elliptic-model-classification}.
Assume that $\ov f$ is exceptional and not of type $Ell$. Then
$\ov f$ is a quotient of a smooth elliptic fibration of type $Ell$
by a cyclic group of order $2$, $3$, $4$, $6$ in cases $\widetilde
D_4$, $\widetilde E_6$, $\widetilde E_7$, and $\widetilde E_8$,
respectively.
\end{corollary1}

\begin{corollary1}
\label{elliptic-model-compl} Let $f\colon X\to Z\ni o$ be a
minimal elliptic fibration. Then there exists a regular complement
of $K_X$.
\end{corollary1}

For convenience we recall Kodaira's classification of singular
elliptic fibers and give a new proof of it using birational
techniques (cf. e.g. \cite[Ch. V, \S 7]{BPV}).

\begin{theorem1}
\label{class}
Let $f\colon X\to Z\ni o$ be a minimal elliptic fibration ($X$ is
smooth) and $F=(f^*o)_{\mt{red}}$, $o\in Z$ the special fiber.
Then there is one of the following possibilities for $F$ (in the
graphs all vertices correspond to $-2$-curves which are components
of $F$):
\begin{description}
\item[$I_b$]
\begin{itemize}
\item[] a smooth elliptic curve ($b=0$);
\item[] a rational curve with one node ($b=1$);
\item[] a wheel of smooth rational curves ($b\ge 2$);
\end{itemize}
\item[${}_mI_b$] multiple $I_b$;
\item[$II$]
a rational curve with a simple cusp;
\item[$III$]
$F=F_1+F_2$ is a pair of smooth rational, tangent each other
curves;
\item[$IV$]
$F=F_1+F_2+F_3$ is a union of three smooth rational curves passing
through one point;
\item[$I^*_b$]
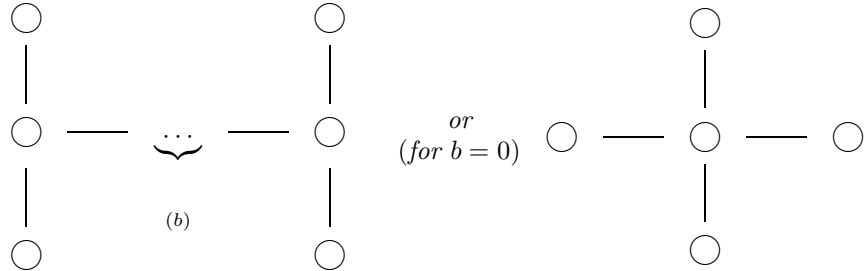
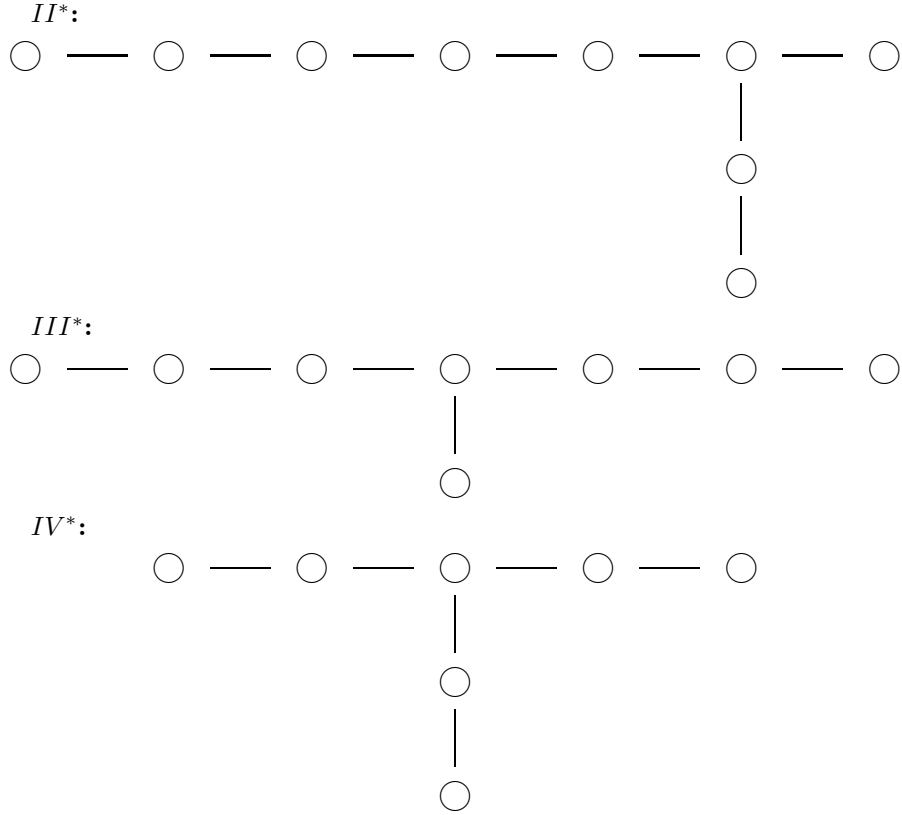
\[
\begin{array}{ccccccc}
&\wci&&&&\wci&\\&\vlin&&&&\vlin&\\
&\wci&\lin&\underbrace{\cdots}&\lin&\wci&\\
&\vlin&&\scriptstyle{(b)}&&\vlin&\\&\wci&&&&\wci&\\
\end{array}
\begin{array}{c}
\text{or}\\
(\text{for $b=0$})\\
\end{array}
\begin{array}{ccccc}
&&\wci&&\\&&\vlin&&\\ \wci&\lin&\wci&\lin&\wci\\&&\vlin&&\\
&&\wci&&\\
\end{array}
\]
\item[$II^*$]
\[
\begin{array}{cccccccccccccc}
\wci&\lin&\wci&\lin&\wci&\lin&
\wci&\lin&\wci&\lin&\wci&\lin&\wci\\
&&&&&&&&&&\vlin&&\\&&&&&&&&&&\wci&&\\
&&&&&&&&&&\vlin&&\\&&&&&&&&&&\wci&&\\
\end{array}
\]
\item[$III^*$]
\[
\begin{array}{ccccccccccccc}
\wci&\lin&\wci&\lin&\wci&\lin&\wci&
\lin&\wci&\lin&\wci&\lin&\wci\\&&&&&&\vlin&&&&&&\\
&&&&&&\wci&&&&&&\\
\end{array}
\]
\item[$IV^*$]
\[
\begin{array}{ccccccccc}
\wci&\lin&\wci&\lin&\wci&\lin&\wci&\lin&\wci\\
&&&&\vlin&&&&\\&&&&\wci&&&&\\&&&&\vlin&&&&\\&&&&\wci&&&&\\
\end{array}
\]
\end{description}
\end{theorem1}

The proof is very similar to that of Theorem~\ref{class_lc}.
\begin{proof}
We are going to apply
Proposition~\ref{elliptic-model-classification}. So we consider a
dlt model $\ov f\colon \ov X\to Z\ni o$ and $\widetilde{h}\colon
\widetilde X\to\ov X$ the minimal resolution of singularities of
$\ov X$. Then we have the following diagram:
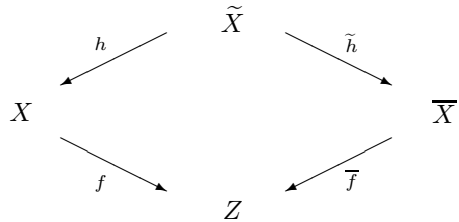
\[
 \mbox{
 \begin{picture}(215,80)
 \put(65,70){\vector(-2,-1) {40}}
 \put(25,30){\vector(2,-1){40}}
 \put(110,70){\vector(2,-1){40}}
 \put(150,30){\vector(-2,-1) {40}}
 \put(90,75){\makebox(0,0){$\widetilde X$}}
 \put(10,40){\makebox(0,0){$X$}}
 \put(170,40){\makebox(0,0){$\ov X$}}
 \put(90,3){\makebox(0,0){$Z$}}
 \put(40,66){\makebox(0,0)[c]{\scriptsize $h$}}
 \put(40,13){\makebox(0,0)[c]{\scriptsize $f$}}
 \put(135,14){\makebox(0,0)[c]{\scriptsize $\ov f$}}
 \put(135,66){\makebox(0,0)[c]{\scriptsize $\widetilde h$}}
 \end{picture}
 }
\]
where $h\colon \widetilde X\to X$ is a sequence of contractions of
$-1$-curves. If $p_a(\ov{C})=1$, then $\ov X=\widetilde X$ and $C$
is a smooth elliptic curve or a wheel of smooth rational curves.
Contracting, if necessary, $-1$-curves we obtain case ${}_mI_b$.
Further, we assume that $p_a(\ov{C})=0$. Then $\ov X$ is singular,
so $\widetilde X\ne\ov X$. Consider the crepant pull back
\[
\widetilde{h}^*(K_{\ov X}+\ov{C})= K_{\widetilde X}+\widetilde
C+\widetilde{B},
\]
where $\widetilde C$ is the proper transform of $\ov{C}$,
$\widetilde{h}_*\widetilde{B}=0$, and $\widetilde{B}\ge 0$. Since
$K_{\ov X}+\ov{C}$, it is easy to see that
$\down{\widetilde{B}}=0$. It is clear also that the set
$\Supp(\ov{C}+\widetilde{B})$ coincides with the fiber over $o$.
By construction, $\Supp \widetilde{B}$ contains no $-1$-curves.
\par
First we consider the case when $\Supp \widetilde C$ also contains
no $-1$-curves. Then $X=\widetilde X$ is exactly the minimal
resolution of $\ov X$. By \ref{Du Val} singular points of $\ov X$
are Du Val. Cases $\widetilde D_n$ ($n\ge 5$), $\widetilde D_4$,
$\widetilde E_6$, $\widetilde E_7$, $\widetilde E_8$ of
Proposition~\ref{elliptic-model-classification} gives cases
$I^*_b$ (with $b\ge 1$), $I^*_0$, $IV^*$, $III^*$, and $II^*$,
respectively. For example, if $\ov C$ is irreducible and there are
exactly three singular points of $\ov X$, then similar to
\ref{cll} the graph of the minimal resolution $\widetilde{h}\colon
\widetilde X\to\ov X$ must be as in Fig.~\ref{figure6}.
\begin{figure}
\vspace{3cm} \qquad\qquad
\begin{picture}(120,80)(100,0)
\unitlength=0.8mm \put(56,70){$-2$} \put(60,60){\circle{6}}
\put(63,60){\line(1,0){32}} \put(57,60){\line(-1,0){32}}
\put(60,57){\line(0,-1){22}} \put(110,60){\oval(30,10)}
\put(105,58){$A_{m_3-1}$} \put(10,60){\oval(30,10)}
\put(5,58){$A_{m_2-1}$} \put(60,20){\oval(10,30)}
\put(69,20){$A_{m_1-1}$}
\end{picture}
\caption{}
\label{figure6}
\end{figure}

By \ref{=} we have the following possibilities for $(m_1, m_2,
m_3)$:
\[
\begin{array}{llll}
\widetilde E_6:\quad& (m_1, m_2, m_3)=(3,3,3)\ &\Longrightarrow\
&\text{case}\ IV^*,\\ \widetilde E_7:\quad& (m_1, m_2,
m_3)=(2,4,4)\ &\Longrightarrow\ &\text{case}\ III^*,\\ \widetilde
E_8:\quad& (m_1, m_2, m_3)=(2,3,6)\ &\Longrightarrow\
&\text{case}\ II^*.\\
\end{array}
\]
Now, we consider the case when $\Supp \widetilde C$ contains a
$-1$-curve. Since $\widetilde{h}\colon \widetilde X\to \ov X$ is a
minimal resolution, all $-1$-curves are contained in $\widetilde
C$, the proper transform of $\ov C$. Using the negative
semidefiniteness for the fiber $\widetilde F\subset \widetilde X$
over $o$ one can show that the dual graph of $\widetilde F$ cannot
contain proper subgraphs of the form
\[
\begin{array}{ccc}
\wcir{-1}&\lin&\wcir{-1}\end{array} \qquad
\text{and}\qquad\begin{array}{ccccc} \wcir{-2}&\lin&\wcir{-1}&
\lin&\wcir{-2}.\end{array}
\]
Suppose that $\ov C$ is irreducible. Then $K_{\ov X}+\ov C$ is plt
and $\widetilde C$ is the only a $-1$-curve. Thus in the case
$\widetilde D_4$ we obtain the dual graph for a fiber of
$\widetilde X\to Z$ as below
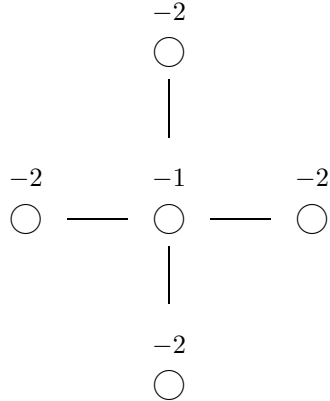
\[
\begin{array}{ccccc}
&&\wcir{-2}&&\\&&\vlin&&\\ \wcir{-2}&\lin&\wcir{-1}
&\lin&\wcir{-2}\\&&\vlin&&\\&&\wcir{-2}&&\\
\end{array}
\]
By the above this is impossible. In other cases we have the dual
graphs as in Fig.~\ref{figure7}.
\begin{figure}
\vspace{3cm} \qquad\qquad
\begin{picture}(120,80)(100,0)
\unitlength=0.8mm \put(56,70){$-1$} \put(60,60){\circle{6}}
\put(63,60){\line(1,0){32}} \put(57,60){\line(-1,0){32}}
\put(60,57){\line(0,-1){22}} \put(110,60){\oval(30,10)}
\put(105,58){$\CC^2/\cyc{m_3}$} \put(10,60){\oval(30,10)}
\put(5,58){$\CC^2/\cyc{m_2}$} \put(60,20){\oval(10,30)}
\put(69,20){$\CC^2/\cyc{m_1}$}
\end{picture}
\caption{}
\label{figure7}
\end{figure}
For $(m_1, m_2, m_3)=(3,3,3)$, $(2,4,4)$ and $(2,3,6)$ we obtain
cases $IV$, $III$ and $II$, respectively. Similarly Case
$\widetilde D_n$, $n\ge 5$ of
Proposition~\ref{elliptic-model-classification} gives Case
$I^*_b$.
\par
Non-simply connected fibers are only of type $I_b$, so only they
can be multiple. This proves the theorem.
\end{proof}
The following table shows correspondence between fibers of minimal
smooth elliptic fibrations and their dlt models:
\medskip
\[
\begin{array}{l||l|l|l|l|l|l|l}
 \ov X&Ell&\widetilde A_n,\, n\ge 1
 &\widetilde D_4&\widetilde D_n,\, n\ge
5&\widetilde E_6&\widetilde E_7&\widetilde E_8\\
 &&&&&&&\\
 \hline
 &&&&&&&\\
 \widetilde X= X&{}_mI_0&{}_mI_n,\, n\ge 2
&I_0^*&I_{n-4}^*&IV^*&III^*&II^*\\ \widetilde X\ne
X&-&{}_mI_b,\,b\le n-1&-&I_b^*,\, b\le n-5&IV&III&II\\
 \compl(\ov X,\ov F)&1&1&2&2&3&4&6\\
\end{array}
\]
\medskip

\begin{theorem}[\cite{Sh1}, cf. {Theorem~\ref{local}}]
\label{bir-nef}
Let $f\colon X\to Z\ni o$ be a contraction from a normal surface
$X$ onto a smooth curve $Z$. Let $D=\sum d_iD_i$ be a boundary on
$X$. Assume that $K_X+D$ is lc and $-(K_X+D)$ is $f$-nef. Then
there exists a regular complement of $K_X+D$. This complement
$K_X+D^+$ can be taken so that $a(E,D)=-1$ implies $a(E,D^+)=-1$
for any divisor $E$ of $\KKK(X)$. Moreover, if there are no $1$,
or $2$-complements, then $(X/Z\ni o,D)$ is exceptional.
\end{theorem}
\begin{proof}
By Corollaries~\ref{two} and \ref{elliptic-model-compl} we may
assume that $K_X+D\equiv 0$ over $Z$ and a general fiber of $f$ is
rational. First, as in the proof of Theorem~\ref{local}, we
replace the boundary $D$ with $D+\alpha f^*o$ so that
$K_X+D+\alpha f^*o$ is maximally lc. Replacing $X$ with its log
terminal modification, we may assume that $X$ is smooth and the
reduced part $C:=\down{D}$ of the boundary is nonempty. Next we
blow up a sufficiently general point on $C:=\down{D}$. We get a
new model such that some component $E$ of $F=f^{-1}(o)$ is
$-1$-curve and it is not contained in $\Supp D$. Moreover,
$E\cap\down{D}$ is a point which is nonsingular for $\Supp D$. Let
$C_1\subset\down{D}$ be a (unique) component passing through
$E\cap\Supp D$. Then the curve $\Supp F\setminus E$ can be
contracted to a point, say $Q$:
\[
f\colon X\stackrel{g}{\longrightarrow} Y\to Z.
\]
The central fiber $g(E)$ of $Y\to Z$ is irreducible. Since
$K_X+D\equiv 0/Y$, the point $Q\in Y$ is lc. Apply
Theorem~\ref{local} to the birational contraction $g\colon X\to
Y$. We get a regular $n$-complement $K_X+D^+$ in a neighborhood of
$g^{-1}(Q)=\Supp(F-E)$. We claim that this complement extends to a
complement in a neighborhood of the whole fiber $F$. We need to
check only that $nD^+\sim{-}nK_X$ in a neighborhood of $F$. But in
our situation the numerical equivalence over $Z$ coincides with
linear one. Therefore the last is equivalent to $D^+\equiv{-}K_X$.
Obviously, both sides have the same intersection numbers with all
components of $F$ different from $E$. For $E$ we have
$1={-}K_X\cdot E$, $E\cdot D^+=E\cdot C_1=1$ (because the
coefficients of $C_1$ in $D$ and $D^+$ are equal to $1$). This
proves the theorem.
\end{proof}

\chapter{Inductive complements}
\label{sect-8}
\section{Examples}
\label{quest-nef}
Roughly speaking the main idea of this chapter is to discuss the
following inductive statement:
\begin{quote}
if a two-dimensional pair $(X/Z,D)$ is lc but not klt and
$-(K_X+D)$ is nef over $Z$, then $K_X+D$ is $1$, $2$, $3$, $4$ or
$6$-complementary.
\end{quote}
It is known that this assertion is true when $-(K_X+D)$ is big
over $Z$ (see Proposition~\ref{Inductive-Theorem-Weak-Form}) and in
the local case. Unfortunately examples~\ref{Atiyah} and
\ref{conterex} below shows that in general, this is false and some
additional assumptions are needed. The main result is the
Inductive Theorem~\ref{Inductive-Theorem} which is a
generalization of \ref{Inductive-Theorem-Weak-Form}.

\begin{example1}[\cite{Sh1}]
\label{Atiyah}
Let $\EEE$ be a indecomposable vector bundle of rank two and
degree $0$ over an elliptic curve $Z$. Then $\EEE$ is a nontrivial
extension
\[
0\longrightarrow \OOO_Z \longrightarrow \EEE \longrightarrow
\OOO_Z \longrightarrow 0
\]
(see e.g., \cite{Ha}). Consider the ruled surface
$X:=\PP_Z(\EEE)$. Let $f\colon X\to Z$ be the projection and $C$ a
section corresponding to the above exact sequence. Then for the
normal bundle of $C$ in $X$ we have $\mathcal{N}_{C/X}=\OOO_C$,
hence $C|_C=0$. In this situation we also have ${-}K_X\sim 2C$
(see \cite{Ha}) and $(K_X+C)|_C=0$. This yields $K_X|_C=0$.

\par
Since $\rho(X)=2$, the Mori cone $\NE(X)$ is generated by two rays
$R_1=\RR_+[F]$, where $F$ is fiber of $X$ and another ray, say
$R$. Since $C^2=0$, $C$ is nef and $C$ generates $R$. In
particular, both ${-}K_X$ and $-(K_X+C)$ are nef and numerically
proportional to $C$.

\par
We claim that $K_X+C$ is not $n$-complementary for any $n$.
Indeed, otherwise we have $L\in |-m(K_X+C)|$ such that $C$ is not
a component of $L$. Then $L\cdot C=0$ and $L\equiv mC$. The
divisor $L-mC$ is trivial on fibers, hence $L-mC=f^*N$ for some
$N\in\Pic(Z)$. Further, $C\cap L=\emptyset$. From this
$(mC-L)|_C\sim 0$ (because $C|_C\sim 0$). Since $f|_C\colon C\to
Z$ is an isomorphism, $f|_C^*N=(mC-L)|_C=0$ gives $N\sim 0$, i.e.,
$L\sim mC$. Then the linear system $|L|$ determines on $X$ a
structure of an elliptic fibration $g\colon X\to\PP^1$ with
multiple fiber $C$. Hence $C|_C$ is an $m$-torsion element in
$\Pic(C)$, a contradiction with $C|_C=0$.
\end{example1}

\begin{example1}
\label{conterex}
Let $X=\PP^1\times\PP^1$. We fix a projection $f\colon X\to
\PP^1$. Let $C$, $H_1$, $H_2$ be different sections of $f$ and
$F_1$, $F_2$, $F_3$ different fibers. Consider the log divisor
$K_X+D$, where
\[
D:=C+\frac1{7}H_1+\frac{6}{7}H_2+\frac12F_1+
\frac23F_2+\frac{5}{6}F_3.
\]
It is clear that $K_X+D$ is lc and numerically trivial. We claim
that there are no regular complements of $K_X+D$. Indeed, assume
that $K_X+D$ has a regular $n$-complement $K_X+D^+$. Then
$K_F+\Diff C(\frac12F_1+\frac23F_2+\frac{5}{6}F_3)$ is also
$n$-complementary. Therefore $n=6$ (see \ref{=}). On the other
hand, by definition we have
\[
D^+\ge C+\frac1{6}H_1+H_2+ \frac12F_1+\frac23F_2+\frac{5}{6}F_3,
\]
a contradiction with $K_X+D^+\equiv 0$.
\end{example1}

The following example shows that under additional assumptions we
can expect some inductive theorems even if $-(K+D)$ is not big.

\begin{example1}
Let $D=C+B$ be a boundary on $X:=\PP^2$ such that $K_X+D$ is lc
and $C:=\down{D}\ne 0$. Assume that $-(K_X+D)$ is nef. Then there
exists a regular complement of $K_X+D$. Indeed, for some
$n\in\RRR_2$ the log divisor $K_C+B|_C$ is $n$-complementary.
Since $H^1(\PP^2,\LLL)=0$ for any invertible sheaf $\LLL$ on
$\PP^2$, $n$-complement (where $n\in\RRR_2$) on $C$ can be
extended to some $\QQ$-divisor $D^+$ on $X$. We write $D^+=C+B^+$.
By Corollary~\ref{inv1}, $K_X+D^+$ is lc near $C$. It is
sufficient to show that $K_X+D^+$ is lc everywhere. But in the
opposite case $K_X+C+\alpha B^+$ is not lc for some $\alpha<1$. By
Connectedness Lemma, $\LCS(X,C+\alpha B^+)$ is connected. This
gives a contradiction.
\end{example1}

\section{Nonrational case}
\label{s-def}
Now we consider the question \ref{quest-nef} for the case when the
surface $X$ is nonrational.
\begin{theorem1}[{\cite{Sh1}}]
\label{Inductive-Theorem-nonrat}
Let $X$ be a normal projective nonrational surface and $D$ a
boundary on $X$ such that $K_X+D$ is lc and $-(K_X+D)$ is nef.
Assume that
\begin{enumerate}
\item
$K_X+D$ is not klt;
\item
there is a boundary $D'$ such that $D'\ge D$ and $K_X+D'$ is lc
and numerically trivial.
\end{enumerate}
Then $K_X+D$ is $n$-complementary for $n\in\RRR_2$.
\end{theorem1}

\begin{proof}
By taking a log terminal modification we may assume that $X$ is
smooth, $K_X+D$ is dlt and $\down{D}\ne 0$ (see \ref{blow-up}).
Since $D\ne 0$, $\kappa(X)=-\infty$. So there is a morphism
$f\colon X\to Z$ onto a curve $Z$ of genus $g\ge 1$.

\begin{lemma1}
\label{lemma1}
Let $f\colon X\to Z$ be a contraction from a projective surface
onto a curve of genus $g\ge 1$. Assume that $K_X+D$ is lc and
$-(K_X+D)$ is nef. Furthermore, assume that the general fiber of
$f$ is a smooth rational curve. Then no components of $\Supp D$
are contained in fibers.
\end{lemma1}
\begin{proof}
Let $L$ be such a component. Replace $(X,D)$ with a dlt
modification. Then we may assume
 $$
\text{$K_X+(1-\ep)D$ is klt for $0<\ep\ll 1$. In particular, $X$
is $\QQ$-factorial.}\leqno{(*)}
 $$
If the fiber $f^{-1}(f(L))$ is reducible, there is its component
$L_1\neq L$ meeting $L$. Then $L_1^2<0$ and $(K_X+D-\ep' L)\cdot
L_1<0$ for $0<\ep'\ll 1$. Thus $L_1$ generates an extremal ray
which is negative with respect to $K_X+(1-\ep'')D-\ep' L$ for
$0<\ep''\ll \ep'$. By Contraction Theorem \cite[3-2-1]{KMM} we can
contract $L_1$ over $Z$. This contraction preserves all
assumptions of the lemma as well as assumption $(*)$ (however, we
can lose the dlt property of $(X,D)$). Continuing the process, we
get the situation when the fiber containing $L$ is irreducible.
Similarly, components of all reducible fibers can be contracted.
We obtain a model $f'\colon X'\to Z$ such that all fibers are
irreducible. Moreover $f'$ is an extremal $K$-negative
contraction. Hence $\rho(X'/Z)=1$. By our construction,
$K_{X'}+D'$ is lc, $X'$ is $\QQ$-factorial and $L'\subset\Supp D'$
is a fiber of $f'$. Let $R$ be an extremal ray on $X'$ other than
that generated by fibers of $f'$. Then $R\cdot L'>0$ and
$(K_{X'}+D')\cdot R\le 0$. Therefore $(K_{X'}+D'-\delta' L')\cdot
R<0$ for $\delta'>0$. Hence there is a curve $M$ on $X'$
generating $R$ (as above, if $K_{X'}+D'$ is not dlt, we can use
Contraction Theorem for $K_{X'}+(1-\delta)D'-\delta' L'$, see also
Appendix~\ref{MMP-in-dimension-2}). In this situation,
$M\simeq\PP^1$ (see Proposition~\ref{MMP-extremal-curves}). But
then the base curve $Z$ also should be rational, a contradiction
with $g\ge 1$.
\end{proof}

\begin{corollary1}
\label{corollary-1-1}
Notation as in Lemma~\ref{lemma1}. Then the pair $(X,D)$ has at
worst canonical singularities.
\end{corollary1}
\begin{outline}
Replace $(X,D)$ with a suitable log terminal modification (see
Proposition~\ref{per-gen}) and apply Lemma~\ref{lemma1}.
\end{outline}

Going back to the proof of Theorem~\ref{Inductive-Theorem-nonrat},
denote $C:=\down{D}$ and $B:=\fr D$. Let $L$ be a component of
$C$. By Lemma~\ref{lemma1}, $f(L)=Z$. If $\rho(X/Z)>1$, there is a
curve $E$ in a (reducible) fiber such that $E^2<0$ and $L\cdot
E>0$. Since $(K_X+D)\cdot E\le 0$, $E$ is a $-1$-curve and
\[
{-}K_X\cdot E=L\cdot E=1,\qquad E\cdot(D-L)=0.
\]
Hence we can contract $E$ and by Lemma~\ref{bir-prop} and
Remark~\ref{bir-prop-comment} we can pull back complements under
this contraction. Repeating the process, we reach the situation
when $\rho(X/Z)=1$. Thus we may assume that $X$ is a (smooth)
ruled surface over a nonrational curve $Z$, i.e., $X=\PP_Z(\EEE)$,
where $\EEE$ is a rank two vector bundle on $Z$. Moreover, $K_X+D$
is dlt (see Proposition~\ref{crepant-singul}). Then
\[
0\ge (K_X+D)\cdot L\ge (K_X+L)\cdot L=2p_a(L)-2\ge 0.
\]
This implies
\begin{equation}\label{eq-new-new}
(K_X+D)\cdot L=(K_X+L)\cdot L=(D-L)\cdot L=0.
\end{equation}
Thus $p_a(L)=1$ and $L$ is a smooth elliptic curve. Hence $g=1$.

\begin{claim*}
We may assume that $(K_X+L)^2\ge 0$.
\end{claim*}
\begin{proof}
Let $F$ be a general fiber of $f$. Since $\rho(X)=2$, there is
exactly one extremal ray $R\neq \RR_+[F]$ on $X$. Assume that
$(K_X+L)^2<0$. Then $-(K_X+L)$ is not nef, $D\neq L$ and
$(K_X+L)\cdot R>0$. If $L^2\le 0$, then $L$ generates an extremal
ray (see Proposition~\ref{MMP-properties}), so $R=\RR_+[L]$. This
contradicts \eqref{eq-new-new}. Therefore, $L^2>0$. In particular,
$L\cdot R\ge 0$.

Further, if $L\cdot R>0$, then $R$ is negative with respect to
$K_X+D-L$. Since $K_X+D-L$ is dlt, the ray $R$ must be generated
by a rational curve. This implies that $Z$ is also rational, a
contradiction. Finally, we have $L\cdot R=0$. Then
\eqref{eq-new-new} implies $[D-L]\in R$ (recall that $D-L$ is
effective). By the Hodge Index Theorem, $R^2<0$. Thus, $D-L=aL'$
and $R=\RR_+[L']$, where $0<a\le 1$ and $L'$ is an irreducible
curve with $L^{\prime 2}<0$. If $a<1$, then $R$ is negative with
respect to $K_X+D+\ep L'$ for $0<\ep\ll 1$. Again we have a
contradiction. Hence $a=1$ and $D=L+L'$. As in \eqref{eq-new-new}
we have
\begin{equation*}
(K_X+D)\cdot L'=(K_X+L')\cdot L'=L\cdot L'=0.
\end{equation*}
Since $\rho(X)=2$, $-(K_X+L')$ is nef, so $(K_X+L')^2\ge 0$.
Replacing $L$ with $L'$ we get our assertion.
\end{proof}

Taking into account the equality $K_X^2=8(1-g(Z))=0$, we obtain
\[
0\le (K_X+L)^2=L\cdot K_X, \qquad L^2\le 0.
\]
Therefore $L$ generates an extremal ray $R$ on $X$ (see
Proposition~\ref{MMP-properties}). It cannot be $K_X$-negative
(otherwise $R$ is generated by a rational curve). Hence $L\cdot
K_X=L^2=0$ and in relations above equalities hold. Thus inequality
\eqref{eq-new-new} gives that $L$ and all the components
$D_i\subset\Supp(D-L)$ are numerically proportional (because
$\rho(X)=2$). In particular,
\begin{equation}
\label{intersections}
D_i^2=0,\quad \forall i\qquad\text{and}\qquad L\cap D_i=D_i\cap
D_j=\emptyset,\quad \forall i,j.
\end{equation}
Let $F$ be a general fiber of $f$. Then $F\cdot(K_X+D)\le 0$ and
$K_X\cdot F=-2$. This yields $1\le C\cdot F\le D\cdot F\le 2$.

\begin{lemma1}
\label{n:1}
Let $f\colon X\to Z$ be a ruled surface over an elliptic curve and
let $C$ be a reduced divisor such that $K_X+C\equiv 0$. Then
$2(K_X+C)\sim 0$. Moreover, if $C$ is reducible, then $K_X+C\sim
0$.
\end{lemma1}

\begin{proof}
Since $K_X+C\equiv 0$, $K_X+C=f^*N$ for some integral divisor of
degree $0$ on $Z$. First we assume that $C$ is irreducible. Then
$C$ is a smooth elliptic curve and we have $f_{C}^*N=
(K_X+C)|_C=K_C=0$. Hence $2N\sim 0$ because $f_C$ is of degree
two. Now we assume that $C=C_1+C_2$, where $C_1$, $C_2$ are
sections. Similarly, $f_{C_1}^*N=(K_X+C)|_C=K_C=0$. Hence $N\sim
0$ because $f_{C_1}$ is an isomorphism.
\end{proof}

Now we finish the proof of Theorem~\ref{Inductive-Theorem-nonrat}.

If $C\cdot F=2$ (i.e., $C$ is a $2$-section of $f$), then $D=C$,
$B=0$ and $(K_X+C)\cdot F=0$. Hence $K_X+C\equiv 0$. By
Lemma~\ref{n:1}, $2(K_X+C)\sim 0$, i.e., we have an $1$ or
$2$-complement (with $C^+=C$).
\par
If $C\cdot F=1$, then $C$ is a section of $f$. Recall that
$X=\PP_Z(\EEE)$, where $\EEE$ is a rank two vector bundle on $Z$.
By assumption (ii) of the theorem, the surface $X$ is not such as
in Example~\ref{Atiyah}. On the other hand, by
\eqref{intersections} $C^2=0$ and the vector bundle $\EEE$ has
even degree. From the classification of rank two vector bundles
over elliptic curves (see e.g., \cite{Ha}), we obtain that $\EEE$
is of splitting type. Hence there is a section $C_1$ such that
$C\cap C_1=\emptyset$. Write $B=\sum b_iB_i$. Then
\begin{equation}
\label{eq-nrat-1}
\sum b_iB_i\cdot F\le 1.
\end{equation}
From \ref{n:1} we have
\begin{equation}
\label{CC1}
{-}K_X\sim C+C_1,
\end{equation}
where $C$, $C_1$ is a pair of disjoint sections.
\par
By \eref{CC1} and by the definition of complements, $D^+=C+C_1$ is
an $1$-complement of $K_X+D$ if $b_i<1/2$ whenever $B_i\ne C_1$.
By \eref{eq-nrat-1} this does not hold only if $B=\frac12B_1$,
where $B_1$ is a $2$-section. Consider this case. As in the proof
of Lemma~\ref{n:1}, $2(K_X+C+B)=f^*N$ and
$f_C^*N=2(K_X+C+B)|_C=2K_C+B_1|_C=0$. Therefore $N=0$ and
$2(K_X+C+B)\sim 0$. This proves our theorem.
\end{proof}

\section{The Main Inductive Theorem}

\index{Inductive Theorem}
\begin{theorem}[Inductive Theorem, \cite{Sh1}]
\label{Inductive-Theorem}
Let $(X,D)$ be a projective log surface such that $D\in\Mm$,
$K_X+D$ is lc, but not klt and $-(K_X+D)$ is nef. Assume
additionally that there exists a boundary $D'\le D$ such that
$K_X+D'$ is klt and $-(K_X+D')$ is nef and big. Then there exists
a regular complement of $K_X+D$.
\end{theorem}

Shokurov proved this theorem under weaker assumptions \cite{Sh1}.
In particular, he showed that if we remove condition $D\in\Mm$ in
the theorem, we obtain a weaker result: there exists an
$n$-complement of $K_X+D$, where
\begin{multline*}
n\in\{1, 2, 3, 4, 5, 6, 7, 8, 9, 10, 11, 12, 13, 14, 15, 16, 17,
18, 19, 20, 21, 22,\\
 23, 24, 25, 26, 27, 28, 29, 30, 31, 35, 36,
40, 41, 42, 43, 56, 57\}.
\end{multline*}

\begin{remark*}
In contrast with \ref{Inductive-Theorem-Weak-Form}, in
Theorem~\ref{Inductive-Theorem} we cannot say that any regular
complement of $K_C+\Diff{C}(B)$ can be extended to $X$. For
example, let $X=\PP^1\times\PP^1$, $C$ is a fiber of
$\mt{pr}_1\colon X\to \PP^1$ and
$B:=\frac12(B'+B'')+\frac23(B_1+B_2+B_3)$, where $B'$, $B''$ are
fibers of $\mt{pr}_1$ and $B_1$, $B_2$, $B_3$ are fibers of
$\mt{pr}_2$. Then $K_X+C+B\equiv 0$, $K_C+\Diff{C}(B)=
K_C+\frac23(P_1+P_2+P_3)$ is $3$-complementary and $K_X+C+B$ has
no $3$-complements. Indeed, otherwise for a $3$-complement, we
have $B^+\ge \frac23(B'+B''+B_1+B_2+B_3)$, a contradiction.
However, $K_X+C+B$ has a $6$-complement with $B^+=B$.
\end{remark*}

\begin{corollary1}
\label{Corr-Inductive-Theorem-discr}
Notation as in Theorem~\ref{Inductive-Theorem}. We can take a
regular complement $K_X+D^+$ so that $a(E,D)=-1$ implies
$a(E,D^+)=-1$ for any $E$.
\end{corollary1}
\begin{proof}
Since $D^+\ge D$ (see Remark~\ref{important_compl}), we have
$a(E,D^+)\le a(E,D)$.
\end{proof}

\begin{proof}[Proof of Theorem~\ref{Inductive-Theorem}]
If $-(K_X+D)$ is big, then by
Proposition~\ref{Inductive-Theorem-Weak-Form} there is a regular
complement. Therefore we assume that $K_X+D\equiv 0$ or
$\kappa(X,-(K_X+D))=1$.
\subsection{}
\label{construction-boundary}
Applying a minimal log terminal modification $f\colon \ov X\to X$
we may assume that $K_{\ov X}+\ov D=f^*(K_X+D)$ is dlt and
$\down{\ov D}\ne 0$. Take the crepant pull back
\[
f^*(K_X+D')=K_{\ov X}+\ov{D'},\quad\text{with}\quad f_*\ov{D'}=D'.
\]
By \ref{crepant-singul}, $K_{\ov X}+\ov{D'}$ is klt, but it is not
necessarily a boundary. Consider the new boundary
\[
\ov{D''}:=\ov{D'}+t(\ov{D}-\ov{D'}),\quad \text{where}\ 0<1-t\ll
1.
\]
Then $K_{\ov X}+\ov{D''}$ is klt and $-(K_{\ov X}+\ov{D''})$ is
nef and big. Further, by \ref{bir-prop-a} we can push-down
complements. So we replace $X$, $D$, $D'$ with $\ov X$, $\ov D$,
$\ov{D''}$. Thus we assume now that $K_X+D$ is dlt, $\down D\ne
0$, $-(K_X+D)$ is nef, and there exists a boundary $D'\le D$ such
that $K_X+D'$ is klt and $-(K_X+D')$ is nef and big. By
Lemma~\ref{lemma_non_rat}, $X$ is rational. Set $C:=\down{D}$ and
$B:=\fr{D}$. The following lemma shows that the Mori cone $\NE(X)$
is polyhedral and generated by contractible extremal curves.

\begin{lemma1}
\label{Kodaira-lemma}
Let $(X/Z,B)$ be a log variety such that $K_X+B$ is klt and
$-(K_X+B)$ is nef and big over $Z$. Then there exists a new
boundary $B'\ge B$ on $X$ such that $K_X+B'$ is again klt and
$-(K_X+B')$ is ample over $Z$.
\end{lemma1}
\begin{outline}
Let $H$ be a very ample divisor on $X$ (over $Z$). By Kodaira's
Lemma, $|-n(K_X+B)-H|\ne\emptyset$ for some $n\gg 0$ (see e.g.,
\cite[0-3-4]{KMM}). Take $L\in |-n(K_X+B)-H|$ and put $B':=B+\ep
L$.
\end{outline}

The lemma shows that we can contract all extremal rays on $X$.
Moreover, if an extremal ray $R$ on $X$ is birational and
generated by a curve $L$ which is not contained in $C$, then the
contraction preserves all assumptions (see \ref{crepant-singul}).
If additionally $R$ is $(K_X+D)$-trivial, we can pull back regular
complements by Proposition~\ref{bir-prop}.

\subsection{Division into cases}
\label{division-into-c}
By Proposition~\ref{dipol2}, $C$ has at most two connected
components. As in proofs of Theorems~\ref{class_lc} and
\ref{class} we distinguish the following cases: $C$ is
disconnected, $C$ is a smooth elliptic curve or a wheel of smooth
rational curves, $C\simeq\PP^1$, $C$ is a tree of rational curves.
Also we should separate cases $K_X+C+B\equiv 0$ and
$K_X+C+B\not\equiv 0$. If $K_X+C+B\not\equiv 0$, then by Base
Point Free Theorem the linear system $|-m(K_X+C+B)|$ determines a
contraction $\nu\colon X\to\PP^1$. By Lemma~\ref{Kodaira-lemma},
$-(K_X+D'')$ is ample for some boundary $D''$, hence a general
fiber is rational. Then we consider cases when $C$ is contained in
a fiber of $\nu$ and $C$ has a horizontal component.

\subsection{Case: $C$ is disconnected}
\label{(aa)}
By Proposition~\ref{dipol2} there exists a contraction $f\colon
X\to Z$ onto a curve such that $C=C_1+C_2$ is a pair of two
disjoint smooth sections (in particular, $K_X+D$ is plt). A
general fiber $F$ of $f$ is $\PP^1$ (see \ref{division-into-c}),
so ${-}K_X\cdot F=C\cdot F=2$, $B\cdot F=0$ and $B$ is contained
in fibers of $f$. Since $X$ is rational, $Z\simeq\PP^1$. In our
case $-(K_X+D)$ is numerically trivial on a general fiber of $f$,
so it is numerically trivial on all fibers. Contracting curves in
fibers we get the situation when $\rho(X/Z)=1$. We can pull back
all complements by Proposition~\ref{bir-prop}. If $C_1^2>0$, then
$-(K_X+(1-\ep)C_1+C_2+B)$ is nef and big for any $\ep>0$. By
Proposition~\ref{Inductive-Theorem-Weak-Form} there exists a
regular complement of $K_X+(1-\ep)C_1+C_2+B$. By
Corollary~\ref{simple-prp-1} $K_X+C_1+C_2+B$ also has a regular
complement. Therefore we may assume that $C_1^2\le 0$ and
$C_2^2\le 0$. Then both $C_1$ and $C_2$ generate extremal rays
which must coincide because $\rho(X)=2$. The ray cannot be
birational, so $C_1^2=C_2^2=0$. This shows that there exists a
nonbirational contraction $g\colon X\to \PP^1$ such that $C_1$ and
$C_2$ are fibers of $g$. If $K_X+C_1+C_2+B\not\equiv 0$, then
again $-(K_X+(1-\ep)C_1+C_2+B)$ is ample for $\ep>0$. As above (by
Proposition~\ref{Inductive-Theorem-Weak-Form} and
Corollary~\ref{simple-prp-1}) there is a regular complement.
Therefore we may assume that $K_X+C_1+C_2+B\equiv 0$. Now it is
sufficient to verify $n(K_X+C+B)\sim 0$ for some $n\in\RRR_2$. By
Lemma~\ref{log-del-Pezzo-Pic} the numerical equivalence in
$\Pic(X)$ coincides with linear one. Therefore it is sufficient to
show only that $n(K_X+C+B)$ is Cartier for some $n\in\RRR_2$. Take
\[
n:=\min\{r\in\NN\mid rB\ \text{is integral.}\}
\]
Since $\down{B}=0$, $n>1$. By Theorem~\ref{bir-nef}, $K_X+C+B$ has
a regular $n_1$-complement $K_X+C+B^+$ near $C_1$ (a fiber of
$g\colon X\to\PP^1$). Then $B^+\ge B$ by
Remark~\ref{important_compl}. This yields $B^+=B$, $n_1B$ is
integral and $n_1(K_X+C+B)\sim 0$ near $C_1$ (because $C_1$
intersects all the components of $B$). Hence $n\mid n_1$.
Similarly, we have a regular $n_2$-complement near $C_2$ and
$n\mid n_2$. Let $n':=\lcm(n_1,n_2)$. Then $n'B$ is integral and
$n'(K_X+C+B)$ is Cartier near $C$. Let $F$ be a fiber of $f$ and
let $P_i:=F\cap C_i$. By Adjunction,
\[
0\ge (K_X+C+F)\cdot F=\deg K_F+\deg\Diff F(C)
\]
and
\[
\Diff F(C)\ge P_1+P_2.
\]
Therefore $\Diff F(C)=P_1+P_2$ and $X$ is smooth outside of $C$.
Further, $n'(K_X+C+B)$ is Cartier everywhere on $X$ and it is
sufficient to show that $n'\in\RRR_2$. Assume the opposite. Then
we have (up to permutations $C_1$ and $C_2$): $n_1=4$, $n_2=6$ and
$n=2$. Since $2B$ is integral, $B\in\Msm$.
Corollary~\ref{coeff-Diff-1} gives $\Diff{C_i}(B)\in\Msm$. By
\ref{=},
\[
\Diff{C_1}(B)=\frac12Q_1+\frac34Q_2+\frac34Q_3,\quad
\Diff{C_2}(B)=\frac12R_1+\frac23R_2+\frac56R_3,
\]
where $Q_1, Q_2, Q_3\in C_1$, $R_1, R_2, R_3\in C_2$ are some
points. On the other hand,
\[
(C_1,\Diff{C_1}{B})\simeq (C_2,\Diff{C_2}{B})
\]
(see \cite[12.3.4]{Ut}) a contradiction. This proves our theorem
in the case when $C$ is disconnected.

\subsection{Case: $C$ is connected and $p_a(C)\ge 1$}
\label{(bb)}
By Lemma~\ref{wheel-1} below there is an $1$-complement. Note that
in this case the assumption $D\in\Mm$ is not needed. First we
claim that $K_X+C+B\equiv 0$. Indeed, by Adjunction we have
\[
0\ge (K_X+C+B)\cdot C\ge (K_X+C)\cdot C\ge\deg K_C\ge 0
\]
If $K_X+C+B\not\equiv 0$, then $C$ is contained in fibers of
$\nu\colon X\to \PP^1$. But $\nu$ has rational fibers (see
\ref{division-into-c}), a contradiction.

\begin{lemma1}
\label{wheel-1}
Let $(X,C+B)$ be a rational projective log surface, where $C$ is
the reduced and $B$ is the fractional part of the boundary. Assume
that $K_X+C+B$ is analytically dlt, $K_X+C+B\equiv 0$, $C$ is
connected and $p_a(C)\ge 1$. Then $B=0$, $K_X+C\sim 0$, $X$ is
smooth along $C$ and has only Du Val singularities outside.
\end{lemma1}
\begin{proof}
By Lemma~\ref{wheel} and Remark~\ref{wheel-rem}, $X$ is smooth and
$B=0$ near $C$. Replace $X$ with a minimal resolution and $C+B$
with the crepant pull back. It is sufficient to show only that
$B=0$ and $K_X+C\sim 0$. Now we contract $-1$-curves on $X$. Since
$C$ is not a tree of rational curves, it cannot be contracted.
This process preserves all the assumptions, so on each step
$C\cap\Supp B=\emptyset$. Since every $-1$-curve $E$ has positive
intersection number with $C+B$, we have either $C\cdot E\ge 1$,
$B\cdot E=0$ or $C\cdot E=0$, $B\cdot E>0$. If $B\ne 0$, then
whole $\Supp B$ also cannot be contracted. At the end we get
$X=\PP^2$ or $\FF_n$ (a Hirzebruch surface). In the case
$X=\PP^2$, $C+B\equiv {-}K_X$ is ample. Hence $\Supp(C+B)$ is
connected. If $B\ne 0$ we derive a contradiction. Consider the
case $X=\FF_n$. Then
\[
0=(K_X+C+B)|_C=(K_X+C)|_C+B|_C=K_C+B_C
\]
and the last two terms are nonnegative. Therefore
$(K_X+C)|_C=B|_C=0$ and $p_a(C)=1$. On the other hand, for a
general fiber $F$ of $\FF_n\to \PP^1$ one has
\[
(K_X+C+B+F)|_F=K_F+(C+B)|_F.
\]
In particular, $(C+B)\cdot F\le 2$. Since $p_a(C)=1$, $C$ is not a
section of $\FF_n\to \PP^1$ at a general point. Hence $C\cdot F=2$
and $B\cdot F=0$. Recall that $B\cdot C=0$. Since $\rho(\FF_n)=2$,
we have $B\equiv 0$ and $B=0$. We proved that $p_a(C)=1$, $B=0$
and $K_X+C\sim 0$ in the case $X=\FF_n$. Therefore on our original
$X$ one also has $B=0$. By Proposition~\ref{bir-prop} we can pull
back an $1$-complement of $K_X+C^+=K_X+C$ under contractions of
$-1$-curves.
\end{proof}

\subsection{Case: $C\simeq\PP^1$ (the exceptional case)
and $K_X+C+B\equiv 0$}
\label{(dd)}
In this case, $K_X+C+B$ is plt. We claim that after a number of
birational contractions $B\in\Msm$. Indeed, otherwise there is a
component $B_i$ of $B$ with coefficient $b_i\notin\Msm$ and
$b_i>6/7$. If $B_i^2>0$, then as in \ref{(aa)} we can reduce $b_i$
a little so that $-(K_X+C+B)$ becomes big and obtain a regular
complement by Proposition~\ref{Inductive-Theorem-Weak-Form} (note
that $b_i^+=1$). This complement is also a complement of our
original $K_X+C+B $ by Proposition~\ref{simple-prp}. If $B_i^2<0$,
then we can contract $B_i$ and pull back complements by
Proposition~\ref{bir-prop} and Remark~\ref{bir-prop-comment}.
Consider the case $B_i^2=0$. By Base Point Free Theorem the linear
system $|mB_i|$ determines a contraction $f\colon X\to\PP^1$ such
that $B_i$ is a fiber of $f$. Contracting curves $\ne C$ in
reducible fibers we get the situation when $\rho(X)=2$, i.e.
$f\colon X\to\PP^1$ is an extremal contraction. Let $g\colon X\to
Z$ be another extremal contraction on $X$ and $F$ a nontrivial
fiber of $g$. Then $F\simeq\PP^1$ and $F\cap B_i\ne\emptyset$
(because $F$ is not a fiber of $f$). Assume that $g$ is
nonbirational and $F$ is sufficiently general. Then $X$ is smooth
along $F$, $\Diff F(C+B)=(C+B)|_F$ and $F$ intersects $\Supp(C+B)$
transversally. Hence $\Diff F(C+B)\in\Mm$ and we can write
\[
\Diff F(C+B)=\sum_j b_jP_{j,l}+\sum_{s=1}^r Q_s,
\]
where $\{P_{j,1},\dots,P_{j,r_j}\}=B_j\cap F$ and
$\{Q_1,\dots,Q_r\}:=C\cap F$. Moreover, the coefficient $b_i$ of
$\Diff F(C+B)$ at points $B_i\cap F$ satisfies $6/7<b_i<1$.
Further,
\[
\deg\Diff F(C+B)=\deg({-}K_F)=2.
\]
Easy computations as in \ref{=} show that this is impossible.
Indeed,
\begin{equation}
\label{eq-Diff--1}
2=\deg\Diff F(C+B)=r+\sum_j b_jr_j.
\end{equation}
Clearly, $\sum_{j\neq i} b_jr_j=2-r-b_ir_i>0$ (otherwise
$b_ir_i=2-r\in\{0,1,2\}$ but $6/7<b_i<1$). Since $b_j\in \Mm$, we
have $\sum_{j\neq i} b_jr_j\ge 1/2$. Thus
\[
2-r=\sum_j b_jr_j\ge 1/2+(6r_i)/7>1.
\]
This gives us $r=0$ and $r_i=1$. Hence $\sum_{j\neq i}
b_jr_j=2-b_i< 8/7$. It is easy to check that the last inequality
has no solutions with $r_j\in \NN$ and $b_j\in\Mm$.
\par
If $g$ is birational and contract $C$ (i.e., $C=F$), then $\Diff
C(B)\in\Mm$ (by Corollary~\ref{coeff-Diff-1}) and has a
coefficient $>6/7$. Moreover, $K_X+C+B$ is plt, so $K_C+\Diff
C(B)$ is klt. As above we derive a contradiction with $\deg\Diff
C(B)=2$. Finally, if $g$ is birational and does not contract $C$,
we can replace $X$ with $Z$. We get the situation when $\rho=1$
and $B_i^2>0$, as above.
\par
Thus we may assume now that $B\in\Msm$. By Base Point Free Theorem
and assumptions of the theorem, there exists $n\in\NN$ such that
$n(K_X+C+B)\sim 0$. Let $n$ be the index of $K_X+C+B$ (i.e., the
minimal positive integer with this property) and $\var\colon X'\to
X$ the log canonical $n$-cover (see \ref{l-c-cover}). It is
sufficient to show that $n\in\RRR_2$. Write
\[
K_{X'}+C'=\var^*(K_X+C+B),
\]
where $C'=\var^*C$. Then $K_{X'}+C'$ is linearly trivial and plt
(see Proposition~\ref{finite-plt}). By Adjunction every connected
component of $C'$ is a smooth elliptic curve.

First we assume that $C'$ is connected. By construction,
$K_{X'}+C'\sim 0$ and we can identify $\mt{Gal}(C'/C)$ with
$\ZZ_n$. We claim that $\mt{Gal}(C'/C)$ contains no translations.
Indeed, let $\xi\in\mt{Gal}(C'/C)$ a translation. Then we put
$X'':=X'/\langle\xi\rangle$ and $C'':=C'/\langle\xi\rangle$. By
Lemma~\ref{wheel-1}, $K_{X''}+C''$ is linearly trivial and plt
(because $p_a(C'')=1$). But then $n''(K_X+C+B)\sim 0$, where $n''$
is the degree of $X''\to X$. By assumptions $n$ is the smallest
positive with this property. The contradiction shows that
$\mt{Gal}(C'/C)$ contains no translations. Then $\mt{Gal}(C'/C)$
is a finite group of order $2$, $3$, $4$, or $6$ (see e.g.,
\cite{Ha}).

If $C'$ is disconnected, then $\mt{Gal}(X'/X)$ interchange
connected components of $C'$. By Proposition~\ref{dipol2} there is
a contraction $X'\to Z'$ onto an elliptic curve with rational
fibers such that components of $C'=C'_1+C'_2$ are sections. This
contraction must be $\mt{Gal}(X'/X)$-equivariant because $X'$ has
a unique structure of a contraction with rational fibers. Set
$G_0:=\mt{Ker}(\mt{Gal}(X'/X)\to \mt{Aut}(Z'))$. Since the
ramification locus of $X'\to X$ does not contain components of
$C'$, $G_0\simeq\cyc{2}$. As above we consider $X'':=X'/G_0$ and
$C'':=C'/G_0$. Then $C''$ is a smooth elliptic curve, hence
$K_{X''}+C''\sim 0$ by Lemma~\ref{wheel-1}. This contradicts our
choice of $n$.

\subsection{Case: $C$ is a tree of rational curves and
$K_X+C+B\equiv 0$}
\label{(ee)}
By Lemma~\ref{chain}, $C$ is a chain and $B$ has coefficients
$=1/2$ near $C$. As in case~\ref{(dd)} we claim that after some
birational contractions $B\in\Msm$. Let $B_k$ be a component with
nonstandard coefficient. If $B_k^2\ne 0$, then we can argue as in
case~\ref{(dd)}. The only nontrivial case is $B_k^2=0$ and
$B_k\cap C=\emptyset$. Then again $|mB_k|$ determines a
contraction $f\colon X\to\PP^1$. Clearly, $B_k$ is a fiber and $C$
is contained in a fiber (because $C\cdot B_k=0$). There is an
extremal rational curve $F$ which is not contained in fibers. Then
$F\cap B_k\ne\emptyset$. If $F^2<0$, we contract $F$ and replace
$X$ with a new birational model. If $F^2=0$ we derive a
contradiction computing $\Diff F(C+B)$ as in \eref{eq-Diff--1} of
case~\ref{(dd)}.
\par
Thus we may assume that $B\in\Msm$. As in case~\ref{(dd)} take the
log canonical $n$-cover $\var\colon X'\to X$. It is sufficient to
show that $n\in\RRR_2$. Again we can write
\[
\var^*(K_X+C+B)=K_{X'}+C'\sim 0,
\]
where $C':=\var^*C$. Obviously, $C'$ is reducible. We claim that
$K_{X'}+C'$ is dlt. By Proposition~\ref{finite-plt} $K_{X'}+C'$ is
plt outside of $\var^{-1}(\Sing C)$. Recall that the ramification
divisor of $\var$ is $\Supp B$. Hence none of irreducible
components of the ramification divisor intersects $\Sing C$. At
points $\Sing C$ the surface $X$ is smooth, so $\var$ is \'etale
over $\Sing C$. Therefore $K_{X'}+C'$ is dlt and $X'$ is smooth at
points $\var^{-1}(\Sing C)$. Since $K_{X'}+C'\sim 0$, $X'$ is
smooth along $C'$ (see \ref{kawamata_1}) and $p_a(C')=1$. By
Lemma~\ref{wheel}, $C'$ is a wheel of smooth rational curves and
by our construction $\mt{Gal}(X'/X)=\cyc{n}$ acts on $C'$
faithfully. Let $C'=\sum_{i=1}^r C'_i$ be the irreducible
decomposition and $P_1,\cdots, P_r$ singular points of $C'$. If
$\mt{Gal}(X'/X)$ contains an element $\xi$ such that $\xi\cdot
C_i'=C_i'$ and $\xi\cdot P_i=P_i$ for all $i$, then
$C'':=C'/\langle\xi\rangle$ is again a wheel of smooth rational
curves. As in case~\ref{(dd)} we derive a contradiction. Therefore
$\mt{Gal}(X'/X)$ acts faithfully on the dual graph of $C'$ and
then it is a subgroup of the dihedral group $\mathfrak{D}_r$. The
same arguments show that every $\xi\in\mt{Gal}(X'/X)$ has a fixed
point on $C'$. This is possible only if
$\mt{Gal}(X'/X)\simeq\cyc{2}$. Therefore $n=2$ and $K_X+C+B$ is
$2$-complementary.

\subsection{Case: $K_X+C+B\not\equiv 0$ and $C$ is contained
in a fiber of $\nu$}
\label{(ff)}
We may contract all components of reducible fibers of $\nu\colon
X\to\PP^1$ which are different from components of $C$. Thus
$C=\nu^{-1}(P)$, $P\in\PP^1$ is a fiber and all other fibers of
$\nu$ are irreducible. Again we may assume that $\nu$-horizontal
components of $B$ have standard coefficients (otherwise some
$b_i>6/7$ and $-(K_X+C+B-\ep B_i)$ is nef and big for $0<\ep \ll
1$, hence we can use Proposition~\ref{Inductive-Theorem-Weak-Form}
and Corollary~\ref{simple-prp-1}). Note that the horizontal part
$B_{\mt{hor}}$ of $B$ is nontrivial (because a general fiber of
$\nu$ is $\PP^1$). Further, there is a regular $n$-semicomplement
of $K_C+\Diff C(B)$ (see Theorem~\ref{1}). For sufficiently small
$\ep>0$ the $\QQ$-divisor $-(K_X+C+B-\ep B_{\mt{hor}})$ is nef and
big. Thus we can extend $n$-semicomplements of $K_C+\Diff C(B-\ep
B_{\mt{hor}})$ from $C$ by
Proposition~\ref{Inductive-Theorem-Weak-Form}. If $n+1$ is not a
denominator of coefficients of $B_{\mt{hor}}$, then by
Proposition~\ref{simple-prp} we obtain a regular complement of
$K_X+C+B$. If $C$ is reducible, then by Theorem~\ref{1} we can
take $n=2$. On the other hand, by Lemma~\ref{chain} coefficients
of $B_{\mt{hor}}$ are equal to $1/2$. Therefore $n+1=3$ is not a
denominator of coefficients of $B_{\mt{hor}}$ in this case and
there is a $2$-semicomplement of $K_X+C+B$. Now we assume that
$C\simeq\PP^1$. Then $\rho(X)=2$. By \ref{=} and by
Corollary~\ref{coeff-Diff-1} we have the following possibilities:
\[\text{
\begin{tabular}{cccc}
 $\Diff S(B)$&\qquad&$n$&\\
 &&\\
 $\frac12P_1+\cdots+\frac12P_4$&&2&\\
 $\frac23P_1+\frac23P_2+\frac23P_3$&&$3$&\\
 $\frac12P_1+\frac34P_2+\frac34P_3$&&$4$&\\
 $\frac12P_1+\frac23P_2+\frac56P_3$&&$6$&\\
\end{tabular}}
\]
By Corollary~\ref{coeff-Diff}, $n+1$ is not a denominator of
coefficients of $B_{\mt{hor}}$ in all cases.

\subsection{Case: $K_X+C+B\not\equiv 0$ and $C$ is not
contained in a fiber of $\nu$} So, we assume that there is a
horizontal component $C_i\subset C$ (i.e., such that
$\nu(C_i)=\PP^1$). It is clear that $-(K_X+C+B-\ep C_i)$ is nef
and big for a sufficiently small positive $\ep$. If $C_i\subsetneq
C$ (i.e., $C$ is reducible), then the same trick as in \ref{(aa)}
(using Corollary~\ref{simple-prp-1}) gives the existence of
regular complements. Thus we may assume that $C=C_i$ and
$\nu(C)=\PP^1$. In particular, $K_X+C+B$ is plt. Contracting
curves $\ne C$ in fibers we get the situation, when fibers are
irreducible, i.e., $\rho(X)=2$. We can pull back complements by
\ref{bir-prop}. There are two subcases: (a) $C$ is a section of
$\nu$, and (b) $C$ is a $2$-section of $\nu$.

\subsection{}
If $C$ is a section of $\nu$, then the horizontal part
$B_{\mt{hor}}$ of $B$ is nontrivial and either $B_{\mt{hor}}
=\frac12B_1+\frac12B_2$ or $B_{\mt{hor}} =\frac12B_0$, where
$B_1$, $B_2$ are sections and $B_0$ is a $2$-section. As above we
can take a regular $n$-complement for $K_X+C+B-\ep B_{\mt{hor}}$.
If $n\ne 1$, then again this gives a regular complement of
$K_X+C+B$ by Proposition~\ref{simple-prp}. But on $C\simeq\PP^1$
there exists a regular $2$, $3$, $4$, or $6$-complement for the
boundary $\De:=\Diff C(B)\ge\Diff C(B-\ep B_{\mt{hor}})$. Indeed,
otherwise by Theorem~\ref{1} there is an $1$-complement
$K_C+\De^+$. By Corollary~\ref{coeff-Diff-1}, $\De\in\Mm$.
Therefore $\De^+\ge \De$ and $\De$ is supported in one or two
points (because $\deg\De^+=2$). Then $K_C+\De^+$ is also
$n$-complement for any $n$. This shows that we have a regular
complement in the case when $C$ is a section.

\subsection{}
\label{last-case-=}
Now let $C$ be a $2$-section of $\nu$. Then $B$ is contained in
the fibers of $X$. Since $C\simeq\PP^1$, the restriction
$\nu\colon C\to\PP^1$ has exactly two ramification points, say
$P_1$, $P_2\in C$. Put $Q_i:=\nu(P_i)$ and
$F_i:=\nu^{-1}(Q_i)_{\red}$, $i=1,2$.

\begin{lemma1}
\label{lemma-2}
Notation as in \ref{last-case-=}. Let $F$ be a fiber of $\nu$ such
that $F\ne F_1,F_2$ and $F\cap C=\{P',P''\}$ (where $P'\ne P''$).
Then
\begin{enumerate}
\item
$\Sing X\subset C\cup F_1\cup F_2$;
\item
$K_X+C+F$ is lc and linearly trivial near $F$;
\item
$K_X+F$ is plt;
\item
$\Diff C(B)$ is invariant under the natural Galois action of
$\ZZ_2$ on $C\to\PP^1$.
\end{enumerate}
\end{lemma1}

\begin{proof}
First we show that $K_X+C+F$ is lc. Assume the converse and regard
$X$ as an analytic germ near $F$. Let $C_1$, $C_2$ be analytic
components of $C$. If $K_X+C_1+C_2+F$ is not lc near $C_1$, then
$K_X+(1-\ep)C_1+C_2+(1-\ep)F$ is not, either. But in this case,
$\LCS(X,(1-\ep)C_1+C_2+(1-\ep)F)$ is not connected. This
contradicts Theorem~\ref{connect}.

Now, by Adjunction
\begin{equation}
\label{eq-K-F-}
{-}K_F\equiv\Diff F(C)\ge P'+P''.
\end{equation}
On the other hand,
\begin{multline}
\label{eq-K-F-1}
\deg K_F+\deg\Diff F(C)=(K_X+C)\cdot F\le 0,\\ \deg\Diff
F(C_1+C_2)\le 2.
\end{multline}
This gives that in \eref{eq-K-F-} and \eref{eq-K-F-1} equalities
hold. Hence $\Diff F(C)=P'+P''$ and $\Sing X=\{P',P''\}$ near $F$
and proves (i). By Theorem~\ref{kawamata_2}, $K_X+C+F$ is Cartier
near $F$. Since $K_X+C+F\equiv 0$, we have $K_X+C+F\sim 0$. This
proves (ii). (iii) easily follows by (ii). Further,
\[
\Diff
F(C)=\left(1-\frac1{m'}\right)P'+\left(1-\frac1{m''}\right)P'',\qquad
m', m''\in \NN.
\]
To show (iv) we just note that $\nu\colon X\to\PP^1$ is of type
$A^*$ of Theorem~\ref{Log-conic-bundles} near $F$ (because $K_X+F$
is plt). In particular, we have $m'=m''$.
\end{proof}

As in Corollary~\ref{m_i} using that $\Diff C(B)\in\Mm$ and
$\deg\Diff C(B)\le 2$, we have the following cases (up to
permutations of $P_1$, $P_2$):
\begin{equation}
\label{eq-case}
\Diff C(B)=\left\{
\begin{array}{l}
\alpha P_1+\beta P_2\\
 \alpha P'+\alpha P''\\
 \alpha P_1+\frac12P'+\frac12P''\\
 \frac12 P_1+\frac23P'+\frac23P'',\\
\end{array}
\right. \qquad
 \Supp B\subset\left\{
\begin{array}{l}
 F_1+F_2\\
 F\\
 F_1+F\\
 F_1+F,\\
\end{array}
\right.
\end{equation}
where $F:=\nu^{-1}(Q)_{\red}$ for some $Q\in\PP^1$, $Q\ne Q_1,
Q_2$, $\{P',P''\}=F\cap C$ and $\alpha,\beta\in\Mm$.
\par
Our strategy is very simple: we construct a boundary $B'\ge B$
such that $B'$ is contained in fibers of $\nu$, $K_X+C+B'$ is lc
and numerically trivial. If $B'\in \Mm$, then we can use proved
cases with $K_X+C+B\equiv 0$. If $B'\notin \Mm$, then we show that
$n(K_X+C+B')\sim 0$ for some $n\in\RRR_2$. The numerical
equivalence in $\Pic(X)$ coincides with linear one (see
Lemma~\ref{log-del-Pezzo-Pic}), it is sufficient to show only that
$n(K_X+C+B')$ is Cartier. Note that $K_X+C+B'$ is trivial on
fibers of $\nu$ (because $B'$ is contained in fibers). On the
other hand,
\[
(K_X+C+B')\cdot C=\deg(K_C+\Diff C(B'))=-2+\deg\Diff C(B').
\]
Since $\rho (X)=2$, we have
\begin{equation}
\label{eq-ravnos}
K_X+C+B'\equiv 0\quad \Leftrightarrow\quad \deg\Diff C(B')=2.
\end{equation}
\begin{lemma1}
\label{lemma-1}
Let $\nu\colon X\to Z\ni o$ be a germ of a contraction from a
surface to a curve, $F_1:=\nu^{-1}(o)_{\red}$, and $C\subset X$ a
germ of a curve such that $C\cap F_1$ is one point. Assume that
$\rho(X/Z)=1$, $K_X+C$ is plt and numerically trivial. Then there
is an $1$ or $2$-complement $K_X+C+\alpha'F_1$ (with
$\alpha'\in\{\frac12, 1\}$) such that $K_X+C+\alpha'F_1$ is not
plt at $C\cap F_1$.
\end{lemma1}

\begin{proof}
Note that a general fiber of $\nu$ is $\PP^1$ and $C$ is a
$2$-section of $\nu$. Take $\alpha$ so that $K_X+C+\alpha F_1$ is
maximally log canonical (i.e., $\alpha$ is maximal such that with
the log canonical property of $K_X+C+\alpha F_1$, see
\ref{def-maximally-lc}). Then $0<\alpha\le 1$. We claim that
$K_X+C+\alpha F_1$ is not plt at $C\cap F_1$. Indeed, otherwise
$\LCS(X,C+\alpha F_1)$ is not connected near $F_1$. This is a
contradiction with Proposition~\ref{dipol1}. In particular,
$(X/Z\ni o,C+\alpha F_1)$ is not exceptional. By
Theorem~\ref{bir-nef}, $K_X+C+\alpha F_1$ has an $1$, or
$2$-complement $K_X+C+R$ which is not plt at $C\cap F_1$. In
particular, $R\ne 0$. Since $C$ is a $2$-section of $\nu$, $R$ has
no horizontal components. Hence $R=F_1$ or $\frac12 F_1$.
\end{proof}

Now we consider possibilities of \eref{eq-case} step by step.

\subsection*{Subcase: $\Diff C(B)=\alpha P_1+\beta P_2$.}
Then $X$ is smooth outside of $F_1\cup F_2$. By
Lemma~\ref{lemma-1} there are $\alpha', \beta'\in\{\frac12,1\}$
such that $\alpha'F_1+\beta'F_2\ge B'$, $K_X+C+\alpha' F_1+\beta'
F_2$ is lc and not plt at $P_1$, $P_2$ and $2(K_X+C+\alpha'
F_1+\beta' F_2)\sim 0$ near $F_1$, $F_2$. By Adjunction (see
\ref{inv} and \ref{inv1}) $\Diff C(\alpha' F_1+\beta'
F_2)=P_1+P_2$. By \eref{eq-ravnos}, $K_X+C+\alpha' P_1+\beta'
P_2\equiv 0$. Moreover, $2(K_X+C+\alpha' F_1+\beta' F_2)$ is
Cartier near $F_1$ and $F_2$. Therefore $2(K_X+C+\alpha'
F_1+\beta' F_2)$ is Cartier on $X$ (because $\Sing X\subset
F_1\cup F_2$).

\subsection*{Subcase: $\Diff C(B)=\alpha P'+\alpha P''$.}
Then $C\cap F=\{P',P''\}$ for some fiber $F$ of $\nu$. By
Lemma~\ref{lemma-2} $K_X+C+F$ is lc. Since $\Diff C(F)=P'+P''$,
$K_X+C+F$ is numerically trivial. By the above cases with
$K_X+C+B\equiv 0$ there is a regular complement of $K_X+C+F$
(actually, $p_a(C+F)=1$ and we can use Case~\ref{(bb)} to show the
existence of an $1$-complement).

\subsection*{Subcase: $\Diff C(B)=\alpha
 P_1+\frac12P'+ \frac12P''$.}
By Lemma~\ref{lemma-1} there is $\alpha_1>0$ such that
$2(K_X+C+\alpha_1F_1+B)\sim 0$ near $F_1$, $K_X+C+\alpha_1F_1+B$
is lc and not plt at $P_1$. Since $\Diff
C(\alpha_1F_1+B)=P_1+\frac12P'+\frac12P''$,
$K_X+C+\alpha_1F_1+B\equiv 0$. Again by the above cases with
$K_X+C+B\equiv 0$ there is a regular complement of $K_X+C+B$ (more
precisely, $K_X+C+\alpha_1F_1+B$ is not plt, so we can use
Case~\ref{(ee)}).

\subsection*{Subcase:
$\Diff C(B)=\frac12P_1+\frac23P'+\frac23P''$.} Then we take $B'=
B+\alpha_1 F_1$ so that $\Diff
C(B')=\frac23P_1+\frac23P'+\frac23P''$. By \eref{eq-ravnos},
$K_X+C+B'\equiv 0$. We show that $6(K_X+C+B')$ is Cartier. First
note that $2(K_X+C)$ (and $2(K_X+C+B')$) is Cartier along $F_2$.
Indeed, $P_2\in X$ is smooth and $C\cdot F_2\in\NN$. On the other
hand, $\nu^{-1}(Q_2)\cdot C=2$. Hence the multiplicity $k$ of the
fiber $\nu^{-1}(Q_2)=kF_2$ is at most $2$ and $2F_2\sim 0$ near
$F_2$ over $\PP^1$. By Lemma~\ref{lemma-1} there is a
$2$-complement $K_X+C+\alpha_1F_2$ near $F_2$ which is not plt at
$P_2$. If $\alpha_1=1$, then $2(K_X+C+F_2)\sim 2(K_X+C)\sim 0$
near $F_2$. Assume that $\alpha_1=1/2$. Then $P_2=\Diff
C(\frac12F_2)=\frac12F_2|_C$. Hence $F_2\cdot C=2$ and the fiber
$\nu^{-1}(Q_2)=F_2$ is not multiple. So $F_2$ is Cartier and $X$
is smooth along $F_2$. This yields $2(K_X+C+\frac12F_2)\sim 0$
near $F_2$.

Write
\[
B=\gamma_1 F_1+\gamma F,\qquad \gamma_1, \gamma\in \Mm
\]
and
\[
\Diff
C(0)=\left(1-\frac1{m_1}\right)P_1+\left(1-\frac1{m}\right)P'+
\left(1-\frac1{m}\right)P'',\quad m_1, m\in\NN.
\]
Then by Corollary~\ref{coeff-Diff},
\[
1-\frac1{m_1}+\frac{\gamma_1n_1}{m_1}=\frac12,\qquad
1-\frac1{m}+\frac{\gamma n}{m}=\frac23,
\]
where $n_1,n\in \NN$. Since $\gamma_1, \gamma\ge 1/2$, we have
$n_1=n=1$. This gives only the following possibilities for $m$:
\begin{enumerate}
\item[]
$m=1$ (i.e., $X$ is smooth along $F$) and $\gamma=2/3$;
\item[]
$m=3$ (i.e., case $A^*$ of Theorem~\ref{Log-conic-bundles} with
$m=3$) and $\gamma=0$.
\end{enumerate}
It is easy to see that $3(K_X+C+B')$ is Cartier near $F$ in both
cases.

Now it is sufficient to show only that $6(K_X+C+B')\sim 0$ near
$F_1$. Similarly, we obtain only the following possibilities for
$m_1$:
\begin{enumerate}
\item[]
$m_1=1$ (i.e., $P_1\in X$ is smooth) and $\gamma_1=1/2$;
\item[]
$m_1=2$ (i.e., $P_1\in X$ is Du Val of type $A_1$) and
$\gamma_1=0$,
\end{enumerate}

Assume that $P_1\in X$ is smooth. Since $K_X+C+\frac12F_1$ is plt,
$F_1$ intersects $C$ transversally (see \ref{exer-odna-comp}).
Hence $B'=B+\frac16F_1$ and the coefficient of $F_1$ in $B'$ is
$2/3$. We claim that $6(K_X+C+\frac23F_1)\sim 0$ near $F_1$.
Indeed, $C\cdot F_1=1$ and $C\cdot \nu^{-1}(Q_1)=2$ implies that
the multiplicity of the fiber $\nu^{-1}(Q_1)$ is at most $2$ and
$2F_1\sim 0$ over $\PP^1$. On the other hand, by
Lemma~\ref{lemma-1} there is a $2$-complement $K_X+C+\alpha_1F_1$
near $F_1$ which is not plt at $P_1$. By our assumptions,
$K_X+C+\frac12F_1$ is plt. Thus $\alpha_1=1$ and $2(K_X+C+F_1)\sim
2(K_X+C)\sim 0$ near $F_1$. This yields
\[
6\left(K_X+C+\frac23F_1\right)\sim 4F_1\sim 0
\]
near $F_1$.

Now we assume that $P_1\in X$ is Du Val of type $A_1$. Then $F_1$
is not a component of $B$. As above, by Lemma~\ref{lemma-1} there
is a $2$-complement $K_X+C+\alpha_1F_1$ near $F_1$ which is not
plt at $P_1$. If $\alpha_1=1$, then $K_X+C+F_1$ is lc and by
Theorem~\ref{kawamata_2}, $C\cdot F_1=1/2$. It is easy to see that
in this case
\[
\Diff C\left(\frac13F_1\right)=\left(\frac12+\frac13C\cdot
F_1\right)P_1=\frac23P_1\quad \text{near}\quad P_1.
\]
Hence we can take $B'=B+\frac13F_1$. Then near $F_1$ we have
\[
6(K_X+C+B')=6(K_X+C+F_1)-4F_1\sim -4F_1.
\]
On the other hand, the multiplicity of the fiber $\nu^{-1}(Q_1)$
divides $4$ (because $\nu^{-1}(Q_1)\cdot C=2$ and $C\cdot
F_1=1/2$). This gives as $4F_1\sim 0$ and $6(K_X+C+B')\sim 0$ near
$F_1$.

\par
Finally, let $\alpha_1=1/2$. By Lemma~\ref{lemma-1}
$K_X+C+\frac12F_1$ is lc but not plt at $P_1$. Then
\[
P_1=\Diff C\left(\frac12F_1\right)=\left(\frac12+\frac12C\cdot
F_1\right)P_1\quad \text{near}\quad P_1.
\]
Hence $C\cdot F_1=1$ and as above, $2F_1\sim 0$. Similarly, we
obtain $B'=B+\frac16F_1$ and \[6(K_X+C+B')=6(K_X+C)+F_1\sim
-3F_1+F_1\sim 0.
\]
near $F_1$.
\par
This finishes the proof of Theorem~\ref{Inductive-Theorem}.
\end{proof}

\section{Corollaries}
The following form of Theorem~\ref{Inductive-Theorem} is very
important for applications.
\begin{corollary1}[see {\cite{Pr1}}]
\label{DelPezzo2a}
Let $(X,D')$ be a log del Pezzo such that $D'\in\Mm$. Assume also
that there exists a boundary $D\ge D'$ such that $-(K_X+D)$ is nef
and $K_X+D$ is not klt. Then $K_X+D'$ has a regular complement
which is not klt.
\end{corollary1}
\begin{proof}
If $K_X+D'$ is not klt, then there is a regular complement by
Proposition~\ref{Inductive-Theorem-Weak-Form}. From now on we
assume that $K_X+D'$ is klt. Replacing $D$ with suitable
$D'+\lambda(D-D')$ we also may assume that $K_X+D$ is lc (and not
klt).

First we consider the case when $\down{D}\ne 0$. Let $D_1$ be a
component of $\down{D}$. Replace $\delta_1$ with $1$:
\[
D'':= D_1+\sum_{i\ne 1}\delta_i D_i,\qquad D'\le D''\le D.
\]
If $-(K_X+ D'')$ is nef, then we can apply \ref{Inductive-Theorem}
(because $\down{D''}\ne 0$ and $D''\in\Mm$). Further, we assume
that $-(K_X+D'')$ is not nef. Then there exists a
$(K_X+D')$-nonpositive extremal ray $R$ such that $(K_X+D'')\cdot
R> 0$. If it is birational, then we contract it. Since $K_X+D$ is
nonpositive on $R$, this preserves the lc property of $K_X+D$ and
$K_X+D''$. We can pull back regular complements of $K_X+D''$
because $D''\in\Mm$ (now we are looking for regular complements of
$K_X+D''$, see Proposition~\ref{bir-prop} and
\ref{bir-prop-comment}). Note also that $(D''-D')\cdot R>0$.
Therefore $D_1$ is not contracted and on each step $K_X+D''$ is
not klt. If on some step $-(K_X+D'')$ is nef, we are done.
Otherwise continuing the process, we obtain a nonbirational
extremal ray $R$ on $X$ such that $(K_X+D'')\cdot R>0$. But on the
other hand,
\[
(K_X+D'')\cdot R\le (K_X+D)\cdot R\le 0,
\]
a contradiction.
\par
Consider now the case $\down{D}=0$. Then $K_X+D$ is lc, but is not
plt. Recall that $K_X+D'$ is klt. As in Proposition~\ref{plt} we
can construct a blowup $f\colon\widetilde X\to X$ with an
irreducible exceptional divisor $E$ such that $a(E, D)=-1$, the
crepant pull back
\[
K_{\widetilde X}+\widetilde D+E=f^*(K_X+D)
\]
is lc and $K_{\widetilde X}+\widetilde D^{\ep}+E$ is plt for any
$\widetilde D^{\ep}:=\widetilde D'+ \ep(\widetilde D-\widetilde
D')$, $\ep<1$. Here $\widetilde D$ and $\widetilde D'$ are proper
transforms of $D$ and $D'$, respectively. Note also that
$\rho(\widetilde X/X)=1$. Write
\[
K_{\widetilde X}+\widetilde D'+\alpha E=f^*(K_X+D'),
\]
where $\alpha<1$. Assume that there exists a curve $C$ such that
$(K_{\widetilde X}+\widetilde D'+E)\cdot C>0$. Then $(\widetilde
D-\widetilde D')\cdot C<0$. Therefore $C$ is a component of
$(\widetilde D-\widetilde D')$ and $C^2<0$. Further, $C\cdot E>0$.
Hence $C\ne E$ and we can choose $\widetilde D^{\ep}<\widetilde D$
so that $(K_{\widetilde X}+\widetilde D^{\ep})\cdot C<0$.
Therefore $C$ is a $(K_{\widetilde X}+\widetilde D')$-negative
extremal curve and its contraction preserves the lc property of
$K_{\widetilde X}+\widetilde D+E$. Again we can pull back
complements of $K_{\widetilde X}+\widetilde D'+E$ (see
Proposition~\ref{bir-prop} and \ref{bir-prop-comment}). Repeating
the process, we get the situation when $-(K_{\widetilde
X}+\widetilde D'+E)$ is nef. All the steps preserve the nef and
big property of $-(K_{\widetilde X}+\widetilde D'+\alpha E)$. If
$\alpha\ge 0$, then we apply Theorem~\ref{Inductive-Theorem} to
$(\widetilde X,\widetilde D'+E, \widetilde D'+\alpha E)$. If
$\alpha<0$, then by the monotonicity, $-(K_{\widetilde
X}+\widetilde D')$ is nef and big. Again apply
Theorem~\ref{Inductive-Theorem} to $(\widetilde X,\widetilde D'+E,
\widetilde D')$. This concludes the proof of the corollary.
\end{proof}

\begin{corollary1}
\label{corr-Inductive}
Let $(X,B)$ be a log del Pezzo surface such that $B\in\Mm$. Then
$(X,B)$ is nonexceptional if and only if there exists a regular
complement $K_X+B^+$ which is not klt.
\end{corollary1}

\begin{corollary1}[cf. {Corollary~\ref{>4}}]
\label{cor-=4}
Let $(X,B)$ be a log del Pezzo surface. Assume that $B\in\Mm$ and
$(K_X+B)^2\ge 4$. Then there exists a regular complement of
$K_X+B$. Moreover, there exists such a complement which is not klt
(in particular, $(X,B)$ is nonexceptional).
\end{corollary1}
\begin{outline}
If $K_X+B$ is not klt, this assertion follows by
Proposition~\ref{Inductive-Theorem-Weak-Form}. Assume that $K_X+B$
is klt. Take $n\in\NN$ so that $H:=-n(K_X+B)$ is an integral
(ample) Cartier divisor. By Riemann-Roch
\begin{multline*}
\dim H^0(X,\OOO_X(H))\ge \frac{H\cdot
(H-K_X)}{2}+1=\frac12n(n+1)(K_X+D)^2-\\ \frac12n(K_X+D)\cdot
D+1\ge 2n(n+1)+1.
\end{multline*}
Pick a smooth point $P\in X$ and let $\mathfrak{m}_P$ be the ideal
sheaf of $P$. Then
\[
\dim H^0(X,\OOO_X/ \mathfrak{m}_P^{2n})=\dim
\CC[x,y]/(x,y)^{2n}=\frac{(2n+1)2n}{2}=(2n+1)n.
\]
From the exact sequence
\begin{multline*}
0\longrightarrow \mathfrak{m}_P^{2n}\otimes\OOO_X(H)
\longrightarrow \OOO_X(H) \longrightarrow \OOO_X/
\mathfrak{m}_P^{2n}\otimes \OOO_X(H)\simeq \OOO_X/
\mathfrak{m}_P^{2n} \longrightarrow 0
\end{multline*}
we see that
\[
H^0(X,\mathfrak{m}_P^{2n}\otimes\OOO_X(H))\neq 0.
\]
Therefore there is $H'\in |H|$ such that $\mult_P(H')\ge 2n$. It
is easy to see that $K_X+B+\frac1nH'$ is not klt. Then
$K_X+B+\alpha H'$ is lc but not klt for some $\alpha\le \frac1n$.
Clearly, $-(K_X+B+\alpha H')$ is nef. Hence we can apply
Corollary~\ref{DelPezzo2a}.
\end{outline}

Note that the above result can be improved: by taking $P\in
\Supp(B)$ or $P\in \Sing(X)$ it is possible to find nonklt
boundary $K_X+B+\alpha H'$ for smaller values of $(K_X+D)^2$. On
the other hand, in the case $B=0$ and $X$ is smooth, it is well
known that $K_X$ is strongly $1$-complementary.

\begin{corollary1}[see {\cite{Pr1}}]
\label{corr-Inductive-sing}
Let $X\ni P$ be a three-dimensional klt singularity, $f\colon
(Y,S)\to X$ a plt blowup, and $K_X+D$ an $n$-complement which is
not klt at $P$. Then one of the following holds
\begin{enumerate}
\item
$a(S,D)=-1$ and $K_Y+S+D_Y:=f^*(K_X+D)$ is an $n$-complement of
$K_Y+S$;
\item
$a(S,D)>-1$ and then there exists a regular complement of $K_Y+S$
which is not plt.
\end{enumerate}
\end{corollary1}
\begin{proof}
(i) is obvious. Assume that $a(S,D)>-1$. Write
\[
K_Y+aS+D_Y:=f^*(K_X+D),
\]
where $a=-a(S,D)<1$ and $D_Y$ is the proper transform of $D$. By
assumptions $K_Y+aS+D_Y$ is lc and not klt (see
\ref{crepant-singul}). Therefore $K_Y+S+D_Y$ is not plt and we can
take $0<b\le 1$ so that $K_Y+S+bD_Y$ is lc but not plt. It is easy
to see that $-(K_Y+S+bD_Y)$ is $f$-ample. If $f(S)=P$, then
$(S,\Diff S(bD_Y))$ is a log Del Pezzo. By \ref{DelPezzo2a} (or by
\ref{Inductive-Theorem-Weak-Form}) there is a regular complement
of $K_S+\Diff S(0)$ and by \ref{prodolj} it can be extended to a
complement of $K_Y+S$. Similarly, in the case when $f(S)$ is a
curve, we can use Theorem~\ref{local}.
\end{proof}

Similar to Theorem~\ref{Inductive-Theorem} one can prove the
following
\begin{proposition1}[cf. {\cite{Bl}}]
Let $(X,D)$ be a log Enriques surface (i.e., $K_X+D$ is lc and
$K_X+D\equiv 0$). Assume that $K_X+D$ is not klt and $D\in \Mm$.
Then $n(K_X+D)\sim 0$ for some $n\in\RRR_2$. In particular, $D\in
\Msm$.
\end{proposition1}
\begin{proof}
By \ref{important_compl} it is sufficient to show the existence of
a regular complement. If $X$ is not rational, then the assertion
follows by Theorem~\ref{Inductive-Theorem-nonrat}. Otherwise we
can apply Theorem~\ref{Inductive-Theorem}. The existence of $D'$
in conditions of the theorem follows by
Proposition~\ref{log-Enriques->log-del-Pezzo} below.
\end{proof}

\begin{proposition}
\label{log-Enriques->log-del-Pezzo}
Let $(X,\Lambda)$ be a log Enriques surface. Assume that $X$ is
rational. Then there is a boundary $\Lambda'\le \Lambda$ such that
\begin{enumerate}
\item
$K_X+\Lambda'$ is klt, and
\item
$-(K_X+\Lambda')$ is nef and big.
\end{enumerate}
\end{proposition}
\begin{proof}
Replace $X$ with its minimal resolution and $\Lambda$ with its
crepant pull back. Then again $(X,\Lambda)$ is a log Enriques
surface. By Corollary~\ref{crepant-contraction} it is sufficient
to construct $\Lambda'$ on this new $X$. Further, there is a
sequence of contractions of $-1$-curves $\varphi\colon X\to
X^{\bullet}$, where $X^{\bullet}\simeq\PP^2$ or
$X^{\bullet}\simeq\FF_n$, $n\ge 0$, $n\neq 1$. Put
$\Lambda^{\bullet}:=\varphi_*\Lambda$. Then
$(X^{\bullet},\Lambda^{\bullet})$ is again a log Enriques surface.
It is sufficient to construct $\Lambda^{\prime\bullet}\le
\Lambda^{\bullet}$ such that
$K_{X^{\bullet}}+\Lambda^{\prime\bullet}$ is klt and
$-(K_{X^{\bullet}}+\Lambda^{\prime\bullet})$ is nef and big.
Indeed, the crepant pull back $\Lambda''$ of
$\Lambda^{\prime\bullet}$ satisfies (i) and (ii). However,
$\Lambda''$ is not necessarily a boundary (i.e., effective). To
avoid this one can take $\Lambda'=\Lambda''+t(\Lambda-\Lambda'')$
for $0<1-t\ll 1$.

Further, if $X^{\bullet}\simeq\PP^2$ or
$X^{\bullet}\simeq\PP^1\times\PP^1$, then we take
$\Lambda^{\prime\bullet}=0$. In the case $X^{\bullet}\simeq \FF_n$
with $n\ge 2$, we write $\Lambda^{\bullet}= \lambda\Sigma_0
+\Lambda^{\circ}$, where $\Sigma_0$ is the negative section of
$X^{\bullet}=\FF_n$, $0\le \lambda\le 1$, $\Lambda^{\circ}\ge 0$,
and $\Sigma_0$ is not a component of $\Lambda^{\circ}$. It is easy
to see
\[
2-n={-}K_{X^{\bullet}}\cdot \Sigma_0=\Lambda^{\prime\bullet}\cdot
\Sigma_0= -n\lambda+ \Lambda^{\circ}\cdot \Sigma_0\ge -n\lambda.
\]
Hence, $\lambda\ge 1-2/n$. Thus we can take
$\Lambda^{\prime\bullet}=(1-2/n)\Sigma_0$.
\end{proof}

\section{Characterization of toric surfaces}
Following Shokurov we prove Conjecture~\ref{conjecture-toric} in
dimension two. Moreover, we prove a generalization of
\ref{conjecture-toric} for $\rhonum$ instead of $\mt{rk}\Weilalg$
(recall that $\rhonum(X)$ is the rank of the quotient of
$\Weil(X)$ modulo numerical equivalence).
\begin{theorem1}[{\cite{Sh1}}]
\label{th-toric}
Let $(X,D=\sum d_iD_i)$ be a projective log surface such that
\begin{enumerate}
\item
$K_X+D$ is lc, and
\item
$-(K_X+D)$ is nef.
\end{enumerate}
Then
\begin{equation}
\label{eq-main-ineq-main}
\sum d_i\le \rhonum(X)+2.
\end{equation}
If the equality holds, then $K_X+D\equiv 0$ and $X$ has only
rational singularities (in particular $X$ is $\QQ$-factorial).
\end{theorem1}

\begin{proof}
Assume that
\begin{equation}
\label{eq-main-ineq}
\sum d_i-\rhonum (X)-2\ge 0.
\end{equation}
\subsection*{}
First we consider the case $K_X+D\equiv 0$.
\subsection*{Step 0}
Apply a minimal log terminal modification as in \ref{log-term}. It
is easy to see that this preserves the left hand side of
\eref{eq-main-ineq}. Thus we may assume that $K_X+D$ is dlt. In
particular, $K_X+\fr{D}$ is klt, $X$ is $\QQ$-factorial and
$\rhonum(X)=\rho(X)$.

\subsection*{Step 1}
Write $D=C+B$, where $C:=\down{D}$ and $B:=\fr{D}$. Then
$-(K_X+B)\equiv C+(\text{nef divisor})$. Hence $K_X+B$ cannot be
nef. Run $(K_X+B)$-MMP, i.e., contract birational extremal rays
$R$ such that $R\cdot (K_X+B)<0$. The left hand side of
\eref{eq-main-ineq} does not decrease. Of course, we can lose the
dlt property of $K_X+D$, but properties (i)--(ii) are preserved.
Moreover, if $C\ne 0$, then on each step we contract a curve $R$
with $R\cdot C>0$. In particular, whole $C$ is not contracted.

At the end we get a nonbirational contraction $\var\colon X\to Z$.

\subsection*{Step 2}
Assume that after Step 1 we get a Fano contraction $\var$ with
$\dim (Z)=1$. Write $D=D^{\mt{vert}}+D^{\mt{hor}}$, where
$D^{\mt{hor}}=\sum_{\mt{hor}} d_iD_i$ is the sum of all components
such that $\var(D_i)=Z$ and $D^{\mt{vert}}=\sum_{\mt{vert}}
d_jD_j$ is the sum of components which are fibers of $\var$. Let
$F$ be a general fiber of $\var$. Then by Adjunction
\begin{multline*}
0\ge (K_X+D)\cdot F=(K_X+D^{\mt{hor}}+F)\cdot F\\ =\deg
K_F+D^{\mt{hor}}\cdot F \ge -2+D^{\mt{hor}}\cdot F.
\end{multline*}
This gives $D^{\mt{hor}}\cdot F\le 2$. In particular,
\begin{equation}
\label{eq-main-ineq-1}
\sum_{\mt{hor}} d_i\le 2,\qquad \sum_{\mt{vert}} d_j\ge 2
\end{equation}
(because $\rho(X)=2$). Now, let $R$ be the extremal ray of
$\NE(X)$ other than $\RR_+[F]$. Then $R\cdot D^{\mt{vert}}>0$.
Hence $R\cdot (K_X+(1-\ep_1)D-\ep D^{\mt{vert}})<0$ for
$0<\ep_1\ll \ep\ll 1$. By Contraction Theorem there is a
contraction $\psi\colon X\to Z_1$ of $R$.

Assume that $\dim Z_1=1$. Then, as above, we have
$\sum_{\mt{vert}} d_j\le 2$ (because components of $D^{\mt{vert}}$
are horizontal with respect to $\psi$). This yields equalities in
\eref{eq-main-ineq-1} and \eref{eq-main-ineq}. If $\dim Z_1=2$,
then $\psi$ is birational. Let $E$ be the $\psi$-exceptional
divisor. If $E$ is a component of $C:=\down D$, then again by
Adjunction we have
\[
0\ge (K_X+D)\cdot E=\deg K_E+\deg\Diff E(D-E),\quad \deg\Diff
E(D-E)\le 2.
\]
Since any component of $D^{\mt{vert}}$ meets $E$, by
Corollary~\ref{coeff-Diff} we obtain $\sum_{\mt{vert}} d_j\le 2$.
This yields equalities in \eref{eq-main-ineq-1} and
\eref{eq-main-ineq}. Finally, if $E$ is not a component of
$C:=\down D$, then we replace $X$ with $Z_1$. Note that in this
case we get strict inequality $\sum d_i-\rho (X)-2>0$ in
\eref{eq-main-ineq} and $\rho(X)=1$. By the next two steps this is
a contradiction.

\subsection*{Step 3}
Assume that $Z$ is a point (and $\rho (X)=1$). Then $\sum d_i\ge
3$. We claim that after perturbation of coefficients one can
obtain the case when $K_X+D$ is not klt. Indeed, assume that
$K_X+D$ is klt. Let $H$ be the ample generator of $\Pic(X)\simeq
\ZZ$ (see Lemma~\ref{log-del-Pezzo-Pic}) and let $D_i\equiv a_iH$,
$a_i>0$. Without loss of generality we may assume that $a_1\le
a_2\le \cdots$. Take $D^t=D+t(D_i-D_j)$, where $i<j$ and $0<t\le
d_j$. Clearly, $-(K_X+D^t)$ is again nef and $D^t$ is effective.
Moreover, for $D^t$ the left hand side of \eref{eq-main-ineq}
remains the same. If $K_X+D^{t_0}$ is lc but not klt for some for
$0\le t_0\le d_j$, then $D^{t_0}$ gives the required boundary. If
$K_X+D^t$ is klt for $t=d_j$, then we replace $D$ with $D^{d_j}$
continue the process with another pair $D_i, D_j$. Since the last
procedure reduces the number of components of $D$, this process
terminates. At the end we get the situation when $K_X+D$ is not
klt.

\subsection*{Step 4}
Now we consider the case when $K_X+D$ is not klt. Apply steps~0--2
again. On Step~2 in \eref{eq-main-ineq} the equality holds. So we
assume that $\rho(X)=1$ and $C:=\down{D}\ne 0$. For any component
$C_i\subset C$ by Adjunction we have
\begin{equation}
\label{eq-toric-step4-1}
2\ge-\deg K_{C_i}=\deg\Diff{C_i}(D-C_i)
\end{equation}
On the other hand, all components of $D$ intersect $C_i$ and
\begin{equation}
\label{eq-toric-step4-2}
\deg\Diff{C_i}(D-C_i)\ge \sum d_j-1\ge 2
\end{equation}
(see Corollary~\ref{coeff-Diff}). Therefore $\sum d_j=3$ and $\deg
K_{C_i}=2$. This completes the proof in the case $K_X+D\equiv 0$.

\subsection*{}
Consider the case $K_X+D\not\equiv 0$. As in Step 0 we may assume
that $K_X+D$ is dlt. Further, similar to Step 1, run
$(K_X+D)$-MMP. This preserves assumption \eqref{eq-main-ineq}. Let
$\var\colon X\to X'$ be a birational extremal contraction,
$D':=\var_*D$, and $E$ the exceptional curve. Clearly,
$K_X+D\equiv \var^*(K_{X'}+D')+aE$, where $a\in \QQ$. Then
$0>(K_X+D)\cdot E=aE^2$. Hence $a>0$. Assume that $K_{X'}+D'\equiv
0$ and $H$ a hyperplane section of $X$. Then $0\ge ( K_X+D)\cdot
H=aE\cdot H>0$, a contradiction. Therefore $K_{X'}+D'\not\equiv
0$. We can replace $(X,D)$ with $(X',D')$ and continue the
process. At the end we get a log surface $(X,D)$ with a
nonbirational $(K_X+D)$-negative extremal contraction $\phi\colon
X\to Z$. In particular $\rho (X)\le 2$. If $Z$ is a point, then
$\rho(X)=1$ and $-(K_X+D)$ is ample. Take $n\ge 0$ so that the
divisor $-n(K_X+D)$ is integral and very ample. Let $G\in
|-n(K_X+D)|$ be a general member. Then $K_X+D+\frac1nG$ is dlt and
numerically trivial. In this case, by the proved inequality
\eqref{eq-main-ineq-main}, $3\le \sum d_i< \sum d_i+1/n\le 3$, a
contradiction. If $Z$ is a curve, then we can use the arguments of
Step 2. Thus $4\le \sum d_i=\sum_{\mt{hor}} d_i+\sum_{\mt{vert}}
d_j< 4$ (because we have strict inequality in
\eqref{eq-main-ineq-1}). The last contradiction proves that
$K_X+D\equiv 0$.
\par
Assume that $X$ has at least one nonrational singularity $P\in X$.
Clearly, $K_X$ is not klt at $P$ and $P\notin \Supp(D)$. Then by
Corollary \ref{cor-rat-nonrat} $P\in X$ is a simple elliptic or
cusp singularity. As in Step 0, let $\var \colon (\tilde X,\tilde
D)\to (X,D)$ be a minimal log terminal modification. If
$\down{\tilde D}$ is connected, then $\down{\tilde
D}=\var^{-1}(P)$ and $p_a\left(\down{\tilde D}\right)=1$. By
Lemma~\ref{wheel-1} $\tilde D=\down{\tilde D}$ and $D=0$, a
contradiction. Therefore $C:=\var^{-1}(P)$ is a connected
component of $\down{\tilde D}$. Denote $\tilde B:=\tilde D-
\down{\tilde D}$ and $C':=\down{\tilde D}-C$. Our assumption
\eqref{eq-main-ineq} implies that $\tilde B\neq 0$. Then $C\cap
C'=C\cap \Supp(\tilde B)=\emptyset$. By Proposition~\ref{dipol2}
there is a contraction $\psi\colon \tilde X\to Z$ with rational
fibers onto a curve $Z$ such that $C$ and $C'$ are (smooth)
disjoint sections. Then $\tilde B$ has no horizontal components.
Let $R$ be a $(K_{\tilde X}+\tilde D-\ep\tilde B)$-negative
extremal rational curve. Since $p_a(Z)=p_a(C)=1$, $R$ cannot be
horizontal. On the other hand, $\rho(\tilde X)\ge 3$ (because
$R\cdot \tilde B>0$). Therefore the contraction of $R$ is
birational. This contraction reduces the left hand side of
\eqref{eq-main-ineq}, a contradiction. This proves
Theorem~\ref{th-toric}.
\end{proof}

\begin{theorem1}[{\cite{Sh1}}]
\label{th-toric-1}
Let $(X,C)$ be a projective log surface with a reduced boundary
$C=\sum_{i=1}^r C_i$ such that $K_X+C$ is lc and $-(K_X+C)$ is
nef. Assume also
\[
r\ge \rhonum(X)+2.
\]
Then
\begin{enumerate}
\item
$r=\rhonum(X)+2=\rho(X)+2$;
\item
$K_X+C\sim 0$ (i.e., $K_X+C$ is $1$-complementary);
\item
$C$ is connected and $p_a(C)=1$;
\item
the pair $(X,C)$ is toric.
\end{enumerate}
\end{theorem1}
\begin{proof}
The assertion (i) follows by Theorem~\ref{th-toric}. This also
shows that $K_X+C\equiv 0$. To prove (ii) we apply steps~0--4 of
the proof of Theorem~\ref{th-toric} to $(X,C)$. At the end we
obtain one of the following:
\par \begin{itemize}
\item
$\rho(X)=1$ and $C$ has exactly three components. Clearly they
intersect each other and does not pass through one point, so
$p_a(C)\ge 1$ (cf. Proof of Corollary~\ref{delta<3}). By
Lemma~\ref{wheel-1}, $K_X+C$ is a $1$-complement. According to
\ref{bir-prop-a} and \ref{bir-prop}, $K_X+C$ on our original $X$
is $1$-complementary.
\item
$\rho(X)=2$ and $C$ has exactly four components. Moreover, there
is an extremal contraction $\var\colon X\to Z$ onto a curve. By
discussions in Step~2 of the proof of Theorem~\ref{th-toric}
(especially, \eref{eq-main-ineq-1}), we have a decomposition
$C=C^{\mt{hor}}+C^{\mt{vert}}$ such that both $C^{\mt{hor}}$ and
$C^{\mt{vert}}$ have two irreducible components. Any component of
$C^{\mt{vert}}$ meets all components of $C^{\mt{hor}}$. As in
\eref{eq-nad-posl-1} we have $p_a(C)\ge 1$. Finally, as above,
$K_X+C$ is $1$-complementary.
\end{itemize}

By Proposition~\ref{dipol2}, $C\subset \LCS(X,C)$ is connected.
Since $K_X+C$ is Cartier, $\Diff C(0)=0$ (see \ref{diff}). Thus
$K_C=0$ and $p_a(C)=1$. This proves (iii). The assertion of (iv)
follows by the lemma below.
\end{proof}

\begin{lemma1}
\label{toric-r=1}
Let $(X,C)$ be a projective log surface such that $C$ is reduced
and connected. Assume that $C$ has exactly $\rhonum(X)+2$
components, $K_X+C$ is lc and linearly trivial. Then $(X,C)$ is a
toric pair.
\end{lemma1}
\begin{proof}
Let $\var\colon \ov X\to X$ be minimal lt modification of $(X,C)$.
Write $\var^*(K_X+C)=K_{\ov X}+\ov C$, where $\ov C$ is reduced
and $\var_*(\ov C)=C$. The exceptional divisor of $\var$ is
contained in $\ov C$. Hence on $(\ov X,\ov C)$ all our conditions
hold. So it is sufficient to prove our assertion for $(\ov X,\ov
C)$. By Lemma~\ref{wheel-1}, $p_a(\ov C)=1$, $\ov C$ a wheel of
smooth rational curves, $\ov X$ is smooth along $\ov C$ and has
only Du Val singularities outside. Run $K_{\ov X}$-MMP. By
\eref{eq-main-ineq-main}, on each step we contracted a component
of $\ov C$ (which is contained into the smooth locus of $\ov X$).
Thus our MMP is a sequence of contractions of $-1$-curves. At the
end we obtain a Fano contraction $\psi\colon (\widehat X,\widehat
C)\to Z$, where $\widehat X$ has only Du Val singularities,
$K_{\widehat X}+\widehat C$ is lc (in fact it is analytically dlt)
and numerically trivial. The sequence of transformations $\ov X\to
\widehat X$ is a sequence of blowups of divisors with
discrepancies $a(\cdot, \widehat C)=-1$. They must preserve the
action of a two-dimensional torus (if $(\widehat X,\widehat C)$ is
a toric pair). Thus it is sufficient to show that the pair
$(\widehat X,\widehat C)$ is toric.
\par
If $\rho(\widehat X)=1$, then $C$ has exactly three components
which are Cartier divisors. Therefore $\widehat X$ is a log del
Pezzo of Fano index $r(\widehat X)\ge 3$ (see \ref{def-Fano-ind}).
By Lemma~\ref{del-Pezzo-index>1}, $\widehat X\simeq\PP^2$ and
$\widehat C=\widehat C_1+\widehat C_2+\widehat C_3$, where the
$\widehat C_i$ are lines. Obviously, $(\widehat X,\widehat C)$ is
toric in this case. Finally, assume that $\dim Z=1$. Then
$\widehat C$ has exactly four components and by
Lemma~\ref{wheel-1} they form a wheel of smooth rational curves.
It is an easy exercise to prove that $Z\simeq\PP^1$ and the fibers
of $\psi$ are rational. Therefore
$C=C^{\mt{hor}}_1+C^{\mt{hor}}_2+
C^{\mt{vert}}_1+C^{\mt{vert}}_2$, where $C^{\mt{hor}}_1$,
$C^{\mt{hor}}_2$ are disjoint sections of $\psi$ and
$C^{\mt{vert}}_1$, $C^{\mt{vert}}_2$ are fibers. We claim that
$\widehat X$ is smooth. Indeed, by construction, $\widehat X$ is
smooth along $\widehat C$. Let $F$ be a fiber of $\psi$ different
from $C^{\mt{vert}}_1$, $C^{\mt{vert}}_2$. Take $c$ so that
$K_{\widehat X}+\widehat C+cF$ is maximally lc. If $c<1$, then
$\LCS(\widehat X,\widehat C+cF)$ has three connected components
near $F$. This contradicts Proposition~\ref{dipol1}. Hence
$K_{\widehat X}+\widehat C+F$ is lc. By Adjunction
$\deg\Diff{F}(\widehat C)=2$. On the other hand,
$\Diff{F}(\widehat C)\ge P_1+P_2$, where $P_i=F\cap
C^{\mt{hor}}_i$. Hence $\Diff{F}(\widehat C)= P_1+P_2$ and
$\widehat X$ is smooth along $F$.
\par
Thus we have shown that $\widehat X$ is smooth. Then $\widehat
X\simeq\FF_n$ and $\psi$ is the natural projection
$\FF_n\to\PP^1$. Since $C^{\mt{hor}}_1$, $C^{\mt{hor}}_2$ are
disjoint sections one of them, say $C^{\mt{hor}}_1$, must be the
minimal section $\Sigma_0$. Now, it is easy to show that the pair
$(\widehat X,\widehat C)$ is toric.
\end{proof}

\chapter{Boundedness of exceptional complements}
\label{sect-9}
\section{The main construction}
\label{notation-boundedness}
In this chapter we discuss some boundedness results for
exceptional log surfaces. The main results are
Theorems~\ref{main_Sh_A-1}, \ref{main_Sh_A-2} and \ref{main_Sh_A}.
Fix the following notation.
\par
Let $(X,B)$ be a projective log surface such that
\begin{enumerate}
\item
$K_X+B$ is lc;
\item
$-(K_X+B)$ is nef;
\item
the coefficients of $B$ are standard or $d_i\ge 6/7$ (i.e.
$B\in\Mm$);
\item
$(X,B)$ is exceptional, i.e., any $\QQ$-complement of $K_X+B$ is
klt;
\item
there is a boundary $B^{\triangledown}\le B$ such that
$K_X+B^{\triangledown}$ klt and $-(K_X+B^{\triangledown})$ is nef
and big.
\end{enumerate}

By Corollary~\ref{corr-Inductive}, (iv) is equivalent to
\begin{enumerate}
\item[(iv)${}'$]
there are no regular nonklt complements.
\end{enumerate}
Note also that by Theorem~\ref{Inductive-Theorem} and
Corollary~\ref{Corr-Inductive-Theorem-discr}, $K_X+B$ is klt.

\subsection{}
\label{1/7-lt}
If $K_X+B$ is $(1/7)$-lt and $b_i<6/7$, $\forall i$, then by
Theorem~\ref{Al}, $X$ belongs to a finite number of families. In
this case $b_i\in\{0, \frac12, \frac23, \frac34, \frac45,
\frac67\}$, so $(X,B)$ is bounded. Therefore we have a finite
number of exceptional values of $\compl(X,B)$ in this case.

\subsection{}
\label{not-1/7-lt}
Now we assume that there is at least one (exceptional or not)
divisor with $a(\cdot,B)\le -1+1/7$. Following Shokurov \cite{Sh1}
we construct a model of $(X,B)$ with $\rho=1$. Similar birational
modifications were used in \cite{KeM} and was called the
\textit{hunt for a tiger}. Let $\mu\colon\widehat X\to X$ be the
blowup of all exceptional divisors with $a(\cdot,B)\le-6/7$ (see
Lemma~\ref{finite-term} and Proposition~\ref{per-gen}). Consider
the crepant pull back
\[
K_{\widehat X}+\widehat{B}=\mu^*(K_X+B),\quad\text{with}\quad
\mu_*\widehat{B}=B.
\]
Then $K_{\widehat X}+\widehat{B}$ is $(1/7)$-lt. By construction,
$\widehat{B}\in\Mm$. As in \ref{construction-boundary} we can
construct a boundary $\widehat{B}^{\triangledown}\le \widehat{B}$
on $\widehat X$ such that $K_{\widehat
X}+\widehat{B}^{\triangledown}$ is klt and $-(K_{\widehat
X}+\widehat{B}^{\triangledown})$ is nef and big. So on $\widehat
X$ all our assumptions (i) -- (v) hold and moreover,
\begin{enumerate}
\item[(i)${}'$]
$K_{\widehat X}+\widehat{B}$ is $(1/7)$-lt and
$\down{\frac76\widehat{B}}\ne 0$.
\end{enumerate}
By Lemma~\ref{Kodaira-lemma} the Mori cone $\NE(\widehat X)$ is
polyhedral and generated by contractible extremal curves.
Moreover, if an extremal ray $R$ on $\widehat X$ is birational,
then the contraction preserves assumptions (i)--(iii), (v) (see
\ref{crepant-singul}).

Write $\widehat{B}=\sum b_i\widehat{B}_i$. Define the boundary
$\widehat{D}=\sum d_i\widehat{B}_i$, where
\begin{equation}
\label{def_D}
d_i= \left\{
\begin{array}{ll}
 1\qquad&\text{if}\ b_i\ge 6/7;\\
 b_i&\text{otherwise.}\\
\end{array}
\right.
\end{equation}
In particular, $\Supp \widehat{D}=\Supp \widehat{B}$. Put
\[
\widehat{C}:=\down{\frac76\widehat{B}}=\sum_{b_i\ge 6/7}
\widehat{B}_i\qquad \text{and}\qquad
\widehat{F}:=\widehat{D}-\widehat{C}=\sum_{b_i<6/7}
b_i\widehat{B}_i.
\]
Then $\widehat{D}=\widehat{C}+\widehat{F}$,
$\widehat{C}=\down{\widehat{D}}$ and $\widehat{F}\in \{0,\frac12,
\frac23,\frac34, \frac45, \frac56\}$. Since $\widehat{D}\in\Msm$,
so $\widehat{D}\in\PPP_n$ for all $n\in\NN$ (see
\ref{def_PPP_n_sect}).

\begin{lemma1}
\label{lc-exc-prop}
$K_{\widehat X}+\widehat{D}$ is lc.
\end{lemma1}
\begin{proof}
By Theorem~\ref{local} the log divisor $K_{\widehat
X}+\widehat{B}$ has a regular $n$-complement $\widehat{B}^+$ near
each point $x\in {\widehat X}$. Then $K_{\widehat
X}+\widehat{B}^+$ is lc and $\widehat{B}^+\ge \widehat{B}$ (see
\ref{important_compl}). Moreover, if $b_i\ge 6/7$, then by
definition of complements $b_i^+\ge \frac1n\down{(n+1)b_i}\ge 1$.
Therefore $\widehat{B}^+\ge \widehat{D}$. This gives that
$K_{\widehat X}+\widehat{D}$ is also lc at $x$.
\end{proof}

Recall that we have assumed $\down{\widehat{D}}\ne 0$. By the
Inductive Theorem~\ref{Inductive-Theorem}, $-(K_{\widehat
X}+\widehat{D})$ is not nef. Assume that $\rho({\widehat X})>1$.
By Lemma~\ref{Kodaira-lemma} the cone $\NE({\widehat X})$ is
polyhedral and generated by contractible curves. Therefore there
exists an extremal ray $R$ which is positive with respect to
$K_{\widehat X}+\widehat{D}$. Recall that $-(K_{\widehat
X}+\widehat{B})$ is nef, so $(K_{\widehat X}+\widehat{B})\cdot
R\le 0$.

\begin{lemma1}
\label{lemma-birat-E}
The contraction of $R$ is birational.
\end{lemma1}
\begin{proof}
Indeed, otherwise $\rho(\widehat X)=2$ and we have an extremal
contraction $\widehat X\to Z$ onto a curve, positive with respect
to $K_{\widehat X}+\widehat{D}$. Let $F$ be a general fiber. Then
\[
\deg(K_F+\widehat{B}|_F)=(K_{\widehat X}+\widehat{B})\cdot F\le 0,
\]
and by \ref{m_i} and by \ref{=}, $K_F+\widehat{B}|_F$ is $1$, $2$,
$3$, $4$, or $6$-complementary. But then (as in the proof of
Lemma~\ref{lc-exc-prop}) $(\widehat{B}|_F)^+\ge \widehat{D}|_F$.
Hence
\[
0=\deg(K_F+\widehat{B}^+|_F)\ge \deg(K_F+\widehat{D}|_F)
=(K_{\widehat X}+\widehat{D})\cdot F> 0,
\]
which is a contradiction.
\end{proof}

Thus $R$ is of birational type and generated by a (rational)
curve, say $E$.

\begin{lemma1}
$E$ is not a component of $\down{\widehat{D}}$.
\end{lemma1}
\begin{proof}
By Theorem~\ref{local} and by Corollary~\ref{two} there exists a
regular complement of $K_{\widehat X}+\widehat{B}$ near $E$. As in
the proof of Lemma~\ref{lemma-birat-E}, $\widehat{B}^+\ge
\widehat{D}$. If $E$ is a component of $\down{\widehat{D}}$, then
$0=(K_{\widehat X}+\widehat{B}^+)\cdot E\ge(K_{\widehat
X}+\widehat{D})\cdot E>0$, a contradiction.
\end{proof}

By Lemma~\ref{ep-lt-contractions} below we have that if the
contraction ${\widehat X}\to {\widehat X}'$ of $E$ is birational,
then $K_{{\widehat X}'}+\widehat{B}'$ is also $(1/7)$-lt. Thus we
can make a birational contraction ${\widehat X}\to {\widehat X}'$.
On ${\widehat X}'$ all our assumptions (i)${}'$, (ii)--(v) hold.
We replace ${\widehat X}$ with ${\widehat X}'$.

\begin{lemma1}[cf. e.g., {\cite[3.38]{KM}}]
\label{ep-lt-contractions}
Let $(X,B)$ be an $\ep$-lt pair and $f\colon X\to X'$ a divisorial
contraction of an $(K_X+B)$-nonpositive extremal ray. Assume that
$X$ is $\QQ$-factorial and $f$ does not contract components of $B$
with coefficients $\ge 1-\ep$. Then $K_{X'}+f(B)$ is again
$\ep$-lt.
\end{lemma1}
\begin{proof}
Denote $B':=f(B)$ and let $E$ be the exceptional divisor. Then we
can write $K_X+B-\alpha E=f^*(K_{X'}+B')$, where $\alpha\ge 0$.
Since $K_X+B-\alpha E$ is $\ep$-lt, by \ref{crepant-singul} it is
sufficient to show only that $f$ does not contract components of
$B-\alpha E$ with coefficients $\ge 1-\ep$.
\end{proof}

After a number of such birational contractions $\var\colon\widehat
X\to \widetilde X$ we obtain the diagram
\begin{equation}
\label{eq-exc-first}
\begin{array}{c}
 \mbox{
 \begin{picture}(80,46)
 \put(28,33){\vector(-1,-1){20}}
 \put(42,33){\vector(1,-1){20}}
 \put(35,42){\makebox(0,0){$\,\widehat X$}}
 \put(3,4){\makebox(0,0){$X$}}
 \put(67,5){\makebox(0,0){$\widetilde X$}}
 \put(13,28){\makebox(0,0){\scriptsize$\mu$}}
 \put(57,28){\makebox(0,0){\scriptsize$\var$}}
 \end{picture}}
 \end{array}
\end{equation}
where $\rho(\widetilde X)=1$.

By construction, $K_{\widetilde X}+\widetilde{B}$ is klt and
$-(K_{\widetilde X}+\widetilde{B})$ is nef. Moreover, there is a
boundary $\widetilde B^{\triangledown}\le \widetilde B$ such that
$K_{\widetilde X}+\widetilde B^{\triangledown}$ is klt and
$-\left(K_{\widetilde X}+\widetilde B^{\triangledown}\right)$ is
nef and big. By the Inductive Theorem~\ref{Inductive-Theorem},
$K_{\widetilde X}+\widetilde{D}$ is ample. Since $\rho(\widetilde
X)=1$, ${-}K_{\widetilde X}$ is ample. Applying Alexeev's
Theorem~\ref{Al} to $(\widetilde X,0)$, we obtain that families of
such ${\widetilde X}$ are bounded. Further, let $\widetilde{B}_i$
be a component of $\widetilde{B}=\sum b_i \widetilde{B}_i$. Then
$H\cdot \widetilde{B}_i\le {-}2K_{\widetilde X}\cdot H$ for any
very ample $H$. This gives that $\widetilde{B}_i$ lies in a finite
number of algebraic families. Therefore $\Supp \widetilde{B}$ is
bounded and we may assume that $(\widetilde X, \Supp
\widetilde{B})$ is fixed. Since $\widetilde{B}\in\Mm$, there is
only a finite number of possibilities for coefficients $b_i\le
6/7$. Therefore we have the following.

\begin{theorem1}[\cite{Sh1}]
\label{main_Sh_A-1}
Let $(X,B)$ be a projective log surface such that $K_X+B$ is lc
and $-(K_X+B)$ is nef. Assume that the coefficients of $B$ are
standard or $d_i\ge 6/7$ (i.e., $B\in\Mm$). Furthermore, assume
that there is a boundary $B^{\triangledown}\le B$ such that
$K_X+B^{\triangledown}$ klt and $-(K_X+B^{\triangledown})$ is nef
and big. Then either
\begin{enumerate}
\item
$(X,B)$ is nonexceptional, and then there exists a regular nonklt
complement of $K_X+B$, or
\item
$(X,\Supp B)$ belongs to a finite number of algebraic families.
\end{enumerate}
\end{theorem1}

\begin{lemma1}
\label{vspomogat-1}
Let $(S\ni o,\De=\sum\delta_i\De_i)$ be a log surface germ. Assume
that $K_S+\De$ is $(1/7)$-lt, $\delta_i\ge 6/7$, $\forall i$ and
$\Supp \De$ is singular at $o$. Then $(S\ni o)$ is smooth, $\De$
has two (analytic) components at $o$ and $\delta_1+\delta_2<13/7$.
\end{lemma1}
\begin{proof}
By Theorem~\ref{local} there is a regular complement $K_S+\De^+$.
Clearly, $\De^+\ge \up{\De}$. Therefore $K_S+\up{\De}$ is lc. By
Theorem~\ref{kawamata_2} $\up{\De}$ has exactly two (analytic)
components at $o$ and we have an analytic isomorphism
\[
(S,\up{\De})\simeq(\CC^2,\{xy=0\})/\cyc{m}(1,q), \qquad
\gcd(q,m)=1,\quad 1\le q\le m-1.
\]
Taking the corresponding weighted blowup, one can compute
\[
-1+\frac17<a(E,\De)=-1+\frac{1+q}{m}- \frac{\delta_1+q\delta_2}{m}
\]
(see Lemma~\ref{discr-tor}). This yields
\begin{equation}
\label{eq-vspom-1-1}
\frac17<\frac{(1-\delta_1)+q(1-\delta_2)}{m}\le \frac{1+q}{7m}\le
\frac17,
\end{equation}
a contradiction. Hence $m=1$ and $S\ni o$ is smooth. The rest is
obvious.
\end{proof}

\begin{remark1}
\label{vspomogat-1-1}
Notation as in Lemma~\ref{vspomogat-1}. If we replace $(1/7)$-lt
condition with $(1/7)$-lc one, then in \eref{eq-vspom-1-1}
equalities may hold. Then we have two cases
\begin{enumerate}
\item
$(S\ni o)$ is smooth, $\De$ has two (analytic) components at $o$
and $\delta_1+\delta_2\le 13/7$;
\item
$(S\ni o)$ is Du Val of type $A_{m-1}$, $\De$ has two (analytic)
components at $o$ and $\delta_1=\delta_2=6/7$.
\end{enumerate}
\end{remark1}

We claim that the coefficients of $\widetilde{B}$ are bounded from
above by an absolute constant $c<1$ (cf. \cite{Ko1}).

\begin{lemma1}
\label{vspomogat}
Let $(S,\De=\sum\delta_i\De_i)$ be a projective log surface with
$\rho(S)=1$ and $\Theta_1$,\dots, $\Theta_m$ irreducible curves on
$X$. Assume that $\delta_i\in\{0,\frac12, \frac23, \frac34,
\frac56, \frac67\}$. Then there is a constant $c<1$ such that
$\theta_i<c$ for all $\theta_1$, \dots, $\theta_m$ whenever
\begin{enumerate}
\item
$\theta_i\ge 6/7$;
\item
$K_S+\De+\sum\theta_i\Theta_i$ is $(1/7)$-lc;
\item
$-(K_S+\De+\sum\theta_i\Theta_i)$ is nef;
\item
$K_S+\De+\sum\theta_i\Theta_i$ has no regular complements.
\end{enumerate}
\end{lemma1}
\begin{proof}
Assume that there is a sequence $\Theta^{(k)}$ of boundaries as in
our lemma such that $\theta^{(k)}_i\longrightarrow 1$ for some
(fixed) $i$. By Lemma~\ref{vspomogat-1},
Remark~\ref{vspomogat-1-1} and because $\rho(S)=1$ we have
$\theta_i+\theta_j\le 13/7$, for all $j\ne i$. Then
$\theta^{(k)}_j\longrightarrow 6/7$ for $j\ne i$. Denote
$\Theta^0:=\lim \Theta^{(k)}$. By the above, $\Theta^0$ is a
$\QQ$-boundary and $\down{\Theta^0}\ne 0$. Moreover
$K_S+\De+\Theta^0$ is lc and $-(K_S+\De+\Theta^0)$ is nef (see
e.g., \cite[2.17]{Ut}). By the Inductive
Theorem~\ref{Inductive-Theorem} there is a regular $n$-complement
$K_S+\De^++\Theta^{0+}$ of $K_S+\De+\Theta^0$ and this is an
$n$-complement of $K_S+\De+\Theta^{(k)}$ if
$\|\Theta^{(k)}-\Theta^0\|\ll 1$ (see
Proposition~\ref{simple-prp}). A contradiction with our
assumptions.
\end{proof}

The above lemma gives that the coefficients of $\widetilde{B}$ are
bounded from above by an absolute constant $c<1$. So are
coefficients of $\widehat{B}$ because $\var\colon\widehat X\to
\widetilde X$ does not contract components of $\widehat{B}$ with
coefficients $\ge 6/7$. Now we can apply Theorem~\ref{Al} to
$(\widehat X,\widehat{B})$. As above we obtain that $\{(\widehat
X, \Supp \widehat{B})\}$ is bounded. Finally, $(\widehat
X,\widehat B^{\triangledown})$ is a log del Pezzo and $K_{\widehat
X}+\widehat B^{\triangledown}$ is klt. By
Lemma~\ref{Kodaira-lemma} the Mori cone $\NE(\widehat X)$ is
polyhedral. Therefore there is only a finite number of
contractions $\widehat X\to X$ and $\{(X,\Supp B)\}$ is bounded.
Moreover, if $B\in\Msm$, then by Lemma~\ref{vspomogat} there are
only a finite number of possibilities for coefficients of $B$.
This shows the following results. Note that Shokurov \cite{Sh1}
proved them under weaker assumptions.

\begin{theorem1}[\cite{Sh1}]
\label{main_Sh_A-2}
Let $(X,B)$ be a projective log surface such that $K_X+B$ is lc
and $-(K_X+B)$ is nef. Assume that the coefficients of $B$ are
standard (i.e., $B\in\Msm$). Furthermore, assume that there is a
boundary $B^{\triangledown}\le B$ such that
$K_X+B^{\triangledown}$ klt and $-(K_X+B^{\triangledown})$ is nef
and big. Then either
\begin{enumerate}
\item
$(X,B)$ is nonexceptional, and then there exists a regular nonklt
complement of $K_X+B$, or
\item
$(X,B)$ belongs to a finite number of algebraic families.
\end{enumerate}
\end{theorem1}

\begin{theorem1}[\cite{Sh1}]
\label{main_Sh_A}
Notation as in \ref{main_Sh_A-2}. There is an absolute constant
$\mt{Const}$ such that $K_X+B$ is $n$-complementary for some
$n\le\mt{Const}$.
\end{theorem1}

\section[Case of log Enriques surfaces]
{Corollaries: Case of log Enriques surfaces}

Note in many applications of the above theorems we do not need
condition (v) of \ref{notation-boundedness}:

\begin{theorem1}[cf. \cite{Bl}, \cite{Z}, \cite{Z1}]
\label{log_Enriques}
Let $(X,B)$ be a log Enriques surface. Assume that $B\in\Mm$. Then
there is an absolute constant $\mt{Const}$ such that $n(K_X+B)\sim
0$ for some $n\le\mt{Const}$.
\end{theorem1}
\begin{proof}
If $X$ is rational, then we can omit ~(v) of
\ref{notation-boundedness} by
Proposition~\ref{log-Enriques->log-del-Pezzo}. Otherwise there is
a regular complement by Proposition \ref{log_Enriques-nonrat}
below.
\end{proof}

\begin{proposition1}[cf. Th.~\ref{Inductive-Theorem-nonrat}]
\label{log_Enriques-nonrat}
Let $(X,B)$ be a log Enriques surface. Assume that $B\in\Mm$ and
$X$ is nonrational. Then $n(K_X+B)\sim 0$ for some $n\in\RRR_2$.
In particular, $nB$ is integral and $B\in\Msm$.
\end{proposition1}
\begin{proof}
By Theorem~\ref{Inductive-Theorem-nonrat} we may assume that
$K_X+B$ is klt. First consider the case when $B\neq 0$. By
Corollary~\ref{corollary-1-1} the pair $(X,B)$ has at worst
canonical singularities. Replace $(X,B)$ with the minimal
resolution and the crepant pull back. It is easy to see that this
preserves all the assumptions. Run $K_X$-MMP. Clearly, whole $B$
cannot be contracted. At the end we get an extremal contraction
$f\colon X\to Z$ to a smooth curve $Z$ with $p_a(Z)\ge 1$.
Moreover $X$ is smooth, so $X$ is a ruled surface. Thus
$h^1(X,\OOO_{X})=p_a(Z)\ge 1$. On each step we can pull back
regular complement, so it is sufficient to prove our statement for
this new $X$. By Lemma~\ref{lemma1} all the components of $B$ are
horizontal. Further, if $B_i^2>0$ for some component
$B_i\subset\Supp{B}$, then $-(K_X+B-\ep B_i)$ is nef and big, a
contradiction (see Lemma~\ref{lemma_non_rat}). On the other hand,
if $B_i^2<0$, then $(K_X+B+\ep B_i)\cdot B_i<0$. Thus $B_i$ is an
extremal rational curve, a contradiction with $p_a(Z)\ge 1$.
Therefore $B_i^2=0$ for all components $B_i\subset\Supp{B}$. Since
$\rho(X)=2$, this gives that all the $B_i$ are numerically
proportional and $B_i\cap B_j=\emptyset$ for $i\neq j$. As in
\eqref{eq-new-new} we have $(K_X+B_i)\cdot B_i=(K_X+B)\cdot
B_i=0$. Hence $p_a(B_i)=p_a(Z)=1$ and all $B_i$ and $Z$ are smooth
elliptic curves. Moreover, $K_X^2=K_X\cdot B_i=0$. Since
$B_i^2=0$, $B_i$ generate an extremal ray of $\NE(X)\subset\RR^2$
(see Proposition~\ref{MMP-properties}). In particular, $X$
contains no curves with negative self-intersections. Restricting
$B=\sum b_iB_i$ to a general fiber $F$ of the rulling $f\colon
X\to Z$ we get a numerically trivial divisor $K_F+B|_F$ with
$B|_F\in \Mm$. Obviously, $B|_F\in \Msm$ and $B\in\Msm$. Thus we
can write
\[
B=\sum_i(1-1/m_i)B_i,\quad\text{where}\quad m_i\in\NN.
\]
\begin{claim*}
There is a structure of an elliptic fibration $g\colon X\to\PP^1$.
All components of $B$ are contained in fibers of $g$.
\end{claim*}
\begin{proof}
If $B$ has at least two components $B_i$ and $B_j$, then $0\equiv
K_X+B \equiv K_X+B+\ep B_i-\ep' B_j$ is klt and numerically
trivial for some small positive $\ep,\ep'\in\QQ$. By the Log
Abundance Theorem \cite[Ch. 11]{Ut}, we have $K_X+B\qq K_X+B+\ep
B_i-\ep' B_j\qq 0$, i.e., $\ep B_i\qq \ep' B_j$. This gives an
elliptic pencil. Assume that $B=(1-1/m_1)B_1$, where $B_1$ is an
irreducible smooth elliptic curve. Clearly, $3\le B_1\cdot F\le
4$. If the rulling $X\to Z$ corresponds to a vector bundle of
splitting type, then there is a section $B_0$ such that $B_0^2=0$.
Again we can apply the Log Abundance Theorem to $K_X+B+\ep
B_0-\ep' B_1$.

Finally, consider the case when $X$ corresponds to an
indecomposable vector bundle $\EEE$. By Example \ref{Atiyah} the
degree of $\EEE$ is odd. Then up multiplication by an invertible
sheaf, $\EEE$ is a nontrivial extension
\[
0\longrightarrow \OOO_Z \longrightarrow \EEE \longrightarrow
\OOO_Z(P) \longrightarrow 0,
\]
where $P$ is a point on $Z$ (see e.g., \cite[Ch. V, \S 2]{Ha}).
Let $C$ be a section corresponding to the above exact sequence.
Denote also $F_P:=f^{-1}(P)$. It is easy to see that ${-}K_X\sim
2C-F_P$ is nef but not big (see \cite[Ch. V, 2.10, 2.21]{Ha}). It
is sufficient to show the existence of an effective divisor
$B_0\equiv {-}K_X\equiv 2C-F$. Indeed, then $B_0$ is irreducible,
$B_0\neq B_1$ and the same arguments as above gives an elliptic
pencil. To find such $B_0$ we consider the linear system $2C$.
Since $C^2=1$, $C$ is ample. By Riemann-Roch and Kodaira
vanishing, we have
\begin{multline*}
h^i(X,\OOO_X(C))=h^i(X,\OOO_X(2C))=0,\quad i>0,\\
h^0(X,\OOO_X(C))=1,\qquad h^0(X,\OOO_X(2C))=3.
\end{multline*}
Thus the surjection
\[
H^0(X,\OOO_X(2C)) \longrightarrow H^0(C,\OOO_C(2C))
\]
is surjective. Therefore the linear system $|2C|$ is base point
free and determines a finite morphism $\phi\colon X\to \PP^2$ of
degree $4$. Since $2C\cdot F=2$, the images of the fibers $F$ form
a pencil of curves of degree $\le 2$. Since any pencil of conics
contains a degenerate member, we derive that $\phi(F_0)$ is a line
for at least one fiber $F_0=f^{-1}(P_0)$. Then the residual curve
$\phi^{-1}(F_0)-F_0$ belongs to the linear system $|2C-F_0|$. The
claim is proved.
\end{proof}

\begin{exercise1}[{\cite[\S 2]{Sh1}}]
As above, let $X=\PP(\EEE)$, where $\EEE$ is an indecomposable
vector bundle of odd degree. Describe the structure of elliptic
fibration on $X$ and multiple fibers explicitly.
\end{exercise1}

Going back to the proof of Proposition~\ref{log_Enriques-nonrat},
assume that $\sum B_i$ contains also all multiple fibers of the
fibration $g$ (we alow $m_i=1$). By the canonical bundle formula
(see e.g. \cite[Ch. V, \S 12]{BPV})
\[
K_X= -2L+\sum (r_i-1)B_i,
\]
where $L$ is a general fiber and $r_i$ is the multiplicity of
$B_i$. This gives us
\begin{equation}
\label{k-k}
K_X+B= -2L+ \sum \left(1-\frac{1}{r_im_i}\right)r_iB_i.
\end{equation}
Since $r_iB_i\sim L$,
\[
\sum \left(1-\frac{1}{r_im_i}\right)=2.
\]
Hence the collection $(r_1m_1,r_2m_2,\dots)$ is one of the
following (see \ref{=}):
\[
(2,2,2,2),\quad (3,3,3),\quad (2,4,4),\quad (2,3,6).
\]
By \eqref{k-k} we see that $n(K_X+B)\sim 0$ for $n=2,3,4$, and
$6$, respectively.

Finally, let $B=0$. Replace $X$ with its minimal resolution and
let $B$ be the crepant pull back of $K_X$. If again $B=0$, then by
the classification of smooth surfaces of Kodaira dimension
$\kappa=0$ we have $nK_X\sim 0$ for $n\in\RRR_2$. If $B\neq 0$,
then we run $K_X$-MMP (i.e. contract $-1$-curves step by step). As
above, at the end we get a contraction $f\colon X\to Z$ to a
smooth curve $Z$ with $p_a(Z)\ge 1$. On the other hand, all
components of $B$ are rational (because so are singularities of
our original $X$) and at least one component of $B$ is horizontal,
a contradiction. The proposition is proved.
\end{proof}

Note that Proposition~\ref{log_Enriques-nonrat} can be proved by
using log canonical covers (cf. \ref{(dd)} and \cite{Ishii1}).

In the case of log Enriques surfaces with a trivial boundary $K_X$
is $n$-complementary if and only if $nK_X\sim 0$. It is known that
this $n$ can be taken $\le 21$ \cite{Bl}, \cite{Z}. Note that here
the lc condition of $K_X$ cannot be removed (as well as in the
case of log del Pezzos): there are examples of series of normal
surfaces with rational singularities and numerically trivial $K_X$
such that their index tends to infinity \cite{Sakai2}.

\section[On the explicit bound]
{On the explicit bound of exceptional complements}
\label{last}
\begin{problem1}
Notation as in Theorem~\ref{main_Sh_A}. Describe the set of all
$\compl(X,B)$. In particular, find the precise value of $\Const$.
\end{problem1}

The proof of Theorem~\ref{main_Sh_A} shows that we can hope to
classify exceptional del Pezzos (at least modulo some birational
modifications). We describe below an explicit method to
reconstruct our original pair $(X,B)$ from a model $(\widetilde
X,\widetilde{B})$ (see \eref{eq-exc-first}).
Let $(X,B)$ be an exceptional log del Pezzo. Assume that $B\in\Mm$
and $K_X+B$ has no regular complements. Shokurov defined the
following invariant (which is finite by Lemma~\ref{finite-term}):
\[
\delta=\delta(X,B):=\#\{\ \text{divisors}\ E\ \text{such that}\
a(E,B)\le -6/7\}.
\]
\index{$\delta(X,B)$} An analog of the invariant $\delta$ was
considered in \cite{Z1}.
\subsection{Case $\delta(X,B)=0$}
\label{delta=0}
In this case, $K_X+B$ is $(1/7)$-lt and coefficients of $B$ are
contained in $\{0,\frac12, \frac23, \frac34, \frac45, \frac56\}$.
By Theorem~\ref{Al} there is only a finite number of families of
such $(X,B)$, but very few classification results are known.

We give alternative proof of the boundedness in this case using
more effective Nikulin's theorem~\ref{Nikulin}. Assume that
${-}K_X$ is not ample. We set
\[
\alpha:=\left\{
\begin{array}{ll}
 0&\text{if}\ {-}K_X\ \text{is nef,}\\
 \min \{t \mid -(K_X+tB)\
 \text{is nef}\}&\text{otherwise.}\\
\end{array}
\right.
\]
Then $0\le\alpha\le 1$, ${-}(K_X+\alpha B)$ is nef and not ample.
Moreover, $K_X+\alpha B$ is $(1/7)$-lt. Further, there is an
extremal ray, say $R$, on $X$ such that $(K_X+\alpha B)\cdot R=0$.
Let $\var\colon X\to X'$ be its contraction. If $\var$ is not
birational, then $0\le (1-\alpha)B\cdot R= (K_X+B)\cdot R\le 0$.
This yields $(K_X+B)\cdot R=0$ and this contradicts to that
$-(K_X+B)$ is big. Therefore $\var$ is birational. Put
$B':=\var_*B$. By Lemma~\ref{ep-lt-contractions}, $K_{X'}+B'$ is
$(1/7)$-lt and all our assumptions hold on $X'$. Note also that
for the exceptional divisor $E$ one has $a(E,B')\le a(E,\alpha
B')=0$. Continuing the process, we obtain a pair $(X'',B'')$ such
that ${-}K_{X''}$ is ample, $K_{X''}+B''$ is $(1/7)$-lt and
${-}(K_{X''}+B'')$ is nef and big. By Theorem~\ref{Nikulin} and by
Lemma~\ref{mult-lemma} below, $\{X''\}$ is bounded and so is $B''$
(because $B''\in\{0, \frac12, \frac23, \frac34, \frac45,
\frac56\}$). By Lemma~~\ref{finite-term} there are only a finite
number of divisors with discrepancies $a(E,B'')\le 0$. Therefore
there is only a finite number of extractions $X\to X''$ (and all
of them are dominated by the terminal modification $X^t\to X$; see
\ref{terminal-mod}). Finally, $B\in\{0, \frac12, \frac23, \frac34,
\frac45, \frac56\}$ gives the boundedness of $\{(X,B)\}$.

\begin{lemma1}
\label{mult-lemma}
Let $X\ni o$ be a germ of two-dimensional $\ep$-lt singularity.
Then the multiplicity $e(X)$ is bounded by $C(\ep)$.
\end{lemma1}
\begin{proof}
By Corollary~\ref{except-bound} it is sufficient to check it for
two series of nonexceptional klt singularities. Consider, for
example, singularities of type $\BA_n$, i.e., it has the following
dual graph of the minimal resolution:
\[
\begin{array}{ccccccccc}
\wcir{-b_1}&\lin&\wcir{-b_2}&\lin&
\cdots&\lin&\wcir{-b_{n-1}}&\lin&\wcir{-b_n}.
\end{array}
\]
Then the discrepancies $a_1,\dots, a_n$ can be found from the
system of linear equations (see \cite{Iliev}, \cite[Ch. 3]{Ut}):

\begin{align*}
 b_1-2&=-b_1a_1+a_2\\
 b_2-2&=a_1-b_2a_2+a_3\\
 \cdots&\phantom{=}\cdots\\
 b_{n-1}-2&=a_{n-2}-b_{n-1}a_{n-1}+a_n\\
 b_n-2&=a_{n-1}-b_na_n.\\
\end{align*}
This yields
\[
\ep\sum(b_i-2)<\sum (b_i-2)(a_i+1)=-a_1-a_n<2-2\ep.
\]
Finally, by \cite{Brieskorn},
\[
e(X)=2+\sum(b_i-2)<2+\frac{2-2\ep}{\ep}=\frac{2}{\ep}.
\]
We left computations in case $\DD_n$ (with the dual graph as in
Fig.~\ref{figure2}) to the reader.
\end{proof}

\subsection{Case $\delta\ge 1$}
We can argue similar to \ref{delta=0}. Notation as in
\eref{eq-exc-first}. Take
\[\alpha:=\max\{ t\mid
-(K_{\widehat X}+\widehat{B}+t(\widehat{D}-\widehat{B}))\ \text{is
nef}\}.
\]
Clearly, $\alpha<1$. Then there is an extremal rational curve $E$
such that $(K_{\widehat X}+\widehat{B}+\alpha(\widehat{D}-
\widehat{B})) \cdot E=0$ and $(\widehat{D}-\widehat{B})\cdot E>0$.
Hence $E\cdot (K_{\widehat X}+\widehat{D})>0$. By
Lemma~\ref{lemma-birat-E} $E$ is of birational type. We contract
it. This contraction is $(K_{\widehat
X}+\widehat{B}+\alpha(\widehat{D}- \widehat{B}))$-crepant.
Repeating this process, we obtain a sequence of contractions
$\widehat X\to\widetilde X$ and the sequence of rational numbers
$\alpha=\alpha_1\le\alpha_2\le\cdots\alpha_n<1$ such that on each
step $K+B$ is $(1/7)$-lt, $-(K+B)$ is nef and big,
$K+B+\alpha(D-B)$ is klt and $-(K+B+\alpha(D-B))$ is nef. On the
last step we have $\rho(\widetilde X)=1$ and $K_{\widetilde
X}+\widetilde{B}+ \alpha_n(\widetilde{D}-\widetilde{B})$ is klt.
Since $\widetilde X$ is $(1/7)$-lt, $\{\widetilde X\}$ is bounded
by Theorem~\ref{Nikulin}. As above $\Supp \widetilde{B}=\Supp
\widetilde{D}$ is also bounded (because $\widetilde{B}$ has
coefficients $\ge 1/2$).
\par

First assume that $\delta=1$, i.e., $\widetilde{C}$ is
irreducible. Since both $-(K_{\widetilde X}+\widetilde{B})$ and
$K_{\widetilde X}+\widetilde{D}$ are ample, there is $\alpha'$,
$\alpha_n<\alpha'<1$ such that $K_{\widetilde X}+\widetilde{B}+
\alpha'(\widetilde{D}-\widetilde{B})$ is numerically trivial (and
klt). Write
\[
K_{\widetilde X}+\widetilde{B}+
\alpha'(\widetilde{D}-\widetilde{B})=K_{\widetilde
X}+\widetilde{B}+ \beta\widetilde{C}, \qquad \beta>0.
\]
Clearly, we can find this $\beta$ even we do not know the
coefficient of $\widetilde{C}$ in $\widetilde{B}$. Thus we may
assume that $\widetilde{B}+ \beta\widetilde{C}$ is fixed. As in
the case $\delta=0$ there are only a finite number of divisors
with discrepancies $a(E,\widetilde{B}+\beta\widetilde{C})\le 0$
(see Lemma~\ref{finite-term}). Therefore there is only a finite
number of extractions $\widehat X\to \widetilde X$.
\par
In the case $\delta\ge 2$ we have to be more careful. Again we can
take $\alpha'$, $\alpha_n\le\alpha'<1$ such that $K_{\widetilde
X}+\widetilde{B}+ \alpha'(\widetilde{D}-\widetilde{B})$ is
numerically trivial (and klt). If this $K_{\widetilde
X}+\widetilde{B}+ \alpha'(\widetilde{D}-\widetilde{B})$ is
$(1/7)$-lc, then by Lemma~\ref{vspomogat} the coefficients of
$\widetilde{B}+ \alpha'(\widetilde{D}-\widetilde{B})$ are bounded
from above by $c<1$ (and we can find this constant on fixed
$\widetilde X$). Then the log divisor $K_{\widetilde
X}+\widetilde{\De}+c\sum \widetilde{C}_i$ is klt, where
$\widetilde{\De}:=\fr{\widetilde{D}}$ (see \ref{first_prop_sing}).
Clearly,
\[
\widetilde{\De}+c\sum \widetilde{C}_i\ge \widetilde{B}+
\alpha_n(\widetilde{D}-\widetilde{B})
\]
Arguing as above, we find only a finite number of extractions
$\widehat X\to \widetilde X$. More precisely, the discrepancy of
any exceptional divisor $E$ of $\widehat X\to \widetilde X$ with
respect to $K_{\widetilde X}+\widetilde{B}+
\alpha_n(\widetilde{D}-\widetilde{B})$ must be nonpositive (i.e.,
$a(E,\widetilde{B}+ \alpha_n(\widetilde{D}-\widetilde{B}))\le 0$).
In particular, the discrepancy of $E$ with respect to
$K_{\widetilde X}+\widetilde{\De}+c\sum \widetilde{C}_i$ must be
nonpositive. By our assumptions, $\widetilde X$,
$\widetilde{\De}$, $\widetilde{C}_i$ and $c$ are fixed. Then
Lemma~\ref{finite-term} gives that there is only a finite number
of extractions $\widehat X\to \widetilde X$.

\par
If $K_{\widetilde X}+\widetilde{B}+
\alpha'(\widetilde{D}-\widetilde{B})$ is not $(1/7)$-lc, then
$K_{\widetilde X}+\widetilde{B}+
\alpha''(\widetilde{D}-\widetilde{B})$ is not $(1/7)$-lc (and
antiample) for some $\alpha''<\alpha'$. Then this new log divisor
has
\[
\delta(\widetilde X,\widetilde{B}+
\alpha''(\widetilde{D}-\widetilde{B}))>\delta(\widetilde X,
\widetilde{B})\ge 2
\]
and $\widetilde{B}+ \alpha''(\widetilde{D}-\widetilde{B})\in\Mm$.
This is impossible by Corollary~\ref{delta<3} below.

\chapter[On classification of exceptional complements]
{On classification of exceptional complements:
case $\delta\ge 1$}
\label{last-1}
Now we study the case $\delta\ge 1$ in details. 
\section{The inequlity $\delta\le 2$}
In this section we show that $\delta\le 2$.
Replace $(X,B)$
with a model $(\widetilde X,\widetilde{B})$. By construction,
$\delta(X,B)=\delta(\widetilde X,\widetilde{B})$. Thus we assume
that $\rho(X)=1$, $B\in\Mm$, $K_X+B$ is $(1/7)$-lt and $-(K_X+B)$
is nef. Moreover, there exists a boundary $D$ defined by
\eref{def_D} such that $K_X+D$ is ample and lc. Let $C:=\down{D}$.
Then $\delta(X,B)$ is the number of components of $C$. Since
$K_X+D$ is lc, $C$ has only nodal singularities. The following is
a very important ingredient in the classification.

\begin{theorem}[{\cite{Sh1}}]
Notation as in \ref{last-1}. Then $p_a(C)\le 1$.
\end{theorem}
\begin{outline}
Assume that $p_a(C)\ge 2$. Consider the following birational
modifications:
\begin{equation}
\label{eq-nad-posl}
\begin{array}{c}
 \mbox{
 \begin{picture}(80,46)
 \put(25,33){\vector(-1,-1){20}}
 \put(45,33){\vector(1,-1){20}}
 \put(35,42){\makebox(0,0){$\phantom{{}^{\min}} X^{\min}$}}
 \put(0,4){\makebox(0,0){$X$}}
 \put(70,5){\makebox(0,0){$X'$}}
 \put(10,28){\makebox(0,0){\scriptsize$\mu$}}
 \put(60,28){\makebox(0,0){\scriptsize$\var$}}
 \end{picture}}
 \end{array}
\end{equation}
where $\mu\colon X^{\min}\to X$ be a minimal resolution and
$\var\colon X^{\min}\to X'$ is a composition of contractions of
$-1$-curves. Since $K_X+C$ is lc, $C$ has only nodal
singularities. By Lemma~\ref{vspomogat-1}, $X$ is smooth at $\Sing
C$. Therefore $C^{\min}\simeq C$. Thus $p_a(C)=p_a(C^{\min})\ge
2$, $C^{\min}$ is not contracted and $p_a(C')\ge 2$. Take the
crepant pull back
\[
\mu^*(K_X+B)=K_{X^{\min}}+B^{\min}, \quad\text{with}\quad
\mu_*B^{\min}=B
\]
and put
\[
B':=\var_*B^{\min}.
\]
Note that both $-(K_{X^{\min}}+B^{\min})$ and $-(K_{X'}+B')$ are
nef and big. Since $\rho(X)=1$ and $C\simeq C^{\min}$, we have
\begin{enumerate}
\item[(*)]
every two irreducible components of $C^{\min}$ intersect each
other.
\end{enumerate}
If $X'\simeq\PP^2$, then $-(K_{X'}+\frac67C')$ is ample. This
gives $\frac67\deg C'<3$, $\deg C'\le 3$ and $p_a(C')\le 1$. Now
we assume that $X'\simeq\FF_n$. We claim that $n\ge 2$. Indeed,
otherwise $X'\simeq\PP^1\times\PP^1$, $X'\ne X^{\min}$ (because
$\rho(X)=1$) and we have at least one blowup $X^{\min}\to X''\to
X'$. Contracting another $-1$-curve on $X''$ we get $\FF_1$
instead of $\PP^1\times\PP^1$ and after the next blowdown we get
$\PP^2$. Thus $n\ge 2$. Let $\Sigma_0$ be a negative section of
$\FF_n$ and $F$ be a general fiber. Since $\frac67C'\cdot F\le
{-}K_{X'}\cdot F=2$, we have $C'\cdot F\le 2$. So $C'$ must be
generically a $2$-section of $\FF_n\to \PP^1$ (otherwise $C'$ is
generically a section and $p_a(C')=0$).
\par
First we consider the case when $\Sigma_0$ is not a component of
$C'$. Then the coefficient of $\Sigma_0$ in $C'$ $\le
2-\frac{2\cdot 6} 7=\frac27$. Thus
\[
0\le -(K_{X'}+B')\cdot \Sigma_0\le
-\left(K_{X'}+\frac27\Sigma_0\right)\cdot \Sigma_0=2-n+\frac{2n}7.
\]
Hence $n=2$, $X'\simeq\FF_2$. If $X^{\min}\ne X'$, then
$X^{\min}\to X'$ contracts at least one $-1$-curve. But then
contracting another $-1$-curve we obtain either $X'=\FF_3$ or
$X'=\FF_1$, a contradiction with our assumptions. Therefore
$X^{\min}=X'$ and $X$ is a quadratic cone in $\PP^3$. Since
$-(K_X+\frac67C)$ is ample, $C\equiv aH$, where $H$ is the ample
generator of $\Pic(X)$ and $a<\frac73$. By Adjunction we have
\[
\deg K_C\le (K_X+C)\cdot C=2(a-2)a<2.
\]
Hence $p_a(C)\le 1$ in this case.
\par
Finally, we consider the case when $\Sigma_0$ is a component of
$C'$. Write $C'=\Sigma_0+\Sigma'$. Then $\Sigma'$ is generically a
section. From $p_a(C')\ge 2$ by genus formula, we have
$\Sigma_0\cdot\Sigma'\ge 3$. But then
\[
0\ge (K_{X'}+B')\cdot \Sigma_0\ge \left(K_{X'} +\Sigma_0+ \frac67
\Sigma'\right) \cdot \Sigma_0\ge -2+\frac67\cdot 3\ge \frac47,
\]
a contradiction.
\end{outline}

\begin{corollary1}[{\cite{Sh1}}]
\label{delta<3}
Notation as in \ref{last-1}. Then $\delta(X,B)\le 2$.
\end{corollary1}
\begin{proof}
Let $C=\sum_{i=1}^{\delta} C_i$. From the exact sequence
\[
0\longrightarrow\OOO_C\longrightarrow\oplus\OOO_{C_i}
\longrightarrow\FFF\longrightarrow 0,
\]
where $\FFF$ is a sheaf with $\Supp \FFF=\Sing C$, we have
\begin{equation}
\label{eq-nad-posl-1}
1\ge p_a(C)=1-\delta+\#\{C_i\cap C_j\mid i\ne j\}+\sum p_a(C_i).
\end{equation}
On the other hand, by (*) we have $\#\{C_i\cap C_j\mid i\ne j\}\ge
\frac12 \delta(\delta-1)$. This yields
\begin{equation}
\label{eq-last-delta}
0\ge \frac12 \delta(\delta-3)+\sum p_a(C_i).
\end{equation}
In particular, $\delta\le 3$. Assume that $\delta=3$. Then $C$ is
a wheel of smooth rational curves and in \eref{eq-last-delta} the
equality holds. Let $H$ be an ample generator of $\Pic(X)$. We
have ${-}K_X\equiv rH$, $C_i\equiv \gamma_iH$ for some positive
rational $r$, $\gamma_1$, $\gamma_2$, $\gamma_3$. Since every
$C_i$ intersects $C_j$ transversally at a (unique) nonsingular
point, $1=C_i\cdot C_j=\gamma_i\gamma_jH^2$. Hence
\[
\gamma_1\gamma_2=\gamma_1\gamma_3=\gamma_2\gamma_3= \frac1{H^2}.
\]
This implies
\begin{equation}
\label{eq-gamma}
\gamma_1=\gamma_2=\gamma_3=\frac1{\sqrt{H^2}}\le 1.
\end{equation}
Since $-(K_X+B)$ is ample,
\[
r>\frac67\gamma_1+\frac67\gamma_2+ \frac67\gamma_3=
\frac{18}7\gamma_1.
\]
Therefore $K_X+C_1+C_2+\frac47C_3\equiv
-\left(r-\frac{18}7\gamma_1\right)H$ is antiample (and lc). We
claim that $X$ is smooth along $C_1$. Indeed, otherwise
$\Diff{C_1}(0)\ge \frac12P$, where $P\notin C_2, C_3$. On the
other hand, by Adjunction we have
\[
2>\deg\Diff{C_1}\left(C_2+\frac47C_3\right)=1+\frac47+\frac12>2.
\]
The contradiction shows that $X$ is smooth along $C_1$, and
similarly $X$ is smooth along $C_2$ and $C_3$. Thus $C_1$, $C_2$,
$C_3$ are Cartier. In particular, $\gamma_i\in\NN$. By
\eref{eq-gamma}, $\gamma_1=1$ and $H^2=1$. Since
$\Pic(X)\simeq\ZZ\cdot H$, $C_1, C_2, C_3\in |H|$. The linear
subsystem of $|H|$ generated by $C_1, C_2, C_3$ is base point free
and determines a morphism $X\to\PP^2$ of degree one (see also
Lemma~\ref{del-Pezzo-index>1} below). Therefore $X\simeq\PP^2$ and
$C_1, C_2, C_3$ are lines in the general position. Simple
computations show that $B$ has no other components. Finally,
$K_X+C$ is an $1$-complements of $K_X+B$, a contradiction proves
the corollary.
\end{proof}

\section{Case $\delta=2$}

Following Shokurov \cite{Sh1} we describe the case $\delta=2$:

\begin{theorem}
\label{delta=2}
Let $(X,B)$ be a log surface such that $K_X+B$ is $(1/7)$-lt,
$-(K_X+B)$ is nef, $B\in \Mm$, $\delta (X,B)=2$ and $\rho (X)=1$.
Assume that $(X,B)$ is exceptional. Let $H$ be a positive
generator of $\Pic(X)$. Write
\begin{eqnarray*}
B &=&b_1C_1+b_2C_2+F,\quad F=\sum (1-1/m_i)F_i,\quad \\
b_1,b_2&\ge& 6/7,\quad m_i\in \{1,2,3,4,5,6\},
\end{eqnarray*}
where $C_1$ and $C_2$ are irreducible curves. Then $C:=C_1+C_2$
has only normal crossings at smooth points of $X$, $\Supp F$ does
not pass $C_1\cap C_2$ and $b_1+b_2<13/7$. We have one of the
following possibilities:
\begin{enumerate}
\item[$\mathrm{(A_2^1)}$]
$X=\PP^2$, $B=b_1C_1+b_2C_2+\frac 12F_1+\frac23F_2$, where $C_1$,
$C_2$, $F_1$, $F_2$ are lines such that no three of them intersect
at a point and $b_1+b_2\le 11/6$;
\item[$\mathrm{(A_2^{1\prime})}$]
$X=\PP^2$, $B=b_1C_1+b_2C_2+\frac 12F_1+\frac34F_2$, where $C_1$,
$C_2$, $F_1$, $F_2$ are lines such that no three of them intersect
at a point and $b_1+b_2\le 7/4$;
\item[$\mathrm{(A_2^2)}$] $X$ is a quadratic cone in $\PP^3$,
$B=b_1C_1+b_2C_2+\frac23F_1$, where $C_1$ is its generator, $C_2$,
$F_1$ are its smooth hyperplane sections, $b_1+2b_2\le 8/3$;
\item[$\mathrm{(A_2^3)}$]
$X$ is a rational cubic cone in $\PP^4$,
$B=b_1C_1+b_2C_2+\frac12F_1$, where $C_1$ is its generator, $C_2$,
$F_1$ are its smooth hyperplane sections, $b_1+3b_2\le 7/2$ and
$\#C_2\cap F_1\ge 2$;
\item[$\mathrm{(A_2^4)}$]
$X=\PP(1,2,3)$, $B=b_1C_1+b_2C_2+\frac 12F_1$, where
$C_1=\{x_2=0\}$, $C_2=\{x_3=0\}$ (i.e., $3C_1\sim H$, $2C_2\sim
H$), $F_1$ is a smooth rational curve $\equiv \frac12H$, $F_1\ne
C_2$ which is given by the equation $x_3=x_1^3+x_1x_2$,
$2b_1+3b_2\le 9/2$;
\item[$\mathrm{(A_2^5)}$]
$X=\PP(1,3,4)$, $B=\frac67(C_1+C_2)+ \frac12F_1$, where
$C_1=\{x_2=0\}$, $C_2=\{x_3=0\}$ (i.e., $4C_1\sim H$, $3C_2\sim
H$), $F_1$ is a smooth rational curve $\equiv \frac13H$, $F_1\ne
C_2$ which is given by the equation $x_3=x_1^4+x_1x_2$, in this
case $14(K_X+B)\sim 0$;
\item[$\mathrm{(A_2^6)}$]
$X=\PP(1,2,3)$, $B=\frac67(C_1+C_2)$, where $C_1$ is a line
$\{x_1=0\}$, $C_2\in |{-}K_X|$ (i.e., $6C_1\sim H$, $C_2\sim H$),
$\Sing X\subset C_1$, in this case $7(K_X+B)\sim 0$;
\item[$\mathrm{(I_2^1)}$]
$X$ is a quadratic cone in $\PP^3$, $B=b_1C_1+b_2C_2+\frac12F_1$,
where $C_1$, $C_2$ are two smooth hyperplane sections, $F_1$ is a
generator of the cone, $b_1+b_2\le 7/4$;
\item[$\mathrm{(I_2^2)}$]
$X=\PP(1,2,3)$, $B=\frac67C_1+\frac67C_2$, where $C_1=\{x_3=0\}$,
$C_2=\{x_2^2=\alpha_1x_1^4+\alpha_2x_1^2x_2+x_1x_3\}$, $\alpha_1,
\alpha_2\in\CC$, $(\alpha_1, \alpha_2)\ne (0,0)$, $2C_1\sim H$,
$3C_2\sim 2H$, in this case $7(K_X+B)\sim 0$.
\end{enumerate}
\end{theorem}
\begin{remark*}
Note that in all cases $\Weillin(X)\simeq\ZZ$. Therefore we can
verify (i) in the definition of complements \ref{ldef}
numerically, i.e., we need to check only that $nB^+$ is integral
and $K_X+B^+\equiv 0$. By the Inductive
Theorem~\ref{Inductive-Theorem}, (ii) of \ref{ldef} holds
automatically whenever $(X,B)$ is exceptional.
\end{remark*}

Shokurov's proof is based on a detailed analysis of the minimal
resolution, cf. \eref{eq-nad-posl}. Our proof uses computations of
Fano indices of $X$ (as in the proof of Corollary~\ref{delta<3}).
We use slightly \ref{Keel-McCern-}. Note that one can avoid using
of \ref{Keel-McCern-}, but then computations become a
little more complicated. 

The important property is that $K_X+D$ is
analytically dlt except for one case:

\begin{lemma1}[{\cite{Sh1}}]
\label{vspomogat-n}
Let $(S\ni o,B=\sum b_iB_i)$ be a log surface germ, where
$B\in\Mm$. Assume that $K_S+B$ is $(1/7)$-lt. As in \eref{def_D},
put
\[
C:=\down{\frac76B}=\sum_{b_i\ge 6/7} B_i,\qquad F:=\sum_{b_i<6/7}
b_iB_i\quad \text{and}\quad D:=C+F.
\]
Then one of the following holds:
\begin{enumerate}
\item
$K_S+D$ is analytically dlt at $o$;
\item
$o\in S$ is smooth and near $o$ we have $D=C+\frac12L$, where
$(S,C+L)\simeq_{\mt{an}} (\CC^2,\{y(y-x^2)=0\})$.
\end{enumerate}
\end{lemma1}
\begin{proof}
Clearly, we may assume that $K_S+D$ is not plt (otherwise we have
case (i)). By Theorem~\ref{local} there is a regular complement
$K_S+B^+$. Since $B\in\Mm$, $B^+\ge D$. In particular, $K_S+D$ is
lc and $C=\down{D}$ has at most two (analytic) components passing
through $o$ (see Theorem~\ref{kawamata_2}). If $C$ has exactly two
components, then $S\ni o$ is smooth by Lemma~\ref{vspomogat-1}.
Obviously, $K_S+D$ is analytically dlt at $o$ in this case. From
now on we assume that $C$ is analytically irreducible at $o$.
Write $B=bC+F$, where $b\ge 6/7$. Recall that $F\in\{0,\frac12,
\frac23,\frac34, \frac45, \frac56\}$.
\par
First we consider the case when $K_S+C$ is not plt. Then $D=C$ and
$(S,C,o)$ is such as in (ii) of \ref{kawamata_2}. In particular,
$2(K_S+C)\sim 0$ and $K_S+C\not \sim 0$. Let $f\colon (\widetilde
S,E)\to S$ be an inductive blowup of $(S,D)$ and $\widetilde C$
the proper transform of $C$. Write
\[
\begin{array}{l}
f^*(K_S+C)=K_{\widetilde S}+\widetilde C +E,\\
f^*(K_S+bC)=K_{\widetilde S}+b\widetilde C +\alpha E,\\
\end{array}
\]
where $\alpha <6/7$. Here $2\left(K_{\widetilde S}+\widetilde C
+E\right)\sim 0$. By Adjunction, $K_E+\Diff{E}(\widetilde C)$ is
not klt and $\deg\Diff{E}(\widetilde C)=2$. Moreover,
$K_E+\Diff{E}(\widetilde C)$ is not $1$-complementary (because
neither is $K_S+C$). Therefore we have (cf. Lemma~\ref{mMmM})
\[
\Diff{E}(\widetilde C)=\frac1{2}P_1+\frac1{2}P_2+P_3,\qquad
\Diff{E}(0)=\frac1{2}P_1+\frac1{2}P_2+\frac{m-1}{m}P_3
\]
for some points $P_1,P_2,P_3\in E$ and some $m\in \NN$. From this
we have
\[
\left(K_{\widetilde S}+E\right)\cdot E+b\widetilde C\cdot E
+(-1+\alpha)E^{2}=0.
\]
By Adjunction
\[
\left(K_{\widetilde S}+E\right)\cdot
E=-2+\frac12+\frac12+1-\frac1m=-\frac1m.
\]
Since $\widetilde C\cap E$ is a point of type $\frac1{m}(1,q)$,
$\widetilde C\cdot E\ge 1/m$. This yields
\[
\frac17\left(-E^{2}\right)<(-1+\alpha)E^{2}\le \frac1{7m}.
\]
Thus $0<-E^{2}<1/m$ and $-1/m<K_{\widetilde S}\cdot E<0$. On the
other hand, $mK_{\widetilde S}$ is Cartier near $E$. Therefore
$mK_{\widetilde S}\cdot E\in\ZZ$, a contradiction.
\par
Now we may assume that $K_S+C$ is plt. By Theorem~\ref{local},
$K_S+D$ is $2$-complementary and $D^+\ge D$, so $2(K_S+D)\sim 0$
and $2F$ is integral. We claim that $(S\ni o)$ is smooth. Assume
the opposite. Then
\[
(S,C)\simeq (\CC^{2},\{y=0\})/\ZZ_m(1,q),\quad \gcd(q,m)=1,\ m\ge
2,\ 1\le q\le m-1.
\]
Consider the weighted blowup with weights $\frac1{m}(1,q)$. By
Lemma~\ref{discr-tor} we get the exceptional divisor $E$ with
\[
a(E,D)=-1+\frac{1+q}{m}-\frac{q}{m}-\frac{\mu}{2}=
-1+\frac1m-\frac{\mu}2,
\]
where $\mu=\mult_E(2F)\in\frac1m\NN$. Since $2(K_S+D)\sim 0$, we
have $a(E,D)=-1$ or $-1/2$. But in the second case $\mu=2/m-1\le
0$, a contradiction. Therefore $a(E,C+F)=-1$ and $\mu=2/m$.
Further,
\[
-1+\frac17<a(E,B)=-1+\frac{1+q}{m}-b\frac{q}{m}-\frac{\mu}2=
-1+\frac{q(1-b)}m<-1+\frac17.
\]
The contradiction shows that $(S\ni o)$ is smooth. Now we claim
that $\up{F}$ is a smooth curve. As above, consider the blowup of
$o\in S$. For the exceptional divisor $E$, we have
\[
-1+\frac17<a(E,B)=1-b-\frac{\mu}2,
\]
where $\mu=\mult_E(2F)\in\NN$. Hence $\mu=1$ and $L=\up{F}$ is
smooth. Finally, $K_S+C+(\frac12-\ep)L$ is plt for any $\ep>0$. By
Adjunction, $\down{\Diff C\bigl((\frac12-\ep)L\bigr)}\le 0$. Hence
$\down{\Diff C(\frac12L)}$ is reduced. This means that $C\cdot
L=2$, i.e., $C$ and $L$ have a simple tangency at $o$. The rest is
obvious.
\end{proof}

We need some (well known) facts about Fano indices of log del
Pezzo surfaces.
\begin{definition}
\label{def-Fano-ind}
Let $(X,D)$ be a log del Pezzo surface. Define the \textit{Fano
index}\index{Fano index} $r(X,D)$\index{$r(X,D)$} of $(X,D)$ by
\[
r(X,D)=\sup\{t \mid -(K_X+D)\equiv tH,\quad\text{for some}\quad
H\in\Pic(X)\}.
\]
\end{definition}
If $K_X+D$ is klt or $K_X+D$ is dlt and $-(K_X+D)$ is ample, then
by Lemma~\ref{log-del-Pezzo-Pic}, $r(X,D)\in\QQ$ and
$-(K_X+D)\equiv r(X,D)H$ for some (primitive and ample) element
$H\in\Pic(X)$ (recall that we consider only $\QQ$-divisors). In
the case $D=0$ we write $r(X)$ instead of $r(X,0)$.

The following is an easy consequece of Riemann-Roch,
Kawamata-Viehweg vanishing and \cite{Fujita}.

\begin{lemma1}
\label{del-Pezzo-index>1}
Let $X$ be a log del Pezzo with klt singularities of Fano index
$r=r(X)$. Assume that ${-}K_X$ is ample and write ${-}K_X\equiv
rH$, where $H$ is a primitive (ample) element of $\Pic(X)$. Then
\begin{enumerate}
\item
$\dim|H|=\frac12(1+r)H^2$, hence $r=\frac{2l}{H^2}-1$, where
$l:=\dim|H|$;
\item
$H^2\ge \dim|H|-1$, hence $r\le 1+\frac2{H^2}$;
\item
if $r>1$, then
\[
\dim|H|=H^2+1,\quad\text{and}\quad r=1+\frac2{H^2}.
\]
Moreover, $X$ is one of the following $X\simeq \PP^2$ ($r=3$),
$X\simeq \PP^1\times \PP^1$ ($r=2$), $X\subset \PP^{d+1}$ is a
cone over a rational normal curve of degree $d=H^2$ ($r=1+2/d$).
\end{enumerate}
\end{lemma1}
\begin{proof}
By Kawamata-Viehweg vanishing \cite[1-2-5]{KMM} one has
$H^i(X,\OOO_X(H))=H^i(X,\OOO_X)=0$ for $i>0$. Therefore by
Riemann-Roch we obtain
\[
\dim|H|=\frac{H\cdot(H{-}K_X)}{2}=\frac{(1+r)H^2}{2}.
\]
This proves (i). Recall (see \cite{Fujita}) that for any polarized
variety $(X,H)$ the following equality holds:
\begin{equation}\label{Fujita-Delta}
\dim X+H^{\dim X}-h^0(X,\OOO_X(H))\ge 0.
\end{equation}
Combining this with (i) we obtain (ii). Finally, assume $r>1$.
Then by (i), $\dim|H|>H^2$. From (ii) we have $H^2= \dim|H|-1$.
Moreover, in \eqref{Fujita-Delta} the equality holds. Such
polarized varieties (of arbitrary dimension) are classified in
\cite{Fujita}. In particular, it is proved that $H$ is very ample
and $X\subset \PP^{\dim|H|}$ are varieties of \textit{minimal
degree}. In the two-dimensional case from \cite{Fujita} we obtain
possibilities as in (iii).
\end{proof}

Log del Pezzo surfaces with $r(X)=1$ are special cases of the
so-called Fujita varieties:
\begin{lemma1}
\label{del-Pezzo-index=1}
Let $X$ be a log del Pezzo with klt singularities of Fano index
$1$. Assume that ${-}K_X$ is ample and $H$ an ample primitive
element of $\Pic(X)$ such that ${-}K_X\equiv H$. Then
\begin{enumerate}
\item
$\dim|H|=H^2$ and $H^2\le 8$;
\item
if $H^2\ge 4$, then $X$ has only DuVal singularities;
\item
if $H^2=6$ and $\rho(X)=1$, then $X$ has exactly two singular
points which are Du Val of types $A_1$ and $A_2$; in this case,
$X$ is isomorphic to the weighted projective plane $\PP(1,2,3)$.
\end{enumerate}
\end{lemma1}
\begin{outline}
Note that by Lemma~\ref{lemma_non_rat}, $X$ is rational. As in
Lemma~\ref{del-Pezzo-index>1}, the first part of (i) follows by
Riemann-Roch and Kawamata-Viehweg vanishing. Set $D:=H+K_X$. If
$D\sim 0$, then $X$ has only DuVal singularities. In this case, by
Noether's formula,
\[
K_{\widetilde X}^2+\rho(\widetilde X)=K_X^2+\rho(\widetilde X)=10,
\]
where $\widetilde X\to X$ is the minimal resolution. This yields
$K_X^2=H^2\le 8$ (because $X\not\simeq\PP^2$).

If $D\not\sim 0$, then by Lemma~\ref{log-del-Pezzo-Pic}, $nD\sim
0$ for some $n\in\NN$. Considering a cyclic cover trick, we get a
cyclic \'etale in codimension one cover $\var\colon X'\to X$.
Moreover, on $X'$ one has ${-}K_{X'}\sim H'$, where $H':=\var^*H$.
Therefore $X'$ is a del Pezzo surface with only DuVal
singularities. Further, by the above arguments,
\[
K_{X'}^2=(\deg\var)K_X^2\le 9.
\]
Hence $K_X^2\le 4$. If $K_X^2=4$, then $K_{X'}^2=8$ and $X$ is a
quotient of $X'$ by an involution $\tau$. In this case, $X'$
cannot be smooth (otherwise $X$ has only singularities of type
$A_1$ and ${-}K_X\sim H$). Let $\widetilde X'\to X'$ be the
minimal resolution. As above, by Noether's formula,
$\rho(\widetilde X')=10{-}K_{\widetilde X'}^2=10{-}K_{X'}^2=2$.
Therefore, $\widetilde X'\to X'$ contracts a single $-2$-curve.
From this, we have only one possibility: $\widetilde
X'\simeq\FF_2$ and $X'$ is a quadratic cone in $\PP^3$. Since
$\Pic(X')=\ZZ\cdot \OOO_{X'}(1)$, one has that $\tau$ acts
linearly in $\PP^3$. Recall that the quotient of the vertex of the
cone is nonGorenstein. The action of $\tau$ on $\PP^3$ is free in
codimension one (because so is the action of $\tau$ on $X'$).
Therefore in some coordinate system,
\[
\tau=
\begin{pmatrix}
 -1&0&0&0\\
 0&-1&0&0\\
 0&0&1&0\\
 0&0&0&1\\
\end{pmatrix}
\]
and $X'$ is given by
\[
q(x_1,x_2)+q'(x_3,x_4)=0,
\]
where $q(x_1,x_2)$ and $q'(x_3,x_4)$ are quadratic forms such that
$\mt{rk}(q+q')=3$. Changing coordinates we may assume that $X'$ is
given by $x_1^2+x_2^2+x_3^2=0$. But then the quotient of the
vertex is a complete intersection singularity $y_1+y_2+x_3^2=0$,
$y_1y_2=y_0^2$, where $y_1=x_1^2$, $y_2=x_2^2$ and $y_0=x_1x_2$.
In particular, it is Gorenstein, a contradiction.
\par
Assume now that $H^2=6$. Then by the above, $X$ is Gorenstein and
$\rho(\widetilde X)=4$, where $\widetilde X\to X$ is the minimal
resolution. Therefore $\widetilde X\to X$ contracts exactly three
$-2$-curves and the configuration of singular points on $X$ is
either $A_3$ or $A_1A_2$. By \cite{Furushima} the only second case
is possible. Moreover, $X$ is unique up to isomorphism (see e.g.,
\cite[3.10]{KeM}). On the other hand, $\PP(1,2,3)$ is a Gorenstein
del Pezzo of degree $6$.
\end{outline}

\begin{remark*}
There is another way to treat the case $H^2=6$: since $\dim|H|=6$,
one can construct a $1$-complement $K_X+C$ such that $C$ has three
components and then use Theorem~\ref{th-toric-1}.
\end{remark*}

\begin{proof}[Proof of Theorem~\ref{delta=2}]
Since $B\ne 0$ and $\rho (X)=1$, ${-}K_X$ is ample. Hence $X$ is
rational. By Lemma~\ref{vspomogat-n} Then $C:=C_1+C_2$ has only
normal crossings at smooth points of $X$, $\Supp F$ does not pass
$C_1\cap C_2$ and $b_1+b_2<13/7$ (by Lemma~\ref{vspomogat-1}).
\par
Write
\begin{equation*}
C_i\equiv d_iH,\quad {-}K_X\equiv rH,\quad F\equiv qH.
\end{equation*}
We assume that $d_1\le d_2$. Since $-(K_X+B)$ is nef,
\begin{equation}
\label{eq-d-d-b-b}
\frac67(d_1+d_2)\le b_1d_1+b_2d_2+q\le r.
\end{equation}
Take $b$ so that $K_X+C_1+bC_2+F\equiv 0$, i.e.
\[
d_1+bd_2+q=r.
\]
Then
\begin{multline}
\label{eq-b}
b=\frac{r-q-d_1}{d_2}\ge \frac{b_1d_1+b_2d_2-d_1}{d_2}=\\
b_2-(1-b_1)\frac{d_1}{d_2}\ge b_1+b_2-1\ge 5/7.
\end{multline}
Since $K_X+C+F$ is ample, $b<1$.
\par
Recall that $K_X+C+F$ is analytically dlt except for the case (ii)
of Lemma~\ref{vspomogat-n}. In particular, $X$ is smooth at points
$C_1\cap C_2$ and $C_1\cap C_2\cap \Supp F=\emptyset$. By
Adjunction,
\begin{equation}
\label{eq-diff-C-1}
K_{C_1}+\Diff{C_1}(bC_2+F)\equiv 0.
\end{equation}
If $p_a(C_1)>0$, then $K_{C_1}=\Diff{C_1}(bC_2+F)=0$. This is
impossible because $C_1\cap C_2\ne \emptyset$. Therefore
$C_1\simeq \PP^1$ and $\deg \Diff{C_1}(bC_2+F)=2$.

\subsection{Case: $X$ is smooth}
Then $X\simeq \PP^2$ and $r=3$. From \eref{eq-d-d-b-b} we obtain
$(d_1,d_2)=(1,2)$ or $(1,1)$. On the other hand, $K_X+C+F$ is
ample. This gives
\[
q>3-d_1-d_2.
\]
If $(d_1,d_2)=(1,2)$, then by \eref{eq-d-d-b-b}, $0<q\le
3-\frac{18}7=\frac37<\frac12$, a contradiction. Therefore $C_1$,
$C_2$ are lines on $X\simeq \PP^2$. Then
\begin{multline}
\label{eq-P2-second}
\frac12\sum \deg F_i\le q=\\
 \sum (1-1/m_i)\deg F_i\le 3-12/7=9/7,
\qquad q>1.
\end{multline}
If $\deg F_1\ge 2$, then $F=\frac12F_1$, $\deg F_1=2$ and $q=1$, a
contradiction. Hence all the components of $F$ are lines. From
\eref{eq-P2-second} we have only two possibilities:
$F=\frac12F_1+\frac23F_2$ and $F=\frac12F_1+\frac34F_2$. These are
cases $\mathrm{(A_2^1)}$ and $\mathrm{(A_2^{1\prime})}$.

From now on we assume that $X$ is singular. Since $p_a(C)\le 1$,
we have two possibilities: $\# C_1\cap C_2 =2$ and $\# C_1\cap C_2
=1$.

\subsection{Case: $\# C_1\cap C_2 =2$}
\label{case-2-points}
Let $C_1\cap C_2 =\{P_1,P_2\}$. Then
\begin{equation*}
2=C_1\cdot C_2=d_1d_2H^2.
\end{equation*}
Equality \eref{eq-diff-C-1} gives
\begin{equation*}
\Diff{C_1}(bC_2+F)=bP_1+bP_2+\Diff{C_1}(F).
\end{equation*}
Hence
\begin{equation*}
\deg \Diff{C_1}(F)=2-2b\le 4/7.
\end{equation*}

By Inversion of Adjunction, $K_X+C_1+F$ is plt near $C_1$. Assume
that $\Diff{C_1}(F)=0$. Then $F=0$ and $b=1$, a contradiction with
$b<1$. Therefore $\Diff{C_1}(F)\ne 0$.

Since $\Diff{C_1}(F)\in \Msm$ (see Corollary~\ref{coeff-Diff}), we
have only one possibility: $\Diff{C_1}(F)=\frac12Q$, where $Q\in
C_1$ is a single point $\ne P_1, P_2$. Moreover, $b=3/4$ and
$d_1+\frac34d_2+q=r$.

If $Q\in X$ is smooth, then $F=\frac12F_1$, where $F_1$ is
irreducible, $F_1\cap C_1=\{Q\}$ and $F_1\cdot C_1=1$. Thus $C_1$
is Cartier (see \ref{diff}), $d_1\in \NN$ and
$r=d_1+\frac34d_2+q>\frac74$. By Lemma~\ref{del-Pezzo-index>1} $X$
is a cone over a rational normal curve of degree $d\ge 2$. In this
case $r=(d+2)/d>7/4$ and $d=2$. Therefore $X\subset \PP^3$ is a
quadratic cone. Further, $d_1=d_2=1$, so $C_1$, $C_2$ are
hyperplane sections (and they do not pass through the vertex of
the cone). Finally, from $F_1\cdot C_1=1$ we see that $F_1$ is a
generator of the cone. This is case $\mathrm{(I_2^1)}$.

Therefore $Q\in X$ is singular. Then it must be DuVal of type
$A_1$. Moreover, $F=0$ and $2C_1$ is Cartier (but $C_1$ is not,
because $C_1$ is smooth at $Q$). Hence $d_1\in \frac12\NN$.
Further, $d_1+\frac34d_2=r$.

If $d_1\ge 1$, then $d_2\ge 1$ and $r\ge 7/4$. By
Lemma~\ref{del-Pezzo-index>1} and our assumption that $X$ is
singular, $r=2$ and $X$ is a quadratic cone. But then $d_2=4/3$, a
contradiction. Hence $d_1=1/2$, $d_2\ge 1/2$. Put $k:=C_1\cdot
H\in \NN$. Then $H^2=2k$, $2=C_1\cdot C_2=\frac12d_2H^2$, so
$d_2=2/k\ge 1/2$, $k\le 4$. This gives
$r=\frac12+\frac34d_2=\frac12+\frac3{2k}$. On the other hand, by
Lemma~\ref{del-Pezzo-index>1}, $r=\frac{l}k-1$, where $l\in\NN$.
Therefore $3k+3=2l$ and $k\in\{1, 3\}$. If $k=1$, then $l=3$,
$r=2$, $d_2=2$. But this contradicts $\frac67(d_1+d_2)\le r$. We
obtain $k=3$, $l=6$, $r=1$, $d_2=2/3$, $H^2=6$. By
Lemma~\ref{del-Pezzo-index=1}, $X\simeq \PP(1,2,3)$. We may assume
that $C_1\in |\OOO_{\PP}(3)|$ and $C_2\in |\OOO_{\PP}(4)|$. Then
$C_1=\{x_3=0\}$ and
$C_2=\{x_2^2=\alpha_1x_1^4+\alpha_2x_1^2x_2+\alpha_3x_1x_3\}$,
$\alpha_1, \alpha_2, \alpha_3\in\CC$. But $\alpha_3\ne 0$
(otherwise $C_2$ is singular at $(0,0,1)$). Moreover, $(\alpha_1,
\alpha_2)\ne (0,0)$, because $C_1\cap C_2$ consists of two points.
This is case $\mathrm{(I_2^2)}$.

\subsection{Case: $p_a(C_2)=1$}
By \ref{case-2-points} we may assume that $C_1\cap C_2$ is a
single point, say $P$. As in \eref{eq-b} take $b'$ so that
$K_X+b'C_1+C_2+F\equiv 0$, i.e.
\begin{equation*}
b'd_1+d_2+q=r.
\end{equation*}
Since $K_{C_2}+\Diff{C_2}(b'C_1+F)\equiv 0$, we have $\deg
(K_{C_2}+\Diff{C_2}(b'C_1))\le 0$ and $K_{C_2}=0$, $b'\le 0$. This
yields
\begin{equation*}
b'=\frac{r-q-d_2}{d_1}=b_1-(1-b_2)\frac{d_2}{d_1}\le 0,
\end{equation*}
\begin{equation}
\label{d/d}
\frac67d_1\le b_1d_1\le (1-b_2)d_2\le \frac17d_2,\quad 6d_1\le
d_2.
\end{equation}
Assume that $r\le 1$. Then
\begin{equation}
\label{eq-pa=1-r}
1\ge r\ge b_1d_1+b_2d_2+q\ge (b_1+6b_2)d_1+q\ge 6d_1.
\end{equation}
On the other hand, by \eref{eq-diff-C-1},
\begin{equation*}
\deg \Diff{C_1}(F)=2-b,\quad \quad
\end{equation*}
where
\begin{equation}
\label{eq-pa=1-two-comp}
1>b\ge b_2-(1-b_1)\frac{d_1}{d_2}\ge b_2+\frac16b_1-\frac16\ge
\frac56.
\end{equation}
(see \eqref{eq-b} and \eqref{d/d}). Hence
\begin{equation*}
1<\deg \Diff{C_1}(F)\le 7/6.
\end{equation*}
Since $\Diff{C_1}(F)\in \Msm$, we have only one possibility
$\Diff{C_1}(F)=\frac12Q_1+\frac23Q_2$ and $b=5/6$. In particular,
$6C_1$ is Cartier (see Theorem~\ref{diff}), so $d_1\ge 1/6$. On
the other hand, $d_1\le 1/6$ (see \eref{eq-pa=1-r}). Hence
$d_1=1/6$ and $kC_1$ is not Cartier for $1\le k\le 5$. This gives
us that $F=0$. Moreover, in \eref{eq-pa=1-two-comp} equalities
hold, so $1=6d_1=d_2$ and $b_1=b_2=6/7$. From \eref{eq-pa=1-r} we
have $r\ge 6d_1=1$. Hence $r=1$. Further, $C_1\cdot
C_2=\frac16H^2=1$, gives $H^2=K_X^2=6$. By
Lemma~\ref{del-Pezzo-index=1}, $X\simeq \PP(1,2,3)$. We get case
$\mathrm{(A_2^6)}$.

Now assume that $r>1$. Then $X$ is a cone. From $2\ge r\ge
b_1d_1+b_2d_2+q\ge (b_1+6b_2)d_1+q\ge 6d_1$ we see that $d_1\le
1/3$ and $ C_1$ is not Cartier. Hence $C_1$ contains the vertex
and $C_2$ does not. Thus $C_2$ is Cartier. Finally, $C_1\cdot
C_2=1$. Therefore $C_1$ is a generator of the cone and $C_2$ is a
smooth hyperplane section. But then $ C_2 $ is rational, a
contradiction.

\subsection{Case $C_1\cap C_2=\{P\}$ and $p_a(C_1)=p_a(C_1)=0$}
Then $C_1\cdot C_2=1$. By \eref{eq-b}, $1>b\ge 5/7$. Hence $1<\deg
(\Diff{C_1}(F))=2-b\le 9/7$. Using $\Diff{C_1}(F)\in \Msm$ we get
the following cases:
\begin{equation}
\label{eq-cases-diff}
\Diff{C_1}(F)=\frac12Q_1+\frac23Q_2,\quad \frac12Q_1+\frac34Q_2.
\end{equation}
By Inversion of Adjunction, $K_X+C_1+F$ is plt near $C_1$. In
particular, either $4C_1$ or $6C_1$ is Cartier (see \ref{diff})
and $F$ has at most two components. Thus $4d_1$ or $6d_1\in \NN$.
Note that
\[
d_1=\frac1{H\cdot C_2}\le 1,\qquad d_2=\frac1{H\cdot C_1}\le 1.
\]

\subsubsection{Subcase $d_2=1$.}
It is easy to see $H\cdot C_1=d_1H^2=1$, so $d_1=1/H^2$. We claim
that $r>1$. Indeed, if $r\le 1$, then
\begin{equation}
1\ge r\ge \frac67(1+d_1)\label{eq-ineq-r-1-2-3}
\end{equation}
and $d_1\le 1/6$. Thus $mC_1$ is not Cartier for $m<6$. By
\eref{eq-cases-diff} we have that $6C_1$ is Cartier,
$\Diff{C_1}(F)=\frac12Q_1+\frac23Q_2$ and $d_1\ge 1/6$. Therefore
$d_1=1/6$ and in \ref{eq-ineq-r-1-2-3} the equality holds. In
particular, $r=1$, $ K_X^2=H^2=6C_1\cdot C_2=6$. By
Lemma~\ref{del-Pezzo-index=1}, $X\simeq \PP(1,2,3)$ and
$\Weillin(X)\simeq \ZZ$. But then $C_2\sim {-}K_X\sim H$ is
Cartier and $p_a(C_2)=1$, a contradiction.

Thus $r>1$ and $X \subset \PP^{d+1}$ is a cone of degree $d:=H^2$
(see \ref{del-Pezzo-index>1}). Hence $C_2$ is a smooth hyperplane
section and $C_1$ is a generator of the cone (i.e., $d_2=1$,
$d_1=1/d$). Write $F_i\equiv \frac{q_i}{d}H$. (Note that $q_i\in
\NN$ and $F_i\sim q_iC_1$ because $\Weillin(X)\simeq \ZZ\cdot C_1$
in our case). We have
\begin{multline}
\label{eq-two-bound-1-2}
1+\frac1d+\sum \left(1-\frac1{m_i}\right)\frac{q_i}d > r=\\
\frac{d+2}d\ge b_2+\frac1db_1+\sum
\left(1-\frac1{m_i}\right)\frac{q_i}d,\quad \\ q_i \in
 \NN,\quad m_i\in \{0,2,3,4,5,6\}\\
\end{multline}
Assume that $F$ has a component $F_1$ which does not pass through
the vertex. Then $q_1\ge d$, so
\[
1+\frac2d\ge b_2+\frac1db_1+1-\frac1{m_1}\ge
\frac67\left(1+\frac1d\right)+1-\frac1{m_1},
\]
\[
8\ge d\left(6-\frac7{m_1}\right)\ge \frac52d.
\]
This gives $d=2$ or $d=3$. If $d=3$, then $m_1=2$. From
\eref{eq-two-bound-1-2} we get $F=\frac12F_1$, i.e., case
$\mathrm{(A_2^3)}$. If $d=2$, then $m_1=2$ or $m_1=3$. In both
cases by \eref{eq-two-bound-1-2} we have
$F=\left(1-\frac1{m_1}\right)F_1$. For $m_1=2$ we derive a
contradiction with the left side of \eref{eq-two-bound-1-2}. We
obtain case $\mathrm{(A_2^2)}$.

Now we assume that all components of $F$ pass through the vertex
$\upsilon$ of the cone (in particular, $F\ne 0$). Since $K_X+C+F$
is plt at $\upsilon$ (see Lemma~\ref{vspomogat-n}), there is at
most one such a component and $F=(1-\frac1{m_1})F_1$. We claim
that either $q_1=1$ or $q_1\ge d+1$. Indeed, assume that $1<q_1\le
d$. Then
\[
F_1\cdot C_1=\frac{q_1}{d^2}H^2=\frac{q_1}d\le 1.
\]
Since $X$ is smooth outside of $\upsilon$, $F_1\cap
C_1=\{\upsilon\}$. By Adjunction, $\down{\Diff{C_1}(F)}=0$ at
$\upsilon$. On the other hand, by \ref{coeff-Diff}, the
coefficient of $\Diff{C_1}(F)$ at $\upsilon$ is
\[
1-\frac1d+\left(1-\frac1{m_1}\right)(F_1\cdot C_1)=1-\frac1d+
\left(1-\frac1{m_1}\right)\frac{q_1}d.
\]
We obtain
\[
\frac1d-\left(1-\frac1{m_1}\right)\frac{q_1}d>0,\quad
1>\left(1-\frac1{m_1}\right)q_1\quad\text{and}\quad
q_1<\frac{m_1}{m_1-1}\le 2,
\]
a contradiction. Therefore $q_1=1$ or $q_1\ge d+1$. But the second
case is impossible by the right side of \eref{eq-two-bound-1-2}.
Hence $q_1=1$. But this contradicts to the left side of
\eref{eq-two-bound-1-2}.
\par
From now on we assume that $d_1\le d_2<1$.

\begin{remark1}
\label{remark-r>1}
If $r>1$, then $X$ is a cone and contains exactly one singular
point, say $P$, and $P\notin C_1\cap C_2$. Hence we may assume
that $P\notin C_1$ and $C_1$ is Cartier. Thus we may assume that
$r\le 1$ and $C_1$, $C_2$ are not Cartier.
\end{remark1}

\subsubsection{Subcase $d_1=1/2$}
Then we have
\begin{equation*}
1=C_1\cdot C_2=d_1H\cdot C_2,\quad H\cdot C_2=2,\quad d_2H^2=2.
\end{equation*}
Since $1>d_2=\frac2{H^2}\ge d_1=\frac12$, $H^2=3$ or $H^2=4$. On
the other hand, $H\cdot C_1=\frac12H^2\in \NN$. Hence $H^2=4$,
$d_2=1/2$ and $\NN\ni {-}K_X\cdot H=rH^2=4r$. By symmetry, taking
into account $d_1=d_2=1/2$, one can see that \eref{eq-cases-diff}
holds also for $C_2$:
\begin{equation*}
\Diff{C_2}(F)=\frac12Q_1'+\frac23Q_2',\quad\text{or}\quad
\frac12Q_1'+\frac 34Q_2'.
\end{equation*}
From $r\ge \frac67(d_1+d_2)=\frac67$ we get $r\ge 1$. Thus $r=1$
and $X$ is Gorenstein by \ref{remark-r>1} and
Lemma~\ref{del-Pezzo-index=1}. By Theorem~\ref{Keel-McCern-} all
singular points are contained in $C$. Since $K_X+C$ is dlt (see
Lemma~\ref{vspomogat-n}), we obtain that $X$ has only DuVal points
of types $A_{n_i}$, $i=1,\dots,s$. Since $\rho(X)=1$,
$\sum_{i=1}^s n_i=10-4-\rho(X)=5$. By \eref{eq-cases-diff},
$n_i\le 3$ and $(n_1,\dots,n_s)\ne (1,1,1,1,1)$. Now we can use
the classification of Gorenstein del Pezzo surfaces with $\rho=1$
(see e.g., \cite{Furushima}). The configuration of singular points
on $X$ is $\{2A_1A_3\}$. We may assume that $C_1$ contains the
point of type $A_3$. Hence $\Diff{C_1}(F)=\frac12Q_1+\frac 34Q_2$
(see \eref{eq-cases-diff}). At least one of points $Q_1$, $Q_2$,
$Q_1'$, $Q_2'$ is smooth. Hence $F\ne 0$ and $\Supp F\cap
C_1=Q_1$. Thus $F=\frac12F_1$, where $F_1\cap C_1=Q_1$ and
$F_1\cdot C_1=1$. This implies $F_1\equiv C_2\equiv\frac12 H$. But
then $1=r<\frac67(d_1+d_2)+q= \frac67+\frac14$, a contradiction.

\subsubsection{Subcase $d_1=1/3$}
Since $4C_1$ is not Cartier, $\Diff{C_1}(F)=\frac
12Q_1+\frac23Q_2$ and $Q_2\in X$ is singular (of type $A_2$ or
$\frac13(1,1)$). Moreover, no components of $F$ pass through
$Q_2$. Further,
\begin{equation*}
1=C_1\cdot C_2=d_1H\cdot C_2,\quad H\cdot C_2=3,\quad d_2H^2=3.
\end{equation*}

Since $1>d_2=\frac3{H^2}\ge d_1=\frac13$, $9\ge H^2\ge 4$. On the
other hand, $H\cdot C_1=\frac13H^2\in \NN$. Thus $H^2=6$ or $9$.
Further, by Lemma~\ref{del-Pezzo-index>1}, $r=\frac{2l}{H^2}-1$,
where $l\in\NN$ and $l\le H^2+1$.

If $H^2=6$, then $d_2=1/2$ and
\[
1\ge r=\frac{l}3-1\ge \frac67\left(\frac13+\frac12\right)=\frac57.
\]
This gives $l=6$ and $r=1$. By Lemma~\ref{del-Pezzo-index=1},
$X\simeq \PP(1,2,3)$. In particular, $\Weillin(X)\simeq \ZZ$.
Since $-(K_X+C)\equiv (1-1/3-1/2)H$ is ample, $F\ne 0$. Therefore
$Q_1=\Supp F\cap C_1$ and moreover $Q_1\in X$ is smooth, $F=\frac
12F_1$ and the intersection of $F_1$ and $C_1$ is transverse. Thus
$1=F_1\cdot C_1=\frac13F_1\cdot H$ and $F_1\equiv \frac12H$. We
may assume that $C_1=\{x_2=0\}$, $C_2=\{x_3=0\}$, and
$F_1=\{x_3=\alpha_1x_1^3+\alpha_2x_1x_2\}$, $\alpha_1,
\alpha_2\in\CC$. But if $F_1=\{x_3=x_1^3\}$, then $K_X+C+F$ is not
lc at $(0,1,0)$. On the other hand, if $F_1=\{x_3=x_1x_2\}$, then
$F_1$ passes through the point $C_1\cap C_2$, a contradiction.
Therefore $\alpha_1, \alpha_2\ne 0$ and we may put
$F_1=\{x_3=x_1^3+x_1x_2\}$. This is case $\mathrm{(A_2^4)}$.

If $H^2=9$, then $d_2=1/3$ and
\[
1\ge r=\frac{2l}9-1\ge \frac67\left(\frac
13+\frac13\right)=\frac47, \quad l\in\ZZ.
\]
This gives $l=9$ or $l=8$. But in the first case $r=1$ which is a
contradiction with $H^2=9$ (see \ref{del-Pezzo-index=1}). Hence
$l=8$ and $r=7/9$. Since $d_1=d_2$, similar to \ref{eq-cases-diff}
we have $\Diff{C_2}(F)= \frac12Q_1'+ \frac23Q_2'$. In particular,
this means that $C$ contains no points of index $>3$. But
$X\setminus (C)$ contains such a point (because $r=7/9$), a
contradiction with \ref{Keel-McCern-}.

\subsubsection{Subcase $d_1=1/4$}
Since $mC_1$ is not Cartier for $m<4$,
$\Diff{C_1}(F)=\frac12Q_1+\frac34Q_2\ge \Diff{C_1}(0)$ and $Q_2\in
X$ is a singular point of type $A_3$ or $\frac14(1,1)$. By
Theorem~\ref{Keel-McCern-}, $Q_1\in X$ is smooth. Thus
$F=\frac12F_1$, where $F_1\cap C_1=Q_1$ and $C_1\cdot F_1=1$. Put
$k:=H\cdot C_1$. Then $H^2=4k$, $ d_2=1/k$. Since $d_2\ge d_1$,
$k\le 4$. If $F_1\equiv q_1H$, then $1=C_1\cdot
F_1=\frac14q_1H^2=q_1k$. Hence $F_1\equiv \frac1kH$. Further, by
Lemma~\ref{del-Pezzo-index>1},
\[
r=\frac{l}{2k}-1\ge \frac67(d_1+d_2)+ \frac12q_1=
\frac67\left(\frac14+ \frac1k\right)+\frac1{2k},\quad l-2k-2\ge
\frac{3k+5}7.
\]
On the other hand, $K_X+C+ \frac12F_1$ is ample, so $0<-r+d_1+d_2+
\frac12q_1$. This gives
\[
0< -\frac{l}{2k}+1+ \frac14+\frac1k+\frac1{2k}=
\frac{-l+2k+k/2+2+1}{2k},\quad l-2k-2<k/2+1.
\]
We get the following case:
\[
k=3,\quad l=10,\quad r=2/3,\quad d_2=1/3,\quad F_1\equiv
\frac13H,\quad H^2=12.
\]
We claim that $K_X+C$ is $1$-complementary. Note that
$-(K_X+C)\equiv (\frac23-\frac14-\frac13)H$ is ample. By
Theorem~\ref{Keel-McCern-} and because $r=2/3$, $C_2$ contains
exactly one singular point of $X$, say $Q'$. Therefore
$\Diff{C}(0)$ is supported at two points $Q'$ and $Q_2$. It is
easy to verify that $K_{C}+Q'+Q_2$ is an $1$-complement. By
Proposition~\ref{prodoljenie} this complement gives an
$1$-complement $K_X+C+\Theta$, where $\Theta$ is reduced and
$\Theta\cap C=\{Q',Q_2\}$. By Theorem~\ref{th-toric-1},
$(X,C+\Theta)$ is a toric pair. Such $X$ is defined by a fan
$\varDelta$ in $\RR^2=\ZZ^2\otimes\RR$. Let $v_1$, $v_2$, $v_3$ be
generators of one-dimensional cones in $\varDelta$. Since
$X\setminus C$ is smooth, we may assume that $v_1$ and $v_2$
generate $\ZZ^2$. Thus we can put $v_1=(1,0,0)$ and $v_2=(0,1,0)$.
Therefore $X$ is a weighted projective space $\PP(1,a_2,a_3)$,
$C_1\sim \OOO_{\PP}(a_2)$, $C_2\sim \OOO_{\PP}(a_3)$ and
${-}K_X\sim \OOO_{\PP}(1+a_2+a_3)$. Since $X\ni Q_2$ is singular
of type $\frac14(1,s)$, where $s=1$ or $3$ and $Q_2\in C_1$, we
can take $a_3=4$. Finally, from
\[
K_X^2=\left(\frac23H\right)^2=\frac{16}3,\quad K_X^2=
\frac{(a_1+a_2+a_3)^2}{a_1a_2a_3}= \frac{(5+a_2)^2}{4a_2}
\]
we obtain $a_2=3$. This is case $\mathrm{(A_2^5)}$.

\subsubsection{Subcase $d_1=1/6$}
Since $mC_1$ is not Cartier for $m<6$,
$\Diff{C_1}(F)=\frac12Q_1+\frac23Q_2\ge \Diff{C_1}(0)$ and
$\Diff{C_1}(0)=\Diff{C_1}(F)=\frac12Q_1+\frac23Q_2$. Hence $F=0$
and points $Q_1,Q_2\in X$ are singular. This contradicts to
Theorem~\ref{Keel-McCern-}.
\par
Theorem~\ref{delta=2} is proved.
\end{proof}

Theorem~\ref{delta=2} completes the classification of log pairs
with $\delta(X,B)=2$. The case $\delta(X,B)=1$ was studied by Abe
\cite{Abe}. In particular, he completely described so called
``elliptic curve case'', i.e., the case $p_a(C)=1$. A different
approach to the classification of exceptional complements was
given in \cite{KeM}.

\section{Examples}
\begin{example1}
\label{P-2-ex-1}
Let $X=\PP^2$ and $B=\sum d_iB_i$, where all $B_i$ are lines on
$\PP^2$ such that no three of them pass through one point, and
$d_i=1-1/m_i$. Assume that $-(K_X+B)$ is ample. By definition,
$K_X+B$ is $n$-complementary if and only if
$\deg({-}nK_X-\down{(n+1)B})\ge 0$ (i.e.,
$\sum\down{(n+1)(1-1/m_i)}\le 3n$). We give the list of all
possibilities for $(m_1,\dots,m_r)$ (with $m_1\le\dots\le m_r$).
These were found by means of a computer program. Here
$n=\compl(X,B)$.
\par\medskip
\textsc{Nonexceptional pairs}\par\nopagebreak
\begin{description}
\item[$n=1$] $(m)$, $(m_1,m_2)$, $(m_1,m_2,m_3)$
(type $\BA1^3_0$, see \ref{types});
\item[$n=2$] $(2,2, m_1, m_2)$, $(2,2,2,2,m)$
(types $\DD 2^2_0$ and $\EE 2^1_0$, respectively);
\item[$n=3$] $(2,3,3,m)$, $(3,3,3,m)$
(type $\EE 3^1_0$);
\item[$n=4$] $(2,3,4,m)$, $(2,4,4,m)$
(type $\EE 4^1_0$);
\item[$n=5$] $(2,3,5,5)$
(there is also a regular $6$-complement of type $\EE 6^1_0$);
\item[$n=6$] $(2,3,5,m)$, $m\ge 6$, $(2,3,6,m)$ (type $\EE 6^1_0$).
\end{description}
\par\medskip
\textsc{Exceptional pairs} \nopagebreak
\begin{description}
\item[$n=4$]
$(3,3,4,4)$, $(3,4,4,4)$, $(2,2,2,3,3)$, $(2,2,2,3,4)$;
\item[$n=5$]
$(2,4,5,5)$, $(2,5,5,5)$ (in these cases there are also regular
$6$-complements);
\item[$n=6$]
$(2,4,5,6)$, $(2,4,6,6)$, $(2,5,5,6)$, $(2,5,6,6)$, $(3,3,4,5)$,
$(3,3,5,5)$, $(3,3,5,6)$, $(3,3,4,6)$, $(2,2,2,3,5)$;
\item[$n=7$]
$(2,3,7,7)$;
\item[$n=8$]
$(2,3,7,8)$, $(2,3,8,8)$, $(2,4,5,7)$, $(2,4,5,8)$, $(2,4,6,7)$,
$(2,4,6,8)$, $(2,4,7,7)$, $(2,4,7,8)$;
\item[$n=9$]
$(2,3,7,9)$, $(2,3,8,9)$, $(2,3,9,9)$, $(3,3,4,7)$, $(3,3,4,8)$,
$(3,3,4,9)$;
\item[$n=10$]
$(2,3,7,10)$, $(2,3,8,10)$, $(2,3,9,10)$, $(2,3,10,10)$,
$(2,4,5,9)$, $(2,4,5,10)$, $(2,5,5,7)$, $(2,5,5,8)$, $(2,5,5,9)$;
\item[$n=12$]
$(2,3,7,11)$, $(2,3,7,12)$, $(2,3,8,11)$, $(2,3,8,12)$,
$(2,3,9,11)$, $(2,3,9,12)$, $(2,3,10,11)$, $(2,3,10,12)$,
$(2,3,11,11)$, $(2,3,11,12)$, $(2,4,5,11)$, $(2,4,5,12)$,
$(2,4,6,9)$, $(2,4,6,10)$, $(2,4,6,11)$, $(3,3,4,10)$,
$(3,3,4,11)$, $(3,4,4,5)$;
\item[$n=14$]
$(2,3,7,13)$, $(2,3,7,14)$;
\item[$n=15$]
$(3,3,5,7)$, $(2,3,7,15)$;
\item[$n=16$]
$(2,3,7,16)$, $(2,3,8,13)$, $(2,3,8,14)$, $(2,3,8,15)$,
$(2,3,8,16)$, $(2,4,5,13)$, $(2,4,5,14)$, $(2,4,5,15)$,
$(2,4,5,16)$;
\item[$n=18$]
$(2,3,7,17)$, $(2,3,7,18)$, $(2,3,8,17)$, $(2,3,8,18)$,
$(2,3,9,13)$, $(2,3,9,14)$, $(2,3,9,15)$, $(2,3,9,16)$,
$(2,3,9,17)$;
\item[$n=20$]
$(2,4,5,17)$, $(2,4,5,18)$, $(2,4,5,19)$;
\item[$n=21$]
$(2,3,7,19)$, $(2,3,7,20)$, $(2,3,7,21)$;
\item[$n=22$]
$(2,3,7,22)$;
\item[$n=24$]
$(2,3,7,23)$, $(2,3,7,24)$, $(2,3,8,19)$, $(2,3,8,20)$,
$(2,3,8,21)$, $(2,3,8,22)$, $(2,3,8,23)$;
\item[$n=28$]
$(2,3,7,25)$, $(2,3,7,26)$, $(2,3,7,27)$, $(2,3,7,28)$,
$(2,4,7,9)$;
\item[$n=30$]
$(2,3,7,29)$, $(2,3,7,30)$, $(2,3,10,13)$, $(2,3,10,14)$,
$(2,5,6,7)$;
\item[$n=36$]
$(2,3,7,31)$, $(2,3,7,32)$, $(2,3,7,33)$, $(2,3,7,34)$,
$(2,3,7,35)$, $(2,3,7,36)$;
\item[$n=42$]
$(2,3,7,37)$, $(2,3,7,38)$, $(2,3,7,39)$, $(2,3,7,40)$,
$(2,3,7,41)$;
\item[$n=66$]
$(2,3,11,13)$.
\end{description}
Thus the set of all $\compl(X,B)$ in this case is
\begin{multline*}
\{1, 2, 3, 4, 5, 6, 7, 8, 9, 10, 12, 14, 15, 16, 18, 20, 21,\\ 22,
24, 28, 30, 36, 42, 66\}.
\end{multline*}
It is easy to see that this set is contained in $\{n\in\NN\mid
\var(n)\le 20, n\ne 60\}$, which is related to automorphisms of
$K3$ surfaces \cite{Ishii} (see also \cite[Sect. 2]{Ts}).
\end{example1}

\begin{example1}
\label{P-2-ex-2}
Replace the condition of the ampleness of $-(K_X+B)$ in
Example~\ref{P-2-ex-1} with numerical triviality. We obtain only
exceptional cases:
\begin{description}
\item[$n=2$]
$(2,2,2,2,2,2)$;
\item[$n=4$]
$(4,4,4,4)$, $(2,2,2,4,4)$;
\item[$n=6$]
$(2,6,6,6)$, $(3,3,6,6)$, $(2,2,2,3,6)$, $(2,2,3,3,3)$;
\item[$n=8$]
$(2,4,8,8)$;
\item[$n=10$]
$(2,5,5,10)$;
\item[$n=12$]
$(2,3,12,12)$, $(2,4,6,12)$, $(3,3,4,12)$, $(3,4,4,6)$;
\item[$n=18$]
$(2,3,9,18)$;
\item[$n=20$]
$(2,4,5,20)$;
\item[$n=24$]
$(2,3,8,24)$;
\item[$n=30$]
$(2,3,10,15)$;
\item[$n=42$]
$(2,3,7,42)$
\end{description}
In these cases $(X,B)$ is a log Enriques surface and $n(K_X+B)\sim
0$. Construction~\ref{l-c-cover} gives a ramified cyclic cover
$\var\colon X'\to\PP^2$ such that $K_{X'}=\var^*(K_X+B)$. Then
$K_{X'}\sim 0$ and is plt, so $X'$ is a surface with Du Val
singularities and $K_{X'}\sim 0$. Note that if we replace the
condition $B\in \Msm$ with $B\in \MMM_{m}$, we can get bigger
values of $\compl(X,B)$. For example, take
$B=\frac12B_1+\frac23B_2+ \frac{18}{19}B_3+ \frac{101}{114}B_4$,
where, as above, $B_i\subset \PP^2$ are lines such that no three
of them pass through one point. Then $\compl(X,B)=78$.
\end{example1}

\begin{example1}
Let $G\subset \PGL_3(\CC)$ be a finite subgroup, $X:=\PP^2/G$, and
$f\colon \PP^2\to X$ the quotient morphism. Define the boundary
$B$ on $X$ by $K_{\PP^2}=f^*(K_X+B)$ (see \eref{eqv-vetvl} and
\eref{eqv-vetvl1_2}). Then $(X,B)$ is exceptional if and only if
$G$ has no semiinvariants of degree $\le 3$ (see \cite{MP}). There
are only four types of such groups up to conjugation in
$\PGL_3(\CC)$.
\end{example1}

\begin{example1}[{\cite{Abe}}]
Let $X:=\PP(1,2,3)$. Take a general member $E\in |{-}K_X|$ (a
smooth elliptic curve) and let $L$ be a line on $X$ (with respect
to ${-}K_X$). Then $E\sim 6L$. Since $(X,L)$ is toric, $K_X+L$ is
plt. Hence $(X,\alpha E+\beta L)$ is a log del Pezzo if and only
if $6\alpha+ \beta<6$, $\alpha\le 1$, $\beta\le 1$. Moreover, if
$\alpha\ge 6/7$ and $\beta\in\Mm$, then $(X,\alpha E+\beta L)$ is
exceptional. Indeed, by Corollary~\ref{corr-Inductive} it is
sufficient to show that there are no regular nonklt complements.
If $K_X+B^+$ is such a complement, then $B^+\ge E+\beta L$, a
contradiction. This gives the following exceptional cases with
$\delta=1$:
\[
\begin{array}{lcl}
 \beta=1/2&\qquad& 6/7\le\alpha<11/12\\
 \beta=2/3&& 6/7\le\alpha<8/9\\
 \beta=3/4&& 6/7\le\alpha<7/8\\
 \beta=4/5&& 6/7\le\alpha<13/15\\
 \beta=5/6&& 6/7\le\alpha<31/36.\\
\end{array}
\]
\end{example1}

\begin{example1}[{\cite{Abe}}]
Let $X\subset\PP^3$ be a quadratic cone, $E\in |{-}K_X|$ a smooth
elliptic curve, and $L$ a generator of the cone. Then
$(X,\frac67E+\frac12L)$ is an exceptional log del Pezzo with
$\delta=1$ and $K_X+\frac67E+\frac47L$ is a $7$-complement.
\end{example1}

\begin{exercise1}
Let $C\subset\PP^2$ be a smooth curve of degree $d$. Assume that
$-(K_X+(1-1/m)C)$ is nef. Prove that $K_X+(1-1/m)C$ is exceptional
only if and only if $(d,m)\in\{(4,3), (4,4), (5,2), (6,2)\}$. For
$(d,m)=(4,3)$, $(5,2)$ such log Del Pezzos can appear as
exceptional divisors of plt blowups of canonical singularities
(see \cite{Pr1}).
\begin{hint}
The nontrivial part is to prove that $K_X+(1-1/m)C$ is exceptional
in these cases. Assuming the opposite we have a regular nonklt
complement $K_X+B$. Then we can use the following simple fact: if
$\sum d_iB_i$ is a boundary on $\CC^2$ such that all the $B_i$ are
smooth curves and $\sum d_i\le 1$, then $(\CC^2,B)$ is canonical.
\end{hint}
\end{exercise1}

\begin{example1}
\label{lc-thresholds-3}
Let $(X\ni o)$ be a three-dimensional klt singularity and $D$ an
effective reduced Weil divisor on $X$. Assume that $D$ is
$\QQ$-Cartier. Let $c_o(X,D)$ be the log canonical threshold.
Assume that $1>c:=c_o(X,D)>6/7$. Let $f\colon Y\to X$ be a plt
blowup of $(X,D)$. Write $K_Y+S+cB=f^*(K_X+cD)$, where $B$ is the
proper transform of $D$. Then $(S,\Diff S(cB))$ is a log Enriques
surface with $\delta\ge 1$. We claim that $K_S+\Diff S(cB)$ is
klt. Indeed, if $K_S+\Diff S(cB)$ is not klt, then by the
Inductive Theorem~\ref{Inductive-Theorem} there is a regular
complement $K_S+\Diff S(cB)^+$. Since $-(K_Y+S+(c-\ep)B)$ is
$f$-ample for $\ep>0$, by Proposition~\ref{prodolj} we have a
regular complement $K_Y+S+(c-\ep)B$. This gives a regular
complement $K_X+A$ of $K_X+(c-\ep)D$. We can take $\ep$ so that
$c-\ep>6/7$. Then $A$ is reduced and $A=D$. Hence $c=1$, a
contradiction. This method can help to describe the set of all lc
thresholds in the interval $[6/7,1]$ (cf. \cite{Ku}). For example,
take $X=\CC^3$ and $D=\{\psi(x,y,z)=0\}$, where
$\psi(x,y,z)=x^3+yz^2+x^2y^2+x^5z$ (see \cite{Ku}). Then
$c_o(\CC^3,D)=11/12$ and $f\colon Y\to \CC^3$ is the weighted
blowup with weights $(4,2,5)$. So $S=\PP(4,2,5)$. It is easy to
compute that $\Diff S(cD)= \frac{11}{12}C+\frac12L$, where
$C:=\{x^3+yz^2+x^2y^2=0\}$ and $L:=\{z=0\}$. Both $C$ and $L$ are
smooth rational curves which intersect each other twice at smooth
points of $S$. Such complements were studied in
\cite[Sect.~2]{Abe} and called there ``sesqui rational curve''
complements.
\end{example1}

\chapter{Appendix}
\section{Existence of complements}
\begin{proposition1}[\cite{Sh1}]
\label{E-C}
Let $f\colon X\to Z\ni o$ be a contraction from a surface and $D$
a boundary on $X$ such that $K_X+D$ is lc and $-(K_X+D)$ is
$f$-nef and $f$-big. Then
\begin{enumerate}
\item
the linear system $|-m(K_X+D)|$ is base point free for some
$m\in\NN$;
\item
$K_X+D$ is $n$-complementary near $f^{-1}(o)$ for some $n\in\NN$;
\item
the Mori cone $\NE(X/Z)$ is polyhedral and generated by
irreducible curves.
\end{enumerate}
\end{proposition1}

We hope that this fact has higher dimensional generalizations (cf.
\cite{Kawamata-models}, see also M.~Reid's Appendix to \cite{Sh}).

\begin{proof}
First we prove (i). We consider only the case of compact $X$. In
the case $\dim Z\ge 1$ there are stronger results (see
Theorem~\ref{local}). Applying a log terminal modification
\ref{blow-up}, we may assume that $K_X+D$ is dlt (and $X$ is
smooth). Set $C:=\down{D}$, $B:=\fr{D}$. Note that $C$ is
connected by Connectedness Lemma. Take sufficiently large and
divisible $n\in\NN$ and consider the exact sequence
\begin{multline*}
0\longrightarrow\OOO_X(-n(K_X+D)-C)
\longrightarrow\OOO_X(-n(K_X+D))\\
\longrightarrow\OOO_C(-n(K_X+D)) \longrightarrow 0.
\end{multline*}
By Kawamata-Viehweg Vanishing \cite[1-2-6]{KMM},
\begin{multline*}
H^1(X,\OOO_X(-n(K_X+D)-C))=\\ H^1(X,\OOO_X(K_X+B-(n+1)(K_X+D)))=0.
\end{multline*}
Therefore $C\cap\Bs |-n(K_X+D)|=\Bs \bigl|-n(K_X+D)|_C\bigr|$.
\par
We claim that $\Bs \bigl|-n(K_X+D)|_C\bigr|=\emptyset$. Indeed, if
$C$ is not a tree of rational curves, then $p_a(C)=1$ and $C$ is
either a smooth elliptic curve or a wheel of smooth rational
curves (see Lemma~\ref{wheel}). Moreover, $\Supp B\cap
C=\emptyset$. But then $(K_X+D)|_C=(K_X+C)|_C=K_C=0$ and $\Bs
\bigl|-n(K_X+D)|_C\bigr|=\emptyset$ in this case. Note also that
here we have an $1$-complement by Lemma~\ref{wheel-1}. Assume now
that $C$ is a tree of smooth rational curves. Then
$\bigl|-n(K_X+D)|_{C_i}\bigr|$ is base point free on each
component $C_i\subset C$ whenever $-n(K_X+D)$ is Cartier. Hence so
is $\bigl|-n(K_X+D)|_C\bigr|$. This proves our claim.
\par
Thus we have shown that $C\cap\Bs |-n(K_X+D)|=\emptyset$. Let
$L\in |-n(K_X+D)|$ be a general member. Then $K_X+D+\frac1{n}L$ is
dlt near $C$ (see \ref{Bertini}). By Connectedness Lemma,
$K_X+D+\frac1{n}L$ is lc everywhere. Hence $K_X+D+\frac1{n}L$ is a
$\QQ$-complement of $K_X+D$. The fact that $|-n(K_X+D)|$ is free
outside of $C$ can be proved in a usual way (see e.g., \cite{R},
\cite{Kawamata-models}). We omit it.

(ii) is obvious. Let us prove (iii). Clearly, we may assume that
$\rho(X)\ge 2$. It follows by \ref{MMP-cone} that any
$(K_X+D)$-negative extremal ray $R$ is generated by an irreducible
curve $C$. By Proposition~\ref{MMP-extremal-curves},
$-(K_X+D)\cdot C\le 2$. Let $\varphi\colon X\to Y\subset \PP^N$ be
the contraction given by the linear system $-m(K_X+D)$ for
sufficiently big and divisible $m\in\NN$. Then $\deg\varphi(C)\le
2$. This implies that $C$ belongs to a finite number of algebraic
families. Thus the cone $\NE(X)$ is polyhedral outside of
$\NE(X)\cap \{z \mid (K_X+D)\cdot z=0\}$. Now consider the
extremal ray $R$ such that $(K_X+D)\cdot R=0$. By the Hodge Index
Theorem, $R^2<0$. Thus, by Proposition~\ref{MMP-properties} $R$ is
generated by an irreducible curve, say $C$. Since $(K_X+D)\cdot
C=0$, we have that $\varphi$ contracts this curve to a point.
Therefore there is a finite number of such curves, so $\NE(X)$ is
polyhedral everywhere.
\end{proof}

\section[MMP in dimension two]{Minimal Model Program 
in dimension two}
\label{MMP-in-dimension-2}
The log Minimal Model Program in dimension two is much easier than
in higher dimensions. Following \cite{A} and \cite{KoKov} (see
also \cite{Sh2}) we present two main theorems \ref{MMP-cone} and
\ref{MMP-contraction} of MMP in the surface case. First we note
that in the surface case it is possible to define the
\textit{numerical pull back} of any $\QQ$-Weil divisor under
birational contractions (see e.g., \cite{Sakai}). Therefore all
definitions of \ref{Singularities-of-pairs} can be given for
arbitrary normal surface (we need not the $\QQ$-Cartier
assumption). It turns out a posteriory that any numerically lc
pair $(X,B)$ satisfies the property that $K_X+B$ is $\QQ$-Cartier
\cite[Sect. 4.1]{KM}, \cite{Masek}. Similarly, the dlt property of
$(X,D)$ implies that the surface $X$ is $\QQ$-factorial
\cite[Sect. 4.1]{KM}. For surfaces there is an alternative way to
define the numerical equivalence: two $1$-cycles
$\Upsilon_1,\Upsilon_2\in Z_1(X/Z)$ are said to be numerically
equivalent if $L\cdot \Upsilon_1=L\cdot \Upsilon_2$ for \emph{all}
Weil divisors $L$ (not only for those, that are $\QQ$-Cartier).
This gives also an alternative way to define $N_1(X/Z)$, $\rho(X/Z)$, and
$\NE(X/Z)$ and leads to a possibly larger dimensional space $N_1(X/Z)$. We
use the standard definition of the numerical equivalence and
$N_1(X/Z)$ \cite{KMM}.

The following properties are well known (see e.g., 
{\cite[1.21--1.22]{KM}} or {\cite[Ch. II, Lemma
4.12]{Ko3}}).
\begin{proposition1}[Properties of the Mori cone]
\label{MMP-properties}
Let $X$ be a normal projective surface.
\begin{enumerate}
\item
Let $z$ be an element of $N_1(X)$ such that $z^2>0$ and $z\cdot
H>0$ for some ample divisor $H$. Then $z$ is contained in the
interior of $\NE(X)$.
\item
Let $C\subset X$ be an irreducible curve. If $C^2\le 0$, then the
class $[C]$ is in the boundary of $\NE(X)$. If $C^2<0$, then the
ray $\RR_+[C]$ is extremal.
\item
Let $R\subset \NE(X)$ be an extremal ray such that $R^2<0$. Then
$R$ is generated by an irreducible curve.
\end{enumerate}
\end{proposition1}
\begin{proof}
We prove only (iii). Take a $1$-cycle $Z$ so that $[Z]\in R$,
$[Z]\neq 0$ and $Z_i$ a sequence of effective $1$-cycles whose
limit is $Z$. Write $Z_i=\sum_j a_{i,j}C_j$, where $C_j$ are
distinct irreducible curves. Since $0>Z^2=\lim Z\cdot Z_i$, there
is at least one curve $C=C_k$ such that $Z\cdot C<0$. Write
$Z_i=c_iC+\sum_{j\neq k} a_{i,j}C_j$, $c_i\ge 0$. Then
\[
0>C\cdot Z=\lim C\cdot Z_i\ge (\lim c_i) C^2.
\]
Thus $C^2<0$ and $\lim c_i>0$. Pick $0<c<\lim c_i$. Then $Z_i-cC$
is effective for $i\gg 0$ and $Z=cC+\lim (Z_i-cC)$ Since $R$ is an
extremal ray, this implies that $[C]\in R$.
\end{proof}

\begin{theorem1}[The Cone Theorem]
\label{MMP-cone}
Let $X$ be a normal projective surface and $K_X+B$ be an effective
$\RR$-Cartier divisor. Let $A$ be an ample divisor on $X$. Then
for any $\ep>0$ the Mori-Kleiman cone of effective curves $\NE(X)$
in $N_1(X)$ can be written as
\[
\NE(X)= \NE_{K+B+\ep A}(X) +\sum R_k
\]
where, as usual, the first part consists of cycles that have
positive intersection with $K+B+\ep A$ and $R_k$ are finitely many
extremal rays. Each of the extremal rays is generated by an
effective curve.
\end{theorem1}

\begin{theorem1}[Contraction Theorem]
\label{MMP-contraction}
Let $X$ be a projective surface with log canonical $K_X+B$. Let
$R$ be a $(K_X+B)$-negative extremal ray. Then there exists a
nontrivial projective morphism $\phi\colon X\to Z$ such that
$\phi_*(\OOO_X)=\OOO_Z$ and $\phi(C)=pt$ if and only if the class
of $C$ belongs to $R$. Moreover, if $\phi$ is birational and
$K_X+B$ is lc (resp. klt) then $K_Z+\phi_*B$ is lc (resp. klt).
\end{theorem1}

\begin{remark1}
Notation as above.
\begin{enumerate}
\item
If $\dim Z=1$, then all fibers of $\phi$ are irreducible smooth
rational curves and $X$ has only rational singularities
\cite{KoKov}.
\item
If $\dim Z=2$, then $C\simeq \PP^1$ and $K_Z+\phi_*B$ is plt at
$\phi(C)$.
\end{enumerate}
\end{remark1}

\begin{proposition1}[Properties of extremal curves]
\label{MMP-extremal-curves}
Let $(X,B)$ be a normal projective log surface and $R$ a
$(K_X+B)$-negative extremal ray on $X$. Assume that $K_X+B$ is lc.
If $R^2\le 0$, then for any irreducible curve $C$ such that
$[C]\in R$ we have $-(K_X+D)\cdot C\le 2$. If $R^2>0$, then $X$ is
covered by a family of rational curves $C_{\lambda}$ such that
$-(K_X+D)\cdot C_{\lambda}\le 3$.
\end{proposition1}
\begin{proof}
Let $\mu\colon \tilde X\to X$ be the minimal resolution and
$K_{\tilde X}+\tilde B=\mu^*(K_X+B)$ the crepant pull back.

Consider the case $R^2\le 0$. Let $\tilde C$ be the proper
transform of $C$. Then
\[
-(K_X+B)\cdot C=-(K_{\tilde X}+\tilde B)\cdot \tilde C\le
-(K_{\tilde X}+\tilde C)\cdot \tilde C\le 2
\]
because $\tilde C^2\le C^2\le 0$ and $\tilde B$ is a boundary.

Now we assume that 
$R^2>0$. Then $-(K_X+B)$ is ample (see \ref{MMP-properties}).
Thus $(X,B)$ is a log del Pezzo surface. By
Corollary~\ref{cor-log-del-Pezzo}, $\tilde X$ is birationally
ruled. It is well known, that in this situation $\tilde X$ is
covered by a family of rational curves $\tilde C_{\lambda}$ such
that ${-}K_{\tilde X}\cdot \tilde C_{\lambda}\le 3$.
Take
$C_{\lambda}=\mu(\tilde C_{\lambda})$. Then
\[
-(K_{X}+B)\cdot C_{\lambda}=-(K_{\tilde X}+\tilde B)\cdot \tilde
C_{\lambda}\le {-}K_{\tilde X}\cdot \tilde C_{\lambda}\le 3.
\]
\end{proof}

\begin{theindex}

  \item $N_1(X/Z)$, 5
  \item $\CC^n/\cyc{m}(a_1,\dots,a_n)$, 12
  \item $\Diff S(B)$, 15
  \item $\FF_n$, 40
  \item $\LCS(X,D)$, 20
  \item $\Mm$, 15
  \item $\Msm$, 15
  \item $\NE(X/Z)$, 5
  \item $\NNN_n(\MMM)$, 32
  \item $\PP(a_1,\dots, a_n)$, 24
  \item $\PPP_n$, 33
  \item $\QQ$-complement, 28
  \item $\RRR_2$, 31
  \item $\Sigma_0$, 40
  \item $\Weil(X)$, 5
  \item $\Weilalg(X)$, 5
  \item $\Weillin(X)$, 5
  \item $\delta(X,B)$, 102
  \item $\down{\cdot}$, 5
  \item $\ep$-log canonical ($\ep$-lc) singularities, 6
  \item $\ep$-log terminal ($\ep$-lt) singularities, 6
  \item $\fr{\cdot}$, 5
  \item $\mathfrak{A}_n$, 55
  \item $\mathfrak{D}_n$, 55
  \item $\mathfrak{S}_n$, 55
  \item $\mt{discrep}(X,D)$, 6
  \item $\rho(X)$, 5
  \item $\rhonum(X/Z)$, 57
  \item $\up{\cdot}$, 5
  \item $a(E,D)$, 6
  \item $c_o(X,D)$, 47
  \item $n$-complement, 29
  \item $n$-semicomplement, 29
  \item $r(X,D)$, 110
  \item ${\operatorname {compl}}'(X,D)$, 32
  \item ${\operatorname {compl}}(X,D)$, 32

  \indexspace

  \item Adjunction, 15

  \indexspace

  \item blowdown, 5
  \item blowup, 5
  \item boundary, 5

  \indexspace

  \item canonical singularities, 6
  \item center of log canonical singularities, 20
  \item complement, 29
  \item Connectedness Lemma, 19
  \item contraction, 5
  \item crepant pull back, 8

  \indexspace

  \item different, 15
  \item discrepancy, 6
  \item divisorial log terminal (dlt), 6
  \item dlt model of elliptic fibration, 65

  \indexspace

  \item elliptic fibration, 65
    \subitem minimal, 65
  \item elliptic singularity, 54
  \item exceptional log variety, 37

  \indexspace

  \item Fano index, 110

  \indexspace

  \item index, 10
  \item inductive blowup, 23
  \item Inductive Theorem, 77
  \item Inversion of Adjunction, 15

  \indexspace

  \item Kawamata log terminal (klt), 6

  \indexspace

  \item locus of log canonical singularities, 20
  \item log Calabi-Yau variety, 40
  \item log canonical (lc) singularities, 6
  \item log canonical threshold, 47
  \item log del Pezzo surface, 40
  \item log Enriques surface, 40
  \item log Fano variety, 40
  \item log resolution, 6
  \item log terminal modification, 21
    \subitem minimal, 22
  \item log variety (log pair), 5

  \indexspace

  \item maximally log canonical, 43

  \indexspace

  \item nodal curve, 29

  \indexspace

  \item purely log terminal (plt) blowup, 23
  \item purely log terminal (plt) singularities, 6

  \indexspace

  \item rational singularity, 54
  \item regular complement, 31
  \item regular log surface, 44

  \indexspace

  \item semilog canonical (slc) singularities, 29
  \item simple $K3$ singularity, 38
  \item simple elliptic singularity, 38
  \item standard coefficient, 15
  \item strong $n$-complement, 29
  \item subboundary, 5

  \indexspace

  \item terminal blowup, 23
  \item terminal singularities, 6
  \item toric pair, 18
  \item trivial complement, 34

  \indexspace

  \item weighted blowup, 24
  \item weighted projective space, 24
  \item weights, 12

\end{theindex}


\begin{thebibliography}{KMM1}

\bibitem[Ab]{Abe}
\textit{Abe T.} Classification of exceptional surface complements,
PhD thesis, Johns Hopkins Univ. (1998)

\bibitem[A]{A}
\textit{Alexeev V.} Boundedness and $K^2$ for log surfaces,
Internat. J. Math. \textbf{5} (1994) 779--810

\bibitem[A1]{A1}
\textit{Alexeev V.} Two two-dimensional terminations, Duke Math.
J. \textbf{69} (1992) 527--545

\bibitem[Ar]{Ar}
\textit{Artin M.} On isolated rational singularities of surfaces,
Amer. J. Math. \textbf{88} (1966) 129--136

\bibitem[BPV]{BPV}
\textit{Barth~W., Peters~C. $\&$ Van~de~Ven~A.} {Compact complex
surfaces.} Springer-Verlag, 1984

\bibitem[Bl]{Bl}
\textit{Blache R.} The structure of l.c. surfaces of Kodaira
dimension zero, J. Algebraic Geom. \textbf{4} (1995) 137--179

\bibitem[Br]{Brieskorn}
\textit{Brieskorn E.} Rationale Singularit\"aten komplexer
Fl\"achen, Invent. Math. \textbf{4} (1968) 336--358

\bibitem[Ca]{Catanese}
\textit{Catanese F.} {Automorphisms of rational double points and
moduli spaces of surfaces of general type,} Compositio Math.
\textbf{61} (1987) 81--102

\bibitem[F]{F}
\textit{Fujino O.} {Abundance theorem for semi log canonical
threefolds}, preprint RIMS-1213

\bibitem[Fuj]{Fujita}
\textit{Fujita T.} On the structure of polarized varieties with
$\Delta$-genera zero. J. Fac. Sci. Univ Tokyo. Sec. IA.
\textbf{22} (1975) 103--115

\bibitem[Fu]{Furushima}
\textit{Furushima M.} Singular del Pezzo surfaces and analytic
compactifications of $3$-dimensional complex affine space
$\CC^{3}$, Nagoya Math. J. \textbf{104} (1986) 1--28

\bibitem[Ha]{Ha}
\textit{Hartshorne R.} Algebraic Geometry, Springer, GTM
\textbf{52}, 1977

\bibitem[Hi]{Hi}
\textit{Hirzebruch F.} \"Uber vierdimensionale Riemannsche
Fl\"achen mehrdeutiger analytischer Funktionen von zwei
Ver\"anderlichen, Math. Ann. \textbf{126} (1953) 1-22

\bibitem[Il]{Iliev}
\textit{Iliev A.} Log-terminal singularities of algebraic
surfaces, Mosc. Univ. Math. Bull. \textbf{41} (1986) 46--53

\bibitem[I]{Ishii}
\textit{Ishii~S.} The quotients of log-canonical singularities by
finite groups, to appear in Adv. Stud. in Pure Math

\bibitem[I1]{Ishii1}
\textit{Ishii~S.} The global induces of log Calabi-Yau varieties
-- a supplement to Fujino's paper: the induces of log canonical
singularities-, e-preprint, math.AG/0003060

\bibitem[IP]{IP}
\textit{Ishii~S. $\&$ Prokhorov Yu. G.} Hypersurface exceptional
singularities, preprint TIT \textbf{92} (1999)

\bibitem[IW]{IW}
\textit{Ishii~S. $\&$ Watanabe~K.} A geometric characterization of
a simple $K3$-singularity, T\^ohoku Math. J. \textbf{44} (1992)
19--24

\bibitem[Is]{Isk}
\textit{Iskovskikh V.~A.} Fano $3$-folds, I $\&$ II, Izv. Akad.
Nauk SSSR Ser. Mat. \textbf{41} (1977) 516--562 $\&$ \textbf{42}
(1978) 504--549 (Russian) [English transl.: Math. USSR, Izv.
\textbf{11}, 485--527 $\&$ \textbf{12}, 469--506]

\bibitem[Ka]{Ka}
\textit{Kachi Y.} Flips from $4$-folds with isolated complete
intersection singularities, Amer. J. Math. \textbf{120} (1998)
43--102

\bibitem[K]{K}
\textit{Kawamata Y.} The crepant blowing-up of $3$-dimensional
canonical singularities and its application to the degeneration of
surfaces, Ann. Math. \textbf{127} (1988) 93--163

\bibitem[K1]{Kawamata-bound}
\textit{Kawamata Y.} Boundednes of $\QQ$-Fano threefolds, Contemp.
Math. AMS \textbf{131} (1992) 439--445

\bibitem[K2]{Kawamata-flip}
\textit{Kawamata Y.} Termination of log-flips for algebraic
threefolds, Internat. J. Math. \textbf{3} (1992) 653--659

\bibitem[K3]{Kawamata-models}
\textit{Kawamata Y.} Log canonical models of algebraic $3$-folds,
Internat. J. Math. \textbf{3} (1992) 351--357

\bibitem[KMM]{KMM}
\textit{Kawamata Y., Matsuda K. $\&$ Matsuki K.} Introduction to
the minimal model program, Adv. Stud. Pure Math. \textbf{10}
(1987) 283--360

\bibitem[KeM]{KeM}
\textit{Keel S. $\&$ McKernan J.} Rational curves on
quasi-projective surfaces, Memoirs AMS \textbf{140} (1999), no.
669

\bibitem[Ko1]{Ko1}
\textit{Koll\'ar J.} Log surfaces of general type; some
conjectures, {Contemporary Math. AMS} \textbf{162} (1994) 261--275

\bibitem[Ko2]{Ko}
\textit{Koll\'ar J.} Singularities of pairs, Proc. Symp. Pure
Math. \textbf{62} (1995) 221--287

\bibitem[Ko3]{Ko3}
\textit{Koll\'ar J.} Rational curves on algebraic varieties,
Ergeb. Math. Grenzgeb. (3), vol \textbf{32}, Springer (1996)

\bibitem[KK]{KoKov}
\textit{Koll\'ar J. $\&$ Kov\'acs S.} Birational geometry of log
surfaces, ftp://ftp.math.utah. edu /u/ma/kollar/

\bibitem[KoM]{KoM}
\textit{Koll\'ar J. $\&$ Mori S.} Classification of three
dimensional flips, J. Amer. Math. Soc. \textbf{5} (1995) 533--703

\bibitem[KM]{KM}
\textit{Koll\'ar J. $\&$ Mori S.} Birational geometry of algebraic
varieties, Cambridge Tracts in Mathematics \textbf{134} (1998)

\bibitem[Ut]{Ut}
\textit{Koll\'ar J. et al.} {Flips and abundance for algebraic
threefolds,} A summer seminar at the University of Utah, Salt Lake
City, 1991. Ast\'erisque. \textbf{211} (1992)

\bibitem[Ku]{Ku}
\textit{Kuwata T.} On log canonical thresholds of surfaces in
$\CC^3$, Tokyo J. Math. \textbf{22} (1999) 245--251

\bibitem[Ma]{Masek}
\textit{Ma\c{s}ek V.} Kawachi's invariant for normal surface
singularities, {Internat. J. Math.} \textbf{9} (1998) 623--640

\bibitem[M]{M}
\textit{Markushevich~D.} Minimal discrepancy for a terminal CDV
singularity is $1$, J. Math. Sci. Univ. Tokyo \textbf{3} (1997)
445--456

\bibitem[MP]{MP}
\textit{Markushevich~D. $\&$ Prokhorov~Yu.~G.} {Exceptional
quotient singularities,} {Amer. J. Math.} \textbf{121} (1999)
1179--1189

\bibitem[Mo]{Mo}
\textit{Mori S.} Flip theorem and the existence of minimal models
for $3$-folds, J. Amer. Math. Soc. \textbf{1} (1988) 117--253

\bibitem[MMM]{MMM} \textit{Mori S., Morrison D.R. $\&$ Morrison I.} On
four-dimensional terminal quotient singularities, Math. Comput.
\textbf{51} (1988) no. 184, 769--786


\bibitem[Mor]{Morrison}
\textit{Morrison D.} The birational geometry of surfaces with
rational double points, Math. Ann. \textbf{271} (1985) 415--438

\bibitem[N1]{N1}
\textit{Nikulin V.~V.} Del Pezzo surfaces with log-terminal
singularities, I. Mat. Sb. \textbf{180} (1989) 226--243 (Russian)
[English transl.: Math. USSR, Sb. \textbf{66} (1990) 231--248]

\bibitem[N2]{N3}
\textit{Nikulin V.~V.} Del Pezzo surfaces with log-terminal
singularities, III. Izv. Akad. Nauk SSSR Ser. Mat. \textbf{53}
(1989) 1316--1334 (Russian) [English transl.: Math. USSR, Izv.
\textbf{35} (1990) 657--675]

\bibitem[P]{Pr}
\textit{Prokhorov Yu.~G.} On extremal contractions from threefolds
to surfaces: The case of one nonsingular point, {Contemporary
Math. AMS} \textbf{207}, 119--141

\bibitem[P1]{Pr1}
\textit{Prokhorov Yu.~G.} {Algebra. Proc. Internat. Algebraic
Conf. on the Occasion of the 90th Birthday of A. G. Kurosh,}
Moscow, Russia, May 25-30, 1998, Yu. Bahturin ed., Walter de
Gruyter, Berlin (2000), 301--317

\bibitem[P2]{Pr2}
\textit{Prokhorov Yu.~G.} Boundedness of nonbirational extremal
contractions, {Internat. J. Math.} \textbf{11} (2000) 393--411

\bibitem[P3]{Pr3}
\textit{Prokhorov Yu.~G.} {Mori conic bundles with a reduced log
terminal boundary}, {J. Math. Sci., New York} \textbf{94} (1999)
no.~1, 1051--1059

\bibitem[PSh]{PSh}
\textit{Prokhorov Yu. G. $\&$ Shokurov V. V.} {The first main
theorem on complements: from global to local}, preprint TIT
\textbf{94} (1999)

\bibitem[R]{R}
\textit{Reid M.} Projective morphism accoding to Kawamata, Warwick
(preprint), 1983

\bibitem[R1]{Pagoda}
\textit{Reid M.} Minimal models of canonical threefolds, Adv.
Stud. Pure Math. \textbf{1} (1983) 131--180

\bibitem[R2]{Reid-canonical}
\textit{Reid M.} Canonical threefolds. G\'eom\'etrie Alg\'ebrique
Angers, A. Beauville ed., Sijthoff and Noordhoff (1980) 273--310

\bibitem[RY]{RYPG}
\textit{Reid M.} Young person's guide to canonical singularities,
Proc. Symp. in Pure Math. \textbf{46} (1987) 343--416

\bibitem[S1]{Sakai}
\textit{Sakai F.} Weil divisors on normal surfaces, Duke Math. J.
\textbf{51} (1984) 877--887

\bibitem[S2]{Sakai1}
\textit{Sakai F.} Classification of normal surfaces, Proc. Symp.
in Pure Math. \textbf{46} (1987) 451--466

\bibitem[S3]{Sakai2}
\textit{Sakai F.} The structure of normal surfaces, Duke Math. J.
\textbf{52} (1985) 627--648

\bibitem[Sh]{Sh0}
\textit{Shokurov V.~V.} Smoothness of the general anticanonical
divisor on Fano $3$-fold, Math USSR, Izv. \textbf{14} (1980)
395--405

\bibitem[Sh1]{Sh-nonvan}
\textit{Shokurov V.~V.} The nonvanishing theorem. Izv. Akad. Nauk
SSSR Ser. Mat. (1985) \textbf{49}, 635--651 (Russian) [English
transl.: Math USSR, Izv. \textbf{26} (1986) 591--604]

\bibitem[Sh2]{Sh}
\textit{Shokurov V.~V.} {$3$-fold log flips,}  
Izv. AN SSSR, Ser.
mat. \textbf{56} (1992), 105--201  (Russian)
[English transl. Russian Acad. Sci. Izv. Math. \textbf{40} (1993)
93--202] 

\bibitem[Sh3]{Sh1}
\textit{Shokurov V.~V.} {Complements on surfaces,} Math. Sci., New
York \textbf{102} (2000) no.~2, 3876--3932

\bibitem[Sh4]{Sh2}
\textit{Shokurov V.~V.} {3-fold log models,} J. Math. Sci., New
York \textbf{81} (1996) 2667--2699

\bibitem[Sz]{Sz}
\textit{Szab\'o E.} {Divisorial log terminal singularities,} {J.
Math. Sci. Univ. Tokyo} \textbf{1} (1994) 631--639

\bibitem[Ts]{Ts}
\textit{Tsunoda S.} Strusture of open algebraic surfaces, I, {J.
Math. Kyoto Univ.} \textbf{23-1} (1983) 95--125

\bibitem[Um]{Um}
\textit{Umezu Y.} On normal surfaces with trivial dualizing sheaf,
Tokyo J. Math. \textbf{4} (1981), 343--354

\bibitem[Z]{Z}
\textit{Zhang D.-Q.} Logarithmic Enriques surfaces, J. Math. Kyoto
Univ. \textbf{31} (1991) 419--466

\bibitem[Z1]{Z1}
\textit{Zhang D.-Q.} Logarithmic Enriques surfaces I, J. Math.
Kyoto Univ. \textbf{33} (1993) 357--397
\end{thebibliography}
\end{document}